\input amstex
\documentstyle{amsppt}
\pageheight{194mm}
\pagewidth{133mm}
\magnification\magstep1

\document


\def\nologo{\let\logo@\empty}

\def\ad{\operatorname{ad}}
\def\Ad{\operatorname{Ad}}

\def\Aut{\operatorname{Aut}}
\def\Aut{\operatorname{Aut}}

\def\BS{\operatorname{BS}}
\def\class{\operatorname{class}}

\def\End{\operatorname{End}}

\def\fil{\operatorname{fil}}

\def\gr{\operatorname{gr}}
\def\Hom{\operatorname{Hom}}
\def\Her{\operatorname{Her}}

\def\Im{\operatorname{Im}}

\def\Int{\operatorname{Int}}

\def\Ker{\operatorname{Ker}}
\def\Lie{\operatorname{Lie}}

\def\modu{\;\operatorname{mod}}
\def\Mor{\operatorname{Mor}}

\def\nspl{\operatorname{nspl}}

\def\rank{\operatorname{rank}}
\def\Re{\operatorname{Re}}

\def\sign{\operatorname{sign}}
\def\SL{\operatorname{SL}}

\def\Spec{\operatorname{Spec}}

\def\spl{\operatorname{spl}}

\def\Sym{\operatorname{Sym}}
\def\toric{\operatorname{toric}}
\def\abtoric{\operatorname{|toric|}}

\def\tsCu{\tsize\bigsqcup}
\def\tOp{\tsize\bigoplus}
\def\tp{\tsize\prod}
\def\ts{\tsize\sum}

\def\bC{\bold C}

\def\bG{\bold G}

\def\bi{\bold i}

\def\bN{\bold N}
\def\bP{\bold P}
\def\bQ{\bold Q}
\def\br{\bold r}
\def\bR{\bold R}
\def\bS{\bold S}

\def\bZ{\bold Z}

\def\cB{{\Cal B}}
\def\cC{{\Cal C}}
\def\cD{{\Cal D}}

\def\cL{{\Cal L}}

\def\cO{{\Cal O}}

\def\cS{{\Cal S}}

\def\cW{{\Cal W}}

\def\fg{{\frak g}}
\def\fh{{\frak h}}

\def\fsl{{\frak s\frak l}}

\def\a{\alpha}
\def\b{\beta}
\def\d{\delta}
\def\g{\gamma}
\def\G{\Gamma}
\def\lam{\lambda}
\def\Lam{\Lambda}

\def\sig{\sigma}
\def\Sig{\Sigma}

\def\ve{\varepsilon}
\def\vf{\varphi}

\def\z{\zeta}

\def\.{$.\;$}
\def\an{{\text{\rm an}}}
\def\add{{\text{\rm add}}}

\def\gp{{\text{\rm gp}}}
\def\loga{{\text{\rm log}}}
\def\mult{{\text{\rm mult}}}
\def\nilp{{\text{\rm nilp}}}

\def\resp.{\text{resp}.\;}

\def\sym{{\text{\rm sym}}}

\def\triv{{\text{\rm triv}}}

\def\val{{\text{\rm val}}}

\def\O^logten{\cO\^log\otimes}

\let\bs=\backslash

\let\da=\downarrow

\let\la=\leftarrow

\let\lan=\langle
\let\lan=\langle

\let\hra=\hookrightarrow

\let\ox=\otimes
\let\op=\oplus
\let\Op=\bigoplus

\let\ran=\rangle

\let\sub=\subset

\let\t=\tilde

\let\x=\times

\def\Dc{\check{D}}


\magnification\magstep1

\NoRunningHeads

\document

\vskip20pt

\topmatter

\title
Classifying spaces of degenerating mixed Hodge structures, II: 
Spaces of $\SL(2)$-orbits
\endtitle

\author
Kazuya Kato\footnote{\text{Partially supported by NFS grant DMS 1001729.}},
Chikara Nakayama\footnote{\text{Partially supported by JSPS Grants-in-Aid for Scientific Research (C) 18540017,
(C) 22540011.}},
Sampei Usui\footnote{\text{Partially supported by JSPS Grant-in-Aid for Scientific Research (B) 19340008.}}
\endauthor

\dedicatory
Dedicated to the memory of Professor Masayoshi Nagata
\enddedicatory


\address
\newline
{\rm Kazuya KATO}
\newline
Department of Mathematics
\newline
University of Chicago
\newline
5734 S.\ University Avenue
\newline
Chicago, Illinois 60637, USA
\newline
{\tt kkato\@math.uchicago.edu}
\endaddress

\address
\newline
{\rm Chikara NAKAYAMA}
\newline
Graduate School of Science and Engineering
\newline
Tokyo Institute of Technology 
\newline
Meguro-ku, Tokyo, 152-8551, Japan
\newline
{\tt cnakayam\@math.titech.ac.jp}
\endaddress

\address
\newline
{\rm Sampei USUI}
\newline
Graduate School of Science
\newline
Osaka University
\newline
Toyonaka, Osaka, 560-0043, Japan
\newline
{\tt usui\@math.sci.osaka-u.ac.jp}
\endaddress

\abstract
We construct an enlargement of the classifying space of mixed Hodge structures with polarized graded quotients, by adding mixed Hodge theoretic version of $\SL(2)$-orbits.  
This space has a real analytic structure and a log structure with sign. 
The $\SL(2)$-orbit theorem in several variables for mixed Hodge 
structures can be understood naturally with this space. 
\endabstract

\footnote"{}"{2000 {\it Mathematics Subject Classification}.
Primary 14C30; Secondary 14D07, 32G20.}

\endtopmatter
\NoRunningHeads

\head
\S0. Introduction
\endhead

\medskip

\hskip80pt
{\it L'impossible voyage aux points \`a l'infini

\hskip80pt
N'a pas fait battre en vain le coeur du g\'eom\`etre}

\hskip184pt
translated by Luc Illusie 

\vskip12pt

This is Part II of our series of papers in which we study degeneration of mixed Hodge structures. 

{\bf 0.1.}
We first review the case of pure Hodge structures. 
Let $D$ be the classifying space of polarized Hodge structures of given weight and given Hodge numbers, defined by Griffiths [G]. 
Let $F_t\in D$ be a variation of polarized Hodge structure with complex analytic parameter $t=(t_1, \dots, t_n)$, $t_1\cdots t_n\neq 0$, which degenerates when $t\to 0=(0, \dots, 0)$. 
It is often asked how $F_t$ and invariants of $F_t$, like Hodge metric of $F_t$ etc., behave when $t\to 0$. 
Usually, $F_t$ diverges in $D$ and invariants of $F_t$ also diverge.

There are two famous theorems concerning the degeneration of $F_t$, which will be roughly reviewed in $0.3$ below. 

\medskip

(1) Nilpotent orbit theorem (\cite{Sc}).

\medskip

(2) $\SL(2)$-orbit theorem (\cite{Sc} and \cite{CKS}). 

\medskip
In \cite{KU2} and \cite{KU3}, we constructed enlargements  $D_{\SL(2)}$ and $D_{\Sig}$ of $D$, respectively. 
Roughly speaking, these theorems (1) and (2) are interpreted as in (1)$'$ and (2)$'$ below, respectively (see \cite{KU3}).

\medskip

(1)$'\;$  $(F_t\mod \G)\in \G\bs D$ converges in $\G \bs D_{\Sig}$, and asymptotic behaviors of invariants of $F_t$ are described by coordinate functions around the limit point on $\G \bs D_{\Sig}$. 

\medskip

(2)$'\;$ $F_t\in D$ converges in $D_{\SL(2)}$, and asymptotic behaviors of invariants of $F_t$ are described by coordinate functions around the limit point on $D_{\SL(2)}$ (see $0.2$ below).
\medskip

Here in (1)$'$, $\Gamma$ is the monodromy group of $F_t$ which acts on $D$ and $\Sig$ is a certain cone decomposition which is chosen suitably for $F_t$. The space $\G \bs D_{\Sig}$ is a kind of toroidal partial compactification of the quotient space $\G \bs D$, and has a kind of complex analytic structure. The space $D_{\SL(2)}$ has a kind of 
real analytic structure. For the study of asymptotic behaviors of real analytic objects like Hodge metrics, $D_{\SL(2)}$ is a nice space to work with.
\medskip

{\bf 0.2.}
Now let $D$ be the classifying space of mixed Hodge structures whose graded quotients for the weight filtrations are polarized, defined in \cite{U1}. The purpose of this paper is to construct an enlargement $D_{\SL(2)}$ of $D$, which is a mixed Hodge theoretic version of
$D_{\SL(2)}$ in \cite{KU2}. A mixed Hodge theoretic version of the $\SL(2)$-orbit theorem of \cite{CKS} was obtained in \cite{KNU1}, and it is also interpreted in the form (2)$'$ above, by using the present $D_{\SL(2)}$ (see \S4.1 of this paper). 

In Part I (\cite{KNU2}) of this series of papers, we constructed the Borel-Serre space $D_{\BS}$ which 
contains $D$ as a dense open subset and which is a real analytic manifold with corners like the original Borel-Serre space in \cite{BS}. 
These spaces $D_{\SL(2)}$ and $D_{\BS}$ belong to  
the following fundamental diagram of eight enlargements of $D$
whose constructions will be given in these series of papers. This fundamental diagram for the pure case $0.1$ was constructed in \cite{KU3}. 
$$
\matrix
&&&&&D_{\SL(2),\val}&\hra
&D_{\BS,\val}&\\
&&&&&&&&\\
&&&&&\da&&\da &\\
&&&&&&&&\\
&D_{\Sig,\val} &\la &D_{\Sig,\val}^\sharp
&\to&D_{\SL(2)}
&&D_{\BS}&\\
&&&&&&&&\\
&\da &&\da&&&&&\\
&&&&&&&&\\
&D_{\Sig}&\la
&D_{\Sig}^\sharp&&&&&
\endmatrix
$$

In the next parts of this series, we 
will construct the rest spaces in this diagram.  
Among them, $D_{\Sig}$ is the space of nilpotent orbits. 
Degenerations of mixed Hodge structures of geometric origin also satisfy a nilpotent orbit theorem (\cite{SZ}, \cite{K}, \cite{Sa}, \cite{P1} etc.;
a review is given in \cite{KNU1}, \S12.10). 
In the next papers in this series, we plan to interpret this in the style (1)$'$ above, by using  $D_{\Sig}$ in this diagram. 
\medskip

{\bf 0.3.}  
We explain the contents of the above $0.1$ and $0.2$ more precisely (but still roughly). 

The nilpotent orbit theorem (in the pure case $0.1$ and in the mixed case $0.2$ also) says roughly that when $t = (t_1,\dots,t_n) \to 0$, we have
$$
(F_t\bmod \G) \;\; \sim\;\;  (\exp(\ts_{j=1}^n z_jN_j)F\bmod \G)
$$
for some fixed Hodge filtration $F$ ($\sim$ expresses \lq\lq very near", but the precise meaning of it is not explained here), where $z_j$ is a branch of $(2\pi i)^{-1}\log(t_j)$ and $N_j$ is the logarithm of the local monodromy of $F_t$ around the divisor $t_j=0$. In \cite{KU3} for the pure case and in the next papers in this series for the mixed case, this is interpreted as the convergence 
$$
(F_t\bmod \G)\;\;\to \;\;((\sig, Z) \bmod \G)\in \G\bs D_{\Sig},
$$
where $\sig$ is the cone $\ts_{j=1}^n \bR_{\geq 0}N_j$ and $Z$ is the orbit
$\exp(\ts_{j=1}^n \bC N_j)F$. 
  (As in the pure case, as a set, $D_{\Sig}$ is a set of such pairs 
$(\sig, Z)$.)

\medskip

The $\SL(2)$-orbit theorem in the pure case $0.1$ obtained in \cite{CKS} says roughly that when $t\to 0$,  $t_j\in \bR_{>0}$, and $y_j/y_{j+1}\to \infty$, 
where $y_j=-(2\pi)^{-1}\log(t_j)$ for $1\le j \le n$ $(y_{n+1} = 1)$, we have
$$
F_t \;\; \sim  \;\; \rho\left(\pmatrix \sqrt{y_1} & 0 \\ 0 &1/\sqrt{y_1}\endpmatrix, \dots, \pmatrix \sqrt{y_n} & 0 \\ 0 & 1/\sqrt{y_n}\endpmatrix\right)\varphi(\bi),
$$
($\sim$ expresses \lq\lq very near'' again) 
where $\rho$ is a homomorphism of algebraic groups $\SL(2,\bR)^n\to \Aut(D)$, $\varphi$ is a complex analytic map $\fh^n\to D$ from the product $\fh^n$ of copies of the upper half plane $\fh$, satisfying $\varphi(gz)=\rho(g)\varphi(z)$ for any $g\in \SL(2,\bR)^n$ and $z\in \fh^n$, and where $\bi=(i, \dots, i)\in \fh^n$. 
In \cite{KU3}, this is interpreted as the convergence 
$$
F_t\;\;\to \;\;\text{class}(\rho, \varphi) \in D_{\SL(2)}.
$$

The $\SL(2)$-orbit theorem in the mixed case $0.2$ obtained in \cite{KNU1} says roughly that when $t\to 0$,  $t_j\in \bR_{>0}$, and $y_j/y_{j+1}\to \infty$, 
where $y_j=-(2\pi)^{-1}\log(t_j)$ for $1\le j \le n$ $(y_{n+1} = 1)$, we have
$$
F_t \;\;\sim  \;\;\text{lift}\left(\tsize\bigoplus_{w\in \bZ} \; y_1^{-w/2}\rho_w\left(\pmatrix \sqrt{y_1} & 0 \\ 0 &1/\sqrt{y_1}\endpmatrix, \dots, \pmatrix 
\sqrt{y_n} & 0 \\ 0 & 1/\sqrt{y_n}\endpmatrix\right)\right)\br,
$$
where $(\rho_w,\varphi_w)$ $(w \in \bZ)$ is the $\SL(2)$-orbit of pure weight $w$ associated to the filtration on $\gr^W_w$ induced from $F_t$, $\br$ is a certain point of $D$ which induces $\varphi_w(\bi)$ on each $\gr^W_w$, and \lq\lq lift" is the lifting to $\Aut(D)$  by the canonical splitting of the weight filtration associated to $\br$ (see 1.2). 
For details, see \cite{KNU1}, and also \S2.4 of this paper. 
By using the space $D_{\SL(2)}$ of this paper, this is interpreted as the convergence 
$$
F_t\;\;\to \;\;\text{class}((\rho_w, \varphi_w)_{w\in \bZ}, \br) \in D_{\SL(2)}.
$$
Since $D_{\SL(2)}$ has a real analytic structure, we can discuss the differential of the extended period map $t\mapsto F_t$ at $t=0$. 
We hope such delicate structure of $D_{\SL(2)}$ is useful for the study of degeneration. 

\medskip

{\bf 0.4.} 
  Precisely, there are two natural spaces $D^I_{\SL(2)}$ and 
$D^{II}_{\SL(2)}$ which can sit in the place of 
$D_{\SL(2)}$ in the fundamental diagram.  
  They coincide in the pure case and coincide always as sets, but 
do not coincide in general. 
  What we wrote in the above 0.3 are valid for both. 
  They have good properties respectively so that we do not choose 
one of them as a standard one.  
  See 3.2.1 for more survey.

\medskip

{\bf 0.5.} 
The organization of this paper is as follows. 
In \S1, we give preliminaries about basic facts on mixed Hodge structures. 
In \S2, we define the space $D_{\SL(2)}$ as a set. 
In \S3, we endow this set with topologies and with real analytic structures (these spaces 
$D^I_{\SL(2)}$ and 
$D^{II}_{\SL(2)}$ 
are not necessarily real analytic spaces, but they have the sheaves of real analytic functions which we call the real analytic structures). 
We study properties of these spaces. 
In \S4, we consider how the degenerations of mixed Hodge structures are related to these spaces. 
\medskip

{\bf 0.6.}  
A large part of this paper was written while one of the authors (K. Kato) was a visitor of University of Cambridge whose hospitality (especially by Professor John Coates) is gratefully appreciated.
The authors thank the referee for careful reading. 

\medskip

The poem at the beginning is a translation by Professor Luc Illusie of a Japanese poem composed by two of the authors (K. Kato and S. Usui).
These poems were placed at the beginning of [KU3]. We put the French version 
here again as it well capture the spirit of this paper.

\vskip20pt

{\bf Notation}
\bigskip

Fix a quadruple
\medskip
\item{}
$\Phi_0=(H_0, W, (\lan\;,\;\ran_w)_{w\in \bZ}, (h^{p,q})_{p, q\in \bZ})$,
\medskip

\noindent
where 
\medskip
\item{}
$H_0$ is a finitely generated free $\bZ$-module, 
\medskip
\item{}
$W$ is an increasing filtration on $H_{0,\bR}:=\bR\otimes_\bZ H_0$ defined over $\bQ$, 
\medskip
\item{}
$\lan\;,\;\ran_w$ is a non-degenerate $\bR$-bilinear form 
$\gr^W_w \times\gr^W_w \to \bR$ defined over $\bQ$ for each $w\in \bZ$ which is symmetric if $w$ is even and anti-symmetric if $w$ is odd, and 
\medskip
\item{}
$h^{p,q}$ is a non-negative integer given for $p, q\in\bZ$ such that 
$h^{p,q}=h^{q,p}$, $\rank_\bZ(H_0) = \tsize\sum_{p,q} h^{p,q}$, and 
$\dim_\bR(\gr^W_w)= \sum_{p+q=w} h^{p,q}$ for all $w$.
\medskip

Let $\Dc$ be the set of all decreasing filtrations $F$ on $H_{0,\bC}:=\bC\otimes_\bZ H_0$ 
satisfying the following two conditions.
\medskip

(1) $\dim (F^p(\gr^W_{p+q})/F^{p+1}(\gr^W_{p+q})) = h^{p,q}$ for any $p, q \in \bZ$. 

\medskip
(2) $\lan\;,\;\ran_w$ kills $F^p(\gr^W_w)\times F^q(\gr^W_w)$ for any 
$p,q,w\in\bZ$ such that $p+q>w$.
\medskip

\noindent
Here $F(\gr^W_w)$ denotes the filtration on $\gr^W_{w,\bC}:= \bC \otimes_\bR\gr^W_w$ induced by $F$.
\medskip

Let $D$ be the set of all decreasing filtrations $F\in\Dc$ which also satisfy the following condition.
\medskip

(3) $i^{p-q}\lan x, \bar x\ran_w > 0$ for any non-zero $x \in F^p(\gr^W_w) \cap \overline{F^q(\gr^W_w)}$ and any $p, q, w \in \bZ$ with $p+q = w$.
\medskip

Then, $D$ is an open subset of $\Dc$ and, for each $F \in D$ and $w \in \bZ$,  $F(\gr^W_w)$ is a Hodge structure on $(H_0\cap W_w)/(H_0\cap W_{w-1})$ of weight $w$ with Hodge number $(h^{p,q})_{p+q=w}$ which is polarized by $\lan\;,\;\ran_w$.
The space $D$ is the classifying space of mixed Hodge structures of type $\Phi_0$ introduced in \cite{U1}, which is a natural generalization to the mixed case of the Griffiths domain in \cite{G}.
These two are related by taking graded quotients by $W$ as follows.
\medskip
\item{}
$D(\gr^W_w)$: the $D$ for $((H_0\cap W_w)/(H_0\cap W_{w-1}), \lan\;,\;\ran_w, (h^{p,q})_{p+q=w})$ for each $w\in \bZ$. 
\medskip
\item{}
$D(\gr^W) = \tp_{w\in \bZ} \; D(\gr^W_w)$.
\medskip
\item{}
$D \to D(\gr^W)$, $F \mapsto F(\gr^W):= (F(\gr^W_w))_{w\in \bZ}$, the canonical surjection.
\medskip

For $A=\bZ, \bQ, \bR,$ or $\bC$, 
\medskip

\item{}
$G_A$: the group of all $A$-automorphisms $g$ of $H_{0,A}:=A \otimes_{\bZ}H_0$ 
compatible with $W$ such that $\gr^W_w(g): \gr^W_w\to \gr^W_w$ are compatible with 
$\lan\;,\;\ran_w$ for all $w$. 
\medskip
\item{}
$G_{A,u}:=\{g\in G_A\;|\;\gr^W_w(g)=1\;\text{for all}\;w\in \bZ\}$, the {\it unipotent radical} of $G_A$.
\medskip
\item{}
$G_A(\gr^W_w)$: the $G_A$ of $((H_0\cap W_w)/(H_0\cap W_{w-1}), \lan\;,\;\ran_w)$ for each $w\in \bZ$.
\medskip
\item{}
$G_A(\gr^W):= \tp_w G_A(\gr^W_w)$.
\medskip

Then, $G_A/G_{A,u}= G_A(\gr^W)$, and $G_A$ is a semi-direct product of $G_{A,u}$ and $G_A(\gr^W)$.
\medskip

The natural action of $G_\bC$ on $\Dc$ is transitive, and $\Dc$ is a complex homogeneous space under the action of $G_{\bC}$. 
Hence $\Dc$ is a complex analytic manifold. 
An open subset $D$ of $\Dc$ is also a complex analytic manifold.
However, the action of $G_\bR$ on $D$ is not transitive in general (see the equivalent conditions (4), (5) below).
The subgroup $G_\bR G_{\bC,u}$ of $G_\bC$ acts always transitively on $D$, and the action of $G_{\bC,u}$ on each fiber of $D\to D(\gr^W)$ is transitive. 

\medskip

\item{}
$\spl(W)$: the set of all isomorphisms 
$s : \gr^W=\tsize\bigoplus_w \gr^W_w @>\sim>> H_{0,\bR}$
of $\bR$-vector spaces such that for any $w\in \bZ$ and $v \in
\gr^W_w$, $s(v) \in W_w$ and $v = (s(v)\bmod W_{w-1})$.
\medskip
\item{}
  We have the action 
$G_{\bR, u}\times \spl(W) \to \spl(W)$, $(g, s) \mapsto gs$.
\medskip

\noindent
For a fixed $s\in\spl(W)$, we have a bijection 
$G_{\bR,u}@>\sim>> \spl(W),\;g\mapsto gs$. 
Via this bijection, we endow $\spl(W)$ with a structure of a real analytic manifold.
\medskip

\item{}
$D_{\spl}:= \{s(F) \;|\; s\in \spl(W),\, F\in D(\gr^W)\} \sub D$,
the subset of {\it $\bR$-split} elements.
\medskip

\noindent
Here $s(F)^p:= s(\bigoplus_w F_{(w)}^p)$ for $F = (F_{(w)})_w \in D(\gr^W)$.
\medskip

\item{}
$D_{\nspl}:=D\smallsetminus D_{\spl}$.
\medskip

Then, $D_{\spl}$ is a closed real analytic submanifold of $D$, and we have a real analytic isomorphism
$\spl(W) \times D(\gr^W)@>\sim>>  D_{\spl}$, $(s, F)\mapsto s(F)$.

\medskip

The following two conditions are equivalent (\cite{KNU2}, Proposition 8.7). 
\medskip

(4) $D$ is $G_\bR$-homogeneous. 
\medskip

(5) $D=D_{\spl}$.
\medskip

\noindent
For example, if there is $w\in \bZ$ such that $W_w=H_{0,\bR}$ and $W_{w-2}=0$, then the above equivalent conditions are satisfied.
But in general these conditions are not satisfied (see Examples I, III, IV in 1.1).
\medskip

For $A=\bQ, \bR, \bC$, 
\medskip

\item{}
$\fg_A:=\Lie(G_A)$ which is identified with 
$\{X \in \End_A(H_{0,A}) \;|\; X(W_w) \sub W_w\;\,\text{for all}\; w,$
\item{}
$\lan \gr^W_w(X)(x), y\ran_w+ \lan x, 
\gr^W_w(X)(y)\ran_w=0\;\,\text{for all}\; w, x, y\}$. 
\medskip
\item{}
$\fg_{A,u}:= \Lie(G_{A,u})=\{X\in\fg_A\;|\; \gr^W_w(X)=0\;\,\text{for all}\;w\}$. 
\medskip
\item{}
$\fg_A(\gr^W_w)$: the $\fg_A$ of $((H_0\cap W_w)/(H_0\cap W_{w-1}),\lan\;,\;\ran_w)$ for each $w\in \bZ$.
\medskip
\item{}
$\fg_A(\gr^W):=\tOp_{w \in \bZ}\fg_A(\gr^W_w)$.

\medskip

  Then, $\fg_A/\fg_{A,u}=\fg_A(\gr^W)$. 

\medskip

\vskip20pt

\head 
\S1. Basic facts
\endhead
\medskip

We examine some examples, review some basic facts, and fix further notation which will be used in this paper. 

\vskip20pt

\head
\S1.1. Examples
\endhead
\medskip

\medskip

{\bf 1.1.1.}
We give six simple examples 0, I, II, $\dots$, V of $D$ 
for which the set $\{w\in \bZ\;|\;\gr^W_w\neq 0\}$ is $\{-1\}$, $\{0,-2\}$, $\{0,-1\}$, $\{0,-3\}$, $\{0,-1,-2\}$, $\{0, 1\}$, respectively. 
Among these, the examples I, II, III are already presented in \cite{KNU2}, 1.10--1.12 to illustrate the results in that paper on each step.
All these examples will be retreated also to illustrate the results in this paper on each step.
\medskip

{\bf Example 0.} 
(This example belongs to the pure case, though Example I--Example V below do not.) 
Let $H_0=\bZ^2=\bZ e_1+\bZ e_2$. 
Let $W$ be the increasing filtration on $H_{0,\bR}$ defined by
$$
W_{-2}=0\subset W_{-1}=H_{0,\bR}.
$$
Let $\langle e_2, e_1\rangle_{-1}=1$. 
Let $h^{-1,0}=h^{0,-1}=1$, and let $h^{p,q}=0$ for all other $(p, q)$. 

For $\tau\in \bC$, let $F(\tau)$ be the decreasing filtration on $H_{0,\bC}$ defined by
$$F(\tau)^1=0\sub F(\tau)^0=\bC(\tau e_1+e_2) \sub F(\tau)^{-1}=H_{0,\bC}.$$
Then we have an isomorphism of complex analytic manifolds
$$
D\simeq \fh, 
$$
where $\fh$ is the upper half plane $\{x+iy\;|\;x, y\in \bR, y>0\}$, in which 
$\tau\in \fh$ corresponds to $F(\tau)\in D$. 
This isomorphism naturally extends to $\Dc \simeq \bP^1(\bC)$.
\medskip

{\bf Example I.} 
Let $H_0=\bZ^2=\bZ e_1+ \bZ e_2$, let $W$ be the increasing filtration on 
$H_{0,\bR}$ defined by 
$$
W_{-3} =0\;\sub\; W_{-2}=W_{-1}= \bR e_1\;\sub \; W_0=H_{0,\bR}.
$$
For $j=1$ (resp. $j=2$), let  $e_j'$  be the image of $e_j$ in $\gr^W_{-2}$ (resp.
$\gr^W_0$). Let $\lan e_2', e_2'\ran_0=1$, $\lan e_1', e_1'\ran_{-2}=1$, and let $h^{0,0}=h^{-1,-1}=1$, $h^{p,q}=0$ for all the other $(p, q)$. 

We have an isomorphism of complex analytic manifolds
$$
D\simeq \bC.
$$  
For $z\in \bC$, the corresponding  $F(z)\in D$ is defined as 
$$
F(z)^1=0 \;\sub \; F(z)^0= \bC(ze_1+e_2) \;\sub \;F(z)^{-1}= H_{0,\bC}.
$$

The group $G_{\bZ, u}$ is isomorphic to $\bZ$ and is generated by 
$\g\in G_\bZ$ which is defined as 
$$
\g(e_1)=e_1, \quad \g(e_2)=e_1+e_2.
$$ 
We have
$$
G_{\bZ,u}\bs D \simeq \bC^\times, 
$$
where $(F(z)\bmod G_{\bZ,u})$ corresponds to $\exp(2\pi iz) \in \bC^\times$. 

This space $G_{\bZ, u}\bs D$ is the classifying space of extensions of mixed Hodge structures of the form $0 \to \bZ(1) \to \;\ast \to \bZ \to 0$.

In this case, $D(\gr^W)$ is a one point set.
\medskip

{\bf Example II.} 
Let $H_0=\bZ^3=\bZ e_1+\bZ e_2+\bZ e_3$, let
$$
W_{-2}=0\; \sub \; W_{-1}=\bR e_1+\bR e_2\;\sub\; W_0=H_{0,\bR}.
$$ 
For $j=1, 2$ (resp. $3$), let $e_j'$ be the image of $e_j$ in $\gr^W_{-1}$ (resp. $\gr^W_0$). Let $\lan e_3', e_3'\ran_0=1$, $\lan e_2', e_1'\ran_{-1}= 1$, and let $h^{0,0}=h^{0,-1}=h^{-1,0}=1$, $h^{p,q}=0$ for all the other $(p,q)$.

Then we have isomorphisms of complex analytic manifolds
$$
D\simeq\fh\times \bC, \quad D(\gr^W)\simeq \fh.
$$ 
Here $(\tau, z)\in \fh \times \bC$ corresponds to $F=F(\tau,z)\in D$ given by 
$$
F^1=0 \;\sub \;F^0 = \bC(\tau e_1+e_2)+\bC (ze_1+ e_3)\;
\sub \; F^{-1}=H_{0,\bC}.
$$ 
The induced isomorphism $D(\gr^W)=D(\gr^W_{-1})\simeq \fh$ is identified with the isomorphism $D\simeq \fh$ in Example 0.

The group $G_{\bZ,u}$ is isomorphic to $\bZ^2$, where $(a, b)\in \bZ^2$
corresponds to the element of $G_\bZ$ which sends $e_j$ to $e_j$ for
$j=1, 2$ and sends $e_3$ to $ae_1+be_2+e_3$.
The quotient space $G_{\bZ,u}\bs D$ is the \lq\lq universal elliptic curve'' over the upper half plane $\frak h$. 
For $\tau \in \fh$, the fiber of $G_{\bZ,u}\bs D \to D(\gr^W)=\fh$ over $\tau$ is identified with the elliptic curve $E_\tau:=\bC /(\bZ \tau +\bZ)$. 
The Hodge structure on $H_0\cap W_{-1}$ corresponding to $\tau$ is 
isomorphic to $H^1(E_\tau)(1)$. 
Here  $H^1(E_\tau)$ denotes the Hodge structure $H^1(E_\tau, \bZ)$ of weight $1$ endowed with the Hodge filtration and $(1)$ here denotes the Tate twist. 
The fiber of $G_{\bZ,u}\bs D\to \fh$ over $\tau$ is the classifying space of extensions of mixed Hodge structures of the form
$$
0 \to H^1(E_\tau)(1) \to \;\ast \to \bZ \to 0.
$$
\medskip

{\bf Example III.}
Let $H_0=\bZ^3=\bZ e_1+\bZ e_2+\bZ e_3$, let
$$
W_{-4}=0 \;\sub \; W_{-3}=W_{-1}=\bR e_1+\bR e_2\;\sub \; W_0=H_{0,\bR}.
$$ 
For $j=1, 2$ (resp. $3$), let $e_j'$ be the image of $e_j$ in $\gr^W_{-3}$ (resp. $\gr^W_0$). Let $\lan e_3', e_3'\ran_0=1$, $\lan e_2', e_1'\ran_{-3}= 1$, and let $h^{0,0}=h^{-1,-2}=h^{-2,-1}=1$, $h^{p,q}=0$ for all the other $(p,q)$.

Then we have isomorphisms of complex analytic manifolds
$$
D\simeq \fh \times \bC^2, \quad D(\gr^W)\simeq\fh.
$$ 
Here $(\tau, z_1, z_2)\in \fh \times \bC^2$ corresponds to $F=F(\tau,z_1,z_2)\in D$ given by 
$$
F^1=0 \;\sub\; F^0 = \bC (z_1 e_1+z_2 e_2+ e_3)\;\sub \;
F^{-1}= F^0 +\bC(\tau e_1+e_2)\;\sub  \;F^{-2}=H_{0,\bC}.
$$ 
The induced isomorphism $D(\gr^W)=D(\gr^W_{-3}) \simeq \fh$ is identified with the isomorphism $D\simeq \fh$ in Example 0 ($F\in D(\gr^W)$ corresponds to the twist $F(-1)$ of $F$, which belongs to the $D$ in Example 0.)

The group $G_{\bZ,u}$ is the same as in Example II. 
The Hodge structure on $H_0\cap W_{-3}$ corresponding to $\tau \in \fh\simeq D(\gr^W_{-3})$ is isomorphic to $H^1(E_\tau)(2)$. 
The fiber of $G_{\bZ,u}\bs D \to D(\gr^W)\simeq\fh$ over $\tau \in \fh$ is the classifying space of extensions of mixed Hodge structures of the form
$$
0 \to H^1(E_\tau)(2) \to \;\ast \to \bZ \to 0.
$$
\medskip

{\bf Example IV.}
Let $H_0=\bZ^4=\bZ e_1+\bZ e_2+\bZ e_3+\bZ e_4$, let
$$
W_{-3}=0 \;\sub \; W_{-2}=\bR e_1\;\sub \; W_{-1}=W_{-2}+\bR e_2+\bR e_3\;\sub \; W_0=H_{0,\bR}.
$$ 
For $j=1$ (resp. $2, 3$, resp. $4$), let $e_j'$ be the image of $e_j$ in $\gr^W_{-2}$ (resp. $\gr^W_{-1}$, resp. $\gr^W_0$). 
Let $\lan e_4', e_4'\ran_0=1$, $\lan e_1', e_1'\ran_{-2}= 1$, and $\lan e_3', e_2'\ran_{-1}=1$, and let $h^{0,0}=h^{0,-1}=h^{-1,0}=h^{-1,-1}=1$, $h^{p,q}=0$ for all the other $(p,q)$.

Then we have isomorphisms of complex analytic manifolds
$$
D=\fh \times \bC^3,\quad
D(\gr^W)=D(\gr^W_{-1})= \fh.
$$
Here $(\tau, z_1,z_2,z_3)\in \fh \times \bC^3$ corresponds to $F=F(\tau, z_1,z_2,z_3)\in D$ given by $F^{-1}=H_{0,\bC}$, $F^1=0$, and
$$
F^0 = \bC (z_1e_1+\tau e_2+e_3)+\bC(z_2e_1+z_3e_2+e_4).
$$ 
The induced isomorphism $D(\gr^W)=D(\gr^W_{-1})\simeq \fh$ is identified with the isomorphism $D\simeq \fh$ in Example 0. 

There is a bijection $G_{\bZ,u} \simeq \bZ^5$ (but not a group isomorphism), where $(a_j)_{1\le j\le5}\in \bZ^5$ corresponds to the element of $G_{\bZ,u}$ which sends $e_1$ to $e_1$, $e_2$ to $a_1e_1+e_2$, $e_3$ to $a_2e_1+e_3$, 
and $e_4$ to $a_3e_1+a_4e_2+a_5e_3+e_4$.
\medskip

{\bf Example V.}
Let $H_0=\bZ^5=\bZ e_1+\bZ e_2+\bZ e_3+\bZ e_4+\bZ e_5$, let
$$
W_{-1}=0 \;\sub \; W_0=\bR e_1+\bR e_2+\bR e_3\;\sub \; W_1=H_{0,\bR}.
$$ 
For $j=1, 2, 3$ (resp. $4, 5$), let $e_j'$ be the image of $e_j$ in $\gr^W_0$ (resp. $\gr^W_1$). Let $\lan e_5', e_4'\ran_1= 1$, $\lan e_1', e_3'\ran_0=2$, $\lan e_2', e_2'\ran_0=-1$, and $\lan e_j', e_k'\ran_0=0$ $(j+k\neq4$, $1\leq j$, $k \leq 3)$, and let $h^{1,-1}=h^{0,0}=h^{-1,1}=h^{1,0}=h^{0,1}=1$, and $h^{p,q}=0$ for all the other $(p,q)$.

Let $\fh^{\pm}=\{x+iy\;|\;x, y\in \bR,\; y\neq 0\}= \fh \sqcup (-\fh)$. 
Then we have isomorphisms of complex analytic manifolds
$$
D\simeq \fh^{\pm} \times \fh \times \bC^3,\quad
D(\gr^W_0)\simeq  \fh^{\pm},\quad 
D(\gr^W_1)\simeq  \fh.
$$
Here $(\tau_0,\tau_1, z_1,z_2,z_3)\in \fh^{\pm} \times\fh \times \bC^3$ corresponds to $F=F(\tau_0,\tau_1, z_1,z_2,z_3)\in D$ given 
by $F^2=0$, $F^{-1}=H_{0,\bC}$, and
$$
\align
&F^1 = \bC(\tau_0^2e_1+2\tau_0 e_2+e_3) 
+ \bC(z_1e_1+z_2e_2+\tau_1 e_4+e_5),\\
&F^0= F^1+\bC(\tau_0 e_1+e_2)
+\bC(z_3e_1+e_4).
\endalign
$$

Let $F(\tau)$ be the filtration in Example 0 corresponding to $\tau\in\fh$.
The induced isomorphism $D(\gr^W_1)\simeq \fh$ sends $\tau\in\fh$ to the Tate twist $F(\tau)(-1)$ of $F(\tau)$.
The induced isomorphism $D(\gr^W_0)\simeq \fh^{\pm}$ sends $\tau \in \fh^{\pm}$ to $\Sym^2(F(\tau))(-1)\in D(\gr^W_0)$ (see 1.1.2 below). 

The group $G_{\bZ,u}$ is isomorphic to $\bZ^6$, where $(a_j)_{1\le j\le6}\in \bZ^6$ corresponds to the element of $G_\bZ$ which sends $e_j$ to $e_j$ for $j=1,2,3$,  $e_4$ to $a_1e_1+a_2e_2+a_3e_3+e_4$, and $e_5$ to $a_4e_1+a_5e_2+a_6e_3+e_5$.
\medskip

{\bf 1.1.2.} {\it Remark.}
For the computations of Example V in 1.1.1 above and in 3.6 and 4.2.4 
later, we describe here the classifying space $D_2$ of polarized Hodge structures of weight $2$ 
underlain by the second symmetric power of the Tate twist $(-1)$ of $(H_0, \lan\;\;,\;\;\ran_{-1})$ 
in Example 0.

The domain $D(\gr^W_0)$ in Example V of 1.1.1 is 
identified with $D_2$ via the Tate twist. 

Let $H_0=\bZ^2=\bZ f_1+\bZ f_2$, let $W_0=0\;\sub \; W_1=H_{0,\bR}$, and let $\lan f_2, f_1\ran_1= 1$.
Then, $\Sym^2(H_0)=\bZ^3=\bZ e_1+\bZ e_2+\bZ e_3$, where $e_1:=f_1^2$, $e_2:=f_1f_2$, $e_3:=f_2^2$, and the induced polarization on $\Sym^2(H_0)$, 
which is defined by 
$$
\langle x_1x_2, y_1y_2\rangle_2= \langle x_1, y_1\rangle_1\langle x_2,y_2\rangle_1 + \langle x_1, y_2\rangle_1\langle x_2,y_1\rangle_1\quad
(x_j, y_j \in H_0, j= 1,2),
$$
is given by
$$
\lan e_1,e_3\ran_2=\lan e_3,e_1\ran_2=2,\quad
\lan e_2,e_2\ran_2=-1,\quad
\lan e_j,e_k\ran_2=0\quad\text{otherwise}.
$$
For $v=\omega_1e_1+\omega_2e_2+\omega_3e_3 \in \bC e_1+\bC e_2+\bC e_3$ to be Hodge type $(2,0)$, the Riemann-Hodge bilinear relations are
$$
\align
&\lan v,v\ran_2 = 4\omega_1\omega_3-\omega_2^2 = 0,\\
&\lan Cv,\bar v\ran_2 = i^2\lan v,\bar v\ran_2 
= -4\Re(\omega_1\bar\omega_3) + |\omega_2|^2 > 0,
\endalign
$$
where $C$ is the Weil operator.
Hence the classifying space $D_2$ and its compact dual $\check D_2$ of the Hodge structures of weight $2$, with Hodge type $h^{2,0}=h^{1,1}=h^{0,2}=1$ and $h^{p,q}=0$ otherwise, and with the polarization $\lan\;\;,\;\;\ran_2$, is as follows.
$$
\align
\check D_2&=\{\bC(\omega_1e_1+\omega_2e_2+\omega_3e_3)
\sub \bC e_1+\bC e_2+\bC e_3 \;|\; 4\omega_1\omega_3-\omega_2^2 = 0\}
\simeq \bP^1(\bC).
\endalign
$$
$$
\align
D_2&=\{\bC(\omega_1e_1+\omega_2e_2+\omega_3e_3)
\in \check D \;|\; -4\Re(\omega_1\bar\omega_3) + |\omega_2|^2 > 0\}
\simeq \fh^{\pm}. \tag1
\endalign
$$
The isomorphism is given by $\omega_1e_1+\omega_2e_2+\omega_3e_3 =\omega^2e_1+2\omega e_2+e_3 \leftrightarrow \omega$.

Assigning $g \in \SL(2,\bR)$ to $\sym^2(g) \in \Aut(H_{0,\bR}, \lan\;\;,\;\ran_2)$, we have an exact sequence
$$
1 \to \{\pm1\} \to \SL(2,\bR) \to \Aut(H_{0,\bR}, \lan\;\;,\;\ran_2) \to \{\pm1\} \to 1. \tag2
$$
The isomorphism (1) is compatible with (2).

\medskip

\vskip20pt

\head
\S1.2. Canonical splittings of weight filtrations for mixed Hodge structures
\endhead
\medskip

Let $W$ and $D$ be as in Notation at the end of Introduction.
In this section, we review the canonical splitting $s=\spl_W(F)\in \spl(W)$ of the weight filtration $W$ associated to $F\in D$, defined by the theory of Cattani-Kaplan-Schmid \cite{CKS}. 
This canonical splitting  $s$ appeared naturally in the $\SL(2)$-orbit theorem for mixed Hodge structures proved in our previous paper \cite{KNU1}.  
The definition of $s$ was reviewed in detail in Section 1 of \cite{KNU1}. The canonical splitting plays important roles in the present series of our papers.
\medskip

{\bf 1.2.1.} 
Let $F=(F_{(w)})_w\in D(\gr^W)$. Regard $F$ as the filtration
$\tsize\bigoplus_w F_{(w)}$ on $\gr^W_{\bC}=\tsize\bigoplus_w\gr^W_{w,\bC}$, and let $H^{p,q}_F= H^{p,q}_{F_{(p+q)}}\sub \gr^W_{p+q,\bC}$. 
Let
$$
L^{-1,-1}_\bR(F)=\{\delta\in 
\End_{\bR}(\gr^W)\;|\;\delta(H^{p,q}_F)\sub \tsize\bigoplus_{p'<p, q'<q}\; 
H^{p',q'}_F\;\text{for all}\;p, q\in\bZ\}.
$$ 
All elements of $L^{-1,-1}_\bR(F)$ are nilpotent.
Let 
$$
\cL=\End_\bR(\gr^W)_{\leq -2}
$$
be the set of all $\bR$-linear maps $\delta:\gr^W\to \gr^W$ such that $\delta(\gr^W_w)\sub \bigoplus_{w'\leq w-2} \gr^W_{w'}$ for any $w\in \bZ$. Denote 
$$
\cL(F)=L^{-1,-1}_\bR(F)\sub \cL.
$$
$\cL(F)$ is sometimes denoted simply by $L$.

In this \S1.2, we review the isomorphism of real analytic manifolds
$$
D\simeq \{(s, F, \delta)\in \spl(W) \times D(\gr^W)\times \cL\;|\; \delta\in \cL(F)\}
$$
obtained in the work \cite{CKS} (see 1.2.5 below). 
For $F'\in D$, the corresponding $(s, F, \delta)$ consists of $F=F'(\gr^W)$, $\delta=\delta(F')\in \cL(F)$ defined in 1.2.2 below, and the canonical splitting $s=\spl_W(F')$ of $W$ associated to $F'$ explained in 1.2.3 below. 
\medskip

{\bf 1.2.2.} 
For $F' \in D$, there is a unique pair $(s', \delta) \in \spl(W)\times \cL(F'(\gr^W))$ such that
$$
F' = s'(\exp(i\delta)F'(\gr^W))
$$
(\cite{CKS}).
This is the definition of $\delta=\d(F')$ associated to $F'$. 
\medskip

{\bf 1.2.3.} 
Let $F' \in D$, and let $s'\in \spl(W)$ and $\d\in \cL(F'(\gr^W))$ be as in 1.2.2. 
Then the canonical splitting $s=\spl_W(F')$ of $W$ associated to $F'$ is defined by
$$
s=s'\exp(\zeta),
$$
where $\zeta=\zeta(F'(\gr^W), \delta)$ is a certain element of $\cL(F'(\gr^W))$ determined by $F'(\gr^W)$ and $\delta$ in the following way. 

Let $\delta_{p,q}$ ($p,q\in \bZ$) be the $(p, q)$-Hodge component of $\delta$ with respect to $F'(\gr^W)$ defined by 
$$
\delta= \ts_{p, q} \;\delta_{p,q} \quad (\delta_{p,q}\in \cL_\bC(F'(\gr^W))=\bC \otimes_{\bR} \cL(F'(\gr^W))),
$$
$$
\delta_{p,q}(H^{k,l}_{F'(\gr^W)}) \sub H^{k+p,l+q}_{F'(\gr^W)} \qquad
\text{for all}\;k, l\in \bZ.
$$
Then the $(p, q)$-Hodge component $\z_{p, q}$ of $\zeta=\zeta(F'(\gr^W), \delta)$ with respect to $F'(\gr^W)$ is 
given as a certain universal Lie polynomial of $\delta_{p', q'}$ ($p', q'\in \bZ$,
$p'\leq -1$, $q'\leq -1$). 
See \cite{CKS} and Section 1 of \cite{KNU1}. 
For example, 
$$
\align
&\z_{-1,-1}=0,\\
&\z_{-1,-2}=-\tfrac{i}2\d_{-1,-2},\\
&\z_{-2,-1}=\tfrac{i}2\d_{-2,-1}.
\endalign
$$
\medskip

{\bf 1.2.4.} 
For $F\in D(\gr^W)$ and $\delta\in\cL(F)$, we define a filtration $\theta(F, \delta)$ on $\gr^W_\bC$ by 
$$
\theta(F, \delta)= \exp(-\zeta)\exp(i\delta)F,
$$ 
where $\zeta=\zeta(F, \delta)$ is the element of $\cL(F)$ associated to the pair $(F, \delta)$ as in 1.2.3. 

\proclaim{Proposition 1.2.5} 
We have an isomorphism of real analytic manifolds
$$
D\simeq \{(s, F, \delta)\in \spl(W) \times D(\gr^W) 
\times \cL\;|\; \delta\in \cL(F)\},\quad
F'\mapsto (\spl_W(F'), F'(\gr^W), \d(F')), 
$$
whose inverse is given by $(s, F, \d)\mapsto s(\theta(F, \delta))$.
\endproclaim

{\bf 1.2.6.} 
For $g=(g_w)_w\in G_\bR(\gr^W)=\tp_w\;G_\bR(\gr^W_w)$, we have
$$
g\theta(F, \delta)=\theta(gF, \Ad(g)\delta),
$$
where $\Ad(g)\delta= g\delta g^{-1}$. 
\medskip

{\bf 1.2.7.} 
For $F\in D(\gr^W)$, $\delta\in \cL(F)$ and $s\in \spl(W)$, the element $s(\theta(F, \delta))$ of $D$ belongs to $D_{\spl}$ if and only if $\delta=0$. 
\medskip

{\bf 1.2.8.} 
{\it Remark.}
The results in this section \S1.2 are valid for 
$W$ defined over $\bR$, i.e., without assuming $W$ being defined over $\bQ$.
\medskip

{\bf 1.2.9.} 
We consider Example I--Example V in 1.1.1.
For these examples, $\cL(F) = L^{-1,-1}_\bR(F)\subset \cL$ in 1.2.1 is independent of the choice of $F \in D(\gr^W)$, and we denote it simply by $L$.
By Proposition 1.2.5, we have a real analytic presentation of $D$ 
$$
D \simeq \spl(W) \times D(\gr^W) \times L. \tag1
$$ 
The relation with the complex analytic presentation of $D$ given in 1.1.1 is as follows. 
We use the notation in 1.1.1.
\medskip

{\bf Example I.} 
We have $\spl(W) \simeq \bR$ by assigning $s \in \bR$ to the splitting of $W$ defined by $e_2' \mapsto se_1 + e_2$, $D(\gr^W)$ is one point, and $L \simeq \bR$, $\delta \leftrightarrow d$, by $\delta(e_2') = de_1'$ (see 1.2.3). 

The relation with the complex analytic presentation $D\simeq \bC$ in Example I in 1.1.1 and the real analytic presentation (1) of $D$ is as follows. 
The composition
$$
\bC\simeq D \simeq \spl(W) \times L \simeq  \bR \times \bR
$$
is given by
$$
z\leftrightarrow (s, d),\quad z=s+id.
$$ 
We have conversely 
$$
s= \text{Re}(z), \quad d=\text{Im}(z).
$$
This is because the $\zeta$ associated to $\delta\in L$ is equal to
$\zeta_{-1,-1}=0$ (1.2.3). 
\medskip

{\bf Example II.} 
We have $\spl(W) \simeq \bR^2$, $s \leftrightarrow (s_1,s_2)$, by $s(e_3') = s_1e_1 + s_2e_2 + e_3$, and $s(e_j')=e_j$ $(j = 1,2)$, and have $L = 0$.

The relation with the complex analytic presentation $D\simeq \fh\times \bC$ in Example II in 1.1.1 and the real analytic presentation (1) of $D$ is as follows. 
The composition 
$$
\fh \times \bC\simeq D \simeq \spl(W) \times D(\gr^W)
\simeq \bR^2\times \fh$$
is given by
$$
(\tau, z)\leftrightarrow ((s_1, s_2), \tau),\quad \text{with} \quad
z=s_1-s_2\tau.
$$ 
We  have conversely 
$$
s_1= \text{Re}(z)-\frac{\text{Im}(z)}{\text{Im}(\tau)}\text{Re}(\tau), \quad 
s_2=-\frac{\text{Im}(z)}{\text{Im}(\tau)}.
$$
\medskip

{\bf Example III.} 
We have $\spl(W) \simeq \bR^2$, $s \leftrightarrow (s_1,s_2)$, by $s(e_3') = s_1e_1 + s_2e_2+e_3$, and $s(e_j')=e_j$ $(j = 1,2)$, and have  $L \simeq \bR^2$, $\delta \leftrightarrow (d_1,d_2)$, by $\delta(e_3') = d_1e_1' + d_2e_2'$.

The relation with the complex analytic presentation $D\simeq \fh\times \bC^2$  in Example III in 1.1.1 and the real analytic presentation (1) of $D$ is as follows. 
The composition
$$
\fh \times \bC^2\simeq D \simeq \spl(W) \times D(\gr^W)\times L
\simeq \bR^2\times \fh\times \bR^2
$$
is given by
$$
(\tau, z_1, z_2)\leftrightarrow ((s_1, s_2), \tau, (d_1, d_2)),\tag2
$$
where 
$$
z_1=s_1+\Big(\frac{\Re(\tau)}{2\Im(\tau)}+i\Big)d_1-\frac{\Re(\tau)^2+\Im(\tau)^2}{2\Im(\tau)}d_2, \tag3
$$
$$
z_2=s_2+\frac{1}{2\Im(\tau)}{d_1}+ \Big(-\frac{\Re(\tau)}{2\Im(\tau)}+i\Big)d_2.
$$
We have conversely 
$$
d_1=\text{Im}(z_1), \quad d_2=\text{Im}(z_2),
$$
$$
s_1=\text{Re}(z_1)-\frac{\Re(\tau)}{2\Im(\tau)}\text{Im}(z_1)+\frac{\Re(\tau)^2+\Im(\tau)^2}{2\Im(\tau)}\text{Im}(z_2),\tag4
$$
$$
s_2=\text{Re}(z_2)-\frac{1}{2\Im(\tau)}\text{Im}(z_1)+\frac{\Re(\tau)}{2\Im(\tau)}\text{Im}(z_2).
$$

We explain that the correspondence (2) is described as in (3) and (4). 
Write $\tau=x+iy$ with $x, y\in \bR$, $y>0$.
We have the Hodge decomposition of $\d(e'_3) \in \gr^W_{-3, \bC}$
$$
d_1e'_1+d_2e'_2= A+B, \quad \text{where}\;\; 
A=\frac{d_1-d_2{\bar \tau}}{2yi}(\tau e'_1+e'_2), \;\; 
B= \frac{-d_1+d_2\tau}{2yi}({\bar \tau}e'_1+e'_2), 
$$
with respect to the element $F\in D(\gr^W_{-3})=D(\gr^W)$ corresponding to $\tau\in \fh$.
This shows that the $(p, q)$-Hodge component $\delta_{p,q}$ of $\delta$ is given as follows. $\delta_{p,q}=0$ for $(p, q)\neq (-1, -2), (-2, -1)$, and $\delta_{-1,-2}$ sends $e'_3$ to $A$, and $\delta_{-2,-1}$ sends $e'_3$ to $B$. 
Since $\zeta(F, \delta)=-\frac{i}{2}\delta_{-1,-2}+\frac{i}{2}\delta_{-2,-1}$ (1.2.3), 
this shows that $\zeta(F, \delta)$ sends $e'_3$ to 
$$
v:=\frac{-d_1x+d_2(x^2+y^2)}{2y}e'_1+ \frac{-d_1+d_2x}{2y}e'_2.
$$
Hence $\theta(F, \delta)=\exp(-\zeta(F, \delta))\exp(i\d)F$ is the decreasing  filtration of $\gr^W_\bC$ characterized by the following properties. 
$\theta(F, \delta)^1=0$, $\theta(F, \delta)^{-2}= \gr^W_{\bC}$,
$\theta(F, \delta)^0$ is generated over $\bC$ by $-v+id_1e'_1+id_2e'_2+e'_3$, and $\theta(F, \delta)^{-1}$ is generated over $\bC$ by $\theta(F, \delta)^0$ and $\tau e'_1+e'_2$. 
(3) follows from this, and (4) follows from (3).
\medskip

{\bf Example IV.} 
We have $\spl(W) \simeq \bR^5$, $s \leftrightarrow (s_j)_{1\le j\le 5}$, by $s(e_1')=e_1$, $s(e_2')=s_1e_1+e_2$, $s(e_3')=s_2e_1+e_3$, and $s(e_4')=s_3e_1+s_4e_2+s_5e_3+e_4$, and have  $L \simeq \bR$, $\delta \leftrightarrow d$, by $\delta(e_4') = de_1'$.

The relation with the complex analytic presentation $D\simeq \fh\times \bC^3$  in Example IV in 1.1.1 and the real analytic presentation (1) of $D$ is as follows. 
The composition
$$
\fh \times \bC^3\simeq D \simeq \spl(W) \times D(\gr^W)\times L\simeq \bR^5\times \fh\times \bR
$$
is given by
$$
(\tau, z_1, z_2, z_3)\leftrightarrow ((s_1, \dots, s_5), \tau, d),
$$
where 
$$
z_1=s_1\tau +s_2, \quad z_2=s_3-s_5(s_1\tau +s_2)+id, \quad z_3=s_4-s_5\tau.
$$
We have conversely
$$
s_1=\frac{\text{Im}(z_1)}{\text{Im}(\tau)}, \quad s_2= \text{Re}(z_1)-\frac{\text{Im}(z_1)}{\text{Im}(\tau)}\text{Re}(\tau),
$$
$$
s_3=\text{Re}(z_2)-\frac{\text{Im}(z_3)}{\text{Im}(\tau)}\Re(z_1), \quad 
s_4=\text{Re}(z_3)-\frac{\text{Im}(z_3)}{\text{Im}(\tau)}\text{Re}(\tau),
$$
$$
s_5= -\frac{\text{Im}(z_3)}{\text{Im}(\tau)}, \quad d= \text{Im}(z_2)-\frac{\text{Im}(z_1)\text{Im}(z_3)}{\text{Im}(\tau)}.
$$

This follows from $\zeta=\zeta_{-1,-1}=0$ (1.2.3).
\medskip

{\bf Example V.} 
We have $\spl(W) \simeq \bR^6$, $s \leftrightarrow (s_j)_{1\le j\le 6}$, by $s(e_4')=s_1e_1+s_2e_2+s_3e_3 + e_4$, $s(e_5')=s_4e_1+s_5e_2+s_6e_3 + e_5$, and $s(e_j')=e_j$ $(j = 1,2,3)$, 
and $L = 0$.

The relation with the complex analytic presentation $D\simeq \fh^{\pm} \times \fh \times \bC^3$ in Example V in 1.1.1 and the real analytic presentation (1) of $D$ is as follows. 
The composition
$$
\fh^{\pm}\times \fh \times \bC^3\simeq D \simeq \spl(W) \times D(\gr^W)\simeq \bR^6\times \fh^{\pm} \times \fh
$$
is given by
$$
(\tau_0, \tau_1, z_1, z_2, z_3)\leftrightarrow ((s_1, \dots, s_6), \tau_0, \tau_1),
$$
where 
$$
z_1= s_1\tau_1-s_3\tau_0^2\tau_1+s_4-s_6\tau_0^2,\quad 
z_2= s_2\tau_1-2s_3\tau_0\tau_1+s_5-2s_6\tau_0, \quad 
z_3=s_1-s_2\tau_0+s_3\tau_0^2.
$$
From this we can obtain presentations of $s_j$ ($1\leq j\leq 6$) in terms of
$\tau_0$, $\tau_1$, $z_1$, $z_2$, $z_3$, but we do not write them down here.

\vskip20pt

\head
\S2. The set $D_{\SL(2)}$
\endhead
\medskip

\head
\S2.1. $\SL(2)$-orbits in pure case
\endhead
\medskip

We review $\SL(2)$-orbits in the case of pure weight. 
We also prove some new results here.

Let $w\in \bZ$ and assume $W_w=H_{0,\bR}$ and 
$W_{w-1}=0$.

\medskip

{\bf  2.1.1.} 
Let $n \geq 0$, and consider a pair  $(\rho, \vf)$ consisting of a homomorphism
$$
\rho: \SL(2,\bC)^n\to G_{\bC}
$$
of algebraic groups which is defined over $\bR$ and a 
holomorphic map $\vf: \bP^1(\bC)^n \to \Dc$ satisfying the following condition.
$$
\vf(gz)=\rho(g)\vf(z)\quad \text{for any}\;g\in \SL(2,\bC)^n, 
\;z\in \bP^1(\bC)^n.
$$
\medskip

{\bf 2.1.2.} 
As in \cite{KU3}, \S5 (see also \cite{KU2}, \S3), we call $(\rho, \vf)$ 
as in 2.1.1 an {\it $\SL(2)$-orbit in $n$ variables} if it further satisfies the
following two conditions (1) and (2). 
$$
\vf(\fh^n)\sub D. \tag1
$$
$$
\rho_*(\fil^p_z (\fsl(2,\bC)^{\op n})) \sub \fil^p_{\vf(z)}(\fg_{\bC}) \quad \text{for any $z \in \bP^1(\bC)^n$ and any $p\in \bZ$}.\tag2
$$
Here in (1), $\fh=\{x+iy\;|\;x,y\in \bR,\; y>0\} \subset \bP^1(\bC)$ as in 1.1. 
In (2), $\rho_*$ denotes the Lie algebra homomorphism
$\fsl(2,\bC)^{\op n}\to \fg_{\bC}$ induced by $\rho$,   
$$
\fil^p_z(\fsl(2,\bC)^{\op n})=\{X\in 
\fsl(2,\bC)^{\op n}\;|\;X(\tsize\bigoplus_{j=1}^n \,F^r_{z_j}(\bC^2))
\sub \tsize\bigoplus_{j=1}^n \,F^{r+p}_{z_j}(\bC^2)\;\,
(\forall r\in \bZ)\},
$$
where for $a\in \bP^1(\bC)$, $F^r_a(\bC^2)=\bC^2$ if $r\leq 0$, 
$F^1_a(\bC^2)=\bC\pmatrix a \\ 1 \endpmatrix$ if $a\in \bC$, 
$F^1_\infty(\bC^2)=\bC\pmatrix 1\\ 0\endpmatrix$, $F^r_a(\bC^2)=0$ 
for $r \geq 2$, and
$$
\fil^p_F(\fg_\bC)=\{X\in \fg_\bC\;|\;XF^r\sub F^{r+p}\;\;\text{for all}\;
r\in \bZ\}\quad \text{for}\;\; F \in \Dc.
$$
\medskip

\proclaim{Proposition 2.1.3} 
Let $(\rho, \vf)$ be as in $2.1.1$. 

\medskip

{\rm(i)} The condition $(1)$ in $2.1.2$ is satisfied if there exists $z\in \fh^n$
such that $\vf(z)\in D$. 

\medskip
{\rm(ii)}
The condition $(2)$ in $2.1.2$ is satisfied if there exists $z\in \bP^1(\bC)^n$ such that $\rho_*(\fil^p_z (\fsl(2,\bC)^{\op n})) \sub 
\fil^p_{\vf(z)}(\fg_{\bC})$ for all $p\in \bZ$. 
\endproclaim

{\it Proof.} 
We prove (i). 
Any element $z'$ of $\fh^n$ is written in the form 
$gz$ with $g\in \SL(2,\bR)^n$. 
Hence $\vf(z')=\rho(g)\vf(z) \in D$. 

We prove (ii).
Any element $z'$ of $\bP^1(\bC)^n$ is written in the form $gz$ with 
$g\in \SL(2,\bC)^n$. 
Hence 
$$
\align
\rho_*(\fil^p_{z'} (\fsl(2,\bC)^{\op n}))
&= \rho_*(\Ad(g)\fil^p_z (\fsl(2,\bC)^{\op n}))
=\Ad(\rho(g))\rho_*(\fil^p_z (\fsl(2,\bC)^{\op n}))\\
&\sub \Ad(\rho(g))\fil^p_{\vf(z)} (\fg_{\bC}) =\fil^p_{\vf(z')} (\fg_{\bC}).
\qed
\endalign
$$
\medskip

{\bf 2.1.4.} 
We fix notation.
Assume that we are given $(\rho, \vf)$ as in 2.1.1. 

Let
$$N_j, Y_j, N_j^+\in \fg_\bR\quad (1\leq j\leq n),$$
$$N_j=\rho_*\pmatrix 0&1\\0&0\endpmatrix_j, \quad Y_j=\rho_*\pmatrix -1&0\\0&1\endpmatrix_j,\quad N_j^+=\rho_*\pmatrix 0&0\\1&0\endpmatrix_j,
$$
where $(\;\;)_j$ means the embedding $\fsl(2)\to \fsl(2)^{\oplus n}$ into the $j$-th factor. 

\medskip

\proclaim{Proposition 2.1.5} 
Let $(\rho, \vf)$ be as in $2.1.1$.
Fix $F\in \vf(\bC^n)$. 
Then the condition $(2)$ in $2.1.2$ is satisfied if and only if
$$
N_jF^p\sub F^{p-1}\quad \text{for any $1\leq j \leq n$ and any $p\in \bZ$}.
\tag{$2^\prime$}
$$ 
\endproclaim

{\it Proof.} 
Since $F=\vf((z_j)_j)=\exp(\ts_{j=1}^n z_jN_j)\vf(\bold0)$ for some
$(z_j)_j\in \bC^n$, where $\bold0=0^n\in \bP^1(\bC)^n$, the condition $(2')$ for $F\in \vf(\bC^n)$ is equivalent to the condition $(2')$ for $F=\vf(\bold0)$. 
Note that $\fil^p_{\bold0}(\fsl(2,\bC)^{\op n})=0$ if $p \geq 2$, 
that $\fil^1_{\bold0}(\fsl(2,\bC)^{\op n})$ is generated as a $\bC$-vector 
space by the matrices $\pmatrix0&0\\1&0\endpmatrix_j$ $(1\leq j\leq n)$, that $\fil^0_{\bold0}(\fsl(2,\bC)^{\op n})$ is generated as a $\bC$-vector space by 
$\fil^1_{\bold0}(\fsl(2,\bC)^{\op n})$ and the matrices
$\pmatrix-1&0\\0&1\endpmatrix_j$  
$(1\le j\le n)$, and that $\fil^p_{\bold0}(\fsl(2,\bC)^{\op n})
=\fsl(2,\bC)^{\op n}$ if $p \leq -1$.
Hence, by Proposition 2.1.3 (ii), the condition (2) in 2.1.2 is equivalent to
$$
N_j\vf(\bold0)^p\sub \vf(\bold0)^{p-1}, \quad Y_j\vf(\bold0)^p\sub 
\vf(\bold0)^p, \quad N_j^+\vf(\bold0)^p\sub \vf(\bold0)^{p+1}\quad 
\text{for any $j$, $p$.}
$$
Hence, if the condition (2) in 2.1.2 is satisfied, then $(2')$ is satisfied for $F=\vf(\bold0)$. 

Assume that the condition $(2')$ is satisfied for $F=\vf(\bold0)$. 
We show that the condition (2) in 2.1.2 is satisfied. 
For any diagonal matrices $g_1, \dots, g_n$ in $\SL(2,\bC)$, we have 
$(g_1, \dots, g_n)\bold0=\bold0$ and hence $\rho(g_1, \dots, g_n)
\vf(\bold0)=\vf(\bold0)$. 
From this, we have $Y_j\vf(\bold0)^p\sub \vf(\bold0)^p$ for all $j$ and all $p\in \bZ$. 
It remains to prove $N_j^+\vf(\bold0)^p\sub \vf(\bold0)^{p+1}$ for all $j$ and all 
$p\in \bZ$. 
The following argument is given in \cite{U2}, \S2, in the case $n=1$. 
By the theory of representations of $\fsl(2, \bR)^{\oplus n}$ and by the property 
$Y_j\vf(\bold0)^p\sub \vf(\bold0)^p$ for any $j$ and any $p$, 
we have a direct sum decomposition as an $\bR$-vector space 
$$
H_{0,\bR}= \tsize\bigoplus_{(a,b)\in S} \; P_{a,b}
$$
with $S=\{(a, b)\in \bZ^n\x\bZ^n\;|\; a\geq b \geq -a,\;a(j)\equiv b(j) 
\bmod 2\;\text{for}\;1\leq j \leq n\}$,
having the following properties (1)--(4).
Here, for $a,b\in\bZ^n$, $a\ge b$ means $a(j)\ge b(j)$ for all $1\le j\le n$.
For $1\leq j \leq n$, let $e_j$ be the element of $\bZ^n$ defined by
$e_j(k)=1$ if $k=j$ and $e_j(k)=0$ if $k\neq j$. 
\medskip

(1) On $P_{a,b}$, $Y_j$ acts as the multiplication by $b(j)$.
\medskip

(2) Let $(a, b)\in S$. If $b(j)\neq -a(j)$, $N_j(P_{a,b})\sub P_{a,b-2e_j}$, 
and the map $N_j: P_{a,b} \to P_{a, b-2e_j}$ is an isomorphism. 
If $b(j)=-a(j)$, $N_j$ annihilates $P_{a,b}$.
\medskip

(3) Let $(a, b) \in S$. If $b(j)\neq a(j)$, $N_j^+(P_{a,b})\sub P_{a,b+2e_j}$, 
and for some non-zero rational number $c$, the map $N_j^+: P_{a,b} 
\to P_{a,b+2e_j}$ is $c$ times the inverse of the isomorphism 
$N_j : P_{a,b+2e_j} @>\sim>> P_{a,b}$. If $b(j)=a(j)$, $N_j^+$ 
annihilates $P_{a,b}$.
\medskip

(4) For any $p\in \bZ$, $\varphi(\bold0)^p= \tsize\bigoplus_{(a, b)\in S} \; 
\vf(\bold0)^p\cap P_{a,b,\bC}$. 
For any $(a, b)\in S$, $P_{a,b}$ with the filtration $(\vf(\bold0)^p
\cap P_{a,b,\bC})_{p\in \bZ}$ is an $\bR$-Hodge structure of weight 
$w+ \ts_{j=1}^n b(j)$.
\medskip

For $(a, b) \in S$ such that $b(j)\neq a(j)$, by (2${}'$) with 
$F=\varphi(\bold 0)$ and (4), the bijection 
$N_j: P_{a,b+2e_j} @>\sim>> P_{a,b}$ in the above (2) 
sends the $(p+1,q+1)$-Hodge component of $P_{a,b+2e_j,\bC}$ with $p+q=w+\ts_{j=1}^n b(j)$ bijectively onto the $(p,q)$-Hodge component of $P_{a,b,\bC}$ for the Hodge structure in (4). 
Hence, by (3), $N_j^+$ sends the $(p, q)$-Hodge component of $P_{a,b,\bC}$ with 
$p+q=w+\ts_{j=1}^n b(j)$ onto the $(p+1,q+1)$-Hodge component of 
$P_{a,b+2e_j,\bC}$. 
This proves $N_j^+\vf(\bold0)^p\sub \vf(\bold0)^{p+1}$ for any $p$. 
\qed
\medskip

{\bf 2.1.6.} 
For $1\leq j \leq n$, define the increasing filtration $W^{(j)}$ on 
$H_{0,\bR}$ as follows. Note that
$$
H_{0,\bR}=\tsize\bigoplus_{m\in \bZ^n} V_m, 
$$
where $Y_j$ acts on $V_m$ as the multiplication by $m(j)$. Let 
$$
W^{(j)}_k= \tsize\bigoplus_{m\in \bZ^n,\, m(1)+\dots+m(j) \leq k-w} \;V_m
$$
$$=(\text{the part of $H_{0,\bR}$ on which eigen values of $Y_1+\dots+Y_j$ are $\leq k-w)$}.
$$
Here $w$ is the integer such that $W_w=H_{0,\bR}$ and $W_{w-1}=0$ as at the beginning of this section.

Let $s^{(j)}$ be the splitting of $W^{(j)}$ given by the eigen spaces of $Y_1+\dots+Y_j$. 
That is, $s^{(j)}$ is the unique splitting of $W^{(j)}$ for which the image of $\gr^{W^{(j)}}_k$ under $s^{(j)}$ is the part of $H_{0,\bR}$ on which $Y_1+\dots+Y_j$ acts as the multiplication by $k-w$ for any $k\in \bZ$. 

\medskip

\proclaim{Proposition 2.1.7} 
An $\SL(2)$-orbit in $n$ variables is determined by $((W^{(j)})_{1\leq j\leq n}, \varphi(\bi))$. 

\endproclaim

This is proved in 3.10 of \cite{KU2}.

\medskip

In 2.1.8 and 2.1.10 below, we characterize the splitting $s^{(j)}$ of $W^{(j)}$ given in 2.1.6 in terms of the canonical splittings and the Borel-Serre 
splittings, respectively. 

\proclaim{Proposition 2.1.8} 
Let $(\rho, \varphi)$ be an $\SL(2)$-orbit in $n$ variables, and take $j$ such that $1\leq j\leq n$. 
Let $y_k\in \bR_{\geq 0}$ $(1\leq k\leq n)$, and assume $y_k> 0$ for $j< k\leq n$. 
Then $(W^{(j)}, \varphi(iy_1, \dots, iy_n))$ is a mixed Hodge structure, and $s^{(j)}$ coincides with the canonical splitting $(1.2.3$, cf.\ $1.2.8)$ 
associated to this mixed Hodge structure. 
\endproclaim 

{\it Proof.} 
Let $F=\varphi(iy_1,\dots, iy_n)$, $F'=\varphi(0, \dots, 0, iy_{j+1}, \dots, iy_n)$. Then, $F = \exp(iy_1N_1+\dots+iy_jN_j)F'$, $(W^{(j)}, F)$ is an $\bR$-mixed Hodge structure, $(W^{(j)}, F')$ is an $\bR$-split $\bR$-mixed Hodge structure, and the canonical splitting of $W^{(j)}$ associated to $F'$ is given by $Y_1+\dots+Y_j$.
We have $\delta(F) = y_1N_1+\dots+y_jN_j$. 
Since this $\delta$ has only $(-1, -1)$-Hodge component, $\zeta=0$ by 1.2.3, and hence $Y_1+\dots+Y_j$ is also the canonical splitting of $W^{(j)}$ associated to $F$. 
\qed

\medskip

{\bf 2.1.9.} 
Let $W'$ be an increasing filtration on $H_{0,\bR}$ 
such that 
there exists a group homomorphism $\alpha: \bG_{m,\bR} \to G_\bR$ such that, for $k \in \bZ$, $W'_k = \Op_{m\le k-w} H(m)$, where $H(m):= \{x \in H_{0,\bR} \;|\; \a(t)x = t^m x \;\;(t \in \bR^\times)\}$.

We define the real analytic map 
$$
\spl_{W'}^{\BS} : D \to \spl(W')
$$
as follows. 
Let $P=(G_{W'}^{\circ})_{\bR}$ 
be the parabolic subgroup of $G_{\bR}$ defined by $W'$ 
($G^{\circ}$ is the connected component of $G$ as an algebraic 
group containing $1$). 
Let $P_u$ be the unipotent radical of $P$, and $S_P$ the maximal $\bR$-split torus of the center of $P/P_u$.
Let $\bG_{m,\bR} \to S_P$, $t \mapsto (t^{k-w}$ on $\gr^{W'}_k)_k$ be the weight map induced by $\a$.
For $F \in D$, let $K_F$ be the maximal compact subgroup of $G_\bR$ consisting of the elements of $G_\bR$ which preserve the Hodge metric 
$\lan C_F(\bullet), \bar\bullet\ran_w$, where $C_F$ is the Weil operator associated to $F$.
Let $S_P \to P$ be the Borel-Serre lifting homomorphism at $F$, which 
assigns $a \in S_P$ to the element $a_F \in P$ uniquely determined by the following condition:
the class of $a_F$ in $P/P_u$ coincides with $a$, 
and $\theta_{K_F}(a_F) = a_F^{-1}$, where $\theta_{K_F}$ is the Cartan involution at $K_F$ which coincides with $\ad(C_F)$ in the present situation 
(\cite{KU3} 5.1.3, \cite{KNU1} 8.1). 
Then, the composite $\bG_{m,\bR} \to S_P \to P$ defines an 
action of $\bG_{m,\bR}$ on $H_{0,\bR}$, and we call 
the corresponding splitting of $W'$ 
the {\it Borel-Serre splitting at $F$}, 
and denote it by $\spl_{W'}^{\BS}(F)$.  

It is easy to see that the map $\spl_{W'}^{\BS} : D \to \spl(W'),\; 
F \mapsto \spl_{W'}^{\BS}(F)$, is real analytic. 

\medskip

\proclaim{Proposition 2.1.10} 
Let $(\rho, \varphi)$ be an $\SL(2)$-orbit in $n$ variables, let $y_j>0$ $(1\leq j\leq n)$, and let $p=\varphi(iy_1,\dots, iy_n)\in D$. 
Then 
$$
s^{(j)}=\spl^{\BS}_{W^{(j)}}(p)\quad (1\leq j \leq n).
$$ 
\endproclaim

See Lemma 3.9 of \cite{KU2} for the proof. 
\medskip

{\bf 2.1.11.} 
Let $E$ be a finite dimensional vector space over a field and let $W'$ be an increasing filtration on $E$ such that $W'_w=E$ for $w \gg 0$ and $W'_w=0$ for 
$w \ll 0$.

Recall (\cite{D}, 1.6) that for a nilpotent endomorphism $N$ of $(E, W')$, 
an increasing filtration $M$ on $E$ is called a {\it relative monodromy 
filtration} of $N$ with respect to $W'$ if the following two conditions are satisfied.
\medskip

(1) $N(M_k) \sub M_{k-2}$ for any $k\in \bZ$.  
\medskip

(2) $N^k$ induces an isomorphism $\gr^M_{w+k}\gr^{W'}_w @>\sim>> \gr^M_{w-k}
\gr^{W'}_w$ for any $w\in \bZ$ and any $k \geq 0$. 
\medskip

If a relative monodromy filtration exists, it is unique and is denoted by $M(N,W')$. 
In the case $W'$ is pure, i.e., $W'_w=E$ and $W'_{w-1}=0$ for some $w$, then $M(N, W')$ exists.

Let $(\rho, \vf)$ be as in 2.1.1.
For the family of filtrations in 2.1.6, we see that, for $0 \leq j \leq k \leq n$, 
$W^{(k)}$ is the relative monodromy filtration of $\ts_{j <l\leq k} \;N_l$ 
with respect to $W^{(j)}$ $(W^{(0)}:=W)$.

\medskip

For an increasing filtration $W'$ on $E$ such that $W'_w=E$ for $w \gg 0$ and $W'_w=0$ for $w \ll 0$,  define the 
{\it mean value of the weights $\mu(W')\in \bQ$ of $W'$} and the {\it variance of the weights $\sigma^2(W')\in \bQ$ of $W'$} by
$$
\align 
\mu(W') &= \ts_{w\in \bZ} \;\dim(\gr^{W'}_w)w/\dim(E),\\
\sigma^2(W')&=\ts_{w\in \bZ}\;\dim(\gr^{W'}_w)(w-\mu(W'))^2/\dim(E).
\endalign
$$

\proclaim{Proposition 2.1.12}
Let $N$ be a nilpotent endomorphism of $(E, W')$ as in $2.1.11$. 
Assume that the relative monodromy filtration $M=M(N, W')$ exists.
Then the following holds. 
\medskip

{\rm(i)} $\mu(M)=\mu(W')$.
\medskip

{\rm(ii)} $\sigma^2(M)>\sigma^2(W')$ unless $M=W'$.
\endproclaim

{\it Proof.}
For each $k$, we have the equality 
$$
\dim(\gr^M_k) =  \ts_w \dim(\gr^{W'}_w\gr^M_k)
=\ts_w \; \dim(\gr^M_k\gr^{W'}_w). \tag1
$$
Taking $\ts_k\;(\cdots)k/\dim(E)$ of (1), and using 2.1.11 (2), we obtain (i).
Let $\mu=\mu(M)=\mu(W')$. 
By taking $\ts_k\;(\cdots)(k-\mu)^2/\dim(E)$ of (1), (ii) is reduced to the
inequality $\ts_k \;d_k(k-\mu)^2> (\ts_k\;d_k)(w-\mu)^2$ unless $d_k=0$ for 
any $k\neq w$, where $d_k=\dim (\gr^M_k\gr^{W'}_w)$ for each $w$.
This inequality is obtained again by using 2.1.11 (2).
\qed

\proclaim{Proposition 2.1.13} 
Let $(\rho, \vf)$ be an $\SL(2)$-orbit in $n$ variables, and let $W^{(j)}$ 
$(1\leq j \leq n)$ be as in $2.1.6$. 
Let $W^{(0)}=W$. 
\medskip

{\rm(i)} Let $1\leq j\leq n$. 
Then, $W^{(j-1)}=W^{(j)}$ if and only if the $j$-th component $\SL(2,\bC) 
\to G_\bC$ of $\rho$ is the trivial homomorphism.
\medskip

{\rm(ii)} For $0 \leq j\leq n$, let $\sigma^2(j)=\sigma^2(W^{(j)})$ be as in $2.1.11$ for the increasing filtration $W^{(j)}$ on the $\bR$-vector space $H_{0,\bR}$.
Then, $\sigma^2(j) \leq \sigma^2(j')$ if $0 \leq j \leq j' \leq n$. 
\medskip

{\rm(iii)} Let $0 \leq j \leq n$, $0 \leq j' \leq n$. 
Then, $W^{(j)}=W^{(j')}$ if and only if $\sigma^2(j)=\sigma^2(j')$.
\endproclaim

The statement (i) was proved in [KU2] \S3. 
The statements (ii) and (iii) follow from Proposition 2.1.12.
\medskip

{\bf 2.1.14.}
Let $(\rho, \vf)$ be an $\SL(2)$-orbit in $n$ variables in pure case. 
Put $W^{(0)}=W$.
We define {\it rank of $(\rho, \vf)$} as the number of the elements of the set $\{j\;|\;1 \leq j \leq n, W^{(j)} \neq W^{(j-1)}\}$.

\medskip

{\bf 2.1.15. Example 0.} 
  Recall that in this case, $D$ is identified with 
the upper-half plane $\fh$. 
  Let $\rho$ be the standard isomorphism $\SL(2,\bC)\to G_{\bC}$, and let 
$\varphi: \bP^1(\bC)\to \Dc$ be the natural isomorphism in 1.1.1. 
Then $(\rho, \vf)$ is an $\SL(2)$-orbit in one variable of rank $1$.
\medskip

\vskip20pt

\head 
\S2.2. Nilpotent orbits and $\SL(2)$-orbits in pure case
\endhead

We consider the pure case. Let  $w\in \bZ$ and assume $W_w=H_{0,\bR}$ and $W_{w-1}=0$.

\medskip

{\bf 2.2.1.} 
Let $F\in \Dc$ and let $N_1, \cdots, N_n$ be elements of $\fg_\bR$ such that $N_jN_k=N_kN_j$ for any $j, k$ and such that $N_j$ is nilpotent as a linear map $H_{0,\bR}\to H_{0,\bR}$ for any $j$. 

We say that the map 
$$
\bC^n\to \Dc, \quad (z_1,\dots, z_n)\mapsto \exp(\ts_{j=1}^n z_jN_j)F
$$
is a {\it nilpotent orbit} if the following conditions (1) and (2) are satisfied.

\medskip

(1) $\exp(\ts_{j=1}^n z_jN_j)F\in D$ if $\text{Im}(z_j)\gg 0$ for all $j$.

\medskip

(2) $N_jF^p\sub F^{p-1}$ for any $j$ and any $p$.

\medskip

In this case, we say also that $(N_1, \dots, N_n, F)$ {\it generates a nilpotent orbit.}

\medskip

{\bf 2.2.2.} 
Assume that 
$(N_1, \dots, N_n, F)$ generates a nilpotent orbit. 
By [CK], for $y_j\in \bR_{\geq 0}$, the filtration $M(y_1N_1+\dots+y_nN_n, W)$ (2.1.11)  depends only on the set $\{j\;|\;y_j\ne 0\}$. 
For $1\leq j\leq n$, let $W^{(j)}=M(N_1+\dots+N_j, W)$.  

\medskip

{\bf 2.2.3.} 
Assume that 
$(N_1, \dots, N_n, F)$ generates a nilpotent orbit. 
Then by Cattani, Kaplan and Schmid \cite{CKS}, an $\SL(2)$-orbit $(\rho, \varphi)$ is  canonically associated to $(N_1,\dots, N_n, F)$. 
(The 
homomorphism $\rho$ is given in \cite{CKS} Theorem 4.20 and $\vf$ is defined
by  $\vf(g\bold0)= \rho(g)\hat F \quad (g\in \SL(2,\bC)^n)$, 
where $\bold0=0^n \in \bP^1(\bC)^n$.)
By \cite{KNU1}, this $\SL(2)$-orbit is characterized in the style of the following theorem. 
\medskip

\proclaim{Theorem 2.2.4} 
Assume that 
$(N_1,\dots, N_n, F)$ generates a nilpotent orbit.

\medskip

{\rm (i)} {\rm(\cite{KNU1}, 8.7)} 
Let $1\leq j\leq n$. 
Then, when $y_k\in \bR_{>0}$ and $y_k/y_{k+1}\to \infty$ $(1\leq k\leq n$,  $y_{n+1}$ means $1)$, the Borel-Serre splitting $\spl_{W^{(j)}}^{\BS}(\exp(\ts_{k=1}^n iy_kN_k)F)$ converges in $\spl(W^{(j)})$. 

Let $s^{(j)}\in \spl(W^{(j)})$ be the limit. 

\medskip

{\rm(ii)} There is a homomorphism $\tau: \bG_{m, \bR}^n \to \Aut_{\bR}(H_{0, \bR})$ of algebraic groups over $\bR$ characterized by the following property. 
For any $1\leq  j\leq n$ and any $k\in \bZ$, we have
$$
s^{(j)}(\gr^{W^{(j)}}_k)=\{v\in H_{0, \bR}\;|\; \tau_j(t)v=t^kv\; \text{for any $t\in \bR^\times$}\}, 
$$
where $\tau_j:\bG_{m,\bR}\to \Aut_{\bR}(H_{0, \bR})$ is the $j$-th component of $\tau$.
\medskip

{\rm(iii)} There exists a unique $\SL(2)$-orbit $(\rho, \vf)$ in $n$ variables
characterized by the following properties $(1)$ and $(2)$. 
\medskip

$(1)$ The associated weight filtrations $W^{(j)}$ $(1\leq j \leq n)$ are the same as $W^{(j)}$ in $2.2.2$.

\medskip

$(2)$ $\varphi(\bi)$ is the limit in $D$ of
$$
\tau\Big(\sqrt{\frac{y_2}{y_1}},\dots, \sqrt{\frac{y_{n+1}}{y_n}}\Big)^{-1}\exp(\ts_{j=1}^n iy_jN_j)F\quad (y_j >0, \;\; y_j/y_{j+1}\to \infty \;\;(1\leq j\leq n))
$$
$(y_{n+1}$ means $1)$,  
where $\tau$ is as in {\rm(ii)}.
\medskip

{\rm(iv)} The associated torus action $\tilde \rho$ {\rm (\cite{KU2} 3.1 (4))}
of $(\rho, \vf)$ and the homomorphism $\tau$ in {\rm(ii)} are related as 
follows{\rom:} $\tau(t_1,\ldots,t_n)=(\prod_{j=1}^nt_j)^w\tilde \rho(t_1,\ldots,t_n)$. 
\endproclaim

\medskip

{\bf 2.2.5. Example 0.} 
Let $(N, F)$ be as follows. 
$N(e_2)=e_1$, $N(e_1)=0$, $F=F(z)$ with $z \in i\cdot \bR$ in the notation 
of 1.1.1. 
Then $(N, F)$ generates a nilpotent orbit, and the associated $\SL(2)$-orbit is the one in 2.1.15.

In fact, $\exp(iyN)F=F(z+iy)$, and $\tau(t)$ in 2.2.4 (ii) sends $e_1$ to $t^{-2}e_1$ and $e_2$ to $e_2$. Hence $\tau(1/\sqrt{y})^{-1}\exp(iyN)F = F((z+iy)/y) \to F(i)$ as $y \to \infty$.

\medskip

{\bf 2.2.6.} 
Assume that $(N, F)$ generates a nilpotent orbit in the pure case in 2.2.1 for $n= 1$. 
Let $W^{(1)}=M(N, W)$ be as in 2.2.2.
Then $(W^{(1)}, F)$ is a mixed Hodge structure, and the splitting $s^{(1)}$ of $W^{(1)}$ given by the $\SL(2)$-orbit (2.1.6) associated to $(N, F)$ coincides with the canonical splitting of $W^{(1)}$ associated to $F$ (1.2).

\medskip

{\bf 2.2.7.} 
More generally, for any mixed Hodge structure, its canonical splitting (1.2) is obtained as in 2.2.6 by embedding the mixed Hodge structure into a pure 
nilpotent orbit.

In fact, let $(M, F)$ be 
a mixed Hodge structure on an $\bR$-vector space $V$.
Let $k$ be an integer such that all the weights of $(M, F)$ are not greater than $k$.
It is shown in \cite{KNU1} that there exist an $\bR$-vector space $V'$, 
an $\bR$-linear injective map $q: V \to V'$, a nilpotent endomorphism $N$ of $V'$, and a decreasing filtration $F'$ on $V'_\bC$ such that the pair $(N, F')$ generates a nilpotent orbit on $V'$ in the pure case of weight $k$ in 2.2.1 for $n= 1$, which satisfy the following conditions.
\medskip

Let $W'$ be the trivial weight filtration on $V'$ of weight $k$, and let $W^{(1)}=M(N, W')$ be as in 2.2.2.
Then, $0@>>>(V, M, F)@>q>>(V', W^{(1)}, F')@>N>>(V', W^{(1)}[-2], F'(-1))$ is an exact sequence of mixed Hodge structures, where $[-2]$ is the shift by $-2$ and $(-1)$ is the Tate twist by $-1$, and the restriction of the splitting $s^{(1)}$ of $W^{(1)}$, given by the $\SL(2)$-orbit  associated to $(N, F')$ on $V'$, to $\Ker(N:\gr^{W^{(1)}}\to \gr^{W^{(1)}[-2]}) \simeq \gr^M$ coincides with the canonical splitting of $M$ on $V$ associated to $F$.
\medskip

For the proof, see \cite{KNU1}, \S3.

\vskip20pt

\head
\S2.3. $\SL(2)$-orbits in mixed case
\endhead
\medskip

Now we consider the mixed version of \S2.1. 
Let $W$ be as in Notation. 

\medskip

{\bf 2.3.1.} 
For $n\geq 0$, let $\cD_{\SL(2),n}'$ be the set of pairs
$((\rho_w, \vf_w)_{w\in \bZ}, \br)$, where $(\rho_w, \vf_w)$ is an 
$\SL(2)$-orbit in $n$ variables for $\gr^W_w$ for each $w \in \bZ$, and 
$\br$ is an element of $D$ such that $\br(\gr^W_w)= \vf_w(\bi)$ for each 
$w\in \bZ$. 
Here $\bi=(i, \dots, i)\in \bC^n\sub \bP^1(\bC)^n$.
\medskip

{\bf 2.3.2.} 
Let $\cD_{\SL(2), n}$ be the set of all triples $((\rho_w, \vf_w)_{w\in \bZ}, \br, J)$, where $((\rho_w, \vf_w)_{w\in \bZ}, \br) \in \cD_{\SL(2),n}'$ and $J$ is a subset of $\{1, \dots, n\}$ satisfying the following conditions (1) and (2). 
Let 
$$
J'=\{j\;|\;1\leq j \leq n,\; \text{there is $w\in \bZ$ such that the $j$-th component}
$$
$$
\hskip60pt
\text{$\SL(2)\to G_\bR(\gr^W_w)$ of $\rho_w$ is a non-trivial homomorphism}\}.
$$ 
\medskip

(1) If $\br\in D_{\spl}$, $J=J'$.
\medskip

(2) If $\br\in D_{\nspl}$, either $J=J'$ or $J=J' \cup \{k\}$ for some 
$k < \min J'$. 
\medskip

Let 
$$
\cD_{\SL(2)}= \tsCu_{n\ge0}\;\cD_{\SL(2),n}.
$$
\medskip

We call an element of $\cD_{\SL(2), n}$ an {\it $\SL(2)$-orbit in $n$ variables}, and an element of $\cD_{\SL(2)}$ an {\it $\SL(2)$-orbit}. 
  Note that, in the pure case, $J$ is determined uniquely by $(\rho_w)_w$ 
since $D=D_{\spl}$. 

We call the cardinality of the set $J$ the {\it rank} of the $\SL(2)$-orbit.
\medskip

{\bf 2.3.3.}  
Let $\cD_{\SL(2),n, \sharp}\sub \cD_{\SL(2), n}$ be the set of all $\SL(2)$-orbits in $n$ variables of rank $n$.

For an element $((\rho_w, \vf_w)_w, \br, J)$ of $\cD_{\SL(2),n, \sharp}$, $J=\{1,\dots, n\}$. 
Hence, by forgetting $J$, the set $\cD_{\SL(2),n, \sharp}$ is identified with the subset of $\cD_{\SL(2),n}'$ (2.3.1) consisting of all elements $((\rho_w, \vf_w)_w, \br)$ satisfying the following conditions (1) and (2).
\medskip

(1) If $2\leq j \leq n$, there exists $w\in \bZ$ such that the $j$-th component of $\rho_w$ is a non-trivial homomorphism.
\medskip

(2) If $\br \in D_{\spl}$ and $n \geq 1$, there exists $w\in \bZ$ such that the 1-st component of $\rho_w$ is a non-trivial homomorphism.
\medskip

As is seen later in \S2.5, for the construction of the space $D_{\SL(2)}$, it is sufficient to consider $\SL(2)$-orbits in $n$ variables of rank $r$ with $r=n$.
We call this type of $\SL(2)$-orbit a {\it non-degenerate $\SL(2)$-orbit of rank $n$}, or for simplicity, an {\it $\SL(2)$-orbit of rank $n$}, and we regard it as an element of $\cD_{\SL(2),n}'$. 

  On the other hand, the generality of the definition in 2.3.2 with the auxiliary data $J$ is natural in 2.4 when we consider the relations with nilpotent orbits.

\medskip

{\bf 2.3.4.} 
If $((\rho_w, \vf_w)_w, \br, J)$ is an $\SL(2)$-orbit in $n$ variables of rank $r$, we have the {\it associated $\SL(2)$-orbit $((\rho'_w, \vf'_w)_w, \br)$ in $r$ variables of rank $r$}, defined as follows.
Let $J=\{a(1), \dots, a(r)\}$ with $a(1)<\dots < a(r)$. 
Then 
$$
\rho'_w(g_{a(1)}, \dots, g_{a(r)}):=\rho_w(g_1, \dots, g_n), 
\quad 
\vf'_w(z_{a(1)}, \dots, z_{a(r)}):=\vf_w(z_1, \dots, z_n).
$$

Note that, for any $w\in \bZ$, $\rho_w$ factors through the projection $\SL(2)^n\to \SL(2)^J$ to the $J$-component, and $\vf_w$ factors through the projection $\bP^1(\bC)^n\to \bP^1(\bC)^J$ to the $J$-component, and hence $(\rho_w, \vf_w)_w$ is essentially the same as $(\rho'_w, \vf'_w)_w$. 
\medskip

{\bf 2.3.5.} 
{\it Associated torus action.}

Assume that we are given an $\SL(2)$-orbit in $n$ variables
$((\rho_w, \vf_w)_w, \br, J)$.

We define the associated homomorphism of algebraic groups over $\bR$
$$
\tau : \bG_{m,\bR}^n \to \Aut_{\bR}(H_{0, \bR}, W)
$$
as follows. 
Let $s_\br: \gr^W @>\sim>> H_{0,\bR}$ be the canonical splitting of $W$ associated to $\br$ (\S1.2). 
Then  
$$
\tau(t_1,\dots, t_n)= s_\br\circ\Big(\tsize\bigoplus_{w\in \bZ}\Big(\tsize\prod_{j=1}^n t_j\Big)^w\rho_w(g_1, \dots, g_n)\;\text{on}\;\gr^W_w\Big)\circ s_\br^{-1}
$$ 
$$
\text{with}\quad g_j=\pmatrix 1/ \tsize\prod_{k=j}^n t_k& 0 \\
0 & \tsize\prod_{k=j}^n t_k\endpmatrix.
$$
For $1\leq j\leq n$, let $\tau_j: \bG_{m,\bR} \to \Aut_{\bR}(H_{0,\bR}, W)$ be the $j$-th component of $\tau$. 
\medskip

{\it Remark.}
The induced action of $\tau(t)$ ($t\in \bR^n_{>0}$) on $D$ is described as follows.
For $s(\theta(F,\d))\in D$ with $s \in \spl(W)$, $F \in D(\gr^W)$, 
$\d\in\cL(F)$ (1.2.5), we have
$$
\tau(t)s(\theta(F, \delta))=s'(\theta(F', \delta'))
$$ 
$$
\text{with}\;\;\; s'=\tau(t)s\gr^W(\tau(t))^{-1}, \;\;
F'=\gr^W(\tau(t))F,\;\; 
\delta'= \Ad(\gr^W(\tau(t)))\delta.
$$
\medskip

{\bf 2.3.6.} 
{\it Associated family of weight filtrations}. 
 
In the situation of 2.3.5, for $1\leq j \leq n$, we define the associated $j$-th weight filtration $W^{(j)}$ on $H_{0,\bR}$ as follows. 
For $k \in \bZ$, $W^{(j)}_k$ is the direct sum of $\{x \in H_{0,\bR}\;|\;\tau_j(t)x=t^\ell x\;(\forall\;t\in\bR^\times)\}$ over all $\ell \leq k$.

By definition, we have $W^{(j)}_k = \ts_{w\in \bZ} s_\br(W^{(j)}_k(\gr^W_w))$, and $W^{(j)}_k(\gr^W_w)$ coincides with the $k$-th filter of the $j$-th weight filtration on $\gr^W_w$ associated to the $\SL(2)$-orbit $(\rho_w, \vf_w)$ in $n$ variables.
\medskip

\proclaim{Proposition 2.3.7} 
{\rm(i)} An $\SL(2)$-orbit in $n$ variables $((\rho_w, \vf_w)_w, \br, J)$ is uniquely determined by $((W^{(j)}(\gr^W))_{1\leq j\leq n}, \br, J)$. 
\medskip

{\rm(ii)} An $\SL(2)$-orbit in $n$ variables $((\rho_w, \vf_w)_w, \br, J)$ is uniquely determined by $(\tau,\br, J)$. 
\endproclaim

{\it Proof.} 
(i) In the pure case, this is Proposition 2.1.7.  
The general case is clearly reduced to the pure case. 

(ii) follows from (i), since the family of weight filtrations $(W^{(j)}(\gr^W))_{1\leq j\leq n}$ is determined by $\tau$.
\qed
\medskip

\proclaim{Proposition 2.3.8} 
Let $((\rho_w, \vf_w)_w, \br, J)$ be an $\SL(2)$-orbit in $n$ variables, and 
let $W^{(j)}$ $(1\leq j \leq n)$ be as in $2.3.6$. 
Let $W^{(0)}=W$. 
\medskip

{\rm(i)} Let $1\leq j\leq n$. 
Then $W^{(j)} = W^{(j-1)}$ if and only if for any $w\in \bZ$, the $j$-th factor
$\SL(2,\bC) \to G_{\bC}(\gr^W_w)$ of $\rho_w$ is the trivial homomorphism. 
\medskip

{\rm(ii)} For $0 \leq j\leq n$, let
$\sigma^2(j) = \ts_{w\in \bZ} \;\sigma^2(W^{(j)}(\gr^W_w))$, 
where $\sigma^2(W^{(j)}(\gr^W_w))$ is the variance $(2.1.11)$
of the increasing filtration $W^{(j)}(\gr^W_w)$ on the $\bR$-vector space 
$\gr^W_w$. 
Then, $\sigma^2(j) \leq \sigma^2(j')$ if $0 \leq j \leq j' \leq n$. 
\medskip

{\rm(iii)} Let $0 \leq j \leq n$, $0 \leq j' \leq n$. 
Then, $W^{(j)}=W^{(j')}$ if and only if $\sigma^2(j)=\sigma^2(j')$.
\endproclaim

{\it Proof.} 
This is also reduced to the pure case Proposition 2.1.13.
\qed
\medskip

{\bf 2.3.9.}
We describe what kind of $\SL(2)$-orbits of positive rank exist in 
Example I--Example V. 
We consider only an $\SL(2)$-orbit in $r$ variables in rank $r$, hence $J = \{1, \dots, r\}$ in the following (2.3.3).
\medskip

{\bf Example I.} 
Any $\SL(2)$-orbit of rank $>0$ is of rank $1$. 
An $\SL(2)$-orbit in one variable of rank $1$ is $((\rho_w, \varphi_w)_w, 
\br)$, 
where $\rho_w$ is the trivial homomorphism from $\SL(2)$ to $G_\bR(\gr^W_w)$ and $\varphi_w$ is the unique map from $\bP^1(\bC)$ onto the one point set $D(\gr^W_w)$, and $\br$ is any element of $D_{\nspl}=\bC\smallsetminus \bR$. 
We have $W^{(1)}=W$. 
Later we refer to the case $\br =F(i) \in D$ (i.e., $\br=i\in \bC=D$) as Example I in 2.3.9. 
\medskip

{\bf Example II.} 
Any $\SL(2)$-orbit of rank $>0$ is of rank $1$. 
An $\SL(2)$-orbit in one variable of rank $1$ is  $((\rho_w, \varphi_w)_w, \br)$, 
where $(\rho_w, \varphi_w)$ is of rank $0$ for $w\neq -1$, and $(\rho_{-1},\varphi_{-1})$ is of rank $1$. 
An example of such $\SL(2)$-orbit is that $(\rho_{-1}, \varphi_{-1})$ is the SL(2)-orbit in 2.1.15, and $\br=F(i, z)$ in the notation of 1.1.1, Example II. 
For this $\SL(2)$-orbit, $W^{(1)}$ is given by 
$$
W^{(1)}_{-3}=0\sub W^{(1)}_{-2}=W^{(1)}_{-1}=\bR e_1 \sub W^{(1)}_0=H_{0,\bR}.
$$ 
\medskip

{\bf Example III.} 
There are three cases for $\SL(2)$-orbits in $r$ variables of rank $r>0$. For any of them, $(\rho_w, \varphi_w)$ is of rank $0$ unless $w=-3$. 
\medskip

{\it Case} 1. 
$r=1$ and $(\rho_{-3}, \varphi_{-3})$ is of rank $1$. 
An example of such $\SL(2)$-orbit is given as follows. $(\rho_{-3}, \varphi_{-3})$ is $(\rho, \varphi(1))$ of 2.1.15 (we identify $\Dc(\gr^W_{-3})$ with $\bP^1(\bC)$ via the Tate twist), 
and $\br=F(i, z_1, i)$ for $z_1 \in \bC$ (1.1.1). 
For this $\SL(2)$-orbit, 
$$
W^{(1)}_{-5}=0\sub W^{(1)}_{-4}=W^{(1)}_{-3}=\bR e_1 \sub W^{(1)}_{-2}=W^{(1)}_{-1}=W^{(1)}_{-3}+\bR e_2\sub W^{(1)}_0=H_{0,\bR}.
$$
\medskip

{\it Case} 2. 
$r=1$ and $(\rho_{-3}, \varphi_{-3})$ is of rank $0$. 
An example of such $\SL(2)$-orbit is given as follows. 
$\rho_{-3}$ is the trivial homomorphism onto $\{1\}$, and $\varphi_{-3}$ is the constant map with value $i\in \fh=D(\gr^W_{-3})$, and $\br=F(i, z_1, z_2)$ with $(z_1, z_2)\in \bC^2\smallsetminus \bR^2$. 
For this $\SL(2)$-orbit, $W^{(1)}=W$. 

\medskip

{\it Case} 3. 
$r=2$ and $(\rho_{-3}, \varphi_{-3})$ is of rank $1$. 
$\rho_{-3}: \SL(2, \bC)^2\to G_{\bC}(\gr^W_{-3})=\SL(2,\bC)$ factors through the second projection onto $\SL(2,\bC)$, and $\varphi_{-3}:\bP^1(\bC)^2\to \Dc(\gr^W_{-3})=\bP^1(\bC)$ factors through the second projection onto $\bP^1(\bC)$. 
An example of such $\SL(2)$-orbit is given as follows. $\rho_{-3}(g_1, g_2)=g_2$, $\varphi_{-3}(p_1, p_2)=p_2$, and $\br=F(i, z_1, z_2)$, 
where $(z_1, z_2)\in \bC^2\smallsetminus \bR^2$. 
For this $\SL(2)$-orbit, $W^{(1)}=W$ and  $W^{(2)}$ is the $W^{(1)}$ in the example in Case 1. 
\medskip

{\bf Example IV.} 
There are three cases for $\SL(2)$-orbits in $r$ variables of rank $r>0$. 
For any of them, $(\rho_w, \varphi_w)$ is of rank $0$ unless $w=-1$. 
\medskip

{\it Case} 1. 
$r=1$ and $(\rho_{-1}, \varphi_{-1})$ is of rank $1$. 
An example of such $\SL(2)$-orbit is given as follows. 
$(\rho_{-1}, \varphi_{-1})$ is the standard one (which is identified with 
$(\rho, \vf)$ in 2.1.15 by the identification of $e'_2, e'_3$ here with 
$e_1, e_2$ there), 
and $\br=F(i, z_1, z_2, z_3)$ for $z_1, z_2, z_3\in \bC$ (1.1.1). 
For this $\SL(2)$-orbit, 
$$
W^{(1)}_{-3}=0\sub W^{(1)}_{-2}=W^{(1)}_{-1}=\bR e_1+\bR e_2 \sub W^{(1)}_0=H_{0,\bR}.
$$
\medskip

{\it Case} 2. 
$r=1$ and $(\rho_{-1}, \varphi_{-1})$ is of rank $0$. 
An example of such $\SL(2)$-orbit is given as follows. 
$\rho_{-1}$ is the trivial homomorphism onto $\{1\}$, and $\varphi_{-1}$ is the constant map with value $i\in \fh=D(\gr^W_{-1})$, and $\br=F(i, z_1, z_2, z_3)$ with $\text{Im}(z_2)\neq \text{Im}(z_1)\text{Im}(z_3)$ (the last condition says $F(i, z_1,z_2,z_3)\in D_{\nspl}$). 
For this $\SL(2)$-orbit, $W^{(1)}=W$. 
\medskip

{\it Case} 3. 
$r=2$ and $(\rho_{-1}, \varphi_{-1})$ is of rank $1$. 
$\rho_{-1}: \SL(2, \bC)^2\to G_{\bC}(\gr^W_{-1})$ factors through the second projection onto $\SL(2,\bC)$, and $\varphi_{-1}:\bP^1(\bC)^2\to \Dc(\gr^W_{-1})=\bP^1(\bC)$ factors through the second projection onto $\bP^1(\bC)$. 
An example of such $\SL(2)$-orbit is given as follows. 
$\rho_{-1}(g_1, g_2)=g_2$, $\varphi_{-1}(p_1, p_2)=p_2$, and $\br=F(i, z_1, z_2, z_3)$ with $\text{Im}(z_2)\neq \text{Im}(z_1)\text{Im}(z_3)$.
For this $\SL(2)$-orbit, $W^{(1)}=W$ and  $W^{(2)}$ is the $W^{(1)}$ in the example in Case 1. 
\medskip

{\bf Example V.} 
There are five cases for $\SL(2)$-orbits in $r$ variables of rank $r>0$. 
For any of them, $(\rho_w, \varphi_w)$ is of rank $0$ if $w\notin\{0, 1\}$. 
\medskip

{\it Case} 1 (resp. {\it Case} 2). 
$r=1$ and $(\rho_0, \varphi_0)$ is of rank $1$ (resp. $0$), and $(\rho_1, \varphi_1)$ is of rank $0$ (resp. $1$). 
An example of such $\SL(2)$-orbit is given as follows. 
$(\rho_0, \varphi_0)$ (resp. $(\rho_1, \varphi_1))$ is $(\Sym^2(\rho), \Sym^2(\varphi)(-1))$ (resp. $(\rho, \varphi(-1)$)) for the standard 
$(\rho, \vf)$ in 2.1.15 via a suitable identification,  where $(-1)$ means the Tate twist, and $\br=F(i, i,  z_1, z_2, z_3)$ for $z_1, z_2, z_3\in \bC$. 
For this $\SL(2)$-orbit, 
$$
\align
W^{(1)}_{-3} = 0 \sub W^{(1)}_{-2} = W^{(1)}_{-1} = \bR e_1 &\sub W^{(1)}_0 = W^{(1)}_{-1} + \bR e_2 \\
&\sub W^{(1)}_1 = W^{(1)}_0 + \bR e_4 + \bR e_5 \sub W^{(1)}_2 = H_{0,\bR}
\endalign
$$
$$
(\text{resp}.\;\; W^{(1)}_{-1} = 0 \sub W^{(1)}_0 =  W^{(1)}_1 = \bR e_1 + \bR e_2 + \bR e_3 + \bR e_4 \sub W^{(1)}_2 = H_{0,\bR}).
$$
\medskip

{\it Case} 3. 
$r=1$, and both $(\rho_0, \varphi_0)$ and $(\rho_1, \varphi_1)$ are of rank $1$. 
An example of such $\SL(2)$-orbit is given as follows. 
$\rho_0=\Sym^2(\rho)$, $\varphi_0=\Sym^2(\varphi)(-1)$, $\rho_1=\rho$, $\vf_1=\vf(-1)$ for the standard $(\rho, \vf)$ in 2.1.15 via a suitable identification, 
and $\br=F(i, i,  z_1, z_2, z_3)$ for $z_1, z_2, z_3\in \bC$. 
For this $\SL(2)$-orbit, 
$$
W^{(1)}_{-3} = 0 \sub W^{(1)}_{-2} = W^{(1)}_{-1} = \bR e_1 \sub W^{(1)}_0 = W^{(1)}_1 =W^{(1)}_{-1} + \bR e_2 + \bR e_4 \sub W^{(1)}_2 = H_{0,\bR}.
$$
\medskip

{\it Case} 4 (resp. {\it Case} 5). 
$r=2$, both $(\rho_0, \varphi_0)$ and $(\rho_1, \varphi_1)$ are of rank $1$, 
$\rho_0: \SL(2, \bC)^2\to G_{\bC}(\gr^W_0)$ factors through the first (resp. second) projection onto $\SL(2,\bC)$, $\varphi_0:\bP^1(\bC)^2\to \Dc(\gr^W_0)$ factors through the first (resp. second) projection onto $\bP^1(\bC)$, $\rho_1: \SL(2, \bC)^2\to G_{\bC}(\gr^W_1)$ factors through the second (resp. first) projection onto $\SL(2,\bC)$, and $\varphi_1:\bP^1(\bC)^2\to \Dc(\gr^W_1)$ factors through the second (resp. first) projection onto $\bP^1(\bC)$. 
An example of such $\SL(2)$-orbit is given as follows. 
For $j=1$ (resp. $2$), $\rho_0(g_1, g_2)=\Sym^2(g_j)$, $\varphi_0(p_1, p_2) = p_j \in \bP^1(\bC) \simeq \Dc(\gr^W_0)$ (cf.\ 1.1.2), $\rho_1(g_1, g_2)=g_{3-j}$, $\varphi_1(p_1, p_2)=p_{3-j}(-1)\in \bP^1(\bC) \simeq \Dc(\gr^W_1)$, and $\br=F(i, i, z_1, z_2, z_3)$ with $z_1, z_2, z_3\in \bC$.
For this $\SL(2)$-orbit, $W^{(1)}$ is the $W^{(1)}$ in the example in Case 1 (resp. Case 2) and  $W^{(2)}$ is the $W^{(1)}$ in the example in Case 3. 
\medskip

\vskip20pt

\head
\S2.4. Nilpotent orbits and $\SL(2)$-orbits in mixed case
\endhead

{\bf 2.4.1.} 
Let $N_j\in \fg_\bR$ ($1\leq j \leq n$) and let $F\in \Dc$. We say $(N_1, \dots, N_n, F)$ {\it generates a nilpotent orbit} if the following conditions (1)--(4) are satisfied.

\medskip

(1) The $\bR$-linear maps $N_j: H_{0,\bR} \to H_{0,\bR}$ are nilpotent for all $j$, and $N_jN_k=N_kN_j$ for all $j, k$.
\medskip

(2) If $y_j \gg 0$ ($1\leq j \leq n$), then $\exp(\ts_{j=1}^n iy_jN_j)F\in D$.
\medskip

(3)  $N_jF^p\subset F^{p-1}$ for all $j$ and $p$ (Griffiths transversality).
\medskip

(4) Let $J$ be any subset of $\{1, \dots, n\}$.
Then for $y_j\in \bR_{>0}$ ($j\in J$), the relative monodromy filtration $M(\ts_{j\in J} y_jN_j, W)$ (2.1.11) exists.
Furthermore, this filtration is independent of the choice of $y_j\in \bR_{>0}$.
\medskip

In the pure case, by 2.2.2, $(N_1, \dots, N_n, F)$ generates a nilpotent orbit in this sense if and only if it does in the sense of 2.2.1. 

Let $\cD_{\nilp,n}$ be the set of all $(N_1, \dots, N_n, F)$ which generate 
nilpotent orbits. 

For $(N_1, \cdots, N_n, F) \in \cD_{\nilp,n}$, we call the map $(z_1,\dots, z_n) \mapsto \exp(\ts_{j=1}^n z_jN_j)F$ a {\it nilpotent orbit in $n$ variables}.

In the terminology of Kashiwara \cite{K}, $\cD_{\nilp, n}$ is the set of all
$(N_1, \dots, N_n, F)$ such that $(H_{0, \bC}; W_{\bC}; F, \bar F; N_1, \dots, N_n)$, with $\bar F$ the complex conjugate of $F$, is an \lq\lq infinitesimal mixed Hodge module''. 
\medskip

We will prove the following results  2.4.2, 2.4.3, and 2.4.5. (i) of 2.4.2  was already proved in Theorem 0.5 of our previous paper \cite{KNU1}.
\medskip

\proclaim{Theorem 2.4.2} 
Let $(N_1, \dots, N_n, F)\in \cD_{\nilp, n}$.
For each $w\in \bZ$, let $(\rho_w, \varphi_w)$ be the  $\SL(2)$-orbit in $n$ variables for $\gr^W_w$ associated to $(\gr^W_w(N_1), \dots, \gr^W_w(N_n), F(\gr^W_w))$, which generates a nilpotent orbit for $\gr^W_w$ $(2.2.3)$. 
Let $k=\min(\{j\;|\;1\leq j \leq n,\; N_j\not=0\} \cup \{n+1\})$. 
\medskip

{\rm(i)} If $y_j\in \bR_{>0}$ and $y_j/y_{j+1}\to \infty$ $(1\leq j\leq n$, $y_{n+1}$ means $1)$, the canonical splitting $\spl_W(\exp(\ts_{j=1}^n iy_jN_j)F)$ of $W$ associated to $\exp(\ts_{j=1}^n iy_jN_j)F$ $(1.2.3)$ converges in $\spl(W)$. 

Let $s\in \spl(W)$ be the limit. 
\medskip

{\rm(ii)} Let $\tau: \bG_{m,\bR}^n \to \Aut_\bR(H_{0,\bR}, W)$ be the homomorphism of algebraic groups defined by 
$$
\tau(t_1,\dots, t_n)=s\circ \big(\tsize\bigoplus_{w\in \bZ}\big(\big(\tsize\prod_{j=1}^n t_j\big)^w\rho_w(g_1, \dots, g_n) \;\text{on}\; \gr^W_w \big)\big)\circ s^{-1}, 
$$
where $g_j$ is as in $2.3.5$.
Then, as $y_j>0$, $y_1=\cdots=y_k$, 
$y_j/y_{j+1}\to \infty $ $(k\leq j\leq n$, $y_{n+1}$ means $1)$,
$$
\tau\left(\sqrt{\frac{y_2}{y_1}},\dots, \sqrt{\frac{y_{n+1}}{y_n}}\right)^{-1}\exp(\ts_{j=1}^n iy_jN_j)F 
$$
converges in $D$. 

Let $\br_1\in D$ be the limit. 
\medskip

{\rm(iii)} Let 
$$
J'=\{j\;|\;1\leq j\leq n, \;\text{the $j$-th component of $\rho_w$ is non-trivial for some $w\in \bZ$}\}.
$$ 
Let $J=J'=\emptyset$ if $k=n+1$, and let $J=J'\cup\{k\}$ if otherwise.
Then
$$
((\rho_w, \varphi_w)_w, \br_1, J)\in \cD_{\SL(2),n}.
$$
\medskip

{\rm(iv)} The family of weight filtrations {\rm(2.3.6)} and the torus action {\rm(2.3.5)} associated to $((\rho_w, \varphi_w)_w, \br_1, J)$ coincide with $(M(N_1+\dots+N_j, W))_{1\leq j\leq n}$ and $\tau$ in {\rm(ii)}, respectively. 
\endproclaim

We prove this theorem later in 2.4.8.

By this theorem, we have a map
$$
\psi: \cD_{\nilp, n} \to \cD_{\SL(2), n}, \quad (N_1, \dots, N_n, F)\mapsto
((\rho_w, \vf_w)_w, \br_1, J), 
$$
(for the notation, see 2.4.1, 2.3.2).
For $p\in\cD_{\nilp, n}$, we call $\psi(p)\in \cD_{\SL(2), n}$ the {\it $\SL(2)$-orbit associated to $p$}. 
  Note that this definition is slightly different from that 
in \cite{KNU1}, 0.2.
  Note also that though in the definition of a nilpotent 
orbit in 2.4.1, the order of $N_1, \dots, N_n$ in $(N_1, \cdots, N_n, F)$ is not important, when we consider 
the SL(2)-orbit associated to $(N_1, \cdots, N_n, F)$, the order of $N_1,\dots, N_n$ becomes essential.  

Even when $k=1$, the $\br_1$ in 2.4.2 (ii) is not $\br$ but 
$\exp(\ve_0)\br$ in the main theorem 0.5 of \cite{KNU1}.
  Note that the $s$ in 2.4.2 (i) coincides with $\spl_W(\br_1)$ (1.2.3).
\medskip

\proclaim{Proposition 2.4.3} 
Let $(N_1, \dots, N_n, F)\in \cD_{\nilp, n}$, and let $W^{(j)}=M(N_1+\dots+N_j, W)$ for $1\leq j\leq n$ (cf. $2.2.2$ in the pure case).
Let $k=\min(\{j\;|\;1\leq j \leq n,\; N_j\not=0\} \cup \{n+1\})$. 
Then the following two conditions {\rm(1)} and {\rm(2)} are equivalent.
\medskip

{\rm(1)} For any $k\leq j \leq n$, 
$(W^{(j)},\exp(\ts_{l=j+1}^n iN_l)F)$ is an $\bR$-split mixed Hodge structure.
\medskip

{\rm (2)} For any $k\leq j\leq n$ and for any $y_l\in \bR_{>0}$ $(j<l \leq n)$, $(W^{(j)}, \exp(\ts_{l=j+1}^n iy_lN_l)F)$ is an $\bR$-split mixed Hodge structure.
\endproclaim

We prove this proposition later in 2.4.9.
\medskip

{\bf 2.4.4.} 
Let $\cD_{\nilp, \SL(2), n}\subset \cD_{\nilp, n}$ be the set of all $(N_1, \dots, N_n, F)\in \cD_{\nilp, n}$ which satisfy the equivalent conditions in 2.4.3. 

For example, $\cD_{\nilp, \SL(2), 1}$ is the set of all $(N, F)\in \cD_{\nilp, 1}$ such that $N = 0$ or 
$(M(N, W), F)$ is an $\bR$-split mixed Hodge structure.

\proclaim{Theorem 2.4.5}
For $p=(N_1, \dots, N_n, F)\in \cD_{\nilp, n}$, 
let $k=\min(\{j\;|\;1\leq j \leq n,\; N_j\not=0\} \cup \{n+1\})$ 
and 
let $\phi(p)=(N_1, \cdots, N_k, N_{k+1}^{\Delta}, \dots, N_n^{\Delta}, F')$, 
where $F'=F$ if $k=n+1$ and $F'=\hat F_{(n)}$ 
if otherwise. 
$(N_j^{\Delta}\in \fg_{\bR}$ $(k< j \leq n)$, $\hat F_{(n)}\in \Dc$ 
are as in {\rm \cite{KNU1}, 10.1--10.2}. 
We review these objects in {\rm 2.4.6--2.4.7} below.$)$
\medskip

{\rm(i)} For $p\in \cD_{\nilp, n}$, we have $\phi(p)\in \cD_{\nilp, n}$ and $\phi(\phi(p))=\phi(p)$.
\medskip

{\rm(ii)} $\cD_{\nilp, \SL(2), n}=\{p\in \cD_{\nilp, n}\;|\;\phi(p)=p\}$. 
\medskip

{\rm(iii)} The map $\psi: \cD_{\nilp, \SL(2), n} \to \cD_{\SL(2),n}$ is injective. 
This map is described via $2.3.7$ as follows. 
For $p=(N_1,\dots, N_n, F)\in \cD_{\nilp, \SL(2), n}$, the  family of weight filtrations associated to $\psi(p)$ is given as in Theorem $2.4.2$ {\rm (iv)}, 
$\br_1= \exp(iN_1+\dots+iN_n)F$, and $J=\{j\;|\;1\leq j\leq n, N_j\neq 0\}$.
If $J=\{a(1), \dots, a(r)\}$ ($a(1)<\dots<a(r)$) and if $p'$ denotes $(N_{a(1)}, \dots, N_{a(r)}, F)$, $\psi(p')$ coincides with the $\SL(2)$-orbit in $r$ variables of rank $r$ associated to $\psi(p)$ $(2.3.4)$.
\medskip

{\rm(iv)} In the pure case, the map $\psi: \cD_{\nilp, \SL(2), n} \to \cD_{\SL(2),n}$ is bijective. 
The converse map is given by $(\rho, \vf)\mapsto (N_1, \dots, N_n, 
\vf(\bold0))$, where $N_j$ is the operator associated to $\rho$ in $2.1.4$.
\medskip

{\rm(v)} The map $\psi : \cD_{\nilp, n}\to \cD_{\SL(2), n}$ factors as
$\cD_{\nilp, n} @>\phi>> \cD_{\nilp, \SL(2), n} \overset \psi \to \hookrightarrow \cD_{\SL(2),n}$.
\medskip

{\rm(vi)} Assume $p = (N_1, \dots, N_n, F)\in \cD_{\nilp, \SL(2), n}$. 
Let 
$$
\psi(p)=((\rho_w, \vf_w)_w, \br_1, J)
$$ 
$(2.4.2)$, and let $(W^{(j)})_{1\leq j\leq n}$ be the family of weight filtrations associated to $\psi(p)$. 
Then $(W^{(j)}, \br_1)$ is a mixed Hodge structure for $1\leq j\leq n$, and 
$p$ is recovered from $\psi(p)$ by the following $(1)$ and $(2)$.
\medskip

$(1)$ Let $k=\min(J \cup \{n+1\})$. 
For $1\leq j< k$, $N_j=0$. 
For $k\leq j\leq n$, $\ts_{l=k}^j N_l=s^{(j)}\delta(W^{(j)}, \br_1)
({s^{(j)}})^{-1}$, where  $s^{(j)}$ is the $s_{\br_1}$-lift (cf.\ $2.4.6$; 
see $2.3.5$ for $s_{\br_1}$) of $(s^{(j)}\;\text{ of }\; (\rho_w, \vf_w))_w$. 

\medskip

$(2)$ 
If $k=n+1$, $F=\br_1$.  If otherwise, 
$(W^{(n)}, F)$ is the $\bR$-split mixed Hodge structure associated to the mixed Hodge structure $(W^{(n)}, \br_1)$.
\endproclaim

We prove this theorem later in 2.4.10.

The injection $\psi: \cD_{\nilp, \SL(2), n}\to \cD_{\SL(2), n}$ need not be surjective though it is bijective in the pure case (2.4.5 (iv)). 
See 2.4.11 Example III.

Some readers may prefer to define an $\SL(2)$-orbit as an element of
$\tsCu_n \cD_{\nilp, \SL(2), n}$, not using $\cD_{\SL(2), n}$. 
The reason why we use  the set $\cD_{\SL(2), n}$ is that the space $D_{\SL(2)}$ of the classes of $\SL(2)$-orbits defined by using $\cD_{\SL(2),n}$  has nice properties (3.5.15, for example).
\medskip

We now give preparations for the proofs of 2.4.2, 2.4.3, 2.4.5.
\medskip

{\bf 2.4.6.} 
Let $(N_1, \dots, N_n, F)\in \cD_{\nilp,n}$. 
In the following, we review an alternative 
construction of $s, \tau$ and $\br_1$ by a finite number of algebraic 
steps, not by a limit.  
In particular, we review the definition of $\hat F_{(n)}$. 

For $0\leq j \leq n$, we denote $M(\ts_{k=1}^j N_k, W)$ by $W^{(j)}$.
In particular, $W^{(0)}=W$. 

For $0 \leq j \leq n$, we define an $\bR$-split mixed Hodge structure $(W^{(j)}, \hat F_{(j)})$ and the associated splitting $s^{(j)}$ of $W^{(j)}$ inductively starting from $j=n$ and ending at $j=0$ (see \cite{KNU1} 10.1; in the pure case,  see \cite{CKS}). 
Note that, in the definition of mixed Hodge structure, we do not assume that the weight filtration is rational (cf. 1.2.8).
First, $(W^{(n)}, F)$ is a mixed Hodge structure as is proved 
by Deligne (see \cite{K}, 5.2.1). 
Let $(W^{(n)}, \hat F_{(n)})$ be the $\bR$-split mixed Hodge
structure associated to the mixed Hodge structure $(W^{(n)}, F)$. 
Then $(W^{(n-1)}, \exp(iN_n)\hat F_{(n)})$ is a mixed Hodge structure. 
Let $(W^{(n-1)}, \hat F_{(n-1)})$ be the $\bR$-split mixed Hodge structure associated to $(W^{(n-1)}, \exp(iN_n)\hat F_{(n)})$. 
Then $(W^{(n-2)}, \exp(iN_{n-1})\hat F_{(n-1)})$ is a mixed Hodge structure. 
This process continues. 
In this way we define $\hat F_{(j)}$ inductively as the $\bR$-split mixed Hodge structure associated to the mixed Hodge structure $(W^{(j)}, \exp(iN_{j+1})\hat F_{(j+1)})$ and define $s^{(j)}$ to be the splitting of $W^{(j)}$ associated to $\hat F_{(j)}$. 
  The splitting $s$ in 2.4.2 (i) is nothing but $s^{(0)}$ (\cite{KNU1} 10.1.2). 

  Thus we have 
$s^{(j)}=\spl_{W^{(j)}}(\exp(iN_{j+1})\hat F_{(j+1)})$,
$\hat F_{(j)}=s^{(j)}((\exp(iN_{j+1})\hat F_{(j+1)})(\gr^{W^{(j)}}))$
($N_{n+1}:=0$, $\hat F_{(n+1)}:=F$).
  We also have $\br_1=\exp(iN_k)\hat F_{(k)}$, where 
$k=\min(\{j\;|\;1\leq j \leq n,\; N_j\not=0\} \cup \{n+1\})$ (cf.\ 2.4.8). 

These  $s^{(j)}$ ($0\leq j\leq n$) are compatible in the sense that we have a direct sum decomposition 
$$
H_{0,\bR} = \tsize\bigoplus_{\theta\in \bZ^{n+1}} \,H_{0,\bR}^{[\theta]}, \quad 
\text{where}\;\; H_{0,\bR}^{[\theta]}= 
 \tsize\bigcap_{j=0}^n \; s^{(j)}(\gr^{W^{(j)}}_{\theta(j)}).
$$ 
This compatibility is expressed also in the following way.  
Let
$$
\tau_j : \bG_{m,\bR} \to \Aut_{\bR}(H_{0,\bR},W)\quad (0 \leq j \leq n)
$$ 
be the homomorphism of algebraic groups over $\bR$ characterized as follows. For $a\in \bR^\times$ and $w\in \bZ$, $\tau_j(a)$ acts on $s^{(j)}(\gr^{W^{(j)}}_w)$ as the multiplication by $a^w$. 
Then the compatibility of $s^{(j)}$ ($0 \leq j \leq n$) in question is expressed as
the fact that $\tau_j(a)\tau_k(b)=\tau_k(b)\tau_j(a)$ for any $j, k$ and any $a, b\in \bR_{>0}$. 
Let 
$$
\tau(a)=\tp_{j=1}^n \tau_j(a_j) \quad \text{for}\;a=(a_j)_j \in \bR^n_{>0}.
$$

  This $\tau$ coincides with the $\tau$ in 2.4.2 (ii). 
  Note also that $s^{(j)}$ is the $s$-lift of
$(s^{(j)}\;\text{ of }\; (\rho_w, \vf_w))_w$, that is, 
it coincides with the composite 
$\gr^{W^{(j)}} \cong \bigoplus_w\gr^{W^{(j)}}(\gr^W_w) \to 
\bigoplus_w\gr^W_w \overset s \to \to H_{0,\bR}$, where the first arrow 
is the sum of the splittings $s^{(j)}$ on $\gr^W_w$ with respect to 
$(\rho_w, \vf_w)$. 
In \cite{KNU1}, 10.3, we denoted $\tau_j(\sqrt{a})^{-1}$ for $a\in \bR_{>0}$ by $t^{(j)}(a)$, and $\tau((\sqrt{a_{j+1}/a_{j}})_j)$ for $a=(a_j)_j\in \bR^n_{>0}$ $(a_{n+1}:=1)$ by $t(a)$.

Any $h \in \fg_{\bR}$ is decomposed uniquely in the form 
$$
h=\tsize\sum_{\theta \in \bZ^{n+1}} \;h^{[\theta]}, \quad h^{[\theta]}\in \fg_\bR, \quad h^{[\theta]}(H_{0,\bR}^{[\theta']})\sub H_{0, \bR}^{[\theta+\theta']} \;\;\;(\forall\; \theta'\in \bZ^{n+1}).
$$
\medskip

\proclaim{Proposition 2.4.7}
Let the notation be as above. 
\medskip

{\rm(i)} Let $1\leq j \leq n$ and let $\theta=(\theta(k))_{0\leq k \leq n}\in \bZ^{n+1}$ $(\theta(k) \in \bZ)$. 
Then $N_j^{[\theta]}=0$ unless $\theta(k)=-2$ for $j\leq k \leq n$. 
\medskip

{\rm(ii)} Let $1\leq j\leq n$, and define $\hat N_j$ $($resp. $N_j^{\Delta}$, resp. $\hat N_j')$ to be the sum of $N_j^{[\theta]}$, where $\theta$ ranges over all elements of $\bZ^{n+1}$ such that $\theta(k)=0$ for $0\leq k \leq j-1$ $($resp. for $1\leq k \leq j-1$, resp. for $k=j-1)$. 
Then 
$$
\hat N_j = \hat N_j'.
$$
Consequently, 
$$
N_j^{\Delta}=\hat N_j \quad \text{for}\;\; 2\leq j\leq n, \quad N^{\Delta}_1=N_1.
$$
\medskip

{\rm(iii)} $N_j \hat N_k= \hat N_k N_j $ if $1\leq j<k \leq n$.
\medskip

{\rm(iv)}
$\hat N_j \hat N_k = \hat N_k \hat N_j$ and $N_j^{\Delta}N_k^{\Delta}=N_k^{\Delta}N_j^{\Delta}$ for all $j, k$.
\medskip

{\rm(v)} Assume $1\leq j\leq k \leq n$. Then $(W^{(k)}, \hat F_{(j)})$ is a
mixed Hodge structure. 
The $\bR$-split mixed Hodge structure associated to $(W^{(k)}, \hat F_{(j)})$ is $(W^{(k)}, \hat F_{(k)})$, $(s^{(k)})^{-1}N_\ell s^{(k)}$ and $(s^{(k)})^{-1}\hat N_{\ell}s^{(k)}$ belong to $L^{-1,-1}_\bR(W^{(k)},\hat F_{(k)})$ {\rm(1.2.1)} for all $\ell \leq k$, and $\delta(W^{(k)}, \hat F_{(j)})= (s^{(k)})^{-1}(\ts_{j<\ell\leq k} \hat N_{\ell})s^{(k)}$. 
\endproclaim

{\it Remark 1.} 
  Thus $(N_1^{\Delta},\ldots, N_n^{\Delta})$ is nothing but 
$(N_1, \hat N_2, \ldots, \hat N_n)$.
  In \cite{KNU1}, (ii) of the above proposition was not recognized 
so that we did not unify the notation $N_j^{\Delta}$ and $\hat N_j$. 

{\it Remark 2.} 
In the case $j \geq k$, $N_j \hat N_k= \hat N_k N_j $ in 2.4.7 (iii) need not be true. For example, in Example III in 1.1.1, 
if we take $N$ in 2.4.11 III below as $N_j$ for $1\leq j\leq n$ and take $F$ 
in 2.4.11 III, then $\hat N_1$ sends $e_1$ and $e_3$ to zero and $e_2$ to 
$e_1$, so 
$N_j \hat N_1=0$ but $\hat N_1 N_j$ is not zero for any 
$j$. 
(On the other hand, in this example, $\hat N_j=0$ for $j \geq 2$ and hence $N_1 \hat N_j=\hat N_j N_1$ is trivially true for $j \geq 2$.)   
\medskip

{\it Proof of Proposition 2.4.7.} 
(i) is explained in \cite{KNU1}, 10.3. We give the proofs of the remaining statements.

Let $1\leq j \leq n$. 
By \cite{KNU1}, 10.1.4, 
  $\hat F_{(j)}= s(\vf(\{0\}^j \times \{i\}^{n-j}))$. 
  Here $s$ is the splitting of $W$ associated to $\br_1$.
  From this, we have 

\medskip

\noindent
(1) $\hat F_{(j)}$ coincides with 
$s^{(k)}(\tsize\bigoplus_w \hat F_{(j)}(\gr^{W^{(k)}}_w))$ if $0 \leq k \leq j$. 
\medskip

\noindent 
By (1) and by $(s^{(j)})^{-1}N_ks^{(j)}
\in L^{-1,-1}_\bR(W^{(j)}, \hat F_{(j)})$ for $1\leq k \leq j$, we have
\medskip

\noindent
(2) $(s^{(j)})^{-1}\hat N_ks^{(j)}$ and $(s^{(j)})^{-1}\hat N_k's^{(j)}$ belong to $L^{-1,-1}_\bR(W^{(j)},\hat F_{(j)})$ for $1\leq k \leq j$. 
\medskip

We prove (ii). 
By (1), we see that 
$$
\hat F_{(j-1)} = \exp(i\hat N_j') \hat F_{(j)}
$$ 
and since $(W^{(j)}, \hat F_{(j)})$ is an $\bR$-split mixed Hodge structure, we have by (2), 
$$
\delta(W^{(j)}, \hat F_{(j-1)})=(s^{(j)})^{-1}\hat N_j's^{(j)}.\tag3
$$
Note that $\zeta=0$ since $\delta$ has only $(-1,-1)$-Hodge component (see 1.2.3).

Next, by \cite{KNU1}, 10.4 (1), 
$$
\hat F_{(j-1)}=\exp(i\hat N_j) \hat F_{(j)}.
$$ 
Hence by (1) and (2), we have 
$$
\delta(W^{(j)}, \hat F_{(j-1)})=(s^{(j)})^{-1}\hat N_js^{(j)}.\tag4
$$
Comparing (3) and (4), we conclude $\hat N_j = \hat N'_j$.

We prove (iii). 
Since $N_j^{[\theta]}=0$ unless $\theta(k-1)=-2$ by (i), and since $\hat N_k = \hat N_k'$ by (ii), $N_j\hat N_k$ (resp. $\hat N_k N_j$) is the sum of $(N_jN_k)^{[\theta]}$ (resp. $(N_kN_j)^{[\theta]}$), where $\theta$ ranges over
all elements of $\bZ^{n+1}$ such that $\theta(k-1)=-2$. 
But $N_jN_k=N_kN_j$. 
(iii) follows.

We prove (iv). 
We may assume $j < k$. 
Then, by (ii), $\hat N_j \hat N_k$ (resp. $\hat N_k \hat N_j$) is the sum of $(N_j\hat N_k)^{[\theta]}$ (resp. $(\hat N_k N_j)^{[\theta]}$), where $\theta$ ranges over all elements of $\bZ^{n+1}$ such that $\theta(j-1)=0$. 
But $N_j \hat N_k= \hat N_k N_j$ by (iii).
The first assertion of (iv) follows, and hence the second follows.

The rest is (v).  
  Again by \cite{KNU1} 10.4 (1), we have 
$\hat F_{(j)} =\exp(\ts_{j<\ell\leq k} i\hat N_{\ell})\hat F_{(k)}$. 
This implies (v) by the same argument in the proof of (ii). 
\qed

\medskip

{\bf 2.4.8.} 
We prove Theorem 2.4.2. (i) is contained in Theorem 0.5 of \cite{KNU1}. 

We prove (ii). 
  It is clear in case when $k=n+1$.
  When $k\le n$, the proof of \cite{KNU1} 10.4 (2) shows (ii) and, 
furthermore,  
$\br_1=\exp(iN_k+iN_{k+1}^{\Delta}+\cdots+iN_{n}^{\Delta})\hat F_{(n)}$. 

We prove (iii). 
We may assume $k \le n$. 
By the pure case, $\br_1(\gr^W_w)=\vf_w(\bi)$. 
By the calculation of $\br_1$ in the above proof of (ii) together 
with 10.4 (1) of \cite{KNU1} and 2.4.7 (ii), we have
\medskip

{\bf Claim.} 
{\it $\br_1$ in {\rm 2.4.2 (ii)} coincides with $\exp(iN_k)\hat F_{(k)}$.}

\medskip

If $\gr^W(N_k)\neq 0$, then, by the pure case, $k \in J'$ and hence there is no problem (2.3.2). 
Assume $\gr^W(N_k)=0$. 
Then $W^{(k)}=W$ and hence  $(W, \hat F_{(k)})$ is an $\bR$-split mixed Hodge structure. 
Since $N_k$ sends the $(p, q)$-Hodge component of $(W, \hat F_{(k)})$ to the $(p-1, q-1)$-Hodge component, we have $\delta(W, \exp(iN_k)\hat F_{(k)})=s^{-1}N_ks$. 
This shows that if $N_k\neq 0$, then  $\br_1=\exp(iN_k)\hat F_{(k)}$ 
(Claim) belongs to $D_{\nspl}$ (1.2.7). 
Hence $((\rho_w, \vf_w)_w, \br_1, J)\in \cD_{\SL(2), n}$ (2.3.2).
 
We prove 2.4.2 (iv). 
Since $s^{(0)}$ in 2.4.6 coincides with $\spl_W(\br_1)$ and also with the $s$ in 2.4.2 (i) (by Claim), 
it is reduced to the pure case that 
$\tau$ in 2.4.6 coincides with the torus action associated to $((\rho_w, \vf_w)_w, \br_1, J)$ in 2.3.5 and also with the $\tau$ in 2.4.2 (ii). 
This also shows the statement for the associated weight filtrations.
\medskip


{\bf 2.4.9.} 
We prove Proposition 2.4.3. 
We may assume $k \le n$. 
It is enough to show that (1) implies (2). 
Assume (1). 
By 2.4.6, we have
$\hat F_{(j)}= \exp(\ts_{l=j+1}^n iN_l)F$ for $k\leq j \leq n$. 
This gives $\hat F^{(n)}=F$ and also
$\delta(W^{(n)}, \hat F_{(j)})= (s^{(n)})^{-1}(\ts_{l=j+1}^n N_l)
s^{(n)}$ for $k\leq j \leq n$.
Comparing this with
$\delta(W^{(n)}, \hat F_{(j)})= (s^{(n)})^{-1}(\ts_{l=j+1}^n \hat N_l)s^{(n)}$ 
($k \leq j \leq n$) obtained in 2.4.7 (v), we have
$\hat N_j=N_j$ for $k < j\leq n$.
  This implies (2).
 
 \medskip
 
{\bf 2.4.10.} 
We prove Theorem 2.4.5. 
 
We prove (i).
We may assume $k \leq n$. 
We show $\phi(p)\in \cD_{\nilp, n}$ by checking the conditions (1)--(4) in 2.4.1. 
(1) is satisfied by 2.4.7 (ii)--(iv). 
2.4.1 (2) is seen by reduction to the pure case. 
2.4.1 (3) (Griffiths transversality) for $N_j$ follows from \cite{KNU1} 5.7, and that for $\hat N_j$ is deduced from it and from (1) in the proof of 2.4.7. 
We show 2.4.1 (4) (concerning relative monodromy filtration). 
By Kashiwara \cite{K}, 4.4.1 and by 2.4.7 (ii), it is sufficient to show that
the relative monodromy filtration exists for $\hat N_j$ $(k\le j \le n)$ and for $N_k$. 
For $N_k$, this is included in the assumption. 
For $\hat N_j$ $(k \leq j \leq n)$, 
this is easy since $\hat N_j$ is of weight $0$ with respect 
to $s^{(0)}$.
Once $\phi(p) \in \cD_{\nilp,n}$ is verified, it is easy to see $\phi\circ\phi =\phi$.

(ii) is essentially proved in 2.4.9. 
(iii) will be proved later. 
(iv) is known as the pure case (see \cite{KU2} \S6). 
(v) is easy. 

We prove (vi).  
  For $k \leq j \leq n$, we have $\br_1=\exp(\ts_{l=k}^j iN_l)\hat F_{(j)}$ 
in the notation in 2.4.6.
  In particular, we have $\br_1=\exp(\ts_{l=k}^n iN_l)F$. 
  (vi) is deduced from these relations. 
 
We prove (iii). The injectivity follows from (vi). 

To prove $J=\{j\;|\;1\leq j\leq n, N_j\neq 0\}$, first we show 

\medskip
{\bf Claim.} 
{\it For $k< j \leq n$, $W^{(j)}\neq W^{(j-1)}$ if and only if 
$N_j\not=0$.}

\demo{Proof}
  Since $N_j$ is of weight $0$ with respect to $s^{(j-1)}$, the $N_j$ is zero 
if and only if $N_j(\gr^{W^{(j-1)}}_w)$ is zero for any $w$.  
The latter condition is equivalent to $W^{(j)}= W^{(j-1)}$.
\qed
\enddemo

By this claim, we have the description of $J$. 
The remaining parts of (iii) are easy. 
\medskip

{\bf 2.4.11.} 
For some examples in Example I--Example III, we describe 
here the map $\psi : \cD_{\nilp,1} \to \cD_{\SL(2),1}$, $(N, F) \mapsto ((\rho_w, \vf_w)_w, \br_1, J)$ in 2.4.2, the torus actions $\tau_j : \bG_{m,\bR} \to \Aut(H_{0,\bR}, W)$ $(j = 0, 1)$, and nilpotent endomorphisms $\hat N$ and $N^\Delta$ (2.4.7).
\medskip

{\bf Example I.}
Let $N(e_1)=0$, $N(e_2)=e_1$, and let $F=F(i)$. 
Then $(N, F)$ generates a nilpotent orbit.
The canonical splitting of $W$ associated to $\exp(iy_1N)F$ $(y_1>0)$ 
sends $e_2'$ to $e_2$. 
From this we have $\tau(t)e_1=t^{-2}e_1$, $\tau(t)e_2=e_2$. 
For $t=1/\sqrt{y_1}$, we have 
$\lim_{t\to 0}\tau(t)^{-1}\exp(iy_1N)F = F(i)$.
Hence the image $((\rho_w, \vf_w)_w, \br_1, J)$ of $(N,F)$ under $\psi$ consists of Example I in 2.3.9 with $J = \{1\}$.

$W^{(1)} = W$, and $\tau_1 = \tau_0 = \tau$.
Hence $\hat N = 0$, $N^\Delta = N$.
\medskip

{\bf Example II.}
We consider the following example $(N, F)$ which generates a nilpotent orbit.
$N(e_2) = e_1$, $N(e_j) = 0$ $(j = 1, 3)$.  
$F = F(i, ia)$ with $a\in \bR$. 

By 1.2.9, $(s_1, s_2) \in \bR^2$ corresponding to the canonical splitting $\spl_W(\exp(iy_1N)F)$ is $s_1 = 0$, $s_2 = - a/(1 + y_1)$.
When $y_1 \to \infty$, $(s_1, s_2)$ converges to $(0, 0)$ in $\spl(W) = \bR^2$.
From this, we have 
$\tau(t)e_1 = t^{-2}e_1$, $\tau(t)e_j = e_j$, $(j = 2, 3)$.
For $t = 1/\sqrt{y_1}$, we have 
$\lim_{t\to 0}\tau(t)^{-1}\exp(iy_1N)F = \br_1 \in D$, where $\br_1^1 := 0$, $\br_1^0:= \bC(ie_1 + e_2) + \bC e_3$, $\br_1^{-1} := H_{0,\bC}$.
Hence the image $((\rho_w, \vf_w)_w, \br_1, J)$ of $(N,F)$ under $\psi$ consists of the example in 2.3.9 II with $z = ia$ 
and $J = \{1\}$.

The torus actions $\tau_0$, $\tau_1$ which induce the splittings of the filtrations $W$, $W^{(1)}$, in 1.1.1 II, 2.3.9 II, respectively, are as follows.
$\tau_0(t)e_j = t^{-1}e_j$ $(j = 1, 2)$, $\tau_0(t)e_3 = e_3$. 
$\tau_1 = \tau$ above.
Hence $\hat N = N^\Delta = N$.
\medskip

{\bf Example III.} 
We consider the following example $(N, F)$ which generates a nilpotent orbit.
$N(e_3)=e_2$, $N(e_2)=e_1$, $N(e_1)=0$. 
$F=F(i, z, i)$ with $z \in \bC$.

By 1.2.9, $(s_1, s_2) \in \bR^2$ corresponding to the canonical splitting $\spl_W(\exp(iy_1N)F)$ is $s_1 = \Re(z)+\frac 12$, $s_2 = - \Im(z)/2(1 + y_1)$.
When $y_1 \to \infty$, $(s_1, s_2)$ converges to $(\Re(z)+\frac 12, 0)$ in $\spl(W) = \bR^2$.
From this, we have 
$\tau(t)e_1 = t^{-4}e_1$, $\tau(t)e_2 = t^{-2}e_2$, $\tau(t)e_3 = e_3
+(1-t^{-4})(\Re(z)+\frac 12)$.
For $t = 1/\sqrt{y_1}$, we have 
$\lim_{t\to 0}\tau(t)^{-1}\exp(iy_1N)F = \br_1 \in D$, where $\br_1^1 := 0$, $\br_1^0:= \bC(\Re(z)e_1 + ie_2 + e_3)$, $\br_1^{-1}:= \br_1^0 + \bC(ie_1 + e_2)$, $\br_1^{-2} := H_{0,\bC}$.
Hence the image $((\rho_w, \vf_w)_w, \br_1, J)$ of $(N,F)$ under $\psi$ consists of the example in Case 1 of 2.3.9 III with 
$z_1=\Re(z)$ and $J = \{1\}$.

There is no nilpotent orbit whose associated SL(2)-orbit is in Case 2 or Case 
3 in 2.3.9 III (cf. the comment after Theorem 2.4.5). 
  (In Examples I, II, IV, and V, all SL(2)-orbits come from nilpotent 
orbits.)

In the following, assume $\Re(z)=-\frac 12$ for simplicity. 
The torus actions $\tau_0$, $\tau_1$ which induce the splittings of the filtrations $W$, $W^{(1)}$, in 1.1.1 III, Case 1 of 2.3.9 III, respectively, are as follows.
$\tau_0(t)e_j = t^{-3}e_j$ $(j = 1, 2)$, $\tau_0(t)e_3 = e_3$. 
$\tau_1 = \tau$ above.
Hence $\hat N = N^{[(0,-2)]}$ is given by $\hat N(e_2) = e_1$, $\hat N(e_j) = 0$ $(j = 1, 3)$.
$N^\Delta = N$.
\medskip

\vskip20pt

\head
\S2.5. Definition of the set $D_{\SL(2)}$
\endhead
\medskip

{\bf 2.5.1.} 
Two non-degenerate $\SL(2)$-orbits $p=((\rho_w, \vf_w)_w, \br)$ and $p'=
((\rho'_w, \vf'_w)_w, \br')$ in $n$ variables of rank $n$ (2.3.3) are said to be {\it equivalent} if there is a $t \in \bR^n_{>0}$ such that
$$
\rho'_w=\Int(\gr^W_w(\tau(t))) \circ \rho_w, \quad \vf'_w=\gr^W_w(\tau(t))\circ \vf_w\quad(\forall\; w\in \bZ), \quad \br'=\tau(t)\br.
$$
Here $\tau : \bG^n_{m,\bR}\to \Aut_{\bR}(H_{0,\bR}, W)$ is the torus action associated to $((\rho_w, \vf_w)_w, \br)$ defined in 2.3.5. 

Note that this is actually an equivalence relation.
We explain this.
For $t=(t_j)_j \in \bR_{>0}^n$, we denote by $\widetilde \rho_w(t)=
\rho_w(g_1, \dots, g_n)$ in 2.3.5. 
Since $\gr^W_w(\tau(t))=(\tp_{j=1}^nt_j)^w \widetilde \rho_w(t)$ 
for $t\in\bR_{>0}^n$ (2.3.5), we have $\widetilde {\rho'}_w
=\widetilde \rho_w$ as homomorphisms $\bG_{m,\bR}^n\to G_\bR(\gr^W_w)$ 
for any $w$.
On the other hand, the splittings of $W$ associated to $\br$ and to $\br' = \tau(t)\br$ coincide by the Remark in 2.3.5.
From these it follows that $\tau$ of $p$ and $\tau$ of $p'$ coincide.
The axioms of equivalence relation can be now easily checked.

An $\SL(2)$-orbit $((\rho_w, \vf_w)_w, \br, J)$ in $n$ variables of rank $r$ and an $\SL(2)$-orbit $((\rho'_w, \vf'_w)_w, \br', J')$ in $n'$ variables of rank $r'$ are said to be {\it equivalent} if $r=r'$ and their associated  $\SL(2)$-orbits in $r$ variables of rank $r$ (2.3.4) are equivalent.

\medskip
The class determines and is determined by 
the associated set of weight filtrations, the associated torus action,
and the associated torus orbit, that is, we have

\proclaim{Proposition 2.5.2} 
Let $p=((\rho_w, \vf_w)_w, \br)$ be a non-degenerate  $\SL(2)$-orbit of rank $n$.
\medskip

{\rm(i)}
The $W^{(j)}$ of $p$, the $\tau$ and the $\tau_j$ of $p$ $(1\leq j \leq n)$, the canonical splitting of $W$ associated to $\br$ $(1.2.3)$, and 
$Z= \tau(\bR^n_{>0})\br$ depend only on the equivalence class of $p$. 
Here $\tau$ is the homomorphism in $2.3.5$ associated to $p$.
$Z$ is called the {\rm torus orbit associated to $p$}.

\medskip

{\rm(ii)} 
The equivalence class of $p$ is determined by $((W^{(j)}(\gr^W))_{1\leq j \leq n}, Z)$, where $Z$ is as above. 

\medskip

{\rm(iii)} 
The equivalence class of $p$ is determined by $(\tau, Z)$, where $\tau$ and $Z$ are as above.
\endproclaim

{\it Proof.} 
We prove (i). 
The statement for $W^{(j)}$ follows from $\tau(t)W^{(j)}=W^{(j)}$ ($t \in (\bR^\times)^n$), the statements for $\tau$ 
and for the splitting were proved in 2.5.1, 
and the rest is clear.

(ii) and (iii) follow from (i) and Proposition 2.3.7.
\qed
\medskip

{\bf 2.5.3.} 
Let $D_{\SL(2)}$ be the set of all equivalence classes of $\SL(2)$-orbits satisfying the following condition (C).  

Take an $\SL(2)$-orbit $((\rho_w, \vf_w)_w, \br, J)$ in $n$ variables which is a representative of the class in question. 
\medskip

\noindent
(C) For each $w\in \bZ$ and for each $1\leq j \leq n$, the weight filtration
$W^{(j)}(\gr^W_w)$ is rational.  
\medskip

\noindent (This condition is independent of the choice of the  
representative by 2.5.2 (i).) 

As a set, we have 
$$
D_{\SL(2)}=\tsCu_{n\ge0}\;D_{\SL(2),n},
$$
where $D_{\SL(2), n}$ is the 
set of equivalence classes of $\SL(2)$-orbits of rank $n$ (2.3.3) with rational associated weight filtrations. 
We identify $D_{\SL(2),0}$ with $D$ in the evident way. 

Let $D_{\SL(2),\spl}$ be the subset of $D_{\SL(2)}$ consisting of the classes of $((\rho_w, \vf_w)_w, \br)$ with $\br\in D_{\spl}$ (see Notation in $\S0$).
(The last condition is independent of the choice of the representative.)
Let $D_{\SL(2),\nspl}=D_{\SL(2)}\smallsetminus D_{\SL(2),\spl}$. 
\medskip

{\bf 2.5.4.} 
We have a canonical projection
$$
D_{\SL(2)}\to D_{\SL(2)}(\gr^W)=\tp_{w\in \bZ}\; D_{\SL(2)}(\gr^W_w),
$$
$$
\class((\rho_w, \vf_w)_w, \br) \mapsto (\class(\rho_w, \vf_w))_w.
$$
Here $D_{\SL(2)}(\gr^W_w)$ is the $D_{\SL(2)}$ for
$((H_0\cap W_w)/(H_0\cap W_{w-1}), \lan\;,\;\ran_w)$.
  Note that in the pure case, the definition of $D_{\SL(2)}$ coincides with that of 
\cite{KU2}.
\medskip

{\bf 2.5.5.} 
As in Notation in $\S0$, \ let $\spl(W)$ be the set of all splittings of $W$. 
We have a canonical map
$$
D_{\SL(2)}\to \spl(W)
$$
as $\class((\rho_w, \vf_w)_w, \br) \mapsto s$, where $s$ denotes the canonical splitting of $W$ associated to $\br$ (see 2.5.2 (i)). 
\medskip

{\bf 2.5.6.} 
For $p\in D_{\SL(2)}$, we denote by $\tau_p$ and $Z_p$ the corresponding $\tau$ and $Z$, respectively (see 2.5.2 (iii)).
\medskip

{\bf 2.5.7.}
Later in \S3.2, we will define two topologies on the set $D_{\SL(2)}$. 
Basic properties of these topologies are the following (3.2, 4.1.1). 
\medskip

(i) If $p\in D_{\SL(2)}$ is the class of $(\tau_p, \br)$, then we have in $D_{\SL(2)}$ 
$$
\tau_p(t)\br \to p \quad \text{when $t\in \bR_{>0}^n$ tends to $\bold0$}.
$$
Here $n$ is the rank of $p$ and $\bold0=(0, \dots, 0)\in \bR_{\geq 0}^n$. 
\medskip

(ii) If $(N_1, \cdots, N_n, F)$ generates a nilpotent orbit and if the monodromy filtration of $\gr^W_w(N_1)+\dots+\gr^W_w(N_j)$ is rational for any $w\in\bZ$ and any $1\leq j\leq n$, then we have in $D_{\SL(2)}$ 
$$
\exp(\ts_{j=1}^n iy_jN_j)F\to p
$$
when $y_j>0, y_j/y_{j+1}\to \infty$ ($1\leq j\leq n$, $y_{n+1}$ denotes $1$), where $p$ denotes the class of the SL(2)-orbit associated to $(N_1, \dots, N_n, F)$ by 2.4.2.
\medskip

This (ii) is the basic principle which lies in our construction of the topologies on $D_{\SL(2)}$. 
Our SL(2)-orbit theorem 0.5 in [KNU1] says roughly that, when $y_j/y_{j+1}\to \infty$ ($1\leq j\leq n$, $y_{n+1}=1$), $\exp(\ts_{j=1}^n iy_jN_j)F$ is near to
$\tau_p(\sqrt{\frac{y_2}{y_1}}, \dots, \sqrt{\frac{y_{n+1}}{y_n}})\br$, 
where $\br \in Z_p$. 
Hence  (i) is natural in view of (ii). 
\medskip

\vskip20pt

\head
\S3. Real analytic structures of $D_{\SL(2)}$
\endhead

\head
\S3.1. Spaces with real analytic structures and log structures with sign
\endhead

\bigskip

We discuss a category $\cB_\bR$ of \lq\lq spaces with real analytic structures'', and its \lq\lq logarithmic version'', a category $\cB_\bR(\log)$. 
In 3.1.11--3.1.13, we consider \lq\lq log modification'' in $\cB_\bR(\log)$ associated to cone decompositions.  

\medskip

{\bf 3.1.1.} {\it The categories $\cB_\bR$, $\cB'_\bR$, and $\cC_\bR$.} 
We define three full subcategories
$$\cB_\bR\sub \cB_\bR'\sub \cC_\bR$$ 
of the category of local ringed spaces over $\bR$. 

We first define $\cB'_\bR$. 
An object of $\cB'_\bR$ is a local ringed space  $(S, \cO_S)$  over $\bR$ such that the following holds locally on $S$. 
There are $n\geq 0$ and a morphism $\iota:S\to \bR^n$ of local ringed spaces over $\bR$ from $S$ to the real analytic manifold $\bR^n$ such that $\iota$ is injective, the topology of $S$ coincides with the one induced from the topology of $\bR^n$ via $\iota$, and the canonical map $\iota^{-1}(\cO_{\bR^n})\to \cO_S$ is surjective. Here $\cO_{\bR^n}$ denotes the sheaf of $\bR$-valued real analytic functions on $\bR^n$ and $\iota^{-1}(\;)$ denotes the inverse image of a sheaf. Morphisms of $\cB'_\bR$ are those of local ringed spaces over $\bR$.

Let $\cB_\bR$ be the full subcategory of $\cB'_\bR$ consisting of all objects for which  locally on $S$, we can take $\iota: S\to \bR^n$ as above such that the kernel of the surjection $\iota^{-1}(\cO_{\bR^n})\to \cO_S$ is a finitely generated ideal. 

Of course a real analytic manifold is an object of $\cB_\bR$. 
An example of an object of $\cB_\bR$ which often appears in this paper is
$\bR_{\geq 0}^n$ with the inverse image of the sheaf of real analytic functions on $\bR^n$.

For an object $(S, \cO_S)$ of $\cB'_\bR$,  we often call $\cO_S$ the sheaf of real analytic functions of $S$ though $(S, \cO_S)$ need not be a real analytic space.

We define another category $\cC_\bR$ as follows.
An object of $\cC_\bR$ is a local ringed space $(S, \cO_S)$ over $\bR$ such that for any open set $U$ of $S$ and for any $n\geq 0$, the canonical map $\Mor(U, \bR^n)\to \cO_S(U)^n,\;\vf\mapsto (\vf_j)_{1\leq j\leq n}$, is bijective, where $\bR^n$ is regarded as a real analytic manifold as usual, $\Mor(U, \bR^n)$ is the set of all morphisms in the category of local ringed spaces over $\bR$, and $\vf_j$ denotes the pull-back of the $j$-th coordinate function of $\bR^n$ via $\vf$.
Morphisms of $\cC_\bR$ are those of local ringed spaces over $\bR$.
 
It is easily seen that real analytic manifolds, $C^\infty$-manifolds (with the sheaves of $C^\infty$-functions), and any topological spaces with the sheaves of real valued continuous functions belong to $\cC_\bR$. 
\medskip

\proclaim{Lemma 3.1.2} 
$$
\cB'_\bR\subset \cC_\bR.
$$
\endproclaim
 
{\it Proof.} 
Let $S$ be an object of $\cB'_\bR$. Let $\Mor_S(-, \bR^n)$ be the sheaf on $S$ of morphisms into $\bR^n$. 
We prove that the map $\Mor_S(-, \bR^n)\to \cO_S^n$ is an isomorphism. 
We first prove the surjectivity. 
A local section of $\cO_S^n$ comes, locally on $S$, from an element of $\cO(V)^n$ for some open set $V$ of $\bR^m$ and for some morphism $S\to V$. Since $\cO(V)^n=\Mor(V, \bR^n)$, a local section of $\cO_S^n$ comes from 
$\Mor_S(-, \bR^n)$ locally on $S$. 
It remains to prove the injectivity of $\Mor_S(-, \bR^n) \to \cO_S^n$. 
We prove

\medskip

{\bf Claim.} 
{\it For any $s\in S$, the local ring $\cO_{S,s}$ is Noetherian.}

\medskip

This is reduced to the fact that the local rings of the real analytic manifold $\bR^n$ are Noetherian. These local rings are the rings of convergent Taylor series. Hence they are Noetherian.

Now we return to the proof of 3.1.2. Assume that two morphisms $f, g:S\to \bR^n$ induce the same element $(\varphi_j)_j$ of $\cO(S)^n$. The underlying map $S\to \bR^n$ of sets induced by $f$ and $g$ are 
given by $s\mapsto (\varphi_j(s))_j$, and hence they coincide. To prove $f=g$, it is sufficient to prove that for any $s\in S$ 
with image $s'=f(s)=g(s)\in \bR^n$ and for any element $h$ of 
$\cO_{\bR^n,s'}$, the pull backs $f^*(h), g^*(h) \in \cO_{S,s}$ coincide. Let $m$ be the maximal ideal of $\cO_{S,s}$ and let $m'$
 be the maximal ideal of $\cO_{\bR^n,s'}$. 
 Let $r\geq 1$. Then $h \bmod (m')^r$ is expressed as a polynomial over $\bR$ in the coordinate functions $t_j$. Hence 
 $f^*(h) \equiv g^*(h)\bmod m^r$. Since $\cO_{S,s}$ is Noetherian, the canonical map $\cO_{S,s}\to \varprojlim_r \cO_{S,s}/m^r$ is injective. Hence $f^*(h)=g^*(h)$ in $\cO_{S,s}$. 
 \qed

\proclaim{Proposition 3.1.3} 
The category $\cB'_\bR$ has fiber products, and $\cB_\bR$ is stable under taking fiber products. 
The underlying topological space of a fiber product in $\cB'_\bR$ is the fiber product of the underlying topological spaces. 
The fiber product in $\cB'_\bR$ is also a fiber product in $\cC_\bR$. 
\endproclaim

{\it Proof.} 
Let $S'\to S$ and $S''\to S$ be morphisms in $\cB'_\bR$.

Working locally on $S$, $S'$ and $S''$, we may assume that there are injective morphisms $\iota:S\to \bR^n$, 
$\iota':S'\to  \bR^{n'}$, $\iota'':S''\to  \bR^{n''}$ such that the topologies of $S$, $S'$, $S''$ are induced from those of
 $\bR^n$, $\bR^{n'}$, and $\bR^{n''}$, respectively, and such that the homomorphisms 
 $\iota^{-1}(\cO_{\bR^n})\to\cO_S$, 
 $(\iota')^{-1}(\cO_{\bR^{n'}})\to \cO_{S'}$ and 
 $(\iota'')^{-1}(\cO_{\bR^{n''}})\to \cO_{S''}$ are surjective. 
  Let $I'$ and $I''$ be the kernels of the last two homomorphisms, 
respectively. 
Let $t_j$ ($1\leq j\leq n$) be the $j$-th coordinate function of $\bR^n$. 
Working locally on $S'$, we may assume that for an open neighborhood  $U'$ of $S'$ in $\bR^{n'}$, there are elements $s'_j\in \cO(U')$ ($1\leq j\leq n$) such that the restriction of $s'_j$ to $S'$ coincides with the pull back of $t_j$ for each $j$. 
Similarly working locally on $S''$, we may assume that for  an open neighborhood $U''$ of $S''$ in $\bR^{n''}$, there are elements  $s''_j\in \cO(U'')$ ($1\leq j\leq n$) such that the restriction of $s''_j$ to $S''$ coincides with the pull back of  $t_j$ for each $j$. 
Let $F:=S'\times_S S''\sub V:=U'\times U''\subset\bR^{n'+n''}$. 
Endow $F$ with the topology as the fiber product, and endow it with the inverse image of 
$$
\cO_V/J \quad \text{with} \;\;J=(I'\cO_V+ I''\cO_V+ (s'_1-s''_1)\cO_V+\dots+(s'_n-s''_n)\cO_V).
$$ 
Here $I'\cO_V+I''\cO_V$ denotes the ideal of $\cO_V$ generated by the inverse images of $I'$ and $I''$. 
  When we regard the diagram $S'\to S\leftarrow S''$ as the one in $\cC_{\bR}$ 
by 3.1.2, we can show that $F$ is the fiber product of it in $\cC_{\bR}$, and 
hence, $F$ is the fiber product also in $\cB'_\bR$. 
  If $S, S', S''$ belong to $\cB_\bR$, we can assume that $I'$ and $I''$ are finitely generated. Then the ideal $J$ is finitely generated. 
\qed
\medskip

We start to discuss log structures.

\medskip

\proclaim{Lemma 3.1.4} 
Let $(S, \cO_S)$ be an object of $\cC_\bR$. 
Let $\cO_{S,>0}^\times$ be the subsheaf of $\cO_S^\times$ consisting of all local sections whose values are $>0$. 
Then $\{\pm 1\} @>\sim>> \cO_S^\times/\cO_{S,>0}^\times$. 
Furthermore, $\cO_{S,>0}^\times$ coincides with the image of $\cO_S^\times \to \cO_S^\times,\;f\mapsto f^2$. 
\endproclaim

{\it Proof.} 
The isomorphisms $\bR_{>0}\times \{\pm 1\} @>\sim>> \bR^\times$ and $\bR_{>0}@>\sim>> \bR_{>0},\;x\mapsto x^2$, of real analytic manifolds induce isomorphisms of sheaves
$$\cO_{S,>0}^\times\times \{\pm 1\} \cong \Mor_S(-, \bR_{>0}\times \{\pm 1\}) @>\sim>> \Mor_S(-, \bR^\times)\cong \cO_S^\times,$$
$$\cO_{S,>0}^\times @>\sim>> \cO_{S,>0}^\times, \;f\mapsto f^2,$$
respectively. This proves 3.1.4. \qed

\medskip

\proclaim{Definition 3.1.5}
For an object $S$ of $\cC_\bR$, a {\rm log structure with sign} on $S$ is an 
integral log structure  $M_S$ on $S$ in the sense of Fontaine-Illusie {\rm(\cite{KU3}, \S2.1)} endowed with a subgroup sheaf $M_{S,>0}^{\gp}$ of $M_S^{\gp}$ satisfying the following three conditions {\rm (1)}--{\rm (3)}. Here $M_S^{\gp}\supset M_S$ denotes the sheaf of commutative groups $\{ab^{-1}\;|\;a, b\in M_S\}$ associated to the sheaf $M_S$ of commutative monoids. 

\medskip

{\rm (1)} $M_{S,>0}^{\gp}\supset \cO^\times_{S,>0}$.

\medskip

{\rm (2)} $\cO_S^\times/\cO_{S,>0}^\times @>\sim>> M_S^{\gp}/M_{S,>0}^{\gp}.$

\medskip

{\rm (3)} Let $M_{S,>0}:= M_S\cap M^{\gp}_{S,>0}\sub M_S^{\gp}$. Then the image of $M_{S,>0}$ in $\cO_S$ under the structural map $M_S\to \cO_S$ of the log structure  has values in $\bR_{\geq 0} \sub \bR$ at any points of $S$. 
 
 \medskip
 
$($We remark $(M_{S,>0})^{\gp}=M_{S,>0}^{\gp}$ and thus $M_{S,>0}^{\gp}$ is recovered from $M_{S,>0}$.$)$
 
  \medskip

Let $\cB_\bR(\log)$ $($resp. $\cB'_\bR(\log)$, resp. $\cC_\bR(\log)$$)$ be the category of objects of $\cB_\bR$ $($resp. $\cB'_\bR$, resp. $\cC_\bR$$)$ endowed with an fs log structure {\rm(\cite{KU3}, \S2.1)} with sign.

\endproclaim

If $S$ is an object of $\cC_\bR(\log)$ such that the structural map $M_S\to \cO_S$ is injective and also the canonical map from $\cO_S$ to the sheaf of real valued functions on $S$ is  injective, then for an object $S'$ of $\cC_\bR(\log)$, a morphism $f:S\to S'$ in $\cC_\bR(\log)$ 
is determined by its underlying map $\bar f$ of sets. For such $S$ and an object $S'$ of $\cC_\bR(\log)$, and for a map $g : S\to S'$ of sets, we sometimes say that $g$ is a morphism of $\cC_\bR(\log)$ if $g=\bar f$ for some morphism $f:S\to S'$ of $\cC_\bR(\log)$. 

\medskip

We introduce some terminologies.

\medskip

{\it Trivial log structure with sign.} 
It is the log structure $M_S=\cO_S^\times$ with $M_{S,>0}^{\gp}=\cO_{S,>0}^\times$. 

\medskip

The {\it inverse image} of a log structure with sign. 
For a morphism $S'\to S$ in $\cC_\bR$ and for a log structure $M_S$ with sign on $S$, the inverse image $M_{S'}$ of $M_S$ on $S'$, which is a log structure with sign on $S'$, is defined as follows. 
As a log structure, $M_{S'}$ is the inverse image of $M_S$ (\cite{KU3}, 2.1.3). $M_{S',>0}^{\gp}$ is the subgroup sheaf of $M_{S'}^\gp$ generated by $\cO_{S',>0}^\times$ and the inverse image of $M_{S,>0}^\gp$. 
\medskip

A {\it chart} of an fs log structure with sign. 
Let $S$ be an object of $\cC_\bR(\log)$.
A chart of  $M_S$ with sign is a pair of an fs monoid $\cS$ and a homomorphism $h:\cS\to M_{S,>0}$ such that $h:\cS\to M_S$ is a chart of the fs log structure $M_S$ (\cite{KU3}, 2.1.5) and such that $M_{S,>0}$ is generated by $\cO_{S,>0}^\times$ and $h(\cS)$ as a sheaf of monoids. A chart of $M_S$ exists locally on $S$.
  This is seen with the fact that $M_{S,>0}/\cO^{\times}_{S,>0} \to 
M_{S}/\cO^{\times}_{S}$ is an isomorphism. 

\medskip

{\bf 3.1.6.} {\it Real toric varieties, real  analytic manifolds with corners.}

As standard examples of objects of $\cB_\bR(\log)$, we have real toric varieties, and also real analytic manifolds with corners.

Let $\cS$ be an fs monoid. We regard $S=\Hom(\cS, \bR_{\geq 0}^{\mult})$ as an object of $\cB_\bR(\log)$ as follows, and call it as a real toric variety associated to $\cS$: 
$\cO_S$ is the sheaf of real valued functions on $S$ which belong to $\cO_X|_S$.  Here $X= \Hom(\cS, \bC^{\mult})=\Spec(\bC[\cS])_{\an}$ and $\cO_X$ denotes the sheaf of complex analytic functions on $X$. 
$M_S$ is the log structure associated to $\cS\to \cO_S$. $M^{\gp}_{S,>0}$ is generated by $\cS^{\gp}$ and $\cO_{S,>0}^\times$. 

For any object $T$ of $\cC_\bR(\log)$, we have 
$$
\Mor(T, \Hom(\cS, \bR_{\geq 0}^{\mult}))= \Hom(\cS, M_{T,>0}).
$$

In the case $\cS=\bN^n$, we have $S=\bR_{\geq 0}^n$. 
We usually regard $\bR_{\geq 0}^n$ as an object of $\cB_\bR(\log)$ in this way.

A real analytic manifold with corners $S$ is a local ringed space over $\bR$ which has an open covering $(U_\lam)_\lam$ such that for each $\lam$, $U_\lam$ is isomorphic to an open set of the object $\bR_{\geq 0}^{n(\lam)}$ of $\cB_\bR(\log)$ for some $n(\lam)\geq 0$. 
The inverse images on $U_\lam$ of the fs log structures with sign of $\bR^{n(\lam)}_{\geq 0}$ glue together to an fs log structure with sign on $S$ canonically. 
Thus a real analytic manifold with corners is regarded canonically as an object of $\cB_\bR(\log)$. 

\medskip

\proclaim{Proposition 3.1.7} 
The category $\cB'_{\bR}(\log)$ has fiber products, and $\cB_\bR(\log)$ is stable under taking fiber products.
A fiber product in $\cB_\bR(\log)$ is a fiber product in $\cC_\bR(\log)$. 
The underlying object of $\cB'_\bR$ $($resp. The underlying topological space$)$ of a fiber product $S'\times_S S''$ in $\cB'_\bR(\log)$ coincides with the fiber product in $\cB'_\bR$ $($resp. fiber product as topological spaces$)$ if one of the following conditions {\rm(1)} and {\rm(2)} is satisfied.

\medskip

{\rm(1)} The log structure of $S$ is trivial.

\medskip

{\rm(2)} The log structure of $S'$ coincides with the inverse image of the log structure of $S$. 
\endproclaim

This is a real analytic version of the complex analytic theory about the category $\cB(\log)$ in \cite{KU3} 2.1.10. 
The proof is given by the same arguments there. 
\medskip

We next consider toric geometry in $\cB_\bR(\log)$ and log modifications in $\cB_\bR(\log)$ and in $\cB'_\bR(\log)$. 
These are real analytic versions of those in $\cB(\log)$ (\cite{KU3}, 3.6).
\medskip

{\bf 3.1.8.} 
Let $N$ be a finitely generated free abelian group whose group law is denoted additively.  
A {\it rational fan} in $N_\bR:=\bR\otimes_{\bZ} N$ is a non-empty set $\Sig$ of sharp rational  finitely generated cones in $N_\bR$ satisfying the following conditions (1) and (2).
\medskip

(1) If $\sig\in \Sig$, any face of $\sig$ belongs to $\Sig$.
\medskip

(2) If $\sig,\tau\in \Sig$, then $\sig\cap \tau$ is a face of $\sig$. 
\medskip

Here a {\it finitely generated cone} in $N_\bR$ is a subset of $N_\bR$ of the form $\{\ts_{j=1}^n a_jN_j\;|\;a_j\in \bR_{\geq 0}\}$ with $N_1,\dots, N_n\in N_\bR$. 

A finitely generated cone in $N_\bR$ is said to be {\it rational} if we can take
$N_1,\dots, N_n\in N_\bQ:=\bQ \otimes_\bZ N$ in the above. 

A finitely generated cone $\sig$ in $N_\bR$ is said to be {\it sharp} if $\sig\cap(-\sig)=\{0\}$. 

For a finitely generated cone $\sig$ in $N_\bR$, a {\it face} of $\sig$ is a non-empty subset $\tau$ of  $\sig$ satisfying the following conditions (3) and (4). 
\medskip

(3) If $x, y\in \tau$ and $a, b\in \bR_{\geq 0}$, then $ax+by\in \tau$.
\medskip

(4) If $x, y\in \sig$ and $x+y\in \tau$, then $x, y\in \tau$. 
\medskip

A face of a finitely generated cone $\sig$ in $N_\bR$ is a finitely generated cone in $N_\bR$. It is rational if $\sig$ is rational.
\medskip

{\bf 3.1.9.}
Let $N$ be as in 3.1.8 and $\Sig$ be a rational fan in $N_\bR$. 
Recalling the definition of the (complex analytic) toric variety $\toric(\Sig)$ corresponding to $\Sig$ ([O] 1.2, see also \cite{KU3} 3.3), we define a subset $\abtoric(\Sig)$ of $\toric(\Sig)$ and a structure of an object of $\cB_\bR(\log)$ on $\abtoric(\Sig)$.

Let $M=\Hom(N, \bZ)$, and denote the group law of $M$ multiplicatively.
 
 For $\sig\in \Sig$, let 
$$
\cS(\sig)=\{\chi\in M\;|\;\chi: N_\bR\to \bR\;\text{sends $\sig$ to}\;\bR_{\geq 0}\}.
$$
Then
$$
\sig=\{x\in N_\bR\;|\;\chi: N_\bR\to \bR\;\text{sends $x$ into}\;\bR_{\geq 0}
\;\text{for any $\chi\in \cS(\sig)$}\}.
$$
We have $\cS(\sig)^{\gp}=M$, where $\cS(\sig)^{\gp}=\{ab^{-1}\;|\;a, b\in \cS(\sig)\}$. 

For $\sig\in \Sig$, let $\toric(\sig)=\Spec(\bC[\cS(\sig)])_\an=\Hom(\cS(\sig),\bC^{\mult})$, where $\bC^{\mult}$ denotes $\bC$ regarded as a multiplicative monoid. Then we have an open covering
$$
\toric(\Sig)=\tsize\bigcup_{\sig\in \Sig}\; \toric(\sig).
$$ 

Let 
$$
\abtoric(\Sig)= \tsize\bigcup_{\sig\in \Sig} \; \abtoric(\sig) \sub \toric(\Sig)= \tsize\bigcup_{\sig\in \Sig} \;\toric(\sig)
$$ 
$$ 
\text{with}\;\;\abtoric(\sig):= \Hom(\cS(\sig), \bR_{\geq 0}^{\mult}).
$$
Then $\abtoric(\Sig)$ has the unique structure of an object of $\cB_\bR(\log)$ whose restriction to each open subsets $\abtoric(\sig)$  coincides with the one given in 3.1.6.

Note that $\abtoric(\Sig) \supset \Hom(M, \bR_{>0})=N \otimes \bR_{>0}$, which is the restriction of $\toric(\Sig) \supset \Hom(M, \bC^\times)=N \otimes \bC^\times$. 
As a subset of $\toric(\Sig)$, $\abtoric(\Sig)$ coincides with the closure of $N \otimes \bR_{>0}$ in $\toric(\Sig)$.

There is a canonical bijection between $\toric(\Sig)$ (resp. $\abtoric(\Sig)$) and the set of all pairs $(\sig, h)$, 
where $\sig\in \Sig$ and $h$ is a homomorphism $\cS(\sig)^\times \to \bC^\times$ (resp. $\cS(\sig)^\times \to \bR_{>0}$). 
Here $\cS(\sig)^\times$ denotes the group of invertible elements of $\cS(\sig)$. 
Indeed, for such a pair $(\sig, h)$, the corresponding element of $\toric(\sig)=\Hom(\cS(\sig), \bC^{\mult})$ (resp. $\abtoric(\sig)=\Hom(\cS(\sig), \bR_{\geq 0}^{\mult})$) is defined to be the homomorphism sending $x\in \cS(\sig)$ to $h(x)$ if $x\in \cS(\sig)^\times$, and to $0$ if $x\notin \cS(\sig)^\times$.  

\medskip

{\bf 3.1.10.} 
Let $\Sig$ and $\Sig'$ be rational fans in $N_\bR$ and assume the following condition (1) is satisfied.
\medskip

(1) For each $\tau\in \Sig'$, there is $\sig\in \Sig$ such that $\tau\sub \sig$.
\medskip

Then, we have a morphism $\toric(\Sig')\to \toric(\Sig)$ of complex analytic spaces (resp. a morphism $\abtoric(\Sig')\to \abtoric(\Sig)$ in $\cB_\bR(\log)$) which induces
the morphisms $\toric(\tau) \to \toric(\sig)$ (resp. $\abtoric(\tau) \to \abtoric(\sig)$)
($\tau\in \Sig'$ $\sig\in \Sig$, $\tau\sub \sig$) induced by the inclusion maps
$\tau\sub \sig$. 

Under the condition (1), let $\Sig'\to \Sig$ be the map which sends $\tau\in\Sig'$ to the smallest $\sig\in \Sig$ with $\tau\sub \sig$. 
Then the map $\toric(\Sig')\to \toric(\Sig)$ (resp. $\abtoric(\Sig')\to \abtoric(\Sig)$) sends the point of $\toric(\Sig')$ (resp. $\abtoric(\Sig')$) corresponding to the pair $(\tau, h')$ ($\tau\in \Sig'$, $h'$ is a homomorphism $\cS(\tau)^\times \to \bC^\times$ (resp. $\cS(\tau)^\times\to\bR_{>0}$) to the point of $\toric(\Sig)$ (resp. $\abtoric(\Sig)$) corresponding to the pair $(\sig, h)$, where $\sig$ is the image of $\tau$ under the map $\Sig'\to \Sig$, and $h$ is the composite of $\cS(\sig)^\times\to \cS(\tau)^\times$ with $h'$. 
\medskip

{\bf 3.1.11.} 
Let $\Sig$ be a finite rational fan in $N_\bR$. 

A {\it finite rational subdivision} of $\Sig$ is a finite rational fan $\Sig'$ in $N_\bR$ satisfying the condition (1) in 3.1.10 and also the following condition (1). 
$$
\tsize\bigcup_{\tau\in \Sig'} \tau = \tsize\bigcup_{\sig\in \Sig} \sig \tag1
$$

For a finite rational subdivision $\Sig'$ of $\Sig$,  the maps $\toric(\Sig')\to \toric(\Sig)$ and $\abtoric(\Sig') \to \abtoric(\Sig)$ are proper. 

\proclaim{Proposition 3.1.12} 
Let $S$ be an object of $\cB_\bR(\log)$ $($resp. $\cB'_\bR(\log))$. Let $\cS$ be an fs monoid and let $\cS\to M_S/\cO_S^\times$ be a homomorphism which lifts locally on $S$ to a chart $\cS\to M_{S,>0}$ of fs log structure with sign $(3.1.5)$. 
Let $\Sig$ be a finite rational subdivision of the cone $\Hom(\cS, \bR_{\geq 0}^{\add})$. 
Then we have an object $S(\Sig)$ of $\cB_\bR(\log)$ $($resp. $\cB'_\bR(\log))$ having the following universal property $(1)$. 

\medskip

{\rm (1)} \; If $T$ is an object of $\cC_\bR(\log)$ over $S$, then there is at most one morphism $T\to S(\Sig)$ over $S$.
We have a criterion for the existence of such a morphism{\rm:}
Such a morphism exists if and only if, for any $t\in T$ and for any homomorphism $h: (M_T/\cO_T^\times)_t
 \to \bN$, there exists $\sig\in \Sig$ such that the composite $\cS \to (M_S/\cO_S^\times)_s \to (M_T/\cO_T^\times)_t \to \bN$ $(s$ is the image of $t$ in $S)$ belongs to $\sig$.

\medskip

The map $S(\Sig)\to  S$ is proper and surjective. 
 \endproclaim

{\it Proof.} 
This $S(\Sig)$ is obtained as follows. 
By taking $N=\Hom(\cS^{\gp}, \bZ)$ and $M=\cS^{\gp}$, define $\abtoric(\Sig)$ as in 3.1.9.
Locally on $S$, take a lift $\cS\to M_{S,>0}$ of $\cS\to M_S/\cO^\times_S$ and consider the corresponding morphism $S\to \Hom(\cS, \bR_{\geq 0}^{\mult})$ (3.1.6). 
Then $S(\Sig)$ is obtained as the fiber product (3.1.7) of $S\to \Hom(\cS, \bR_{\geq 0}^{\mult})\leftarrow \abtoric(\Sig)$.
The universal property is proved similarly to the complex analytic case (\cite{KU3}, 3.6.1, 3.6.11).
\qed

\medskip

The object $S(\Sig)$ is called the {\it log modification} of $S$ associated to the subdivision $\Sig$ of the cone $\Hom(\cS, \bR_{\ge0}^\add)$. 
It is the real analytic version of the complex analytic log modification in the category $\cB(\log)$ in \cite{KU3}, 3.6.12.
\medskip

{\bf 3.1.13.} 
We use the notation in 3.1.12.
As a set, the log modification $S(\Sig)$ is identified with the set of all triples $(s, \sig, h)$, where $s\in S$, $\sig\in \Sig$, and if $P(\sig)$ denotes the image of $\cS(\sig)$ (3.1.9 for $N=\Hom(\cS^{\gp}, \bZ)$ and $M= \cS^{\gp}$) in $(M_S/\cO_S^\times)_s^{\gp}$ and $P'(\sig)$ denotes the inverse image of 
$P(\sig)$ in $M^{\gp}_{S,>0, s}$, then $h$ is a homomorphism $P'(\sig)^\times\to \bR_{>0}$, satisfying the following conditions (1) and (2). 
\medskip

(1) $P(\sig)^{\times} \cap (M_S/\cO_S^\times)_s=\{1\}$. 
\medskip

(2) The restriction of $h$ to $\cO_{S,>0,s}^\times$ ($\sub P'(\sigma)^\times$) is the evaluation map at $s$. 
\medskip

This is the real analytic version of the complex analytic theory \cite{KU3}, 3.6.15.

\vskip20pt

\head
\S3.2. Real analytic structures of $D_{\SL(2)}$ 
\endhead

{\bf 3.2.1.} 
We will define two structures on the set $D_{\SL(2)}$ as an object of $\cB_\bR(\log)$. We will denote $D_{\SL(2)}$ with these structures by $D_{\SL(2)}^I$ and $D_{\SL(2)}^{II}$. 
There is a morphism $D_{\SL(2)}^I\to D_{\SL(2)}^{II}$ whose underlying map is the identity map of  $D_{\SL(2)}$.  
The log structure with sign of $D_{\SL(2)}^I$ coincides with the inverse image (3.1.5) of that of $D_{\SL(2)}^{II}$. 

In the pure case, these two structures coincide, and the topology of $D_{\SL(2)}$ given by these structures coincides with the one defined in \cite{KU2}.

$D_{\SL(2)}^{II}$ is proper over $\spl(W)\times D_{\SL(2)}(\gr^W)$ (3.5.16 below). 
This shows that our definition of $D_{\SL(2)}$ in the mixed case provides sufficiently many points at infinity. 
This properness is  a good property of $D_{\SL(2)}^{II}$ which $D_{\SL(2)}^I$ need not have. 
On the other hand, $D_{\SL(2)}^I$ is nice for norm estimate (4.2.2 below), but $D_{\SL(2)}^{II}$ need not be. 

The sheaf of rings on $D_{\SL(2)}^I$ is called the {\it sheaf of real analytic functions} 
(or the {\it real analytic structure}) {\it on $D_{\SL(2)}$ in the first sense}, and that on $D_{\SL(2)}^{II}$ is called the {\it sheaf of real analytic functions} 
(or the {\it real analytic structure}) {\it 
on $D_{\SL(2)}$ in the second sense}.
The topology of $D_{\SL(2)}^I$ is called the {\it stronger topology of $D_{\SL(2)}$}, and that of $D_{\SL(2)}^{II}$ is called the {\it weaker topology of $D_{\SL(2)}$}. 
These two topologies often differ. 

In this \S3.2, we characterize the structures of $D_{\SL(2)}^I$ and $D_{\SL(2)}^{II}$ as objects of $\cB_{\bR}(\log)$ by certain nice properties of them 
(Theorem 3.2.10). The existences of such structures will be proved in \S3.3 and \S3.4. 
\medskip

{\bf 3.2.2.} 
We define sets $\cW$, $\overline{\cW}$, a subset $D_{\SL(2)}^I(\Psi)$ of $D_{\SL(2)}$ for $\Psi\in \cW$, and a subset $D_{\SL(2)}^{II}(\Phi)$ of $D_{\SL(2)}$ for $\Phi\in \overline{\cW}$, as follows. 

For $p\in D_{\SL(2)}$, let $\cW(p)$ be the set of weight filtrations associated to $p$. 

By an {\it admissible set of weight filtrations on $H_{0,\bR}$}, we mean a finite set $\Psi$ of increasing filtrations on $H_{0, \bR}$ such that $\Psi=\cW(p)$ for some element $p$ of $D_{\SL(2)}$. 
We denote by $\cW$ the set of all admissible sets of weight filtrations on $H_{0,\bR}$. 

For $\Psi \in \cW$, we define a subset $D_{\SL(2)}^I(\Psi)$ of $D_{\SL(2)}$ by
$$
D_{\SL(2)}^I(\Psi)= \{p \in D_{\SL(2)}\;|\; \cW(p) \sub \Psi\}.
$$
Note that $D_{\SL(2)}$ is covered by the subsets $D_{\SL(2)}^I(\Psi)$  for $\Psi\in \cW$. Furthermore, $D_{\SL(2)}$ is covered by the subsets $D_{\SL(2)}^I(\Psi)$  for $\Psi\in \cW$ with $W\notin \Psi$ and the subsets 
${D_{\SL(2)}^I(\Psi)}{}_{\nspl}:=D_{\SL(2)}^I(\Psi)\cap D_{\SL(2),\nspl}$ for $\Psi\in \cW$ with $W\in \Psi$. 
As is stated in Theorem 3.2.10 below, these are open coverings of $D_{\SL(2)}$ for the topology of $D_{\SL(2)}^I$. 
  
For $p\in D_{\SL(2)}$, let 
$$
\overline{\cW}(p)=\{W'(\gr^W)\;|\; W'\in \cW(p),\, W'\neq W\},
$$ 
where $W'(\gr^W)$ is the filtration on $\gr^W=\bigoplus_w \gr^W_w$induced by $W'$, i.e., $W'(\gr^W)_k := \bigoplus_w W'_k(\gr^W_w)\subset \bigoplus_w \gr^W_w$. 

By an {\it admissible set of weight filtrations on $\gr^W$}, we mean a finite set $\Phi$ 
of increasing filtrations on $\gr^W$ such that $\Phi=\overline{\cW}(p)$ for some element $p$ of $D_{\SL(2)}$. 
 We denote by $\overline{\cW}$ the set of all admissible sets of weight filtrations on $\gr^W$. 

For $\Phi\in \overline{\cW}$, we define a subset $D_{\SL(2)}^{II}(\Phi)$ of $D_{\SL(2)}$ by
$$
D_{\SL(2)}^{II}(\Phi)= \{p\in D_{\SL(2)}\;|\;\overline{\cW}(p)\sub \Phi\}.
$$
As a set, $D_{\SL(2)}$ is covered by  $D_{\SL(2)}^{II}(\Phi)$ ($\Phi\in \overline{\cW}$). As is stated in Theorem 3.2.10 below, this is an open covering for the topology of $D_{\SL(2)}^{II}$.

We have a canonical map 
$$
\cW\to \overline{\cW}
$$ 
which sends $\Psi\in \cW$ to $\bar \Psi:=\{W'(\gr^W)\;|\;W'\in \Psi,\, W'\neq W\}\in \overline{\cW}$. 
For $\Psi\in \cW$, we have $D_{\SL(2)}^I(\Psi)\subset D_{\SL(2)}^{II}(\bar \Psi)$.

\medskip

{\bf 3.2.3.} 
Let $\Psi\in \cW$. 
A homomorphism $\alpha\;:\;\bG_{m,\bR}^{\Psi} \to \Aut_{\bR} (H_{0,\bR}, W)$ of algebraic groups over $\bR$ is called a {\it splitting of $\Psi$} if it satisfies the following conditions (1) and (2).

\medskip

(1) The corresponding direct sum decomposition 
$$
H_{0,\bR}=\tsize\bigoplus_{\mu\in X} \; S_\mu \quad
(X:= \bZ^\Psi)
$$
into eigen $\bR$-subspaces $S_\mu$ satisfies 
$$
W'_{w'}= \ts_{\mu\in X,\, \mu(W')\leq w'} \;S_\mu
$$
for all $W'\in \Psi$ and for all $w'\in \bZ$. 

\medskip

(2) For all $w\in \bZ$ and all $t \in \bG_{m,\bR}^\Psi$, $\iota(t)^{-w}\a_w(t)$ is contained in $G_{\bR}(\gr^W_w)$, 
where $\alpha_w\;:\;\bG_{m,\bR}^{\Psi} \to \Aut_{\bR} (\gr^W_w)$ is the 
homomorphism induced by $\a$, and 
$\iota$ is the composite of the multiplication 
$\bG_{m,\bR}^{\Psi} \to \bG_{m,\bR}$ and the canonical 
map $\bG_{m,\bR} \rightarrow \Aut_{\bR}(\gr_w^W)$, $a \mapsto \text{(multiplication by $a$)}$.

\medskip 

A splitting of $\Psi$ exists: 
If $\Psi$ is associated to $p\in D_{\SL(2)}$, the torus action $\tau_p$ associated to $p$ (2.5.6, 2.3.5) is a splitting of $\Psi$. 
  Here and hereafter, we identify 
$\{1, \ldots, n\}$ ($n$ is the rank of $p$) with $\Psi$ via the bijection 
$j \mapsto W^{(j)}$, which is independent of 
the choice of $p$ by 2.3.8. 

Let $\Phi\in \overline{\cW}$.
A homomorphism $\alpha\;:\;\bG_{m,\bR}^{\Phi} \to \prod_w \Aut_{\bR}
(\gr^W_w)$ of algebraic groups over $\bR$ is called a {\it splitting of $\Phi$} if it satisfies the  following conditions ($\bar 1$) and ($\bar 2$).

\medskip

($\bar 1$) The corresponding direct sum decomposition 
$$
\gr^W=\tsize\bigoplus_{\mu\in X} \; S_\mu \quad
(X:= \bZ^\Phi)
$$
into eigen $\bR$-subspaces $S_\mu$ satisfies 
$$
W'_{w'}= \ts_{\mu\in X,\, \mu(W')\leq w'} \;S_\mu
$$
for all  $W'\in \Phi$ and for all $w'\in \bZ$. 

\medskip

($\bar 2$) For all $w\in \bZ$ and all $t \in \bG_{m,\bR}^\Phi$, $\iota(t)^{-w}\a_w(t)$ is contained in $G_{\bR}(\gr^W_w)$, 
where $\alpha_w\;:\;\bG_{m,\bR}^{\Phi} \to \Aut_{\bR} (\gr^W_w)$ is the 
$w$-component of $\a$, and 
$\iota$ is the composite of the multiplication 
$\bG_{m,\bR}^{\Phi} \to \bG_{m,\bR}$ and the canonical 
map $\bG_{m,\bR} \rightarrow \Aut_{\bR}(\gr_w^W)$. 

\medskip 

A splitting of $\Phi$ exists: 
For $p\in D_{\SL(2)}$, let $\bar \tau_p$ be $\gr^W(\tau_p)$ in the case $W\notin \cW(p)$, and in the case $W\in \cW(p)$, let $\bar \tau_p$ be the restriction of $\gr^W(\tau_p)$ to $\bG_{m,\bR}^{\overline{\cW}(p)}$ which we identify with the part of $\bG_{m,\bR}^{\cW(p)}$ with the $W$-component removed. 
Then if $\Phi=\overline{\cW}(p)$, $\bar \tau_p$ is a splitting of $\Phi$. 

\medskip

{\it Remark.}
  Under the condition $(\bar 1)$, the condition $(\bar2)$ is equivalent to 
the following condition: 
For all $w\in \bZ$, the direct sum decomposition 
$$\gr^W_w=\tsize\bigoplus_{\mu\in X} \; S_{w,\mu}$$
corresponding to $\a_w$ satisfies 
$$
\langle S_{w,\mu}, S_{w,\mu'} \rangle =0
$$
unless $\mu + \mu' = (2w,\ldots,2w)$.

\medskip

{\bf 3.2.4.}
Let $\Psi\in \cW$. 
Assume $W\notin \Psi$ (resp. $W\in \Psi$).
If a real analytic map $\beta : D \to \bR_{>0}^{\Psi}\; (\text{resp.}\; D_{\nspl} \to \bR_{>0}^{\Psi}$) satisfies the following (1) for any splitting $\a$ of $\Psi$, then we call $\beta$ a {\it distance to $\Psi$-boundary}.
$$
\beta(\alpha(t)p)=t\beta(p)\quad (t \in \bR_{>0}^{\Psi}, \ p \in D \; (\text{resp}.\;  D_{\nspl})).\tag1
$$

Let $\Phi\in \overline{\cW}$. 
If a real analytic map $\beta : D(\gr^W) \to \bR_{>0}^{\Phi}$ satisfies the following ($\bar 1$) for any splitting $\a$ of $\Phi$, then we call $\beta$ a {\it distance to $\Phi$-boundary}.
$$
\beta(\alpha(t)p)=t\beta(p)\quad (t \in \bR_{>0}^{\Phi}, \ p \in D(\gr^W)).\tag
{$\bar 1$}
$$

The proofs of the following propositions 3.2.5--3.2.7 and 3.2.9  will be given in \S3.3.

\proclaim{Proposition 3.2.5}
{\rm (i)} Let $\Psi\in \cW$. 
Then a distance to $\Psi$-boundary exists. 
\medskip
 
{\rm (ii)} Let $\Phi\in \overline{\cW}$. 
Then a distance to $\Phi$-boundary exists. 
\endproclaim

\proclaim{Proposition 3.2.6} 
{\rm(i)} Let $\Psi\in \cW$, let $\a$ be a splitting of $\Psi$, and let $\b$ be a distance to $\Psi$-boundary. 
Assume $W\not\in \Psi$ $($resp. $W \in \Psi)$ and consider the map 
$$
\nu_{\a,\b}:D\ 
(resp.\  D_{\nspl})
\to \bR_{>0}^{\Psi}\times D
\times
\spl(W)
\times \tp_{W'\in \Psi} \spl(W'(\gr^W)),
$$ 
$$
p\mapsto 
(\b(p),\; \a\b(p)^{-1}p,\; \spl_W(p),\; (\spl^{\BS}_{W'(\gr^W)}(p(\gr^W)))_{W'\in \Psi}).
$$ 
Here $\spl_W(p)$ is the canonical splitting of $W$ associated to $p$ in $\S1.2$, and $\spl_{W'(\gr^W)}^{\BS}(p(\gr^W))$ is the Borel-Serre splitting of $W'(\gr^W)$ associated to $p(\gr^W)$ in $2.1.9$. Let $p\in D_{\SL(2)}^I(\Psi)$
$($resp. $D_{\SL(2)}^I(\Psi)_{\nspl})$,  $J$ the set of weight filtrations associated to $p$ $(2.3.6)$, $\tau_p:\bG_{m,\bR}^J\to \Aut_\bR(H_{0, \bR}, W)$ the associated torus action $(2.5.6, 2.3.5)$, and  $\br\in D$ a point on the torus orbit $(2.5.2)$ associated to $p$.
Then, when $t\in \bR_{>0}^J$ tends to $0^J$ in $\bR_{\geq 0}^J$, $\nu_{\a,\b}(\tau_p(t)\br)$ converges in $\bR_{\geq0}^{\Psi}\times D \times \spl(W) \times \tp_{W'\in \Psi} \spl(W'(\gr^W))$. 
  This limit depends only on $p$ and is independent of the choice of $\br$. 

\medskip

{\rm(ii)} Let $\Phi\in \overline{\cW}$, let $\a$ be a splitting of $\Phi$, and let $\b$ be a distance to $\Phi$-boundary. 
Consider the map 
$$
\nu_{\a,\b}:D\to \bR_{>0}^{\Phi}\times D(\gr^W)\times \cL \times \spl(W) \times \tp_{W'\in \Phi} \spl(W'),
$$ 
$$
p\mapsto (\b(p(\gr^W)), \;\a\b(p(\gr^W))^{-1}p(\gr^W), \;\Ad(\a\b(p(\gr^W)))^{-1}\delta(p), \;
$$ 
$$
\hskip140pt
\spl_W(p), \; (\spl^{\BS}_{W'}(p(\gr^W)))_{W'\in \Phi}).
$$
Here $\cL$ is in $1.2.1$ and $\delta(p)$ denotes $\delta$ of $p$. 
Let $p\in D_{\SL(2)}^{II}(\Phi)$,  $J$ the set of weight filtrations associated to $p$, $\tau_p:\bG_{m,\bR}^J\to \Aut_\bR(H_{0, \bR}, W)$ the associated torus action, and  $\br\in D$ a point on the torus orbit 
associated to $p$.  
Then, when $t\in \bR_{>0}^J$ tends to $0^J$ in $\bR_{\geq 0}^J$, 
$\nu_{\a,\b}(\tau_p(t)\br)$ converges in $\bR_{\geq0}^{\Phi}\times D(\gr^W)\times 
\bar \cL \times \spl(W) \times \tp_{W'\in \Phi} \spl(W')
$. 
  This limit depends only on $p$ and is independent of the choice of $\br$. 
\endproclaim

We recall the compactified vector space $\bar V$ associated to a weightened finite dimensional $\bR$-vector space $V=\bigoplus_{w\in \bZ} V_w$ such that $V_w=0$ unless $w\leq -1$. 
It is a compact real analytic manifold with boundary. 
For $t\in \bR_{>0}$ and $v=\sum_{w\in \bZ} v_w\neq 0$ ($v_w\in V_w$), let $t\circ v= \sum_w t^wv_w$. 
Then as a set, $\bar V$ is the disjoint union of $V$ and the points $0\circ v$ ($v\in V\smallsetminus \{0\})$, where $0\circ v$ is the limit point in $\bar V$ of $t\circ v$ with $t\in \bR_{>0}$, $t\to 0$. 
$0\circ v=0\circ v'$ if and only if $v'=t\circ v$ for some $t\in \bR_{>0}$. 

Since $\bar V$ is a real analytic manifold with boundary (a special case of a real analytic manifold with corners),  $\bar V$ is regarded as an object of $\cB(\log)$ (3.1.6). 

Since $\cL$ is a finite dimensional weightened $\bR$-vector space of weights $\leq -2$, we have the associated compactified vector space $\bar \cL\supset \cL$. 

\medskip

In Proposition 3.2.6, in both (i) and (ii), we denote the limit of $\nu_{\a, \b}(\tau_p(t)\br)$ by $\nu_{\a,\b}(p)$. 

\medskip

As we will see in 3.3.10, in 3.2.6 (ii), the $\bar \cL$-component of $\nu_{\a,\b}(p)$ belongs to $\cL$ (resp. $\bar \cL\setminus \cL$) if and only if $W\notin\cW(p)$ (resp. $W\in \cW(p)$).

\proclaim{Proposition 3.2.7} 
{\rm(i)} Let $\Psi\in \cW$, let $\a$ be a splitting of $\Psi$, and let $\b$ be a distance to $\Psi$-boundary. 
Then, in the case $W\not\in \Psi$ $($resp.  $W \in \Psi)$ the map 
$$
\multline
\nu_{\a,\b}:D_{\SL(2)}^I(\Psi)\ 
(resp.\  D_{\SL(2)}^I(\Psi)_{\nspl})
\\ \to \bR_{\geq0}^{\Psi}\times D
\times
\spl(W)
\times \tp_{W'\in \Psi} \spl(W'(\gr^W))
\endmultline
$$ 
is injective.
\medskip

{\rm(ii)} Let $\Phi\in \overline{\cW}$, let $\a$ be a splitting of $\Phi$, and let $\b$ be a distance to $\Phi$-boundary. 
Then the map 
$$
\nu_{\a,\b}:D_{\SL(2)}^{II}(\Phi)\to \bR_{\geq0}^{\Phi}\times D(\gr^W)\times \bar \cL \times \spl(W) \times \tp_{W'\in \Phi} \spl(W')
$$ 
is injective.
\endproclaim

\medskip

{\bf 3.2.8.}  
Here, for $\Psi\in \cW$, we define a structure of an object of $\cB_\bR'(\log)$ on the set $D_{\SL(2)}^I(\Psi)$ (resp. $D_{\SL(2)}^I(\Psi)_{\nspl}$) in the case $W\notin \Psi$ (resp. $W\in \Psi$), depending on choices of a splitting $\a$ of $\Psi$ and a distance to $\Psi$-boundary $\b$. 
Also, for $\Phi\in \overline{\cW}$, we define a structure of an object of $\cB_\bR'(\log)$ on the set $D_{\SL(2)}^{II}(\Phi)$  depending on choices of a splitting $\a$ of $\Phi$ and a distance to $\Phi$-boundary $\b$.
 
 Let $\Psi\in \cW$.  
Assume $W\notin \Psi$ $($resp. $W\in \Psi)$. Let $A = D_{\SL(2)}^I(\Psi)$ $($resp. $A = D_{\SL(2)}^I(\Psi)_{\nspl})$, let $B=\bR_{\geq0}^{\Psi}\times D \times \spl(W) \times \tp_{W'\in \Psi} \spl(W'(\gr^W))$, and regard $B$ as an object of $\cB_\bR(\log)$.  Define the topology of $A$ to be the one 
 as a subspace of $B$ in which $A$ is embedded by $\nu_{\a,\b}$ in {\rm 3.2.7 (i)}. We define the sheaf of real analytic functions on $A$ as follows.
For an open set $U$ of $A$ and a function $f : U\to \bR$,  we say $f$ is a real analytic function if and only if, for each $p\in U$, there are an open neighborhood $U'$ of
$p$ in $U$, an open neighborhood $U''$ of $U'$ in $B$,  and a real analytic function $g : U''\to \bR$ such that the restrictions to $U'$ of $f$ and $g$ coincide. Then $A$ belongs to $\cB'_\bR$. Define the log structure with sign on $A$ as  the inverse image $(3.1.5)$ of the log structure with sign of $B$. 
 
Let $\Phi\in \overline{\cW}$. 
Let $A = D_{\SL(2)}^{II}(\Phi)$, let $B=\bR_{\geq0}^{\Phi}\times D(\gr^W) \times \bar\cL\times \spl(W) \times \tp_{W'\in \Phi} \spl(W')$, and regard $B$ as an object of $\cB_\bR(\log)$.
Define the topology of $A$ to be the one 
 as a subspace of $B$ in which $A$ is embedded by $\nu_{\a,\b}$ in {\rm 3.2.7 (ii)}. We define the sheaf of real analytic functions on $A$ as follows.
For an open set $U$ of $A$ and a function $f : U\to \bR$,  we say $f$ is a real analytic function if and only if for each $p\in U$, there are an open neighborhood $U'$ of
$p$ in $U$, an open neighborhood $U''$ of $U'$ in $B$,  and a real analytic function $g : U''\to \bR$ such that the restrictions to $U'$ of $f$ and $g$ coincide. Then $A$ belongs to $\cB'_\bR$. Define the log structure with sign on $A$ as  the inverse image of the log structure with sign of $B$.

\proclaim{Proposition 3.2.9} {\rm(i)} 
Let $\Psi\in \cW$. Assume $W\notin \Psi$ $($resp. $W\in \Psi)$. 
Then the structure of an object of $\cB_\bR'(\log)$ on $D_{\SL(2)}^I(\Psi)$ $($resp. $D_{\SL(2)}^I(\Psi)_{\nspl})$ in $3.2.8$ is independent of the choices of $\a$ and $\b$.

\medskip

{\rm (ii)} Let $\Phi\in \overline{\cW}$. 
Then the structure of an object of $\cB_\bR'(\log)$ on $D_{\SL(2)}^{II}(\Phi)$  in $3.2.8$ is independent of the choices of $\a$ and $\b$.

\endproclaim

The following theorem will be proved in \S3.4. 

\proclaim{Theorem 3.2.10} 
{\rm(i)} There exists a unique structure $D_{\SL(2)}^I$ of an object of $\cB_\bR(\log)$ on the set $D_{\SL(2)}$ having the following property{\rm:} 
For any $\Psi\in \cW$, $D_{\SL(2)}^I(\Psi)$ and $D_{\SL(2)}^I(\Psi)_{\nspl}$ are open in $D_{\SL(2)}^I$, and if $W\notin\Psi$ $($resp. $W\in \Psi)$, the induced structure on $D_{\SL(2)}^I(\Psi)$ $($resp. $D_{\SL(2)}^I(\Psi)_{\nspl})$ coincides with the structure in $3.2.9$ as objects of $\cB'_\bR(\log)$. 
 
\medskip

{\rm(ii)} There exists a unique structure $D_{\SL(2)}^{II}$ of an object of $\cB_\bR(\log)$ on the set $D_{\SL(2)}$ having the following property{\rm:} 
For any $\Phi\in \overline{\cW}$, $D_{\SL(2)}^{II}(\Phi)$ is open in $D_{\SL(2)}^{II}$, and the induced structure on $D_{\SL(2)}(\Phi)$ coincides with the structure in $3.2.9$ as objects of $\cB'_\bR(\log)$. 
 
 \medskip

{\rm (iii)} The topology of $D_{\SL(2)}^{II}$ is coarser than or equal to 
that of $D_{\SL(2)}^{I}$, and 
the sheaf of real analytic functions on $D_{\SL(2)}^{II}$ is contained in the sheaf of real analytic functions on $D_{\SL(2)}^{I}$. 
Thus we have a morphism $D_{\SL(2)}^I\to D_{\SL(2)}^{II}$ of local ringed spaces over $\bR$. 
The log structure with sign on $D_{\SL(2)}^I$ coincides with the inverse image of that of $D_{\SL(2)}^{II}$. 
Thus we have a morphism $D_{\SL(2)}^I\to D_{\SL(2)}^{II}$ in $\cB_\bR(\log)$ whose underlying map of sets is the identity map of $D_{\SL(2)}$. 
In the pure case $($i.e., in the case $W_w=H_{0,\bR}$ and $W_{w-1}=0$ for some $w\in \bZ)$, the last morphism is an isomorphism, and the topology of $D_{\SL(2)}$ given by these structures coincides with the one defined in \cite{KU2}.

\endproclaim

{\bf 3.2.11.} 
In 3.2.12 below, we give characterizations of the topologies of $D_{\SL(2)}^I$ and $D_{\SL(2)}^{II}$. 
Recall [Bn], Ch.1, \S8, no.4, that a topological space $X$ is said to be {\it regular}
if it is Hausdorff and if for any point $x$ of $X$ and any neighborhood $U$ of $x$, there is a closed neighborhood of $x$ contained in $U$. 

Recall [Bn], Ch.1, \S8, no.5 that the topology of a regular topological space $X$ is determined by the restrictions of neighborhoods of each point to a dense subset $X'$ of $X$. 
Precisely speaking, if $T_1$ and $T_2$ are topologies on a set $X$ and if $X'$ is a subset of $X$, then $T_1$ and $T_2$ coincide if the following conditions (1) and (2) are satisfied. 
\medskip

(1) The space $X$ is regular for $T_1$ and also for $T_2$, and the subset $X'$
is dense in $X$ for $T_1$ and also for $T_2$.
\medskip

(2) Let $x\in X$, and for $j=1, 2$, let $S_j$ be the set  $\{X'\cap U\;|\;U \;\text{is a neighborhood of $x$ in}$ $\text{$X$ for $T_j$}\}$ of subsets of $X'$. 
Then $S_1=S_2$. 
\medskip

This condition (2) is equivalent to the following condition (2$'$).
\medskip

(2$'$) For any $x\in X$ and for any directed family $(x_\lam)_\lam$ of elements of $X'$, $(x_\lam)_\lam$ converges to $x$ for $T_1$ if and only if it converges to $x$ for $T_2$.
\medskip

The topologies of $D_{\SL(2)}^I$ and that of $D_{\SL(2)}^{II}$ have the following characterizations. 

\proclaim{Proposition 3.2.12} 
{\rm (i)} The topology of $D_{\SL(2)}^{I}$ is the unique topology which satisfies the following conditions $(1)$ and $(2)$.
\medskip

$(1)$ For any admissible set $\Psi$ of weight filtrations on $H_{0,\bR}$, $D_{\SL(2)}^I(\Psi)$ $(3.2.2)$ is open and regular, and $D$ is dense in it.
\medskip

$(2)$ For any $p \in D_{\SL(2)}$ and for any family $(p_\lam)_{\lam\in \Lam}$ of points of $D$ with a directed ordered set $\Lam$, $(p_\lam)$ converges to $p$ 
in $D_{\SL(2)}^I$ if and only if the following {\rm(a)}, {\rm(b)}, and 
{\rm(c.I)} are satisfied.
Let $n$ be the rank of $p$ $(2.5.1$, $2.3.2$--$2.3.3)$, let $((\rho_w, \vf_w)_w, \br)$ be an $\SL(2)$-orbit in $n$ variables which represents $p$, let $\Psi = \cW(p)$, and let $\tau: \bG^\Psi_m\to \Aut_\bR(H_{0,\bR}, W)$ be the homomorphism of algebraic groups associated to $p$ $(2.3.5)$. 
\medskip

{\rm(a)} The canonical splitting of $W$ associated to $p_\lam$ converges to the canonical splitting of $W$ associated to $\br$.
\medskip

{\rm(b)} For each $1\leq j \leq n$ and $w\in \bZ$, the Borel-Serre splitting $\spl^{\BS}_{W^{(j)}(\gr^W_w)}(p_\lam(\gr^W_w))$ of $W^{(j)}(\gr^W_w)$ at $p_\lam(\gr^W_w)$ {\rm (2.1.9)} converges to the Borel-Serre splitting of $W^{(j)}(\gr^W_w)$ at $\br(\gr^W_w)$.
\medskip

{\rm(c.I)} There is a family $(t_\lam)_{\lam\in \Lam}$  of elements of $\bR^n_{>0}$  such that $t_\lam\to \bold0$ in $\bR^n_{\ge0}$ and such that $\tau(t_\lam)^{-1}p_\lam\to \br$.  
\medskip

{\rm (ii)} The topology of $D_{\SL(2)}^{II}$ is the unique topology which satisfies the following conditions $(1)$ and $(2)$.
\medskip

$(1)$ For any admissible set $\Phi$ of weight filtrations on $\gr^W$, $D_{\SL(2)}^{II}(\Phi)$ $(3.2.2)$ is open and regular, and $D$ is dense in it.
\medskip

$(2)$ For any $p \in D_{\SL(2)}$ and for any family $(p_\lam)_{\lam\in \Lam}$ of points of $D$ with a directed ordered set $\Lam$, $(p_\lam)$ converges to $p$ 
in $D_{\SL(2)}^{II}$ if and only if {\rm(a)} and {\rm(b)} in {\rm(i)} and the following {\rm(c.II)} are satisfied.
Let $n$, $((\rho_w, \vf_w)_w, \br)$, $\Psi$ and $\tau$ be as in {\rm(2)} of  {\rm (i)}.
Let $\Phi = \overline\cW(p) = \bar \Psi$. 
\medskip

{\rm(c.II)} There is a family $(t_\lam)_{\lam\in \Lam}$  of elements of $\bR_{>0}^\Phi \sub \bR_{>0}^\Psi$ such that $t_\lam\to \bold0$ in $\bR_{\ge0}^\Phi \sub \bR_{\ge0}^\Psi$ and such that 
$(\tau(t_\lam)^{-1}p_\lam)(\gr^W)@>>> \br(\gr^W)$ and $\delta(\tau(t_\lam)^{-1}p_\lam) @>>>\delta(\br)$.
\endproclaim

The proof of 3.2.12 is given in \S3.4. 

\medskip

{\bf 3.2.13.} 
{\bf Example 0.} 
Consider the pure case Example 0 in 1.1.1. 
Let $\Psi=\{W'\}$, 
where $W'_{-3}=0\sub W'_{-2}=W'_{-1}=\bR e_1\sub W'_0=H_{0,\bR}$. 
Then we have a splitting $\a$ of $\Psi$ defined by $\a(t)e_1=t^{-2}e_1$, $\a(t)e_2=e_2$, and we have a distance $\b$ to $\Psi$-boundary defined by $\b(x+iy)=y^{-1/2}$ ($x+iy\in \fh=D,\; x, y\in \bR, y>0$). 
Then the map
$$
\nu_{\a,\b}\,:\, D\to \bR_{>0}\times D\times \spl(W'),\quad
p\mapsto (\b(p), \;\a\b(p)^{-1}p,\; \spl_{W'}^{\BS}(p)),
$$
is described as 
$$
x+iy\mapsto (y^{-1/2},\; \tfrac{x}{y}+i,\; x)\quad (x, y\in \bR, y>0),
$$
where we identify $\spl(W')$ with $\bR$ in the standard way.  
We can identify $D_{\SL(2)}^I(\Psi)$ with $\{x+iy\;|\;x, y\in \bR,\, 0<y\leq \infty\}$ (see 3.6.1 below). 
The extended map $\nu_{\a,\b}: D_{\SL(2)}^I(\Psi)\to \bR_{\geq 0}\times D\times \spl(W')$ sends $x+i\infty$ to $(0, i, x)$. 

\vskip20pt

\head
\S3.3. Proofs of propositions 3.2.5--3.2.7 and 3.2.9
\endhead

{\bf 3.3.1.} Let $\overline{\cW}$ be as in 3.2.2. For each $w\in \bZ$, let 
$\cW(\gr^W_w)$ be the set of all admissible sets of weight filtrations on $\gr^W_w$.  We have a canonical map
$$\overline{\cW}\to \cW(\gr^W_w),\quad \Phi\mapsto \{W'(\gr^W_w)\;|\; W'\in \Phi, \;W'(\gr^W_w)\neq W(\gr^W_w)\}.$$
This map sends $\overline{\cW}(p)$ for $p\in D_{\SL(2)}$ to $\cW(p(\gr^W_w))$. 

For $\Phi\in \overline{\cW}$ and $w\in \bZ$, let $\Phi(w)\in \cW(\gr^W_w)$ 
be the image of $\Phi$ under the above map. 

We will sometimes denote elements of $\Phi$ and elements of $\Phi(w)$ by the small letters $j, k$ etc.

Note that $\Phi$ is a totally ordered 
set by 2.3.8 (for $j, k\in \Phi$, $j\leq k$ means 
$\sig^2(j)\leq \sig^2(k)$), and $\{W(\gr^W_w)\}\cup \Phi(w)$ is also a totally ordered set by 2.1.13 with respect to $\sig^2$ 
(note $W(\gr^W_w)\leq j$ for any $j\in \Phi(w)$). 
The canonical map $\Phi \to \{W(\gr^W_w)\}\cup 
\Phi(w),\;W'\mapsto W'(\gr^W_w)$, preserves the ordering. 
\medskip

\proclaim{Lemma 3.3.2} 
We use the notation in $3.3.1$.
\medskip 

{\rm (i)} For $\Phi\in \overline{\cW}$ and $w\in \bZ$, the map 
$\Phi\to \prod_{w\in \bZ} (\{W(\gr^W_w)\}\cup \Phi(w))$, $W'\mapsto (W'(\gr^W_w))_{w\in \bZ}$, is injective. 

By this injection, we identify $\Phi$ and its image, and denote the latter also by $\Phi$.
\medskip 

{\rm (ii)} We have the bijection from $\overline{\cW}$ onto the set of pairs $(\Phi', (\Phi'(w))_{w\in \bZ})$, 
where $\Phi'(w)$ is an element of $\cW(\gr^W_w)$ for each $w \in \bZ$ and 
$\Phi'$ is a subset of $\prod_{w\in \bZ} (\{W(\gr^W_w)\}\cup \Phi'(w))$ satisfying the conditions $(1)$--$(3)$ below. 
  The bijection sends $\Phi \in \overline{\cW}$ to $(\Phi, (\Phi(w))_{w \in \bZ})$. 

\medskip

$(1)$ For each $w\in \bZ$, the image of the projection $\Phi'\to \{W(\gr^W_w)\} \cup \Phi'(w)$, which we denote by $j\mapsto j(w)$, contains $\Phi'(w)$.

\medskip

$(2)$ For each $j\in \Phi'$, there is  $w\in \bZ$ such that $j(w)\in \{W(\gr^W_w)\}\cup \Phi'(w)$ belongs to $\Phi'(w)$. 
\medskip
$(3)$ For any $j,k\in \Phi'$, one of the following $(a)$, $(b)$ holds. 

\medskip
$(a)$ $j(w)\leq k(w)$ for all $w\in \bZ$. 

\medskip
$(b)$ $j(w)\geq k(w)$ for all $w\in \bZ$. 

\endproclaim

{\it Proof.} 
(i) is clear. 

We prove (ii).  
The injectivity of the map $\Phi\mapsto (\Phi, (\Phi(w))_w)$ follows from (i). 
We prove the surjectivity. 
Let $(\Phi', (\Phi'(w))_w)$ be a pair satisfying (1)--(3). 
For $w\in \bZ$, let $n(w)$ be the cardinality 
of $\Phi'(w)$, let $(\rho'_w, \varphi'_w)$ be an $\SL(2)$-orbit on $\gr^W_w$ in $n(w)$ variables of rank $n(w)$ whose associated set of weight filtrations is $\Phi'(w)$, and let $\br(w) \in D(\gr^W_w)$ be a point on the torus orbit associated to $(\rho'_w,\vf_w')$. 
Take a point $\br$ of $D_{\SL(2)}$ such that $\br(\gr^W)=(\br(w))_w$. 
Let $n$ be the cardinality of $\Phi'$, write $\Phi'=\{\phi_1,\dots, \phi_n\}$ ($\phi_1(w)\leq \dots \leq \phi_n(w)$ for all $w\in \bZ$), write 
$\Phi'(w)=\{\phi_{w,1},\dots, \phi_{w,n(w)}\}$ ($\phi_{w,1}<\dots<\phi_{w,n(w)}$), and let $e_w:\{1,\dots, n(w)\}\to \{1,\dots, n\}$ be the injection defined by $e_w(k)=\min\{j\;|\;\phi_j(w)=\phi_{w,k}\}$. 
Let $p\in D_{\SL(2)}$ be the class of the $\SL(2)$-orbit $((\rho_w, \vf_w)_w, \br)$ in $n$ variables of rank $n$, where 
$$
\rho_w(g_1,\dots, g_n)=\rho'_w(g_{e_w(1)},\dots, g_{e_w(n(w))}), \quad 
\vf_w(z_1,\dots, z_n)=\vf'_w(z_{e_w(1)},\dots, z_{e_w(n(w))}).
$$ 
Then the pair $(\Phi', (\Phi'(w))_w)$ is the image of $\cW(p)\in \overline{\cW}$. 
\qed
\medskip

\proclaim{Lemma 3.3.3} 
Let $\Phi\in \overline{\cW}$ and let $(\Phi(w))_w$ be the image of $\Phi$ in $\prod_w \cW(\gr^W_w)$. 
Then there is a bijection between the set of all splittings of $\Phi$ and the set of all families $(\a_w)_{w\in\bZ}$, where $\a_w$ is a splitting of $\Phi(w)$ for each $w$. 
This bijection sends a splitting 
$\a$ of $\Phi$ to the following family $(\a_w)_w$. 
For $w\in \bZ$, let $e_w: \Phi(w)\to \Phi$ be the map defined by $e_w(k)=\min\{j\in \Phi\;|\;j(\gr^W_w)=k\}$. 
Then $\a_w$ is the composite $\bG_{m,\bR}^{\Phi(w)} \to \bG_{m,\bR}^{\Phi}\to \Aut(\gr^W_w)$, where the first arrow is induced from $e_w$ and the second arrow is given by $\a$. 

\endproclaim

{\it Proof.} 
From a family $(\a_w)_w$ of splittings $\a_w$ of $\Phi(w)$, the corresponding splitting $\a$ of $\Phi$ is recovered as follows. 
For $w\in \bZ$, let $R(w)=\{W(\gr^W_w)\}\cup \Phi(w)$. 
Let $\bG_{m,\bR}^{\Phi}\to \bG_{m,\bR}^{R(w)}= \bG_{m,\bR}\times \bG_{m,\bR}^{\Phi(w)}$ be the homomorphism induced by the map $\Phi\to R(w)$, $W'\mapsto W'(\gr^W_w)$. 
Then the action of $\bG_m^{\Phi}$ on $\gr^W_w$ by $\a$ is defined to be the composite 
$\bG_{m,\bR}^{\Phi}\to \bG_{m,\bR}\times \bG_{m,\bR}^{\Phi(w)}\to 
\Aut(\gr^W_w)$, where the last arrow is $(t, t')\mapsto t^w\a_w(t')$.
\qed

\medskip

\proclaim{Lemma 3.3.4} 
Let $\Phi\in \overline{\cW}$. 
For each $w\in \bZ$, let $\b_w:D(\gr^W_w)\to \bR_{>0}^{\Phi(w)}$ be a distance  to $\Phi(w)$-boundary. 
Let $h:\bZ^{\Phi}\to \prod_{w\in \bZ} \bZ^{\Phi(w)}$ be an injective 
homomorphism induced by the map $\Phi \to \prod_{w\in \bZ}(\{W(\gr^W_w)\}\cup \Phi(w))$.
Then there is a homomorphism $h':\prod_{w\in \bZ} \bZ^{\Phi(w)}\to  \bZ^{\Phi}$ such that the composite $\bZ^{\Phi} @>h>> \prod_{w\in \bZ} \bZ^{\Phi(w)} @>{h'}>> \bZ^{\Phi}$ is the identity map, and, for such an $h'$, the composite $D(\gr^W) \to \prod_{w\in \bZ} \bR_{>0}^{\Phi(w)} \to \bR_{>0}^{\Phi}$, where the first arrow is $(\b_w)_w$ and the second arrow is induced by $h'$, is a distance to $\Phi$-boundary.
\endproclaim 

{\it Proof.} 
Since the cokernel of $h$ is torsion free, there is such an $h'$. 
The rest follows from 3.3.3. \qed

\medskip

\proclaim{Lemma 3.3.5} Let $\Psi\in \cW$ and let $\Phi\in \overline{\cW}$ be the image of $\Psi$ under the canonical map $\cW\to \overline{\cW}$ 
$(3.2.2)$. Let $\b: D(\gr^W) \to \bR_{\geq 0}^{\Phi}$ be a distance to $\Phi$-boundary. 

\medskip

{\rm (i)} Assume $W\notin \Psi$.
Then the map $$D\to D(\gr^W)@>\b>> \bR_{> 0}^{\Phi}\simeq \bR_{> 0}^{\Psi}, \quad x\mapsto \b(x(\gr^W)),$$ is a distance to $\Psi$-boundary, where the last isomorphism is induced from the canonical bijection $\Psi\to \Phi,\;W'\mapsto W'(\gr^W)$.

\medskip

{\rm (ii)} Assume $W\in \Psi$. Let $\g: D_{\nspl}\to \bR_{>0}$ be a real analytic map such that $\g(\a(t)x)=t_W\g(x)$ for any $t\in \bR_{>0}^{\Psi}$ and $x\in D_{\nspl}$, where
$t_W$ denotes the $W$-component of $t$.
Then the map
$$D_{\nspl} \to \bR_{>0}\times 
\bR_{>0}^{\Phi}\simeq \bR_{>0}^{\Psi}, \quad x\mapsto (\g(x),\; \b(x(\gr^W))),$$ is a distance to $\Psi$-boundary.

\endproclaim

This is proved easily. 
\medskip

{\bf 3.3.6.} We prove Proposition 3.2.5 (the existence of $\b$). 

\medskip

{\it Proof.} Assume first that 
we are in the pure case. In this case, the existence of $\b$ is proved in 
 \cite{KU2}, 4.12. 

In fact, there is a mistake in \cite{KU2}, for, loc.\ cit.\ 4.12 does not hold for a general compatible family of $\bQ$-rational increasing filtrations in the sense of \cite{KU2}. 
The proof for 4.12 there assumed the injectivity of the splitting (denoted $\nu$ there), but, for a general compatible family, a splitting is not necessarily injective. 
On the other hand, for an admissible set of weight filtrations, any splitting 
is injective, and for such 
a family, the proof there is correct and hence the conclusion of \cite{KU2}, 4.12 holds.

\medskip
The existence of a distance to $\Phi$-boundary $\b$ for $\Phi\in \overline{\cW}$ follows from the pure case by Lemma 3.3.4.
\medskip

We prove the existence of a distance to $\Psi$-boundary $\b$ for $\Psi\in \cW$.
 Let $\Phi$ be the image $\bar\Psi$ of $\Psi$ in $\overline{\cW}$ as in 3.2.2.

If $W\notin \Psi$, the existence  of $\b$ follows from Lemma 3.3.5 (i). 
Assume $W\in \Psi$. It is sufficient to construct 
a real analytic map $\gamma:D_{\nspl}\to \bR_{>0}$ having the property stated in 3.3.5 (ii). 
Fix $\bar \br=(\bar \br_w)_w\in D(\gr^W)$ and, for each $w\leq -1$, fix a $K'_{\bar \br_w}$-invariant positive definite symmetric $\bR$-bilinear form $(\;,\;)_w$ on the 
component $L_w$ of $L:=\cL(\bar \br)$ (1.2.1) of weight $w$. 
  Here $K'_{\bar \br_w}$ is the isotropy subgroup of $G_{\bR}(\gr^W_w)$ at $\bar 
\br_w$, which is compact so that there is such a form. 
Let $f : L-\{0\} \to \bR_{>0}, \;f(v) := (\ts_{w\leq -1} \;(v_w,v_w)_w^{-1/w})^{-1/2}$, where $v_w$ denotes the component of $v$ of weight $w$. 
For $F\in D(\gr^W)$, if $g$ is an element of $G_{\bR}(\gr^W)$ such that $F=g\bar \br$, then we have an isomorphism $\Ad(g)^{-1}:\cL(F)@>\sim>> L$.  
The map $f_F : \cL(F)-\{0\}\to \bR_{>0}, \;v\mapsto f(\Ad(g)^{-1}v)$, is independent of the choice of $g$. 
This is because $(g')^{-1}g\in \prod_w K'_{\bar \br_w}$ if $g, g'\in G_\bR(\gr^W)$ and $g\bar \br = g'\bar \br$. 
Define $\g': D_{\nspl}\to\bR_{>0}$ by $\g'(s(\theta(F, \d)))= f_F(\d)$. 
  Let $\a$ be any splitting of $\Psi$. 
  Then 
$\g'(\a(t)x)=(\tp_{W'\in \Psi} t_{W'})\g'(x)$ for $t\in \bR_{>0}^{\Psi}$ and $x\in D_{\nspl}$, where $t_{W'}\in \bR_{>0}$ denotes the $W'$-component of $t$. For $x\in D_{\nspl}$, define 
$\g(x) = \g'(x)\cdot\prod_{W'\in \Phi} \b(x(\gr^W))_{W'}^{-1}$, 
where $\b(x(\gr^W))_{W'}$ denotes the $W'$-component of $\b(x(\gr^W))$. Then $\g$ has the property stated in 3.3.5 (ii). 
\qed
\medskip 

{\bf 3.3.7.} 
We start to prove Proposition 3.2.6. 
The last assertions of (i) and (ii) are clear once the preceding convergences are shown.  
We will prove the convergences in 3.3.7--3.3.12. 

In this 3.3.7, we prove the following part of Proposition 
3.2.6 (i): 
\medskip

{\it Let $\Psi\in \cW$ and assume $W\notin \Psi$ $($resp. $W\in \Psi)$, let $\b$ be a distance to $\Psi$-boundary, $p\in D_{\SL(2)}^I(\Psi)$ $($resp. 
$D_{\SL(2)}^I(\Psi)_{\nspl})$,and let $\br \in D$ be a point on the torus orbit associated to $p$. 
Let $J$ be the set of weight filtrations associated to $p$. 
Then $\b(\tau_p(t)\br)$ $(t\in \bR_{>0}^J)$ converges in $\bR_{\geq 0}^{\Psi}$ when $t$ tends to $0^J$.}

\medskip

{\it Proof.}
Take a splitting $\a$ of $\Psi$ and 
let $\alpha_J:\bG_{m,\bR}^J\to \Aut(H_{0, \bR})$  be the restriction of $\alpha$ to the $J$-component $\bG_{m,\bR}^J$ of $\bG_{m,\bR}^{\Psi}$. 
Let $H_{0,\bR}=\bigoplus_{m\in \bZ^J} S(J,m)$ be the decomposition associated to $\a_J$. Since both $\tau_p$ and $\a_J$ split $J$, there is a unique element $u$ of $G_{\bR}$ such that $\tau_p=\Int(u)(\alpha_J)$ and such that $(1-u)S(J, m) \sub \bigoplus_{m'<m} S(J, m')$ for any $m\in \bZ^J$. We have
$$
\b(\tau_p(t)\br)=\b(u\a_J(t)u^{-1}\br)=\b(\a_J(t)u_tu^{-1}\br)=\iota_J(t)
\b(u_tu^{-1}\br),
$$
where $u_t=\Int(\a_J(t))^{-1}(u)$, and $\iota_J: \bR_{>0}^J \to \bR_{>0}^{\Psi}$ is the canonical injective homomorphism from the $J$-component. 
When $t\to 0^J$, $u_t$ converges to $1$ as is easily seen. 
Hence $\b(\tau_p(t)\br)$ converges to $0^J\b(u^{-1}\br)$ in $\bR_{\geq 0}^{\Psi}$, where $0^J$ denotes the element of $\bR_{\geq 0}^{\Psi}$ whose $j$-th component for $j\in \Psi$ is $0$ if $j\in J$ and is $1$ if $j\notin J$. 
\qed
\medskip

{\it Remark.} 
In \cite{KU2}, 4.12, the corresponding statement in the pure case was treated, but, the second line after the proof of it, the factor corresponding to \lq\lq$0^J$'' here is missing. 

\medskip

{\bf 3.3.8.} 
We prove the following part of Proposition 3.2.6 (ii): 
\medskip

{\it Let $\Phi\in\overline{\cW}$, $\b$  a distance to $\Phi$-boundary, $p\in D_{\SL(2)}^{II}(\Phi)$, and $\br$  a point on the torus orbit associated to $p$. Let $J$ be the set of weight filtrations associated to $p$. 
Then $\b(\tau_p(t)\br(\gr^W))$ $(t\in \bR_{>0}^J)$ converges in $\bR_{\geq 0}^{\Phi}$ when $t$ tends to $0^J$.}

\medskip

{\it Proof.}
Let $\bar J\in \overline{\cW}$ be the image of $J$, and let $\bar \tau_p$ be 
 as in 3.2.3. Take a splitting $\a$ of $\Phi$,
let $\alpha_{\bar J}:\bG_{m,\bR}^{\bar J}\to \Aut_{\bR}(\gr^W)$ be the restriction of $\alpha$ to the $\bar J$-component $\bG_{m,\bR}^{\bar J}$ of $\bG_{m,\bR}^{\Phi}$, and let $\gr^W=\bigoplus_{m\in \bZ^{\bar J}} \bar S(\bar J, m)$ be the decomposition associated to $\a_{\bar J}$.  Since both $\bar \tau_p$ and $\a_{\bar J}$ split $\bar J$, there is a unique element $u$ of $G_{\bR}(\gr^W)$ such that $\bar \tau_p=\Int(u)(\alpha_{\bar J})$ and such that $(1-u)\bar S(\bar J, m) \sub \bigoplus_{m'<m} \bar S(\bar J, m')$ for any $m\in \bZ^{\bar J}$. 
We have
$$\align
\b(\tau_p(t)\br(\gr^W))&=\b(u\a_{\bar J}(t_{\bar J})u^{-1}\br(\gr^W)) \\
&=\b(\a_{\bar J}(t_{\bar J})u_tu^{-1}\br(\gr^W))=
\iota_{\bar J}(t_{\bar J})\b(u_tu^{-1}\br(\gr^W)),
\endalign
$$
where $u_t=\Int(\a_{\bar J}(t_{\bar J}))^{-1}(u)$, $\iota_{\bar J}: \bR_{>0}^{\bar J} \to \bR_{>0}^{\Phi}$ is the canonical injective homomorphism from the $\bar J$-component, and $t_{\bar J}$ is the $\bar J$-component of $t$.
Here we identify $\bar J$ with $J$ (resp. $J\smallsetminus \{W\}$) if $W\notin J$ (resp. $W\in J$). 
When $t\to 0^J$,  $u_t$ converges to $1$ as is easily seen. 
Hence $\b(\tau_p(t)\br(\gr^W))$ converges to $0^{\bar J}\b(u^{-1}\br(\gr^W))$, where $0^{\bar J}$ denotes the element of $\bR_{\geq 0}^{\Phi}$ whose $j$-th component for $j\in \Phi$ is $0$ if $j\in \bar J$ and is $1$ if $j \not\in \bar J$.
\qed

\medskip

{\bf 3.3.9.} 
We prove the following part of Proposition 3.2.6 (i): 
\medskip

{\it Let the notation be as in $3.3.7$, let $\a$ be a splitting of $\Psi$, and let $\mu: D\to D$ be the map $x\mapsto \a\b(x)^{-1}x$. 
Then, $\mu(\tau_p(t)\br)$ converges in $D$ when $t\in \bR_{>0}^J$ tends to $0^J$ in $\bR_{\geq 0}^J$.}

\medskip

{\it Proof.}
We have
$$
\mu(\tau_p(t)\br)=\mu(u\a_J(t)u^{-1}\br)=\mu(\a_J(t)u_tu^{-1}\br)=\mu(u_tu^{-1}\br)\to \mu(u^{-1}\br),
$$
when $t\to 0^J$.
\qed
\medskip

{\bf 3.3.10.} 
We prove the following part of Proposition 3.2.6 (ii): 
\medskip

{\it Let the notation be as in $3.3.8$, let $\a$ be a splitting of $\Phi$, and let $\mu=(\mu_1, \mu_2): D\to D(\gr^W)\times \cL$ be the map $x\mapsto (\a\b(x(\gr^W))^{-1}x(\gr^W),\, \Ad(\a\b(x(\gr^W)))^{-1}\delta(x))$. 
Then, $\mu(\tau_p(t)\br)$ converges in $D(\gr^W)\times \bar \cL$ when $t\in \bR_{>0}^J$ tends to $0^J$ in $\bR_{\geq 0}^J$.}

\medskip

{\it Proof.}
This $\mu_1$ factors through the projection $D\to D(\gr^W)$ and 
$$\align
\mu_1(\tau_p(t)\br)&=\mu_1(u\a_{\bar J}(t_{\bar J})u^{-1}\br(\gr^W))
=\mu_1(\a_{\bar J}(t_{\bar J})u_tu^{-1}\br(\gr^W))\\ 
&=\mu_1(u_tu^{-1}\br(\gr^W))\to \mu(u^{-1}\br(\gr^W))
\endalign$$
when $t\to 0^J$. 
Assume $W\notin J$, and identify $J$ and $\bar J$ via the canonical bijection. 
Then  
$$\align
\mu_2(\tau_p(t)\br) &= (\Ad\a\b(\bar \tau_p(t)\br(\gr^W)))^{-1}\Ad(\bar \tau_p(t))\delta(\br) \\
&=\Ad(\a\b(u_tu^{-1}\br(\gr^W)))^{-1}\Ad(u_tu^{-1})\delta(\br) \\
&\to \Ad(\a\b(u^{-1}\br(\gr^W)))^{-1}\Ad(u^{-1})\delta(\br)
\endalign
$$
when $t\to 0^J$. Next assume $W\in J$ and identify $J\smallsetminus \{W\}$ with $\bar J$ via the canonical bijection. For $t\in \bR_{>0}^J$, write $t=(t', t_{\bar J})$, 
where $t'\in \bR_{>0}$ denotes the $W$-component of $t$ and $t_{\bar J}$ denotes the $\bar J$-component of $t$. Then 
$$\align
\mu_2(\tau_p(t)\br) &= \Ad(\a\b(\bar \tau_p(t)\br(\gr^W)))^{-1}(t'\circ \Ad(\bar \tau_p(t))\delta(\br))\\
&=t'\circ \Ad(\a\b(u_tu^{-1}\br(\gr^W)))^{-1}\Ad(u_tu^{-1})\delta(\br)\\
&\to 0\circ \Ad(\a\b(u^{-1}\br(\gr^W)))^{-1}\Ad(u^{-1})\delta(\br) 
\endalign
$$
when $t\to 0^J$.
Here, for $t \in \bR_{>0}$ and $\delta=\sum_{w \le -2}\delta_w \in 
\cL$, we write $t \circ \delta = \sum t^w \delta_w$, and 
$0 \circ \delta = \lim\limits_{t \to 0}t \circ \delta$ in $\bar \cL$. 
\qed
\medskip

{\bf 3.3.11.} 
Since the convergences of the canonical splittings are trivial (cf. 2.4.6, 2.5.5), to prove 3.2.6, the rest is the convergences of BS-splittings. 
To see the latter, we may and do assume that we are in the pure case. 

Let $\Psi$ be an admissible set of weight filtrations and let $W'\in \Psi$. 
Fix an $\SL(2)$-orbit $q$ whose associated set of weight filtrations is $\Psi$. 
Let $X=\bZ^{\Psi}$ and let $\fg_\bR=\bigoplus_{m\in X} \; \fg_{\bR,m}$ be the direct sum decomposition, where $t\in (\bR^\times)^{\Psi}$ acts via $\tau_q$ 
on $\fg_{\bR,m}$ as the multiplication by $t^m$. 

In this paragraph, we prove the following: 
\medskip

{\it Let $\br$ be a point on the torus orbit associated to $q$.
Let $J$ be a subset of $\Psi$ and let $\tau_J$ be the restriction of $\tau_q$ to the $J$-component $\bG_{m,\bR}^J$ of $\bG_{m,\bR}^{\Psi}$. 
Let $h \in \fg_{\bR}$ be an element whose $m$-component is zero $(m \in X)$ unless $m(j) < 0$ for all $j \in J$. 
Then, there are an open neighborhood $U$ of $0^J$ in $\bR_{\geq 0}^J$ and real analytic maps $f_1:U\to G_{W', \bR}$ and $f_2: U\to K_{\br}$ such that $\Int(\tau_{J}(t))^{-1}(\exp(h))= f_1(t)f_2(t)$ for any $t\in U\cap \bR_{>0}^J$, and, furthermore, $\Int(\tau_{J}(t))(f_1(t))$ extends to a real analytic map 
on $U$.}
\medskip

  To prove this, first, we take 
an $\bR$-subspace of $\fg_\bR$ satisfying the following (1)--(3).

\medskip

\noindent
(1) $\fg_\bR=V \oplus \Lie(K_\br)$.

\smallskip

\noindent
(2) $V$ is the sum of $V_{\pm m}:=V\cap (\fg_{\bR,m}+\fg_{\bR,-m})$ for $m\in X$. 

\smallskip  

\noindent
(3) $\Lie(G_{W',u,\bR})\sub V \sub \Lie(G_{W',\bR})$. 

\medskip

  Then 
there exist an open neighborhood $O$ of
$0$ in $\fg_\bR$ and a real analytic function
$a=(a_1, a_2):O\to 
V\oplus\Lie(K_{\br})$ having the following properties
{\rm(4)--(7)}.
\medskip

\noindent
{\rm(4)} For any $x\in O$,
$\exp(x)=\exp(a_1(x))\exp(a_2(x))$.
\smallskip

\noindent
{\rm(5)} $a(0)=(0,0)$.
\smallskip

\noindent
{\rm(6)} $\exp:O\to G_\bR$ is an injective open
map. 
\smallskip

\noindent
{\rm(7)} 
For $k=1,2$, $a_k$ has the form of
absolutely convergent series $a_k=\ts_{r=0}^\infty
a_{k,r}$, where $a_{k,r}$ is the part of degree $r$
in the Taylor expansion of $a_k$ at $0$, such that
$a_{k,r}(x)=l_{k,r}(x\ox\cdots\ox x)$ for some
linear map $l_{k,r}:\fg_\bR^{\ox r}\to\fg_\bR$
having the following property$:$
If $m_1, \dots, m_r\in X$ and $x_j\in\fg_{\bR,{m_j}}$ for $1\le j\le r$, then $l_{k,r}(x_1\ox\cdots\ox x_r)
\in\ts_m\fg_{\bR,m}$, where $m$ ranges over
all elements of $X$ satisfying 
$|m| \le |m_1|+\cdots+|m_r|$. 
  Here $|\quad |: \bZ^{\Psi} \to \bN^{\Psi}$ is the map sending 
$(m(j))_j$ to $(|m(j)|)_j$. 

\medskip

  This is proved similarly as \cite{KU3}, 10.3.4. 
  Or, if we choose $V$ such that 
$V=\Lie(\t\rho(\bR^n_{>0})) \oplus L$ for some $L$ 
as in \cite{KU3}, 10.1.2 (such a choice is always possible), this is 
seen by loc.\ cit.\ 10.3.4 just by taking $a_1(x)=H(f_1(x), f_2(x))$, 
$a_2(x)=f_3(x)$, where $H(x,y)=x+y+\frac 12[x,y]+\cdots$ is a Hausdorff 
series. 

  Now consider the decomposition $h=\sum_{m\in X} h_m$ $(h_m \in \fg_{\bR,m})$. 
  By assumption, $h_m=0$ unless $m(j) < 0$ for any $j \in J$. 
  Then, $\Ad(\tau_{J}(t))^{-1}(h)= \sum_{m\in X} t^{-m_J}h_m$ 
$(t\in \bR_{>0}^J)$ extends to a real analytic map $g:\bR_{\ge0}^J
\to \fg_{\bR}$ sending $0^J$ to $0$, where $m_J \in \bZ^J$ is the 
$J$-component of $m$. 
  Let $U=g^{-1}(O)$, $f_j=\exp\circ a_j \circ g$ $(j=1,2)$. 
  It is enough to show that $\Ad(\tau_{J}(t))(a_1(g(t)))$ extends to a real 
analytic map around $0$.
  This is a consequence of the property (7) of $a_1$.  
  In fact, in the notation in (7), 
$a_1(g(t))=a_1(\sum t^{-m_J}h_m)$ is the infinite formal sum of 
$t^{-((m_1)_J+\cdots +(m_r)_J)}l_{1,r}(h_{m_1}\otimes \cdots \otimes h_{m_r})$ 
$(m_j \in X, h_{m_j} \in \fg_{\bR, m_j}$ $(1 \le j \le r))$. 
  Since the weights $m$ of $l_{1,r}(h_{m_1}\otimes \cdots \otimes h_{m_r})$ 
satisfies $|m| \le |m_1|+\cdots+ |m_r|$, we conclude that 
$\Ad(\tau_{J}(t))(a_1(g(t)))$ extends to a real analytic map over $0^J$, as 
desired.

\medskip

{\bf 3.3.12.}
  We continue to assume that we are in the pure situation. 

Let $\Psi$ be an admissible set of weight filtrations and let $W'\in \Psi$. 
We prove the following, which completes the proof of Proposition 3.2.6:
\medskip

{\it Let $p\in D_{\SL(2)}(\Psi)$, and let $\br$ be a point on the torus orbit associated to $p$. 
Let $J$ be the set of weight filtrations associated to $p$. 
Then, $\spl^{\BS}_{W'}(\tau_p(t)\br)$ $(t\in \bR_{>0}^J)$ converges in 
$\spl(W')$ when $t$ tends to $0^J$.}

\medskip

{\it Remark 1.} 
The proof is easy when $W'\in J$ (BS-splitting is then constant on the torus orbit) but is not when $W'\notin J$. 
  See Remark 3 after the proof. 

\medskip

{\it Proof.}
 Since $\Psi$ is admissible, $\Psi$ is the set of weight filtrations associated to some $q\in D_{\SL(2)}$. Let $\br_q$ be a point on the torus orbit associated to $q$. 
Then, by Claim 1 in \cite{KU3} 6.4.4, there exist $v\in G_{J, \bR}$ and $k\in K_{\br_q}$ such that $\tau_p=\Int(v)(\tau_{q,J})$ and $\br=vk\br_q$. Here 
$G_{J,\bR}=\{g \in G_{\bR}\;|\; g W'' = W''$ for any $W'' \in J\}$, and 
$\tau_{q,J}$ denotes the restriction of $\tau_q$ to the $J$-component $\bG_{m,\bR}^J$ of $\bG_{m,\bR}^{\Psi}$. 

Let $G_{\bR}(J)$ be the $\bR$-algebraic subgroup of $G_{J,\bR}$ consisting of all elements of $G_\bR$ 
which commute with any element of $\tau_{q,J}(\bG_{m,\bR}^J)$. 
Then we have the projection $G_{J, \bR}\to G_\bR(J), \;a\mapsto a(J)$, 
where $a(J)$ on $S(J, m)$ ($m\in \bZ^J$) (3.3.7) 
is defined to be the ($S(J, m) \to S(J, m)$)-component of $a:S(J, m) \to \bigoplus_{m'\leq m} S(J, m')$. The composite $G_{\bR}(J)\to G_{J, \bR}\to G_\bR(J)$ 
is the identity map. 
  Since $G_{\bR}(J)$ is reductive, 
any element of $G_\bR(J)$ is expressed in the form $bc$, where $b\in G_{\bR}(J)\cap G_{W', \bR}$ and 
$c\in G_\bR(J)\cap K_{\br_q}$. 
Write the image of $v$ in $G_\bR(J)$ as $bc$ using such $b$ and $c$. 
Then $v=bv_uc$ with $v_u\in G_{J, \bR}$ satisfying
$(v_u-1)S(J, m)\sub \bigoplus_{m'<m} S(J, m)$ for any $m\in \bZ^J$.
We have  $\Int(\tau_{q,J}(t))^{-1}(v_u)\to 1$ 
when $t\to 0^J$ in $\bR_{\geq 0}^J$. Hence by 3.3.11, there are an open neighborhood $U$ of $0^J$ in $\bR_{\geq 0}^J$ and real analytic maps $b_u:U\to G_{W', \bR}$ and $c_u: U\to K_{\br_q}$ such that $\Int(\tau_{q,J}(t))^{-1}(v_u)= \Int(\tau_{q,J}(t))^{-1}(b_u(t))c_u(t)$ for any $t\in U\cap \bR_{>0}^J$. 
We have, for $t\in U\cap \bR_{>0}^J$, 
$$
\tau_p(t)\br=v\tau_{q,J}(t)k\br_q= bb_u(t)\tau_{q,J}(t)c_u(t)ck\br_q,
$$
and hence 
$$
\align
\spl_{W'}^{\BS}(\tau_p(t)\br)
&=\Int(bb_u(t))\Int(\tau_{q,J}(t))(\spl_{W'}^{\BS}(c_u(t)ck\br_q))\\
&=\Int(bb_u(t))\Int(\tau_{q,J}(t))(\spl_{W'}^{\BS}(\br_q))
= \Int(bb_u(t))(\spl_{W'}^{\BS}(\br_q))\\
&\to \Int(bb_u(0^J))(\spl_{W'}^{\BS}(\br_q)).
\qed
\endalign
$$ 
\medskip

{\it Remark 2.} 
In the above proof, $\spl_{W'}^{\BS}(\br_q)$ coincides with the splitting of $W'$ associated to $q$.
\medskip

{\it Remark 3.} 
In the case $W'\in J$, $\spl_{W'}^{\BS}(\tau_p(t)\br)$ 
constantly coincides with
$\Int(v)\spl_{W'}^{\BS}(\br_q)$ with $v$ as in the above proof.

\medskip

{\bf 3.3.13.}  
We prove Proposition 3.2.7 (injectivity of $\nu_{\a,\b}$). Recall that a point of $D_{\SL(2)}$ is determined by the associated weight filtrations and the associated torus orbit (2.5.2 (ii)). 

First, let $\Psi\in\cW$. Assume $W\notin \Psi$ (resp. $W\in \Psi$). We prove that the map
$$
\hskip-60pt
\nu_{\a,\b} : D_{\SL(2)}^I(\Psi) \;(\text{resp.}\;\; D_{\SL(2)}^I(\Psi)_{\nspl})
$$
$$
\hskip30pt
\to \bR_{\geq 0}^{\Psi} \times D \times \spl(W) \times \tsize\prod_{W' \in \Psi} \spl(W'(\gr^W))
$$
is injective. 
  Denote $\nu_{\a,\b}(p)$ by $(\b(p), \mu(p), \spl_W(p), (\spl^{\BS}_{W'(\gr^W)}(p(\gr^W)))_{W' \in \Psi})$ (note that the symbol $\mu$ was introduced 
in 3.3.9). 

Let $p\in D_{\SL(2)}^I(\Psi)$ (resp. $D_{\SL(2)}^I(\Psi)_{\nspl}$). Then 
the set $J\sub \Psi$ of weight filtrations associated to $p$ is recovered from $\b(p)$ as
$$
J=\{j\in \Psi\;|\;\b(p)_j=0\}.
$$
Let $\a_J$ be the restriction of $\a$ to the $J$-component $\bG_{m,\bR}^J$ of $\bG_{m,\bR}^{\Psi}$. Since both $\gr^W(\tau_p)$ and $\gr^W(\a_J)$ split 
$W'(\gr^W)$ for all $W'\in J$, there is a unique element $u$ of $G_{\bR}(\gr^W)$ such that $\gr^W(\tau_p)=\Int(u)(\gr^W(\a_J))$ and such that $(1-u)\bar 
S(\bar J, m)\sub \bigoplus_{m'<m} \bar S(\bar J, m')$ for any 
$m\in \bZ^{\bar J}$ (cf.\ 3.3.8). 
This $u$ is characterized by the following property (1). 

\medskip 

(1) For any $W'\in J$, $u^{-1}\spl_{W'(\gr^W)}^{\BS}(p(\gr^W))$ coincides with the splitting of $W'(\gr^W)$ defined by the $W'$-component of $\alpha$. 

\medskip

The torus orbit associated to $p$ is recovered as
$$
\{\spl_W(p)\theta(u\gr^W(\a(t))(\mu(p)(\gr^W)),\, \Ad(u \a(t))(\delta(\mu(p))))\;|\; t\in \bR_{>0}^{\Psi},\, \b(p)=0^Jt\}.
$$

Next, let $\Phi\in\overline{\cW}$.  
We prove that the map
$$
\nu_{\a,\b} : D_{\SL(2)}^{II}(\Phi) \to \bR_{\geq 0}^{\Phi} \times D(\gr^W) \times \bar\cL\times \spl(W) \times \tsize\prod_{W' \in \Phi} \spl(W')
$$
is injective. 
Denote $\nu_{\a,\b}(p)$ by $(\b(p(\gr^W)),\, \mu(p(\gr^W)),\, \spl_W(p),\, (\spl^{\BS}_{W'}
(p(\gr^W)))_{W' \in \Phi})$. 

Let $p\in D_{\SL(2)}^{II}(\Phi)$. 
Let $J$ be the set of weight filtrations associated to $p$.
Let $\bar J=\{W'(\gr^W)\;|\; W'\in J, W'\neq W\}\sub \Phi$. Then $\bar J$ is recovered from $\b(p(\gr^W))$ as
$$
\bar J=\{j\in \Phi\;|\; \b(p(\gr^W))_j =0\}.
$$
Let $\mu(p(\gr^W))=(x, y)$ with $x\in D(\gr^W)$ and $y\in \bar \cL$ (3.3.10). 
If $y\in \cL$, $J$ is the lifting of $\bar J$ on $H_{0,\bR}$ by $\spl_W(p)$. 
If $y\in \bar \cL\smallsetminus \cL$, $J$ is the union of $\{W\}$ and the  lifting of $\bar J$ on $H_{0,\bR}$ by $\spl_W(p)$.

Let $\a_{\bar J}$ be the restriction of $\a$ to the $\bar J$-component $\bG_{m,\bR}^{\bar J}$ of $\bG_{m,\bR}^{\Phi}$. Since both $\bar \tau_p$ and $\a_{\bar J}$ 
split 
all $W'\in \bar J$, 
there is a unique element $u$ of $G_{\bR}(\gr^W)$ such that $\gr^W(\bar \tau_p)=\Int(u)(\a_{\bar J})$ and such that $(1-u)\bar S(\bar J, m)\sub \bigoplus_{m'<m} \bar S(\bar J, m')$ for any $m\in \bZ^{\bar J }$. This $u$ is characterized by the following property (1). 

\medskip 
(1) For any $W'\in \bar J$, $u^{-1}\spl_{W'}^{\BS}(p(\gr^W))$ coincides with the splitting of $W'$ defined by the $W'$-component of $\alpha$. 

\medskip

If $W\notin J$ (note $y\in \cL$ in this case), 
the torus orbit associated to $p$ is recovered as 
$$
\{\spl_W(p)\theta(u\a(t)(x), \Ad(u\a(t))(y))\;|\; t\in \bR_{>0}^{\Phi},\, \b(p)=0^Jt\}.
$$
If $W\in J$, $y$ has the shape $0\circ z$ with $z\in \cL \smallsetminus \{0\}$ (3.3.10), and the torus orbit associated to $p$ is recovered as  
$$
\{\spl_W(p)\theta(u\a(t)(x), t'\circ \Ad(u\a(t))(z))\;|\; t\in \bR_{>0}^{\Phi},\, \b(p)=0^Jt,\, t'\in \bR_{>0}\}.
$$

Proposition 3.2.7 is proved.

\medskip

{\bf 3.3.14.} 
We prove Proposition 3.2.9. The proofs of (i) and (ii) are similar. We give here the proof of (ii). 

To prove that another choice $(\a',\b')$ gives the same structure as $(\a, \b)$, we may assume either $\a=\a'$ or $\b=\b'$.

Assume first $\a=\a'$. 
Then we have a commutative diagram in which the right vertical arrow is a morphism of $\cB_\bR(\log)$.
$$
\matrix D_{\SL(2)}^{II}(\Phi)& @>{\text{by $\nu_{\a,\b}$}}>>& \bR_{\geq 0}^{\Phi} \times D(\gr^W) \times \bar \cL&& (t, y, \delta)\\
\Vert && \downarrow && \downarrow \\
D_{\SL(2)}^{II}(\Phi) &@>{\text{by $\nu_{\a,\b'}$}}>> & \bR_{\geq 0}^{\Phi} \times D(\gr^W) \times \bar \cL && \;(t\b'(y), \a\b'(y)^{-1}y, \Ad(\a\b'(y))^{-1}\delta).\endmatrix 
$$
Assume $\b=\b'$. 
Then $\a'=\Int(u)\a$ for some $u\in G_{\bR}$ such that $(u-1)W'_w\sub W'_{w-1}$ for any $W' \in \Phi$ and for any 
$w\in \bZ$. 
For $t\in \bR_{>0}^{\Phi}$, let  $u_t=\a(t)^{-1}u\a(t)$. 
Then as is easily seen, the map $\bR_{>0}^{\Phi}\to G_\bR,\;t\mapsto u_t$, extends to a real analytic map $\bR_{\geq 0}^{\Phi}\to G_\bR$, which we still denote by $t\mapsto u_t$. 
We have a commutative diagram in which the right vertical arrow is a morphism in $\cB_\bR(\log)$.  
$$
\matrix D_{\SL(2)}^{II}(\Phi)& @>{\text{by $\nu_{\a,\b}$}}>>& \bR_{\geq 0}^{\Phi} \times D(\gr^W) \times \bar \cL&& (t, y, \delta)\\
\Vert && \downarrow && \downarrow \\
D_{\SL(2)}^{II}(\Phi) &@>{\text{by $\nu_{\a',\b}$}}>> & \bR_{\geq 0}^{\Phi} \times D(\gr^W) \times \bar \cL && \;(t, uu_t^{-1}y, \Ad(uu_t^{-1})\delta).\endmatrix 
$$
These commutative diagrams prove (ii) of Proposition 3.2.9.

\hskip20pt

\head
\S3.4. Local properties of $D_{\SL(2)}$
\endhead

In this subsection, we prove Theorem 3.2.10 and Proposition 3.2.12, 
give  local descriptions of $D_{\SL(2)}^I$ and $D_{\SL(2)}^{II}$ 
(Theorem 3.4.4, Theorem 3.4.6), and 
 prove a criterion (Proposition 3.4.29) for the coincidence of 
$D_{\SL(2)}^I$ and $D_{\SL(2)}^{II}$.

\medskip

{\bf 3.4.1.}
Let $p\in D_{\SL(2)}$, let $\Phi=\overline{\cW}(p)$ (3.2.2), 
let $\br$ be a point on the torus orbit associated to $p$, and let $\bar \br=\br(\gr^W)$. 
Fix $\bR$-subspaces $$R\sub \fg_\bR(\gr^W), \quad S\sub\Lie(K_{\bar \br})$$ satisfying the following conditions (a), (b), and (c).
Here $K_{\bar \br}=\prod_w K_{\bar \br_w}$ with $K_{\bar \br_w}$ the maximal compact subgroup of $G_\bR(\gr^W_w)$ corresponding to $\bar\br_w$ (see \cite{KU3}, 5.1.2), 
where we write $\bar \br = (\bar \br_w)_w$ as in 3.3.6. Note $K_{\bar \br_w}
\supset K'_{\bar \br_w}$ for all $w$. 

\medskip

\noindent
(a) $\fg_\bR(\gr^W)=R\op\Lie(\t\rho(\bR_{>0}^\Phi)) \op\Lie(K_{\bar \br})$.

\medskip

Here $\t\rho$ is  the homomorphism $\bG_{m,\bR}^\Phi \to G_\bR(\gr^W)$ defined by
$$
\tilde\rho(t_1,\dots, t_n) = \tsize\bigoplus_{w\in \bZ}(\rho_w(g_1, \dots, g_n)\;\text{on}\;\gr^W_w) \quad 
\text{with} \quad 
g_j=\pmatrix 1/\prod_{k=j}^n t_k& 0 \\
0 & \prod_{k=j}^n t_k\endpmatrix,
$$
where $n$ is the number of the elements of $\Phi$ and  $((\rho_w,\vf_w)_w, \br)$ is the $\SL(2)$-orbit in $n$ variables of rank $n$ with class $p$ (cf. 2.3.5).
\medskip

\noindent
(b) $\Lie(K_{\bar \br}) = S\op\Lie(K'_{\bar \br})$, 
where $K'_{\bar \br}=\prod_w K'_{\bar \br_w}$ (cf. 3.3.6). 

\medskip

We introduce notation to state the condition (c). Let $\fg_\bR(\gr^W)=
\bigoplus_{m\in \bZ^{\Phi}} \;\fg_{\bR}(\gr^W)_m$ be the direct decomposition associated to the adjoint action of $\bG_{m,\bR}^{\Phi}$ via $\t\rho$. Note that this action coincides with the adjoint action of $\bG_{m,\bR}^{\Phi}$ via $\bar \tau_p$ (3.2.3). Thus
$$
\fg_{\bR}(\gr^W)_m := \{x\in \fg_\bR(\gr^W)\;|\;\Ad(\bar \tau_p(t))x=t^mx \;\;\text{for all}\; t\in (\bR^\times)^{\Phi}\}.
$$

The condition (c) is the following.
\medskip
\noindent
(c) $R=\ts_{m\in \bZ^{\Phi}} \; R\cap (\fg_\bR(\gr^W)_m+\fg_\bR(\gr^W)_{-m}).$
\medskip

Such $R$ and $S$ exist. 
The proof of the existence for the pure case is in \cite{KU3}, 10.1.2, and 
the general case is similar to it.  
  We remark that 
when we are given a parabolic subgroup $P$ of $G_\bR(\gr^W)$, we can take $R\sub \Lie(P)$. 

\medskip

{\bf 3.4.2.} 
Let the notation be as in 3.4.1. 
 We define objects $Y^{II}(p,\br,S)$ and $Y^{II}(p, \br, R, S)$ 
 of $\cB_\bR(\log)$.

 Let $L=\cL(\bar \br)$ (1.2.1). 

We define sets $Z(p)$ and $Z(p,R)$. Let  
$$
Z(p) \subset \bR_{\geq 0}^{\Phi} \times \fg_\bR(\gr^W) \times \fg_\bR(\gr^W) \times \fg_\bR(\gr^W)
$$
be the set  of all $(t, f, g, h)$ satisfying the following conditions (1) and (2). 
Let $J=J(t):=\{j\in \Phi\;|\; t_j=0\}$.
\medskip
(1) For $m\in \bZ^{\Phi}$, $g_m=0$ 
unless $m(j)=0$ for all $j\in J$, $f_m=0$ unless $m(j)\leq 0$ for all $j\in J$, and $h_m=0$ unless $m(j)\geq 0$ for all $j\in J$.

Here $(\;)_m$ for $m\in \bZ^{\Phi}$ denotes the $m$-component for the adjoint action of $\bG_{m,\bR}^{\Phi}$ under $\bar \tau_p$.
\medskip

(2) Let $t'$ be any element of $\bR_{>0}^{\Phi}$ such that $t'_j=t_j$ for any $j\in \Phi\smallsetminus J$. If $m\in \bZ^{\Phi}$ and $m(j)=0$ for any $j\in J$, then $g_m=\Ad(\bar\tau_p(t'))^{-1}(f_m)$ and $g_m= \Ad(\bar\tau_p(t'))(h_m)$. 

\medskip
Let
$$Z(p,R)\subset Z(p)$$
be the subset consisting of all elements $(t,f,g,h)$ satisfying the following condition (3).

\medskip

(3) $g\in R$ and $f_m+h_{-m}\in R$ for all $m\in \bZ^\Phi$. 
\medskip

Let 
$$
Y^{II}(p,\br,S)\sub Z(p)\times S\times \bar L\times \fg_{\bR,u}
$$
$$
\text{(resp}.\quad Y^{II}(p,\br,R,S)\sub Z(p,R)\times S\times \bar L\times \fg_{\bR,u})
$$ 
be the set consisting of all elements $(t,f,g,h,k, \delta,u)$ ($(t,f,g,h)\in Z(p)$ $(\text{resp.}\;Z(p,R))$, $k\in S$, $\delta\in \bar L$, $u\in \fg_{\bR,u}$) satisfying the following condition (4). 

\medskip

(4) $\exp(k)\bar \br \in (K_{\bar \br}\cap G_{\bR}(\gr^W)_{J})\cdot \bar \br$ with $J=J(t)$, where $G_{\bR}(\gr^W)_{J}=\{g \in G_{\bR}(\gr^W)\;|\;g W' = W'$ 
for any $W' \in J\}$. 
\medskip

We endow $Y^{II}(p,\br,S)$ (resp.\ $Y^{II}(p, \br, R, S))$ with the following structure as an object of $\cB_{\bR}(\log)$. 

Let $E=\bR_{\geq 0}^{\Phi}\times \fg_\bR(\gr^W) \times \fg_\bR(\gr^W)\times \fg_\bR(\gr^W)\times S \times \bar L \times \fg_{\bR,u}$. 
Let $A=Y^{II}(p, \br,S)$ 
(resp.\ $A=Y^{II}(p, \br, R, S))$.

We endow $A$ with the topology as a subspace of $E$. 

We define the sheaf of real analytic functions on $A$ as follows. 
For an open set $U$ of $A$ and for a map $f: U\to \bR$, we say that 
$f$ is real analytic if and only if, for any $p\in U$, there are an open neighborhood $U'$ of $p$ in $U$, an open neighborhood $U''$ of $U'$ in $E$, and a real analytic function $g$ on $U''$, such that the restrictions to $U'$ of $f$ and $g$  coincide. 

We show that with this sheaf of rings over $\bR$, $A$ is an object of $\cB_\bR$. Let $\cO_E$ be the sheaf of real analytic functions on $E$. 
Let $I$ be the ideal of $\cO_E$ generated by the following sections $a_{m,l}$ and $b_{m,l}$ given for 
elements $m$ of $\bZ^{\Phi}$ and for $\bR$-linear maps $l:\fg_\bR(\gr^W)\to \bR$:
$$
a_{m,l}(t,f,g,h,k,\delta, u)= (\tsize\prod_{j\in \Phi,\, m(j)\leq 0} t_j^{-m(j)})l(f_m)-(\tsize\prod_{j\in \Phi,\, m(j)\geq 0} t_j^{m(j)})l(g_m),
$$
$$
b_{m,l}(t,f,g,h,k,\delta, u)= (\tsize\prod_{j\in \Phi,\, m(j)\leq 0} t_j^{-m(j)})l(g_m)-(\tsize\prod_{j\in \Phi,\, m(j)\geq 0} t_j^{m(j)})l(h_m).
$$
Here $(\;\;)_m$ denotes the $m$-th component with respect to the adjoint action of $\bG_{m,\bR}^{\Phi}$ by $\bar \tau_p$,  $\prod_{j\in \Phi, m(j)\leq 0}$ means the product over all $j\in \Phi$ such that $m(j)\leq 0$, and  $\prod_{j\in \Phi, m(j)\geq 0}$ is defined in the similar way. 
Then $I$ is a finitely generated ideal. 
Indeed, if $l_1, \dots, l_r$ form a basis of the dual $\bR$-vector space of $\fg_\bR(\gr^W)$, $a_{m,l_j}$ and $b_{m,l_j}$ ($1\leq j\leq r$) such that $\fg_\bR(\gr^W)_m\neq 0$ (there are only finitely many such $m$) generate $I$. 
Furthermore, the inverse image of $\cO_E/I$ 
on $Y^{II}(p,\br,S)$ coincides with the sheaf of real analytic functions on $Y^{II}(p,\br,S)$. 
Hence $Y^{II}(p,\br,S)$ is an object of $\cB_\bR$. 
Let $I'$ be the ideal of $\cO_E$ generated by $I$ and by the following sections $c_l$ 
and $d_{m,l}$ given for elements $m$ of $\bZ^{\Phi}$ and $\bR$-linear maps $
l: \fg_\bR(\gr^W)\to  \bR$ which kill $R$:
$$
c_l(t, f,g,h,k,\delta,u)=l(g),
$$
$$
d_{m,l}(t,f,g,h,k,\delta,u)=l(f_m+h_{-m}).
$$
As is easily seen, $I'$ is a finitely generated ideal. Furthermore, the inverse image of $\cO_E/I'$ 
on $Y^{II}(p,\br,R,S)$ coincides with the sheaf of real analytic functions on $Y^{II}(p,\br,R,S)$. Hence $Y^{II}(p,\br, R, S)$ is also an object of $\cB_\bR$.

We define the log structures with sign of $Y^{II}(p, \br,S)$ and $Y^{II}(p, \br,R,S)$ to be the inverse images of the  log structure with sign of $\bR_{\geq 0}^{\Phi}$. 
This endows $Y^{II}(p,\br, S)$ and $Y^{II}(p, \br,R,S)$ with structures of  objects of $\cB_\bR(\log)$.

\medskip

{\bf 3.4.3.} 
Define an open subset $Y_0^{II}(p,\br,S)$ of $Y^{II}(p,\br,S)$ by
$$Y_0^{II}(p,\br,S)=\{(t,f,g,h,k,\delta,u) \in Y^{II}(p,\br,S)\;|\;t\in \bR_{>0}^{\Phi}, \delta\in L\}.$$ We define 
an open subset $Y_0^{II}(p,\br,R,S)$ of $Y^{II}(p,\br, R, S)$ by
$$
Y_0^{II}(p,\br, R, S)=Y^{II}(p, \br,R,S)\cap Y_0^{II}(p,\br,S).
$$

We have isomorphisms of real analytic manifolds
$$
Y_0^{II}(p, \br,S)@>\simeq>>
\bR_{>0}^{\Phi}\times \fg_\bR(\gr^W)\times S \times L \times \fg_{\bR,u},
$$
$$
Y_0^{II}(p, \br, R, S)@>\simeq>>
\bR_{>0}^{\Phi}\times R\times S \times L \times \fg_{\bR,u},
$$
given by
$$
(t, f,g,h,k,\delta, u)\mapsto (t,g,k,\delta,u),
$$
whose inverse maps are given by $$f=\Ad(\bar \tau_p(t))(g),\quad h=\Ad(\bar\tau_p(t))^{-1}(g).$$

We have a morphism of real analytic manifolds
$$\eta_{p,\br,S}^{II}: Y_0^{II}(p,\br,S)\to D,\quad (t, f,g,h,k,\delta, u)\mapsto \exp(u)s_{\br}\theta(d\bar \br, \Ad(d)\delta)$$ $$\text{with}\quad s_\br=\spl_W(\br), \quad d=\bar \tau_p(t)\exp(g)\exp(k)=\exp(f)\bar \tau_p(t)\exp(k).$$
Let $$\eta_{p,\br,R,S}^{II}: Y_0^{II}(p,\br,R, S)\to D$$
be the induced morphism.
\medskip

\proclaim{Theorem 3.4.4}  
Let the notation be as above. 
If $U$ is a sufficiently small open neighborhood of $0:=(0,0,0,0)$ in 
$\fg_\bR(\gr^W)\times\fg_\bR(\gr^W)\times\fg_\bR(\gr^W)\times S$ and if $Y^{II}(p,\br, S,U)$ $($resp. $Y^{II}(p,\br, R,S,U)$$)$ denotes the open set of 
$Y^{II}(p,\br, S)$ $($resp. $Y^{II}(p,\br,R,S)$$)$ consisting of all elements $(t,f,g,h,k,\delta, u)$ such that $(f,g,h,k)\in U$, we have{\rm :}

\medskip

{\rm (i)} There is a unique morphism $Y^{II}(p,\br,S,U)\to D_{\SL(2)}^{II}(\Phi)$ in the category $\cB_\bR'(\log)$ 
  whose restriction to 
$Y_0^{II}(p,\br,S,U)=Y_0^{II}(p,\br,S)\cap Y^{II}(p,\br,S,U)$ coincides with the restriction of $\eta_{p,\br,S}^{II}$ $(3.4.3)$. 

\medskip

{\rm (ii)} The restriction of the morphism in {\rm (i)} 
induces an open immersion
$Y^{II}(p,\br,R, S,U)\to D_{\SL(2)}^{II}(\Phi)$ in the category $\cB_\bR'(\log)$ 
which sends
$(0^{\Phi}, 0,0,0,0,\delta(\br), 0)\in Y^{II}(p,\br, R, S, U)$ to $p$.

\endproclaim
The proof of this theorem will be given later in 3.4.18--3.4.19.

\medskip

{\it Remark.} 
From the proof of 3.4.4 given below, we see that if $q$ is  the image of $(t,f,g,h,k,\delta,u)\in Y^{II}(p, \br, S, U)$ in $D_{\SL(2)}(\Phi)$,  then $q\in
D_{\SL(2),\spl}$ if and only if $\delta=0$, and $W\in \cW(q)$  if and only if $\delta\in \bar L\smallsetminus L$. 
\medskip

{\bf 3.4.5.} 
Next we consider $D_{\SL(2)}^I$. 

Let $\Psi=\cW(p)$. Let $\Phi, \br, R, S$ be as before in 3.4.1.

We define an object $Y^I(p, \br, R, S)$ of $\cB_\bR(\log)$ first in the case 
$W\notin \cW(p)$. Let $$Y^I(p, \br,R,S)\sub Y^{II}(p,\br,R,S)\times \fg_{\bR,u}
\tag{$*$}$$ be the set consisting of all elements $(t,f,g,h,k,\delta,u,v)$ ($(t,f,g,h,k,\delta,u)\in Y^{II}(p,\br,R,S)$, $v\in \fg_{\bR,u}$) satisfying the following conditions (5)--(7).
Via the bijection $\Psi\to \Phi$, we regard $\tau_p$ as a homomorphism $\bG_{m,\bR}^{\Phi}\to \Aut(H_{0,\bR}, W)$.  Let $\fg_{\bR,u}=
\bigoplus_{m\in \bZ^{\Phi}}\fg_{\bR,u,m}$ be the corresponding direct sum 
decomposition. 
  Denote by $u_m$ the $m$-component of $u \in \fg_{\bR,u}$. 
  
\medskip

(5) For $m\in \bZ^{\Phi}$, $u_m=0$ unless $m(j)\leq 0$ for all $j\in J=J(t)$, and $v_m=0$ unless $m(j)=0$ for all $j\in J$. 

\medskip

(6) Let $t'$ be any element of $\bR_{>0}^{\Phi}$ such that $t'_j=t_j$ for any $j\in \Phi\smallsetminus J$. 
If $m\in \bZ^{\Phi}$ and $m(j)=0$ for any $j\in J$, then $v_m=\Ad(\tau_p(t'))^{-1}(u_m)$. 

\medskip
(7) $\delta\in L$ in $\bar L$. 
\medskip

  We endow $Y^I(p, \br, R, S)$ with a structure of an object of 
$\cB_\bR(\log)$ via the injection $Y^I(p, \br, R, S)\hookrightarrow E \times \fg_{\bR, u}$, 
just as we endowed $Y^{II}(p, \br, R, S)$ with it via the injection 
$Y^{II}(p, \br, R, S)\hookrightarrow E$ in 3.4.2.

Next, in the case $W\in\cW(p)$, we define an object $Y^I(p, \br, R, S)$ of $\cB_\bR(\log)$ by fixing a closed real analytic subspace $L^{(1)}$ of $L\smallsetminus \{0\}$ such that $\bR_{>0}\times L^{(1)}\to L\smallsetminus \{0\}, (a,x)\mapsto a\circ x$, is an isomorphism of real analytic manifolds. Via the evident bijection between $\Psi$ and the disjoint union of $\{W\}$ and $\Phi$, we regard $\tau_p$ as a homomorphism $\bG_{m,\bR}\times \bG_{m,\bR}^{\Phi}\to \Aut(H_{0,\bR}, W)$. 
 Let $$Y^I(p, \br,R,S)\sub \bR_{\geq 0}\times Y^{II}(p,\br,R, S)\times \fg_{\bR,u}\tag{$*$}$$ be the set consisting of all elements $(t_0,t,f,g,h,k,\delta,u,v)$ ($t_0\in \bR_{\geq 0}$, $(t,f,g,h,k,\delta,u)\in Y^{II}(p,\br,R,S)$, $v\in \fg_{\bR,u})$  satisfying the 
 following conditions (5${}'$)--(7${}'$).

\medskip

(5${}^\prime$)  (5) holds and furthermore, in the case $t_0=0$, we have $\exp(v)s_\br=s_\br$.

\medskip

(6${}^\prime$) Let $t'$ be any element of $\bR_{>0}^{\Phi}$ such that $t'_j=t_j$ for any $j\in \Phi\smallsetminus J$. Let $m\in \bZ^{\Phi}$ and assume $m(j)=0$ for any $j\in J$. If $t_0\neq 0$, then $v_m=\Ad(\tau_p(t_0, t'))^{-1}(u_m)$. 
If $t_0=0$, then $v_m=\Ad(\tau_p(1, t'))^{-1}(u_m)$.

\medskip

(7${}^\prime$) $\delta\in L^{(1)}$. 

\medskip

  We endow $Y^I(p,\br,R,S)$ with a structure of an object in $\cB_{\bR}(\log)$ 
via the injection $Y^I(p,\br,R,S) \hookrightarrow \bR_{\ge0}\times B \times 
\fg_{\bR,u}$.

We define a canonical morphism $Y^I(p,\br, R,S) \to Y^{II}(p, \br, R, S)$. In the case $W\notin \cW(p)$, it is just the canonical projection. In the case $W\in \cW(p)$, it is the morphism $(t_0,t',f,g,h,k,\delta,u,v) \mapsto
(t', f,g,h,k,t_0\circ \delta, u)$. In both cases, this morphism is injective. 

  Define an open subset $Y_0^I(p,\br,R,S)$ of $Y^I(p,\br,R,S)$ by the inverse image of $Y_0^{II}(p,\br,R,S)$ (3.4.3).
Then we have an isomorphism of real analytic manifolds 
$Y_0^I(p,\br,R,S) @>\simeq>> Y_0^{II}(p,\br, R, S)$.

Combining this with $\eta^{II}_{p,\br,R,S}$ (3.4.3), we have a morphism of real analytic manifolds 
$$
\eta^I_{p,\br,R,S}: Y^I_0(p,\br,R,S)\to D.
$$

\proclaim{Theorem 3.4.6}  
Let the notation be as above. 
Assume $W\notin \Psi$ $($resp. $W\in \Psi$$)$. 
Then if $U$ is a sufficiently small open neighborhood of $0:=(0,0,0,0)$ in 
$\fg_{\bR}(\gr^W) \times R\times \fg_{\bR}(\gr^W) \times S$ and if 
$Y^{I}(p,\br,R, S,U)$ denotes the open set of 
$Y^{I}(p,\br,R, S)$ defined as the inverse image of $U$ 
by the canonical map $Y^{I}(p,\br,R, S)\to 
\fg_{\bR}(\gr^W) \times R\times \fg_{\bR}(\gr^W) \times S$, 
then there is an open immersion $Y^{I}(p,\br,R, S,U)\to D_{\SL(2)}^{I}(\Psi)$ 
in the category $\cB_\bR'(\log)$ 
 which sends $(0^{\Phi}, 0,0,0,0,\delta(\br), 0, 0)$ $($resp. 
$(0^{\Psi}, 0,0,0,0,\delta(\br)^{(1)}, 0, 0)$, 
where $\delta(\br)^{(1)}\in L^{(1)}$  $(3.4.5)$ such that $\delta(\br)=0\circ \delta(\br)^{(1)}$$)$ to $p$ and whose restriction to 
$Y^{I}(p,\br,R, S,U) \cap Y_0^{I}(p,\br,R,S)$ coincides with the restriction of $\eta_{p,\br,R,S}^{I}$ $(3.4.5)$. 
\endproclaim
The proof will be given later in 3.4.20.

\medskip

{\bf 3.4.7.} 
  Before we start to prove 3.4.4 and 3.4.6, we make some 
preparation. 

  Let the notation be as in 3.4.1. 
  Then, there exist an open neighborhood $O$ of
$0$ in $\fg_\bR(\gr^W)$ and a real analytic function
$c=(c_1, c_2, c_3):O\to 
\bR^{\Phi}_{>0} \times R \times S$ having the following properties
{\rm(1)--(4)}.
\medskip

\noindent
{\rm(1)} For any $x\in O$,
$\exp(x)\bar \br=\bar \tau_p(c_1(x))\exp(c_2(x))\exp(c_3(x))\bar \br$.
\medskip

\noindent
{\rm(2)} $c(0)=(1,0,0)$.
\medskip

\noindent
{\rm(3)} $\exp:O\to G_\bR(\gr^W)$ is an injective open map. 
\medskip

\noindent
{\rm(4)} 
For $k=2,3$, $c_k$ has the form of
absolutely convergent series $c_k=\ts_{r=0}^\infty
c_{k,r}$, where $c_{k,r}$ is the part of degree $r$
in the Taylor expansion of $c_k$ at $0$, such that
$c_{k,r}(x)=l_{k,r}(x\ox\cdots\ox x)$ for some
linear map $l_{k,r}:\fg_\bR(\gr^W)^{\ox r}\to\fg_\bR(\gr^W)$
having the following property$:$
If $m_1, \dots, m_r\in \bZ^{\Phi}$ and $x_j\in\fg_{\bR}(\gr^W)_{m_j}$ for $1\le j\le r$, then $l_{k,r}(x_1\ox\cdots\ox x_r)
\in\ts_m\fg_{\bR}(\gr^W)_{m}$, where $m$ ranges over
all elements of $\bZ^{\Phi}$ of the form 
$\sum_{1 \le j \le r} e_jm_j$ with $e_j \in \{1, -1\}$ for 
each $j$. 

\medskip

  This is proved similarly as \cite{KU3}, 10.3.4. 
  Cf.\ also 3.3.11. 
  It is clear that there is a real analytic $c$ satisfying (1)--(3) 
unique up to restrictions of domains of definitions. 
  The property (4) of Taylor expansion can be checked formally as follows. 

Consider the following formal calculation. 

\medskip

$$
\align
\exp(x)&=\exp(t^{(1)}+b^{(1)}+k^{(1)})
=\exp(t^{(1)})\exp(b^{(1)}+x^{(1)})\exp(k^{(1)})\\
&= \exp(t^{(1)})\exp(b^{(1)}+t^{(2)}+b^{(2)}+k^{(2)})\exp(k^{(1)})\\
&= \exp(t^{(1)})\exp(t^{(2)})\exp(b^{(1)}+b^{(2)}+x^{(2)})\exp(k^{(2)})\exp(k^{(1)})=\cdots
\endalign
$$

\medskip

\noindent 
  Here $x \in O$, $t^{(j)} \in \Lie(\t\rho(\bR^{\Phi}_{>0}))$ with $\t\rho$ being 
as in 3.4.1 (note that the actions of $\t\rho(t)$ and $\bar\tau_p(t)$ for $t\in \bR_{>0}^{\Phi}$ on $D(\gr^W)$ coincide), 
$b^{(j)} \in R$, $k^{(j)} \in S$, $x^{(j)} \in \fg_{\bR}(\gr^W)$ for any $j$. 
  Then, we have $\t\rho(c_1(x))=\exp(t^{(1)})\exp(t^{(2)})\cdots$, 
$c_2(x)=b^{(1)}+b^{(2)}+\cdots$, 
and 
$\exp(c_3(x))=\cdots\exp(k^{(2)})\exp(k^{(1)})$, formally.
  From these, we can prove the property (4) formally. 

\medskip

{\bf 3.4.8.}  
We prove 3.4.4 till 3.4.19.  After that, we prove 3.4.6. 
Let $p, \Phi$ and $\br$ be as in 3.4.1.
  In the notation in 3.4.7, let $U=\exp(O)\bar \br$ which is 
an open neighborhood of $\bar \br$ in $D(\gr^W)$.
  By 3.4.7, there is a real analytic map 
$$a=(a_1,a_2,a_3): U\to \bR_{>0}^{\Phi} \times R \times S$$
such that for any $y \in U$, we have 
$y=\bar\tau_p(a_1(y))\exp(a_2(y))\exp(a_3(y))\bar \br.$
(Just put $a_j(\exp(x)\bar \br)=c_j(x)$ for $x \in O$.) 

\medskip

{\bf 3.4.9.} 
Fix a distance $\b$ to $\Phi$-boundary such that $\b(\bar \br)=1$. 
Here we denote $\b(x)=\b(x(\gr^W))$ $(x \in D)$ by abuse of notation. 
Let $\mu: D(\gr^W) \to D(\gr^W)$ be the real analytic map defined by 
$\mu(x)=\bar\tau_p(\b(x))^{-1}x$.
Denote the composite 
$D \to D(\gr^W) @>\mu>> D(\gr^W)$ also by $\mu$ by abuse of notation. 
Let $D(U) \subset D$ be the inverse image of $U$ by $\mu$. 

Let 
$$
b=b_{R,S}:D(U) \to Y^{II}_0(p, \br, R,S)
$$ 
be the real analytic map 
$x\mapsto (t,f,g,h,k,\delta,u)$, 
where 
$t=\b(x)a_1(\mu(x))$, $f=\Ad(\bar\tau_p(t))(a_2(\mu(x)))$, $g=a_2(\mu(x))$, $h=\Ad(\bar\tau_p(t))^{-1}(a_2(\mu(x)))$, $k=a_3(\mu(x))$, $\delta=\Ad(\bar \tau_p(t)\exp(g)\exp(k))^{-1}(\delta(x))$, and 
$u$ is characterized by $\spl_W(x)=\exp(u)\spl_W(\br)$. 

\medskip

  Recall that, in 3.4.4, for an 
open neighborhood $U'$ of $0$ in $\fg_\bR(\gr^W)\times \fg_\bR(\gr^W)\times 
\fg_\bR(\gr^W)\times S$, 
we denote by $Y^{II}(p,\br, S,U')$ the subset of $Y^{II}(p, \br,S)$ 
consisting of all elements $(t,f,g,h,k,\delta,u)$ such that $(f,g,h,k)\in U'$.  
  We also defined $Y_0^{II}(p,\br, S,U')$ and 
$Y^{II}(p, \br, R, S, U')$ there. 
  Now, we define 
$Y_0^{II}(p, \br, R, S, U')=Y^{II}(p, \br, R, S, U')\cap 
Y_0^{II}(p,\br,S)$. 

\medskip

The next two lemmas are easily seen. 

\proclaim{Lemma 3.4.10}  The composite $D(U) \overset b \to 
\to Y^{II}_0(p,\br,R,S)\to D$ is 
the canonical inclusion.
\endproclaim

\proclaim{Lemma 3.4.11}
 If $U'$ is sufficiently small, the image of $Y^{II}_0(p,\br,S,U')$ 
in $D$ is contained in $D(U)$ and 
the map $Y^{II}_0(p,\br,R,S,U') \to D(U) \to Y^{II}_0(p,\br,R,S)$ is the 
canonical inclusion. 
\endproclaim 

\medskip

{\bf 3.4.12.} 
We define 
 $$p(J, \br, z, \delta, u)\in D_{\SL(2)}^{II}(\Phi)$$ as follows, for a subset $J$ of $\Phi$, a point $\br$ on the torus orbit associated to $p$ (2.5.2),  an element $z$ of $G_\bR(\gr^W)$ which satisfies
\medskip

(1) $z\in G_\bR(\gr^W)_J$,

\medskip
\noindent 
an element $\delta$ of $\bar L$, and an element $u$ of $\fg_{\bR,u}$.

\medskip

This $p(J, \br, z,\delta,u)$ is the unique element of $D_{\SL(2)}$ which satisfies the following (2)--(5).

\medskip

(2) The set of weight filtrations on $\gr^W$ associated to $p(J, \br, z, \delta,u)$ is $J$.

\medskip

(3) The torus action $\bar \tau$ associated to $p(J, \br, z, \delta,u)$ is 
$\Int(z)(\bar \tau_{p,J}):\bG_{m,\bR}^J\to \Aut_\bR(\gr^W)$, where $\bar \tau_{p,J}$ denotes the restriction of $\bar \tau_p:\bG_{m,\bR}^{\Phi}\to \Aut_\bR
(\gr^W)$ (2.5.6, 2.3.5) to the $J$-component of $\bG_{m,\bR}^{\Phi}$.

\medskip

(4) $\delta\in L$ in $\bar L$ if and only if $W$ does not belong to the set of weight filtrations associated to $p(J, \br,z,\delta,u)$. 

\medskip

(5) The torus orbit associated to $p(J, \br, z, \delta,u)$ (2.5.2) 
contains $\exp(u)s_{\br}\theta(z(\br(\gr^W)), \Ad(z)(\delta))$ if $\delta\in L$, and contains $\exp(u)s_{\br}\theta(z(\br(\gr^W)), \Ad(z)(\delta'))$ if $\delta\in \bar L\smallsetminus L$ and $\delta=0\circ \delta'$ with $\delta'\in L\smallsetminus \{0\}$. 
\medskip

This $p(J,\br, z,\delta,u)$ is constructed as follows.  
Let $n$ be the cardinality of $\Psi=\cW(p)$, and identify $\Psi$ with $\{1,\dots, n\}$ as a totally ordered set for the ordering in 2.3.8. 
In the case $W\notin \Psi$, consider the bijection $\Psi\to \Phi$. 
In the case $W\in \Psi$, consider the bijection $\Psi\smallsetminus \{W\}\to \Phi$. 
Via these bijections, embed $J\sub \Phi$ into $\Psi$. 
In the case $\delta\in L$ (resp. $\delta\in \bar L\smallsetminus L$),
let $m=\sharp(J)$ (resp. $m=\sharp(J)+1$)  
 and write  $J=\{j_1, \dots, j_m\}\sub \Psi$ with $j_1<\dots <j_m$ (resp. 
 $J=\{j_2,\dots, j_m\}\sub \Psi$ with $j_2<\dots<j_m$). 
Let $((\rho_w, \vf_w)_w, \br)$ be an $\SL(2)$-orbit in $n$ variables of rank $n$ whose class in $D_{\SL(2)}$ is $p$. 
Then, in the case $\delta\in L$ (resp. $\delta\in \bar L\smallsetminus L$), 
the $p(J,\br,z,\delta,u)$ is the class of the following SL(2)-orbit 
$((\rho',\vf')=(\rho'_w, \vf'_w)_w, \br')$ in $m$ variables of rank $m$.
$$
\rho'(g_1,\dots,g_m):=\Int(z)(\rho(g'_1,\dots,g'_n)),
$$
$$
\vf'(z_1,\dots,z_m):=z\vf(z'_1,\dots,z'_n),
$$
$$
\br':= \exp(u)s_{\br}\theta(z(\br(\gr^W)), \Ad(z)(\delta))
$$ 
$$
(\text{resp.}\;\;
\br':=\exp(u)s_{\br}\theta(z(\br(\gr^W)), \Ad(z)(\delta'))\;\;\text{with}\;\;\delta'\in  L\smallsetminus\{0\},\;\;\delta=0\circ\delta'),
$$
where 
$g'_j$ and $z'_j$ $(1\leq j\leq n)$ are as follows.
If $j\le j_k$ for some $k$, define $g'_j:=g_k$ and $z'_j:=z_k$ for the smallest integer $k$ with $j\le j_k$.
Otherwise, $g'_j:=1$ and $z'_j:=i$. 

\medskip
Let $Y_1:=Y_1^{II}(p,\br,S)$ 
be the subset of 
$Y^{II}(p, \br,S)$ consisting of all elements 
$(t,f,g,h,k,\delta,u)$ such that $h_m =0$ unless $m(j)=0$ for all $j \in J(t)$. We have $Y_1\supset Y_0:=Y_0^{II}(p,\br,S)$. We have

\medskip

(6) A point $(t,f,g,h,k,\delta,u)\in Y^{II}_1(p,\br, S)$ is the limit of 
$y(t', \delta')\in Y_0^{II}(p,\br,S)$ defined by $y(t',\delta')=
(t', f,\; \Ad(\tau_p(t'))^{-1}(f),\; \Ad(\tau_p(t'))^{-2}(f),\; k,\delta', u)$,
where $t'\in \bR_{>0}^{\Phi}$, $\delta'\in L$, and 
$t'$ tends to $t$ and $\delta'$ tends to $\delta$. Write $\exp(k)\cdot \bar \br=k'\cdot\bar \br$ with $k'\in K_{\bar \br}\cap G_\bR(\gr^W)_J$. Note that $k'$ commutes with $\bar \tau_p(t')$. 
The image of $y(t', \delta')$ in $D$ is
$\exp(u)s_{\br}\theta(z(\bar \br), \Ad(z)(\delta''))$, where $z=\exp(f)k'\bar \tau_p(t') $ and $\delta''=\Ad((k')^{-1}\exp(k))(\delta')$. 

\medskip

 We extend the map $\eta^{II}_{p,\br,S}: Y_0\to D$ in 
3.4.3 to a map 
$$
\eta^{II}_{p,\br,S}: Y_1\to D_{\SL(2)}^{II}(\Phi),$$ $$\eta^{II}_{p,\br,S}(t,f,g,h, k, \delta,u) = p(J,\br,z, \delta',u),
$$  
where $J$, $z$, and $\delta'$ are defined as follows. 
  Let $J=\{j\in \Phi\;|\;t_j=0\}$.
  Let $t'$ be an element of $\bR_{>0}^{\Phi}$ such that $t'_j=t_j$ for any $j\in \Phi\smallsetminus J$, and 
let $k'$ be an element of $K_{\bar \br}\cap G_\bR(\gr^W)_J$ 
such that $\exp(k)\cdot \bar\br=k'\cdot\bar \br$.
  Let $z=\exp(f)k'\bar \tau_p(t')$ 
and $\delta'=\Ad((k')^{-1}\exp(k))\delta$.

\medskip

We will use the following fact (7) which is deduced from 10.2.16 of \cite{KU3}. 

\medskip

(7) Let $\mu: D^{II}_{\SL(2)}(\Phi)\to D(\gr^W)$ be the extension of $\a\b(x(\gr^W))^{-1}x(\gr^W)$ $(x \in D)$ given in 3.2.6 (ii). 
Then, if $p'\in D^{II}_{\SL(2)}(\Phi)$ and if $\mu(p')$ is sufficiently near to $\mu(p)$, $p'$ is expressed as $p(J, \br, z,\delta', u)$ as above. 

\proclaim{Lemma 3.4.13} 
There are an open neighborhood $U'$ of $0$ in $\fg_{\bR}(\gr^W)\times \fg_{\bR}(\gr^W)\times \fg_{\bR}(\gr^W)\times S$ and a morphism $\xi : Y^{II}(p,\br,S,U')\to Y^{II}(p,\br, R, S)$ which satisfy the following conditions: 
$\eta^{II}_{p,\br,S}$ sends $Y^{II}_0(p,\br,S,U')$ into $D(U)$ and the restriction of $\xi$ to $Y^{II}_0(p,\br,S,U')$ coincides with the composite 
$Y^{II}_0(p,\br,S,U') @>\eta^{II}_{p,\br,S}
>> D(U) @>b>> Y^{II}_0(p,\br, R, S)$, where $b$ is in $3.4.9$.
\endproclaim

\demo{Proof} 
  Let $x=\eta^{II}_{p,\br,S}(t,f,g,h,k,\delta,u)$ and write 
$b(x)$ as $(t',f',g',h',k',\delta',u')$. 

  First we show that each component $t',f', g', \ldots$ extends real analytically 
over the boundary of $Y^{II}(p,\br,S,U')$ for some $U'$.  
  Since $\mu(x)=\bar \tau_p(\b(\exp(g)\exp(k)\bar \br))^{-1}
\exp(g)\exp(k)\bar \br$, this extends over the boundary.
  Hence, so does $a_j\mu(x)$ for each $j=1,2,3$ (3.4.8). 
  On the other hand, $\b(x)=t\b(\exp(g)\exp(k)\bar \br)$, and this 
is also real analytic over the boundary because $\b$ is so. 
  Thus, $t', g', k'$ extend. 
  Further, $u'=u$ trivially extends. 
  We have $\delta'=\Ad(\bar \tau_p(t')\exp(g')\exp(k'))^{-1}\Ad(\bar \tau_p(t)
\exp(g)\exp(k))(\delta)$. 
  Since $g'$ and $k'$ already extend and since $t't^{-1}=\b(\exp(g)\exp(k)\bar \br) a_1\mu(x)$ also extends, so does $\delta'$. 

  The rest are $f'$ and $h'$, that is, to see that 
$\Ad(\bar \tau_p(t'))^{\pm 1}a_2\mu(x)$ extend real analytically. 
  We can replace $t'$ in the last formula with $t$ because $t'=t \b
(\exp(g)\exp(k)\bar \br)a_1(\mu(x))$. 
  Further, by 3.4.7 with the formal construction there, $a_2(\mu(x))
=c_2(g)$.  
  Hence, it is enough to show that $\Ad(\bar \tau_p(t))^{\pm 1}c_2(g)$ 
extend. 
 
   Consider the decomposition $g=\sum_{m\in \bZ^{\Phi}}g_m$ 
$(g_m \in \fg_{\bR}(\gr^W)_m)$. 
  Then, by the property (4) of $c_2$ in 3.4.7, 
$c_2(g)=c_2(\sum g_m)$ is the infinite formal sum of 
$l_{2,r}(g_{m_1}\otimes \cdots \otimes g_{m_r})$ 
$(m_j \in \bZ^{\Phi}$, $g_{m_j} \in \fg_{\bR}(\gr^W)_{m_j}$ $(1 \le j \le r))$. 
  Now the weights $m$ of $l_{2,r}(g_{m_1}\otimes \cdots \otimes g_{m_r})$ 
satisfy $m=\sum e_j m_j$ with $e_j \in \{1, -1\}$. 
  Decompose 
$l_{2,r}(g_{m_1}\otimes \cdots \otimes g_{m_r})$ into 
$\sum_m l_{2,r,m}(g_{m_1}\otimes \cdots \otimes g_{m_r})$ according to 
the weights, where $m$ ranges over such $\sum e_j m_j$. 
  We will see that, for each $m$ and $j \in \{1, -1\}$, 
$\bar \tau_p(t)^jl_{2,r,m}(g_{m_1}\otimes \cdots \otimes g_{m_r})$ extends over 
the boundary.  
  We explain the proof for $j=1$.  The other case is similar. 
In this case, we observe $\bar \tau_p(t)l_{2,r,m}(g_{m_1}\otimes \cdots \otimes g_{m_r})$ is 
$(\prod (t^{m_j})^{e_j})l_{2,r,m}(g_{m_1}\otimes \cdots \otimes g_{m_r})
= l_{2,r,m}((t^{m_1})^{e_1}g_{m_1}\otimes \cdots \otimes (t^{m_r})^{e_r}g_{m_r})$. 
  Since $t^mg_m=f_m$ and $t^{-m}g_m=h_m$, 
the last function extends to a real analytic map over the boundary. 
  Shrinking $U'$ if necessary, we may assume that $f$ and $h$ are sufficiently 
near to $0$, and the above infinite sum converges, as desired. 

  Next we show that 
in the ambient product space containing $Y^{II}(p,\br,R,S)$, the image 
of each element 
$y=(t,f,g,h,k,\delta, u)$ of 
$Y^{II}(p,\br,S,U')$ by the 
extended coordinate functions in fact belongs to $Y^{II}(p,\br,R,S)$, 
which completes the proof. 
  For $t'\in \bR_{>0}^{\Phi}$ such that $t'_j=t_j$ for any $j\in \Phi\smallsetminus J$ with $J=J(t)$  
and for $\delta'\in L$, let $y(t',\delta')=(t', f', g', h', k, \delta', u)\in Y_0^{II}(p,\br,S)$, where $f', g', h'\in \fg_{\bR}(\gr^W)$ are defined as follows. 
Let $m\in \bZ^{\Phi}$. 
Then $f'_m=(t')^{2m}h_m$, $g'_m=(t')^mh_m$, $h'_m=h_m$ if $m(j)\geq 0$ for any $j\in J$, $f_m'=f_m$, $g'_m=(t')^{-m}f_m$, $h'_m=(t')^{-2m}f_m$ if $m(j)\leq 0$ for any $j\in J$ and $m(j)<0$ for some $j\in J$, 
and $f'_m=g'_m=h'_m=0$ if otherwise.  
Here $(t')^m:=\prod_{j\in \Phi} (t'_j)^{m(j)}$ etc. 
Then, $y(t',\delta')\to y$ in $Y^{II}(p,\br,S)$ 
when $t'\to t$ and $\delta'\to \delta$. 
  
  We have to prove that the limit 
$(t_0,f_0,g_0,\ldots)$ of the image $(t'',f'',g'',\ldots)$ of $y(t', \delta')$ 
in the ambient product space 
satisfies the conditions (1)--(4) in 3.4.2. 
  First, it is easy to see $J:=J(t_0)=J(t)$. 
  (2) and (3) are deduced from the corresponding conditions on 
$(t'',f'',g'',\ldots)$. 
  (1) is also seen from the condition (2) on $(t'',f'',g'',\ldots)$. 
  For example, we show $(f_0)_m = 0$ unless $m(j) \leq 0$ for any $j \in J$. 
  We have $f''_m = (t'')^m g''_m$ for any $m \in \bZ^{\Phi}$. 
Since 
$t''=t'\b(\exp(g')\exp(k)\bar \br)a_1\mu(y(t',\delta'))$, 
if there is some $j \in J$ such that $m(j)>0$, 
the above equality implies $f''_m \to 0\cdot(\lim g''_m)=0$.  Hence 
we have $(f_0)_m = 0$. 
  Finally, (4) is seen as follows.  
  Let $k'$ be the element of $\Lie(K_{\bar \br})$ such that 
$\exp(g)=\exp(g_0)\exp(k')$ and $k'_m = 0$ unless $m(j)=0$ for any 
$j \in J$.
  Then, we have $\exp(k_0)=\exp(k')\exp(k)$. 
  Hence $k_0$ satisfies (4). 
\qed
\enddemo

\medskip
\proclaim{Lemma 3.4.14} 
There are an open neighborhood $U'$ 
of $0$ in $\fg_{\bR}(\gr^W)\times \fg_{\bR}(\gr^W)\times \fg_{\bR}(\gr^W)\times S$ and a morphism $Y^{II}(p,\br,S,U')\to 
B:= \bR_{\geq 0}^{\Phi}\times D(\gr^W)\times \bar \cL\times \spl(W) \times \prod_{W'\in \Phi} \spl(W')$ whose restriction to $Y^{II}_0(p,\br,S,U')$ coincides with the composite $\nu_{\bar \tau_p,\b}\circ \eta^{II}_{p, \br,S}$ 
$(3.2.6, 3.4.3)$. 
\endproclaim

\demo{Proof}
  It is enough to show that the composite map from 
$Y^{II}_0(p,\br,S,U')$ extends componentwise over the boundary. 
  The components except the last ones (BS-splittings) are easily treated. 
  For example, the first two were already treated in the proof of 3.4.13.
  The extendability of BS-splittings is reduced to 3.4.13. 
  In fact, let $W' \in \Phi$. 
  Then, by 3.4.10 and 3.4.13, it is sufficient to prove that the composite 
$Y^{II}_0(p,\br,R,S,U') \to Y^{II}_0(p,\br,S,U') \to \spl(W')$ extends to a real analytic map on   $Y^{II}(p,\br, R, S,U')$ under the assumption 
$R \sub \Lie(G_{\bR}(\gr^W)_{W'})$. 
Assuming this, we prove $f_m\in \Lie(G_{\bR}(\gr^W)_{W'})$ for  
any $(t,f,g,h,k,\delta,u)\in Y^{II}(p,\br,R,S,U')$ and any 
$m\in \bZ^{\Phi}$. This is clear if $m(W')\leq 0$. If $m(W')\geq 0$, since $f_m+h_{-m}\in R\sub \Lie(G_{\bR}(\gr^W)_{W'})$ and $h_{-m}\in \Lie(G_{\bR}(\gr^W)_{W'})$, we have $f_m\in \Lie(G_{\bR}(\gr^W)_{W'})$.  
  Thus, $\exp(f)$ belongs to $G_{\bR}(\gr^W)_{W'}$ so that the concerned component is 
$\spl^{\BS}_{W'}(\exp(f)\bar \tau_p(t)\exp(k)\bar \br)=\exp(f)\spl^{\BS}_{W'}(\bar \br)\gr^{W'}
\exp(f)^{-1}$, which is real analytically extends over the boundary. 
\qed\enddemo

\proclaim{Lemma 3.4.15} 
There exist open neighborhoods $U'' \sub U'$ of $0$ in 
$\fg_{\bR}(\gr^W) \times \fg_{\bR}(\gr^W) \times \fg_{\bR}(\gr^W) \times S$ 
such that, for any $y\in Y^{II}(p,\br,S,U'')$, there exists 
$y_1\in Y_1^{II}(p,\br,S) \cap Y^{II}(p,\br,S, U')$ such 
that $(y_1,y)$ belongs to the closure of $Y^{II}_0(p,\br,S)\times_D 
Y^{II}_0(p,\br,S)$ in $Y_1^{II}(p,\br,S) \times Y^{II}(p,\br,S)$. 
\endproclaim

\demo{Proof}  
  For any subset $J$ of $\Phi$, 
take $R=R_J$ as in 3.4.1 such that 
$\Lie(G_{\bR}(\gr^W)_{J, u}) \sub R_J$. Here $G_{\bR}(\gr^W)_{J, u}$ denotes the unipotent part of $G_\bR(\gr^W)_J$. 
  For this $R=R_J$, let $U_J$ be the neighborhood $U$ in 3.4.8, and 
let $U'$ be a neighborhood of 
$0$ in 
$\fg_{\bR}(\gr^W) \times \fg_{\bR}(\gr^W) \times \fg_{\bR}(\gr^W) \times S$ 
such that $Y^{II}(p,\br,S, U')$ 
is contained in $(\eta_{p,\br,S})^{-1}(\bigcap_J D(U_J))$.

  Let $y=(t,f,g,h,k,\delta, u) \in Y^{II}(p,\br,S, U')$. 
For $t'\in \bR_{>0}^{\Phi}$ such that $t'_j=t_j$ for any $j\in \Phi\smallsetminus J$ with $J=J(t)$  
and for $\delta'\in L$, consider $y(t',\delta')$ 
in the proof of 3.4.13.

  Let $R=R_{J(t)}$.
Then, for any $(t',\d')$ which is sufficiently near to $(t,\d)$, 
the point $y_1(t',\delta'):=b_{R,S}(\eta^{II}_{p,\br,S}(y(t',\delta')))$ 
is well-defined and 
$(y_1(t',\delta'), y(t',\delta'))\in Y^{II}_0(p,\br,S)\times_D Y^{II}_0(p,\br,S)$. 
  Furthermore, $y_1(t',\delta')$ converges to an element $y_1$ 
of $Y^{II}(p,\br,S)$ 
when $t'\to t$ and $\delta'\to \delta$ by 3.4.13. 
  We show that the limit $y_1=(t_0, f_0,g_0,h_0,\ldots)$ belongs 
to $Y_1^{II}(p,\br,S)$, that is, $(h_0)_m=0$ if $m(j)\geq0$ for any $j \in J(t)=
J(t_0)$ and if $m(j)>0$ for some $j \in J(t_0)$. 
  Fix such an $m$. 
  Then, we have $\fg_{\bR}(\gr^W)_{-m} \sub R_J$ and 
$\fg_{\bR}(\gr^W)_m \cap R_J=\{0\}$.
  Hence the property $(f_0)_{-m} + (h_0)_m \in R_J$ 
implies $(h_0)_m=0$.

  Finally, for a sufficiently small $U'' \sub U'$, the above correspondence 
$y \mapsto y_1$ sends $Y^{II}(p,\br,S,U'')$ into $Y^{II}(p,\br,S,U')$. 
\qed
\enddemo

\proclaim{Lemma 3.4.16} 
{\rm (i)} On the intersection of $Y_1=Y_1^{II}(p,\br,S)$ 
and $Y(U'):=
Y^{II}(p,\br,S,U')$, 
the map $Y(U')\to B$ in $3.4.14$ coincides with the restriction of the composite 
$Y_1\to D_{\SL(2)}^{II}(\Phi)\to B$.

\medskip

{\rm (ii)}
  For a sufficiently small $U'$, 
the image of $Y(U')\to B$ in $3.4.14$ is contained in the image of $D_{\SL(2)}^{II}(\Phi)$.
\endproclaim

\demo{Proof} (i) 
follows from (6) of 3.4.12. 
(ii) follows from (i) and Lemma 3.4.15. \qed
\enddemo

\medskip

\proclaim{Lemma 3.4.17} 
Let $U$ be a sufficiently small open neighborhood of $\bar \br$ in $D(\gr^W)$ and let $D_{\SL(2)}^{II}(U)$ be the inverse image of $U$ under $D_{\SL(2)}^{II}(\Phi)@>\mu>>D(\gr^W)$. Let $q \in D_{\SL(2)}^{II}(U)$, and 
let $\br_q$ be a point on the torus orbit associated to $q$. 
Then, the limit 
$\lim\limits_{t \to 0^{\cW(q)}}b(\tau_q(t)\br_q)$ exists 
in $Y^{II}(p,\br,R,S)$ and independent 
of the choice of $\br_q$.  
\endproclaim
\medskip

\demo{Proof} 
  We reduce this to 3.4.13.
  First, by (7) in 3.4.12, we may assume that $q$ has the form 
$p(J,\br,z,\d,u)$ such that $\br_q$ is the point in 3.4.12 (5). 
  Hence it is the image of some 
$y_1=(s,f,g,h,k,\d,u) \in Y_1$ by $\eta^{II}_{p,\br,S}$ in 3.4.12.
  Then, $\tau_q(t)\br_q$ 
is the image of $y_1(t):=
(t',\,f,\, \Ad(\bar \tau_p(t'))^{-1}f,\, 
\Ad(\bar \tau_p(t'))^{-2}f,\, 
k,\d'',u)$, 
where $t' \in \bR_{>0}^{\Phi}$ such that $t'_j = t_j$ for any $j \in J$ 
and $t'_j =s_j$ for any $j \in \Phi  
\smallsetminus J$ and $\d''=\d$ if $\d \in L$ and 
$\d''=t_{W} \circ \d'$ for $\d' \in L$ in 3.4.12 (5) if $\d \in \bar L
\smallsetminus L$. 
  Since $y_1(t)$ converges to $y_1$, the sequence 
$b(\tau_q(t)\br_q)$ converges to the image of $y_1$ by $\xi$ in 3.4.13.
    The last independency is clear. 
\qed
\enddemo

Denote this limit by $b(q)$.  
Thus, $b$ in 3.4.9 is extended to a map $D_{\SL(2)}^{II}(U) \to Y^{II}(p,\br,R,S)$.   

\medskip

{\bf 3.4.18.} 
We prove 3.4.4. (i) of 3.4.4  follows from (ii) of 3.4.16.
We prove (ii) of 3.4.4. We first describe the idea of the proof. 

Locally on $Y(R,S):=Y^{II}(p,\br, R, S)$, we define an object $X$ of $\cB_\bR(\log)$ which contains $Y(R,S)$ having the following properties.  

\medskip

(1) The morphism $Y(R,S)\to B$ (defined locally) extends to some 
explicit morphism $X\to B$ (locally). 
(It will be explained in 3.4.19.)

\medskip
(2) As an object of $\cB_\bR(\log)$, $X$ is isomorphic to the product $\bR_{\geq 0}^{\Phi}\times$  (a real analytic manifold) $\times \bar L$. 
Hence, for any $x\in X$, 
the local ring $\cO_{X,x}$ is isomorphic to the ring of convergent power series in $n$ variables over $\bR$ for some $n$. 
Note that $Y(R,S)$ need not have this last property (because $Y(R,S)$ can have a singularity of the style $t_1^2x=t_2y$), and this is the reason why we use $X$ here.  

\medskip

(3) $\cO_X|_{Y(R,S)}\to \cO_{Y(R,S)}$ is surjective. Here $\cO_X|_{Y(R,S)}$ 
is the inverse image of $\cO_X$ on $Y(R,S)$.

\medskip

(4) $\cO_B|_X\to \cO_X$ is surjective. Here $\cO_B|_X$ denotes the inverse image of $\cO_B$ on $X$. 

\medskip

Though (3) is shown easily, (4) is not. 
But by the property of the local rings explained in (2), 
the property (4) is reduced to the surjectivity 
of $m_{B,y}/m_{B,y}^2\to m_{X,x}/m_{X,x}^2$, where $x\in X$ and $y$ is the image of $x$ in $B$. 
This is the injectivity of the map of tangent spaces $T_x(X)\to T_y(B)$, where $T_x(X)$ and $T_y(B)$ are $\bR$-linear duals of $m_{X,x}/m_{X,x}^2$ and $m_{B,y}/m_{B,y}^2$, respectively,  and this injectivity will be explained in 3.4.19. 

By (3) and (4), we have the surjectivity of $\cO_B|_{Y(R,S)}\to \cO_{Y(R,S)}$. Since $Y(R,S)\to B$ factors (locally) as $Y(R,S)\to A\to B$ by 3.4.16 (ii), 
where $A:= D_{\SL(2)}^{II}(\Phi)$, we see that the map $\cO_A|_{Y(R,S)}\to \cO_{Y(R,S)}$ is surjective. 

Since the map $Y(R,S)\to A$ has the inverse map $A\to Y(R,S)$ (locally) by 3.4.17, $Y(R,S)\to A$ is bijective locally. 

Since $\cO_A|_{Y(R,S)}\to \cO_{Y(R,S)}$ is injective (they are sub-sheaves of the sheaves of functions), we have $(Y(R,S), \cO_{Y(R,S)})\cong (A, \cO_A)$ locally. It is easy to see that this isomorphism preserves the log structures with sign. 

\medskip

{\bf 3.4.19.} 
We give the definition of $X$ and the proof of (4) in 3.4.18.

\medskip
Actually $X$ is constructed at each point of $Y(R,S)$. We give the construction at $\tilde p=(0^{\Phi},0,0,0,0,\delta(\br),0)\in Y(R,S) $ and the proof of (4) for the tangent space at $\tilde p$. The general case is similar.

We define the set $X$ to be the subset of 
 $E:=\bR^\Phi_{\ge0} \times \fg_\bR(\gr^W)\times \fg_{\bR}(\gr^W)\times \fg_\bR(\gr^W) \times S \times \bar L \times \fg_{\bR,u}$ consisting of 
  all elements $(t, f, g, h, k, \delta,u)$ satisfying the following conditions (1)--(3). 

\medskip

(1) If $m\in \bZ^{\Phi}$ and $m(j)\geq 0$ for any $j\in \Phi$, then 
$f_m=t^mg_m$ and $g_m=t^mh_m$. 
Here $t^m:=\prod_{j\in \Phi} t_j^{m(j)}$. 
\medskip

(2) If $m\in \bZ^{\Phi}$ and $m(j)\leq 0$ for any $j\in \Phi$, then 
$h_m=t^{-m}g_m$ and $g_{m}=t^{-m}f_{m}$. 

\medskip

(3) $g\in R$ and $f_m+h_{-m}\in R$ for all $m\in \bZ^{\Phi}$. 
\medskip

Define the structure on $X$ as an object of $\cB_\bR(\log)$ by using the embedding $X\subset E$ just as we defined the structure of $Y(R,S)$ as an object of $\cB_\bR(\log)$ by using the embedding $Y(R,S)\sub E$ in 3.4.2. 
Then it is clear that $X$ is isomorphic to a product $\bR_{\geq 0}^{\Phi}\times $ (a real analytic manifold) $\times \bar L$ as an object of $\cB_\bR(\log)$. 

  We give a morphism $X \to B$ which extends $Y(R,S) \to B$ and prove 
the property (4) in 3.4.18 for it.
  We define the morphism componentwise. 
  Let $X_0$ be the inverse image of $\bR_{>0}^{\Phi} \times L$ by 
the natural map $X \to \bR_{\ge0}^{\Phi} \times \bar L$.
  First, we define $X_0 \to 
B':=\bR_{\geq 0}^{\Phi}\times D(\gr^W)\times \bar \cL\times \spl(W)$ 
as the projection after $\nu_{\bar \tau_p,\b}\circ \eta$, where 
$\eta$ sends $(t,f,g,h,k,\delta,u)$ to 
$\exp(u)s_{\br}\theta(d\bar \br, \Ad(d)\delta)$ with $d=\bar \tau_p(t)\exp(g)\exp(k)$.
  Then this map $X_0 \to B'$ extends to $X \to B'$, as seen 
easily in the same way as in 3.4.14.
  Next, for each $j=W' \in \Phi$, we give an extension to $\spl(W')$. 
  Define $X_0 \to \spl(W')$ as follows. 
  Consider the decomposition $\fg_{\bR}(\gr^W)=
\fg_{\le} \oplus \fg_{>0}$, where 
$\fg_{\le} = \sum_{m(j)\le0}\fg_{\bR}(\gr^W)_m$ and 
$\fg_{>} = \sum_{m(j)>0}\fg_{\bR}(\gr^W)_m$. 
  Then there are a neighborhood $V_1$ of $0$ in $\fg_{\bR}(\gr^W)$ and 
a real analytic map $(c_{\le},\, c_{>}): 
V_1 \to \fg_{\le} \times \fg_{>}$ such that 
for any $g \in V_1$, we have $\exp(g)=\exp(c_{\le}(g))\exp(c_{>}(g))$. 
  Further, let $M$ be an $\bR$-subspace of 
$\sum_{m(j)\ge0}\fg_{\bR}(\gr^W)_m$ containing $\fg_{>}$ 
such that $\fg_{\bR}(\gr^W)=M \oplus \Lie(K_{\bar \br})$.
  Then, there are a neighborhood $V_2$ of $0$ in $\fg_{\bR}(\gr^W)$ and 
a real analytic map $(c'_1, c'_2): V_2 \to M \times \Lie(K_{\bar \br})$ 
such that for any $g' \in V_2$, we have 
$\exp(g')=\exp(-Cc'_1(g'))\exp(c_2'(g'))$, 
where $C$ is the Cartan involution at $\bar \br$. 
  We define $X_0 \to \spl(W')$ (locally) as 
$\spl_{W'}^{\BS}(\exp(c_{\le}(f))\exp
(-C(c'_1(c_{>}(h))))
\bar\br)$. 
  This extends to $X \to \spl(W')$ and gives an extension of $Y(R,S) \to \spl(W')$, 
since $\Int(\bar \tau_p(t))\exp(g)=\exp(f)$ etc. on $Y(R, S)$.

We prove the surjectivity of $\cO_B|_X\to \cO_X$. We write the proof of the surjectivity for the stalk at $\tilde p$. (The general case is similar.) 
It is sufficient to prove the injectivity of $T_{\tilde p}(X)\to T_q(B)$, 
where $q$ denotes the image of $\tilde p$ in $B$.

 The first tangent space is identified with the vector subspace $V$ of $\bR^{\Phi} \times \fg_\bR(\gr^W) \times \fg_\bR(\gr^W) \times \fg_\bR(\gr^W) \times S \times L\times \fg_{\bR,u}$
consisting of all elements
$(t, f, g, h, k, \delta,u)$ satisfying the following conditions (1) and (2).

\medskip

(1) $f_m=g_m=0$ if $m(j)\geq 0$ for any $j\in \Phi$, and $g_m=h_m=0$ if $m(j)\leq 0$ for any $j\in \Phi$.

\medskip 

(2)  $g\in R$ and $f_m+h_{-m}\in R$ for all $m\in \bZ^\Phi$. 

\medskip

The injectivity of the map of tangent spaces in problem is reduced 
to the injectivity of the following map.

\medskip

$$
V\to \bR^\Phi \times R  \times S\times L\times \fg_{\bR,u}\times  (\tsize\prod_{j\in \Phi} \fg_\bR(\gr^W)),
$$
$$
(t, f,g,h,k,\delta,u) \mapsto (t, g,  k, \delta, u, (v_j)_{j\in \Phi}),
$$
$$
\text{where}\quad
v_j=\tsize\sum_{m(j)<0}  (f_m-C(h_{-m})).
$$ 

Assume that the image of $(t, f, g, h, k, \delta,u)\in V$ under this map is zero. Then clearly we have $t=g=k=\delta=u=0$. We have also 

\medskip

(i) If $m(j)<0$ for some $j\in \Phi$, then $f_m=h_{-m}=0$.

\medskip

Indeed, if $m(j)<0$ for some $j\in \Phi$, then
$f_m-C(h_{-m})=0.$ 
Since $h_{-m}+C(h_{-m})\in \Lie(K_{\bar \br})$,  $f_m+h_{-m}\in R\cap \Lie(K_{\bar\br})=0$, and consequently we have (i).

This shows that if $m(j)<0$ and $m(j')>0$ for some $j, j'\in \Phi$, then $f_m=f_{-m}=h_m=h_{-m}=0$. If $m(j) \leq 0$ for any $j\in \Phi$ and if $m(j)<0$ for some $j\in \Phi$, then $f_m=h_{-m}=0$ by (i) and $f_{-m}=h_m=0$ by the definition of $V$. If $m(j)\geq 0$ for any $j\in \Phi$, we have similarly $f_m=h_m=f_{-m}=h_{-m}=0$. 
\qed

\medskip

{\bf 3.4.20.} 
We prove 3.4.6. We deduce it from 3.4.4 (ii) as follows.

Let $\Psi\in \cW$ and let $\Phi$ be the image of $\Psi$ in $\overline{\cW}$ (3.2.2). 
Take a distance to $\Phi$-boundary $\b$. 

Let $E$ be the subset of $\bR_{\geq 0}^{\Psi} \times \fg_{\bR,u} \times \fg_{\bR,u}$ consisting of all elements $(t, u, v)$ satisfying the conditions 
(5) and (6) (resp.\ 
(5${}^\prime$) and  (6${}^\prime$)) in 3.4.5 
in the case where $W \not\in \Psi$ 
(resp.\ $W \in \Psi$). 
We regard $E$ as an object of $\cB_\bR(\log)$ in the similar way as in the case of $Y^{II}(p,\br, R,S)$ (3.4.2). 

Assume first $W\notin \Psi$. Let $D_{\SL(2)}^{II}(\Phi)'$ be the open set of $D_{\SL(2)}^{II}(\Phi)$ consisting of all elements $q$ such that $W\notin \cW(q)$ (this condition is equivalent to the condition that the $\bar \cL$-component of $\nu_{\bar \tau_p,\b}(q)$ (3.2.6 (ii)) is contained in $\cL$). 
Then $D_{\SL(2)}^I(\Psi)$ is the fiber product of
$$
D_{\SL(2)}^{II}(\Phi)'\to \bR_{\geq 0}^{\Phi}\times \fg_{\bR,u}\leftarrow E
$$
in $\cB'_\bR(\log)$, where the first arrow is given by $x\mapsto (\b(x), u)$ with $\spl_W(x)=\exp(u)s_{\br}$, the second arrow sends $(t,u,v)$ to $(t,u)$, and the morphism $D_{\SL(2)}^I(\Psi)\to E$ is given by $x\mapsto (\b(x), u, v)$ with $\spl_{W}(x)=\exp(u)s_{\br}$ and $\spl_W(y)=\exp(v)s_\br$ for the $D$-component $y$ of $\nu_{\tau_p, \b}$ (3.2.6 (i)). 
Since $Y^{I}(p,\br, R, S)$ is 
the fiber product of $Y^{II}(p,\br, R, S)\to \bR_{\geq 0}^{\Phi}\times \fg_{\bR,u}\leftarrow E$, 3.4.6 is reduced to 3.4.4.

Next assume $W\in \Psi$. 
Let $\b_0:\bar \cL\smallsetminus \{0\}\to \bR_{>0}$ be a real analytic function such that  $\b_0(a\circ \delta)=a\b_0(\delta)$ for any $a\in \bR_{>0}$ and $\delta\in \bar\cL\smallsetminus \{0\}$. 
Denote the composite $D_{\SL(2)}^{II}(\Phi)_{\nspl}\to \bar \cL\smallsetminus \{0\}\to \bR_{\geq 0}$ also by $\b_0$, where the first arrow is the $\bar \cL$-component of $\nu_{\tau_p,\b}$ (3.2.6 (ii)). 
Then $(\b_0,\b): D\to \bR_{>0}^{\Psi}=\bR_{>0}\times \bR_{>0}^{\Phi}$ is a distance to $\Psi$-boundary. 
As an object of $\cB'_\bR(\log)$, $D_{\SL(2),\nspl}^I(\Psi)$ is the fiber product of 
$$
D_{\SL(2),\nspl}^{II}(\Phi)\to \bR_{\geq 0}^{\Psi}\times \fg_{\bR,u}\leftarrow E,
$$
where the first arrow is given by $x\mapsto ((\b_0,\b)(x), u)$ with $\spl_W(x)=\exp(u)s_{\br}$. On the other hand, if we denote by $Y^*(p,\br, R, S)_{\nspl}$
($*=I, II$) the open set of $Y^*(p, \br, R, S)$ consisting of all elements satisfying $\delta\neq 0$, $Y^I(p,\br, R, S)_{\nspl}$ is the fiber product of
$$
Y^{II}(p,\br,R,S)_{\nspl}\to \bR_{\geq 0}^{\Psi}\times \fg_{\bR,u}\leftarrow 
E
$$
in $\cB'_\bR(\log)$, where the first arrow is given by $(t,f,g,h,k,\delta,u)\mapsto ((a, t), u)$ for $\delta=a\circ \delta^{(1)}$ with $\delta^{(1)}\in L^{(1)}$ (3.4.5).
From these facts, 3.4.6 is reduced to 3.4.4 also in the case $W\in \Psi$. 

3.4.6 is proved.

\medskip

{\bf 3.4.21.} 
We prove 3.2.10.

We first prove (ii) of 3.2.10. Let $\Phi\in \overline \cW$. We prove 

\medskip

{\bf Claim 1.}  
{\it For $\Phi'\sub \Phi$, the inclusion map $D_{\SL(2)}^{II}(\Phi')\to D_{\SL(2)}^{II}(\Phi)$ is an open immersion in $\cB'_\bR(\log)$.}

\medskip

 Let $\a$ be a splitting of $\Phi$ and let $\b$ be a distance to $\Phi$-boundary. Since $D_{\SL(2)}^{II}(\Phi')$ is the inverse image of $\{t\in \bR_{\geq 0}^{\Phi}\;|\; t_j\neq 0\;\text{if}\;j\in \Phi\smallsetminus \Phi'\}$ under the map $\b: D_{\SL(2)}^{II}(\Phi)\to \bR_{\geq 0}^{\Phi}$, it is an open subset of $D_{\SL(2)}^{II}(\Phi)$. 
Let $\alpha': \bG_{m,\bR}^{\Phi'}\to 
\Aut(\gr^W)$ be the $\Phi'$-component of $\alpha$ and let $\b': D(\gr^W)\to \bR_{>0}^{\Phi'}$ be the $\Phi'$-component of $\b: D(\gr^W)\to \bR_{> 0}^{\Phi}$. 
Then we have a commutative diagram
$$
\matrix D_{\SL(2)}^{II}(\Phi') &\to &\bR_{\geq 0}^{\Phi'} \times D(\gr^W)'\times \bar \cL\\
&&\\
\cap & & \downarrow\\
&&\\
D_{\SL(2)}^{II}(\Phi) & \to & \bR_{\geq 0}^{\Phi} \times D(\gr^W) \times \bar \cL,\endmatrix
$$
where $D(\gr^W)'=\{x\in D(\gr^W)\;|\;\b'(x)=1\}$, the upper horizontal arrow is induced by $(\a',\b')$ as in 3.2.6, the lower horizontal arrow is induced by $(\a,\b)$ as in 3.2.6, and the right vertical arrow sends 
$(t, x,\delta)\in \bR_{\geq 0}^{\Phi'}\times D(\gr^W)'\times \bar \cL$ to 
$((t,\b(x)), \a\b(x)^{-1}x, \Ad(\a\b(x))^{-1}\delta)$. Here by the fact $\b(x)_j=1$ for any $j\in \Phi'$, we regard $(t, \b(x))$ as an element of $\bR_{\geq 0}^{\Phi'}\times \bR_{>0}^{\Phi\smallsetminus\Phi'}\subset \bR_{\geq 0}^{\Phi}$. 
From this, we obtain
\medskip

{\bf Claim 2.}  
{\it Let $D_{\SL(2)}^{II,\Phi}(\Phi')$ be the set $D_{\SL(2)}^{II}(\Phi')$  endowed with the structure of an object of $\cB_\bR'(\log)$ as an open set of $D_{\SL(2)}^{II}(\Phi)$. Then the canonical inclusion map 
$D_{\SL(2)}^{II,\Phi}(\Phi') \to D_{\SL(2)}^{II}(\Phi')$ is a morphism
 in $\cB'_\bR(\log)$. 
 This morphism is an isomorphism if and only if, for any $W'\in \Phi$, the composite $D_{\SL(2)}^{II}(\Phi') \to D_{\SL(2)}^{II}(\Phi) \to \spl(W')$, where the last arrow is induced by $\spl_{W'}^{\BS}$,  is a morphism in $\cB'_\bR(\log)$.} 
\medskip

By Claim 2 and  by 3.4.4, for the proof of Claim 1, it is sufficient to prove

\medskip

{\bf Claim 3.} 
{\it Let $p'\in D_{\SL(2)}^{II}(\Phi)$ and let $\Phi'=\overline{\cW}(p')\sub \Phi$. Let $\br'$ be a point on the torus orbit associated to $p'$. 
Then, for a sufficiently small open neighborhood $U$ of $(0^{\Phi'},0,0,0,0,\delta(\br'), 0)$ in $Y^{II}(p',\br', S)$ ($S$ is taken for $\br'$), the composite $U\to D_{\SL(2)}^{II}(\Phi') \to D^{II}_{\SL(2)}(\Phi)\to \spl(W')$ is a morphism of $\cB_\bR'(\log)$.} 
\medskip

We prove Claim 3. 
Take $p\in D_{\SL(2)}^{II}(\Phi)$ such that $\Phi=\overline{\cW}(p)$. 
Let $\a=\bar \tau_p$ and take a distance to $\Phi$-boundary $\b$ such that $\b(K_{\bar \br}\cdot\bar\br)=1$. 
  Note that such a $\b$ exists (cf.\ \cite{KU2}, 4.12). 
For each $w \in \bZ$, let $Q(w) \in \cW(\gr^W_w)$ be the image of $\Phi$, and let $Q = (Q(w))_w$.
Let $D_{\SL(2)}(\gr^W)(Q) = \prod_w D_{\SL(2)}(\gr^W_w)(Q(w))$.
Let $\bar\mu: D_{\SL(2)}(\gr^W)(Q) \to D(\gr^W)$ be the extension of $D(\gr^W) \to D(\gr^W)$, $x\mapsto \a\b(x)^{-1}x$, induced by 3.2.6 (ii). 
Let $\a'=\bar \tau_{p'}$. 
We first prove

\medskip

\noindent
{\bf Claim 4.} 
{\it There exists $y\in G_\bR(\gr^W)_{W'}$ such that $\bar\mu(y^{-1}\bar p')\in K_{\bar \br}\cdot \bar \br$, where $\bar p' = p'(\gr^W)$.} 

\medskip

In fact, by Claim 1 in 6.4.4 of \cite{KU3}, there are $z\in G_\bR(\gr^W)_{\Phi'}$ and $k\in K_{\bar \br}$ such that
$\a'=\Int(z)(\a_{\Phi'})$ and $\bar\br'=zk\bar\br$, where $\a_{\Phi'}$ is the restriction of $\a$ to $\Phi'$. Write $z=z_0z_u$, where $z_0$ commutes with $\a_{\Phi'}(t)$ ($t\in (\bR^\times)^{\Phi'}$) and $z_u\in G_{\bR}(\gr^W)_{\Phi',u}$. 
We can write $z_0=yk_0$, where $y$ and $k_0$ commute with $\a_{\Phi'}(t)$ ($t\in (\bR^\times)^{\Phi'}) $, $y\in G_\bR(\gr^W)_{W'}$, and $k_0\in K_{\bar\br}$. 
We have $\bar\mu(y^{-1}\bar p')=k_0k\bar\br$. 
In fact, since $\bar p'=\lim \a'(t)\bar \br'=\lim z\a(t)k\bar \br$, 
$\bar\mu(y^{-1}\bar p')$ is the limit of
$\bar\mu(y^{-1}z\a(t)k\bar \br)= \bar\mu(\a(t)y^{-1}z_tk\bar \br) = 
\bar\mu(y^{-1}z_tk\bar \br)$, where $z_t=\bar \tau_p(t)^{-1}z\bar \tau_p(t)$, which converges to $\bar\mu(y^{-1}z_0k\bar \br)=\bar\mu(k_0k\bar \br)=k_0k\bar \br \in 
K_{\bar\br}\cdot \bar\br$. 

\medskip

Let  $y$ be as in Claim 4. 
Then, for $q\in D$ near $p'$ in $D_{\SL(2)}^{II}(\Phi')$, $\spl^{\BS}_{W'}(\bar q)= y \spl^{\BS}_{W'}(y^{-1}\bar q) y(\gr^{W'})^{-1}$, 
where $\bar q = q(\gr^W)$.
We denote the right-hand side of the last equation by $\Int(y)\spl^{\BS}_{W'}(y^{-1}\bar q)$.
From this, we may replace $\bar p'$ by $y^{-1}\bar p'$ and hence we may assume $\bar\mu(\bar p')\in K_{\bar \br} \cdot \bar\br$.

Take an $\bR$-subspace $V$ of $\Lie(G_{\bR}(\gr^W)_{W'})$ such that $\fg_\bR(\gr^W)= V\oplus \Lie(K_{\bar \br})$.  
For $q\in D$ near $p'$ in $D_{\SL(2)}^{II}(\Phi')$, write
$\bar\mu(\bar q)\in \exp(v(\bar q))\cdot K_{\bar \br}\cdot \bar \br$ with $v(\bar q)\in V$ and write $f(\bar q)=\Int(\a\b(\bar q))(\exp(v(\bar q)))\in G_\bR(\gr^W)_{W'}$. 
Then, since $\bar q=\a\b(\bar q)\bar\mu(\bar q)$, we have
$$
\spl_{W'}^{\BS}(\bar q)=\Int(f(\bar q))(\spl_{W'}^{\BS}(\a\b(\bar q)\bar \br))= \Int(f(\bar q))(\spl_{W'}^{\BS}(\bar \br)).
$$ 
Here the last equality follows from  $\Int(\a(t))\spl^{\BS}_{W'}(\bar \br)=
\spl^{\BS}_{W'}(\bar \br)$ for any $t$. 
By 3.4.4 and the real analycity of $a_1$ in 3.3.11, $v(\bar q)$ extends 
over the boundary, and hence so does $f(\bar q)$, that is, 
for a sufficiently small open neighborhood $U$ of
$(0^{\Phi'},0,0,0,0,\delta(\br'),0)$ in 
$Y^{II}(p',\br', S)$, there is a morphism $U\to G_{\bR}(\gr^W)_{W'}$ which is compatible with the map 
$Y^{II}_0(p',\br', S) \to G_{\bR}(\gr^W)_{W'}$ induced by $f$.
  Hence $\spl^{\BS}_{W'}$ extends over the boundary. 
This completes the proof of Claim 3, and hence the proof of Claim 1.

 By Claim 1, on $D_{\SL(2)}$, there is a unique structure as an object of $\cB'_\bR(\log)$ for which each $D_{\SL(2)}^{II}(\Phi)$ ($\Phi\in \overline \cW$) 
is open and whose restriction to $D_{\SL(2)}^{II}(\Phi)$ coincides with the 
structure of $D_{\SL(2)}^{II}(\Phi)$ as an object of $\cB'_\bR(\log)$. 
 By 3.4.4, this object $D_{\SL(2)}^{II}$ of $\cB_\bR'(\log)$ belongs to  $\cB_\bR(\log)$.

Next, (i) of 3.2.10 follows from (ii) of 3.2.10 and 3.4.6.

We prove (iii) of 3.2.10. It is clear that the identity map of $D_{\SL(2)}$ is a morphism
$D_{\SL(2)}^I\to D_{\SL(2)}^{II}$ in $\cB_\bR(\log)$ 
and that the log structure with sign on $D^I_{\SL(2)}$ is the pull-back of that 
of $D^{II}_{\SL(2)}$. 
It is also clear that, in the pure case, this morphism $D_{\SL(2)}^I\to D_{\SL(2)}^{II}$ is an isomorphism.

It remains to prove that
 in the pure case, the topology of $D_{\SL(2)}$ defined in \cite{KU2} coincides with the topology defined in this paper. 
 
  Assume that we are in the pure case.

The topology of $D_{\SL(2)}$ defined in  \cite{KU2} is characterized by the following properties (1) and (2) (\cite{KU3}). 

\medskip

(1) For any $\Psi\in \cW$, $D^I_{\SL(2)}(\Psi)$ is open and is a regular space.

\medskip

(2) Let $p\in D_{\SL(2)}$, let $\br$ be a point on the torus orbit associated to $p$, and let $\Psi=\cW(p)$. 
Then, for a directed family $(p_\lam)_\lam$ of points of $D$, $(p_\lam)_\lam$  converges to $p$ in $D_{\SL(2)}(\Psi)$ if and only if there exist $t_\lam \in \bR_{>0}^{\Psi}$, $g_\lam\in \fg_\bR$, $k_\lam\in \Lie(K_\br)$ such that $p_\lam=\tau_p(t_\lam)\exp(g_\lam)\exp(k_\lam)\br$, $t_\lam\to 0^{\Psi}$ in $\bR_{\geq 0}^{\Psi}$, $\Ad(\tau_p(t_\lam))^j(g_\lam)\to 0$ for $j=\pm 1, 0$, and $k_\lam\to 0$.
\medskip

It is sufficient to prove that the topology of $D_{\SL(2)}^{II}$ ($=$ the topology of $D_{\SL(2)}^I$) in this paper satisfies these (1) and (2). (1) is clearly satisfied. We prove (2). 

Assume $p_\lam\to p$ for the topology of this paper. 
By 3.4.4 (ii), for some 
$\tilde p_\lam=(t_\lam, f_\lam, g_\lam, h_\lam, k_\lam)\in Y_0(p,\br,R,S)\subset  \bR_{>0}^{\Psi}\times \fg_\bR\times \fg_\bR \times \fg_\bR\times \Lie(K_\br)$
such that
$p_\lam=\tau_p(t_\lam)\exp(g_\lam)\exp(k_\lam)\br$, we have
$\tilde p_\lam \to (0^{\Psi},0,0,0,0)$ in 
$Y(p,\br, R,S)$. Since $f_\lam=\Ad(\tau_p(t_\lam))(g_\lam)$ and $h_\lam=\Ad(\tau_p(t_\lam))^{-1}(g_\lam)$, 
we have $t_\lam\to 0^{\Psi}$, $\Ad(\tau_p(t_\lam))^j(g_\lam)\to 0$ for $j=\pm 1, 0$, and $k_\lam\to 0$. 
Conversely, assume $p_\lam=\tau_p(t_\lam)\exp(g_\lam)\exp(k_\lam)\br$ for some $t_\lam \in \bR_{>0}^{\Psi}$, $g_\lam\in \fg_\bR$, $k_\lam\in \Lie(K_\br)$ such that $t_\lam\to 0^{\Psi}$, $\Ad(\tau_p(t_\lam))^j(g_\lam)\to 0$ for $j=\pm 1, 0$, and $k_\lam\to 0$.
Then if we put $f_\lam=\Ad(\tau_p(t_\lam))(g_\lam)$ and $h_\lam=\Ad(\tau_p(t_\lam))^{-1}(g_\lam)$, then $(t_\lam,f_\lam,g_\lam,h_\lam,k_\lam)$ converges to $(0^{\Psi},0,0,0,0)$ in $Y(p, \br, S)$. By Theorem 3.4.4 (i), this shows that $\tau_p(t_\lam)\exp(g_\lam)\exp(k_\lam)\br$ converges to $p$ for the topology of this paper. 
\medskip

{\bf 3.4.22.} 
In 3.4.23 and 3.4.27 below, we give local descriptions of $D_{\SL(2)}^{II}$ and $D_{\SL(2)}^I$ as topological spaces, respectively. 
Compared with the real analytic local descriptions in 3.4.4 and 3.4.6, we have simpler descriptions here.

We define a topological space
$Z^{II}_{\text{top}}(p,R)$ as the subspace of
$\bR_{\geq 0}^{\Phi} \times R$ consisting of all elements $(t,a)$ satisfying the following condition (1).

\medskip
(1) Let $m\in \bZ^{\Phi}$. 
Then $a_m=0$ unless either $m(j)\geq 0$ for all $j\in J$ or $m(j)\leq 0$ for all $j\in J$. 
 
 \medskip
We define a topological space  
$Y_{\text{top}}^{II}(p, \br, R,S)$ as the subspace  of $Z^{II}_{\text{top}}
(p, R) \times S\times \bar L\times \fg_{\bR,u}$ consisting of all elements $(t, a, k, \delta,u)$ ($(t,a) \in Z_{\text{top}}^{II}(p,R)$, $k\in S$, $\delta\in \bar L$, $u\in \fg_{\bR,u}$) such that $(t, k)$ satisfies the condition (4) in 3.4.2.
Let $Y_{0,\text{top}}^{II}(p,\br,R, S)$ be the open set 
$\bR_{>0}^{\Phi}\times R\times S\times L\times \fg_{\bR,u}$ 
of $Y_{\text{top}}^{II}(p,\br, R, S)$, 
and let 
$$
\eta_{p,\br,R,S,\text{top}}^{II}: Y_{0,\text{top}}^{II}(p,\br,R, S)\to D
$$ 
be the continuous map 
$$
(t, a, k, \delta, u)\mapsto
\exp(u)s_{\br}\theta(d\bar \br, \Ad(d)\delta)
$$
$$
\text{with}\;\; d=\bar \tau_p(t)\exp(\tsize\sum_{m\in \bZ^{\Phi}}\; g_m/(t^m+t^{-m}))\exp(k).
$$
Here $t^m=\prod_{j\in \Phi} t_j^{m(j)}$.

\proclaim{Proposition 3.4.23}  
Let the notation be as in $3.4.4$.
Then there are an open neighborhood $V$ of
$(0^{\Phi},0,0,\delta(\br),0)$ in $Y^{II}_{\text{top}}(p,\br,R, S)$ and an 
 open immersion $V\to D_{\SL(2)}^{II}(\Phi)$ of topological spaces 
 which sends
$(0^{\Phi}, 0,0,\delta(\br), 0)$ to $p$ and whose restriction to 
$V \cap Y_{0,\text{top}}^{II}(p,\br,R,S)$ coincides with the restriction of $\eta^{II}_{p,\br,R,S,\text{top}}$ $(3.4.22)$. 
\endproclaim

{\bf 3.4.24.} 
This 3.4.23 follows from 3.4.4, because we have a homeomorphism 
$$
Y^{II}(p,\br,R,S)\cong Y^{II}_{\text{top}}(p,\br,R,S),\quad (t,f,g,h,k,\delta,u)\leftrightarrow (t,a, k, \delta,u)
$$
$$
\text{with}\quad a=f+h,
$$
$$ 
f=\tsize\sum_m (1+t^{-2m})^{-1}a_m, \quad 
g=\tsize\sum_m (t^m+t^{-m})^{-1}a_m, \quad 
h=\tsize\sum_m (t^{2m}+1)^{-1}a_m,
$$ 
where, in $\sum_m$, $m$ ranges over all elements of $\bZ^{\Phi}$ such that either 
$m(j)\geq 0$ for any $j\in J(t)$ or $m(j)\leq 0$ for any $j\in J(t)$ (note that $(1+t^{-2m})^{-1}, (t^m+t^{-m})^{-1}, (t^{2m}+1)^{-1}\in \bR$ are naturally defined for such $m$). 

\medskip

{\bf 3.4.25.} {\it Remark.}
In the pure case, at the beginning of \cite{KU3} \S10, it is suggested that the local homeomorphism with 
$Y^{II}_{\text{top}}(p, \br, R, S)$ 
in 3.4.23 may be used to define a real analytic structure of $D_{\SL(2)}$. 
  If we do so, we regard $Y^{II}_{\text{top}}(p,\br,R,S)$ 
as an object of $\cB_\bR(\log)$ by using the embedding 
$Y^{II}_{\text{top}}(p,\br,R,S) \hookrightarrow\bR_{\geq 0}^{\Phi} 
\times R \times S \times \bar L\times \fg_{\bR,u}$
 in the same way as we did so for $Y^{II}(p,\br, R,S)$ by using the injection $Y^{II}(p,\br,R,S)\hookrightarrow E$ (3.4.2). 
 However the definition of the real analytic structure of $D_{\SL(2)}$ in this paper, which is given by the local homeomorphism with 
$Y^{II}(p, \br, R, S)$, is slightly different from this suggested one in \cite{KU3} \S10. The above map 
$Y^{II}(p, \br, R, S)\to Y^{II}_{\text{top}}(p, \br, R, S)$ is real analytic and is a homeomorphism but the inverse map 
need not be real analytic at $(0^{\Phi}, 0, 0, \delta(\br),0)$.

\medskip

{\bf 3.4.26.} 
  We define the topological space $Y^I_{\text{top}}(p,\br,R,S)$ as 
follows. 

In the case $W\notin \Psi$, 
let $Y^I_{\text{top}}(p, \br,R,S)$ be the subset of $Y^{II}_{\text{top}}(p,\br,R,S)\times \fg_{\bR,u}$ consisting of all elements $(t,a,k,\delta,u,v)$ ($(t,a,k,\delta,u)\in Y^{II}_{\text{top}}(p,\br,R,S)$, $v \in \fg_{\bR,u}$) such that $(t,\d,u,v)$ 
satisfies the conditions (5)--(7) in 3.4.5. 

  Similarly, in the case $W \in \Psi$, 
let $Y^I_{\text{top}}(p, \br,R,S)$ be the subset of $\bR_{\ge0} \times 
Y^{II}_{\text{top}}(p,\br,R,S)\times \fg_{\bR,u}$ consisting of all elements 
$(t_0, t,a,k,\delta,u,v)$ ($t_0 \in \bR_{\ge0}$, $(t,a,k,\delta,u)\in Y^{II}_{\text{top}}(p,\br,R,S)$, $v \in \fg_{\bR,u}$) such that $(t_0, t, \delta, u, v)$ 
satisfies the conditions (5${}'$)--(7${}'$) in 3.4.5.

\medskip
  We define a canonical map
$Y^I_{\text{top}}(p,\br, R,S) \to Y^{II}_{\text{top}}(p,\br, R,S)$. 
  If $W \not\in \Psi$, it is the canonical projection. 
  If otherwise, it is $(t_0,t',a,k,\d,u,v) \mapsto (t',a,k,t_0\circ \d,u)$. 
Let $Y^I_{0,\text{top}}(p,\br, R,S)$ be the open set 
of $Y^I_{\text{top}}(p,\br, R,S)$ defined by the inverse image of 
$Y^{II}_{0,\text{top}}(p,\br, R,S)$ by this canonical map. 
  Then, $Y^I_{0,\text{top}}(p,\br, R,S) \to 
Y^{II}_{0,\text{top}}(p,\br, R,S)$ is a homeomorphism. 
  Let 
$\eta_{p,\br,R,S,\text{top}}^{I}: Y_{0,\text{top}}^{I}(p,\br, R, S)\to D$ 
be the continuous map obtained from $\eta_{p,\br,R,S,\text{top}}^{II}$ and the 
last homeomorphism.

\medskip

\proclaim{Proposition 3.4.27}  Let the notation be as in $3.4.6$.
Assume $W\notin \Psi$ $($resp. $W\in \Psi$$)$. 
Then there is an open neighborhood $V$ of
v:=$(0^{\Psi},0,0,\delta(\br),0, 0)$ 
$($resp. $(0^{\Psi}, 0,0,\delta(\br)^{(1)}, 0, 0)$, 
where $\delta(\br)^{(1)}\in L^{(1)}$ $(3.4.5)$ such that $\delta(\br)=0\circ \delta(\br)^{(1)}$$)$ 
in $Y^{I}_{\text{top}}(p,\br,R, S)$ and an 
 open immersion $V\to D_{\SL(2)}^{I}(\Psi)$ of topological spaces 
 which sends
$v$ to $p$ and whose restriction to 
$V \cap Y_{0,\text{top}}^{I}(p,\br,R,S)$ coincides with the restriction of $\eta^{I}_{p,\br,R,S,\text{top}}$ $(3.4.26)$. 
\endproclaim

 This follows from Theorem 3.4.6, just as that 3.4.23 follows from 3.4.4 
in 3.4.24. 

\medskip

{\bf 3.4.28.}
We prove 3.2.12. 

We prove (i). It is sufficient to prove that the topology of $D_{\SL(2)}^I$ has the property (2). Let $p\in D_{\SL(2)}$ and let $\Psi$ be the set of weight filtrations associated to $p$. 
  In the following, we assume $W \not \in \Psi$.  
  The case where $W \in \Psi$ is similar. 
Assume first that 
$(p_\lam)_\lam$ ($p_\lam \in D$) converges to $p$. Then clearly (a) and (b) are satisfied. 
Take a distance to $\Psi$-boundary $\b$ 
such that $\b(\br)=1$ 
and let $\mu: D_{\SL(2)}^I(\Psi)\to D$ be the extension of $x\mapsto \tau_p\b(x)^{-1}x$ given in 3.2.6 (i). 
We show that (c.I) is satisfied for $t_\lam:=\b(p_\lam)$. 
We have $t_\lam=\b(p_\lam)\to \b(p)=0^{\Psi}$, and 
 $\tau_p(t_\lam)^{-1}p_\lam=\mu(p_\lam) \to \mu(p)=\br$. 
Next assume (a), (b), and (c.I) are satisfied. Take 
$\alpha=\tau_p$ and take $\b$ such that $\b(\br)=1$. We prove $p_\lam\to p$. 
It is sufficient to prove that $\nu_{\a,\b}(p_\lam) $
converges to $\nu_{\a,\b}(p)=(0^{\Psi},\, \br,\, \spl_W(\br),\, (\spl^{\BS}_{W'(\gr^W)}(\br(\gr^W)))_{W'\in \Psi})$ in $\bR_{\geq 0}^{\Psi}\times D\times \spl(W)\times \prod_{W'\in \Psi} \spl(W'(\gr^W))$.
The $\spl(W)$-component and the $\spl(W'(\gr^W))$-component of $\nu_{\a,\b}(p_\lam)$ converge to $\spl_W(\br)$ and to $\spl^{\BS}_{W'(\gr^W)}(\br(\gr^W))$ by (a) and (b), respectively. 
Let $a_\lam = t_\lam^{-1}\b(p_\lam)\in \bR_{>0}^{\Psi}$. 
By taking $\b$ of  $\tau_p(t_\lam)^{-1}p_\lam\to \br$, we have $a_\lam\to 1$.  Since $t_\lam\to 0^{\Psi}$,  $\b(p_\lam)=t_\lam a_\lam$ converges to $0^{\Psi}$. Finally, $\a\b(p_\lam)^{-1}p_\lam =\tau_p(a_\lam)^{-1}\tau_p(t_\lam)^{-1}p_\lam \to \br$.  

The proof of (ii) is similar to that of (i). 

3.2.12 is proved.
\medskip

\proclaim{Proposition 3.4.29}
The following conditions {\rm(1)--(3)} are equivalent. 

\medskip

$(1)$ The topology of $D_{\SL(2)}^I$ coincides with that of $D_{\SL(2)}^{II}$.

\medskip

$(2)$ $D_{\SL(2)}^I$ and $D_{\SL(2)}^{II}$ coincide in $\cB_\bR(\log)$.

\medskip

$(3)$ For any $p\in D_{\SL(2)}$, for any $w, w'\in \bZ$ such that $w>w'$, for any member $W'$ of the set of weight filtrations associated to $p$, and for any $a, b\in \bZ$ such that $\gr^{W'}_a(\gr^W_w)\neq 0$ and $\gr^{W'}_b(\gr^W_{w'})\neq 0$, we have $a\geq b$.
\endproclaim

\medskip

{\it Remarks.} 
(i) Assume that the equivalent conditions of Proposition 3.4.29 are satisfied. 
Then, for any $\Psi\in \cW$ and for 
$\bar \Psi=\{W'(\gr^W)\;|\;W'\in \Psi, W'\neq W\}\in \overline{\cW}$, $D_{\SL(2)}^I(\Psi)=D_{\SL(2)}^{II}(\bar \Psi)$ in $\cB_{\bR}(\log)$ if $W\in \Psi$, and $D_{\SL(2)}^I(\Psi)$ is an open subobject of $D_{\SL(2)}^{II}(\bar \Psi)$ in general.
\medskip

(ii) As is easily seen from 2.3.9, Examples I--IV in 1.1.1 satisfy the above condition (3), but Example V does not. See 3.6.2.  

\medskip

{\bf 3.4.30.}
{\it Proof of Proposition 3.4.29.}

We first prove that (1) implies (3). 
Assume that (3) is not satisfied. 
Then for some  $p\in D_{\SL(2)}$,  there exists $W'\in \cW(p)$ 
having the following property. 
There are $w, w', a, b\in \bZ$ such that $\gr^{W'}_a(\gr^W_w)$ and $\gr^{W'}_b(\gr^W_{w'})$ are not zero, and $w>w'$ and $a<b$. 
There is a non-zero element $u$ of $\fg_{\bR,u}$ such that the $W'$-component $\tau_{p,W'}$ satisfies $\Ad(\tau_{p,W'}(t))u= t^{b-a}u$ for all $t\in \bR^\times$. 
Take any real number $c$ such that $0<c<b-a$. We have $W'\neq W$. 
For $t\geq 0$, let $\epsilon(t)$ be the element of $\bR_{\geq 0}^{\Psi}$ whose $W'$-component is $t$ and whose all the 
other components are $1$. Let $\Phi$ be the image of $\Psi$ in $\overline{\cW}$ (3.2.2). 
Take a point $\br\in D$ on the torus orbit associated to $p$,  consider an element $p'$ of $Y^I(p,\br, R, S)$ of the form $p'=(\epsilon(0), 0,0,0,0,\delta, 0,0)\in Y^I(p,\br, R, S)$, let $\bar \epsilon(t)$ be the image of $\epsilon(t)$ in $\bR_{\geq 0}^{\Phi}$, and let $p''=(\bar\epsilon(0), 0,0,0,0,\delta, 0)\in Y^{II}(p,\br, R, S)$ be the image of $p'$. 
When $t \in \bR_{>0}$ tends to $0$, 
$(\bar\epsilon(t), 0,0,0,0,\delta, t^cu)\in Y^{II}(p,\br, R, S)$ converges to $p''$. But this element of $Y^{II}(p,\br, R, S)$ is the image of
$(\epsilon(t), 0,0,0,0,\delta, 
t^cu, t^{c+a-b}u)\in Y^I(p,\br, R,S)$ which does not converge to 
$p'$ when $t\to 0$ 
because $c+a-b<0$. 
By 3.4.4 and 3.4.6, 
this proves that the topology of $D_{\SL(2)}^I$ and that of $D_{\SL(2)}^{II}$ are different. 

It is clear that (2) implies (1). 

It remains to prove that (3) implies (2). 
Assume that (3) is satisfied. As in (i) of Remarks after 3.4.29, 
$D_{\SL(2)}^I(\Psi)$ is an open set of $D_{\SL(2)}^{II}(\Phi)$. 

 By 3.4.4 and 3.4.6, it is sufficient to prove

 \medskip
 
{\bf Claim.}
{\it For a splitting $\alpha$ of $\Psi$, the map
 $\bR_{>0}^{\Psi}\times \fg_{\bR,u}\to \fg_{\bR,u},\;(t,u)\mapsto \Ad(\alpha(t))^{-1}(u)$
 extends to a real analytic map $\bR_{\geq 0}^{\Psi}\times \fg_{\bR,u}\to \fg_{\bR,u}$.} 
 
 \medskip 
 
 By (3), for the adjoint action of $\bG_{m,\bR}^{\Psi}$ by $\alpha$, $\fg_{\bR,u}$ is the sum of the eigen spaces $(\fg_{\bR,u})_m$ for all $m\in \bZ^{\Psi}$ such that $m\leq 0$. 
This proves the claim.
\qed

\medskip

\vskip 20pt

\head
\S 3.5. Global properties of $D_{\SL(2)}$
\endhead

\medskip

In this \S3.5, we prove that the projection $D_{\SL(2)}^{II}\to \spl(W)\times D_{\SL(2)}(\gr^W)$ is proper (Theorem 3.5.16). 
We prove also results on the actions of a subgroup $\G$ of $G_\bZ$ on $D_{\SL(2)}^I$ and on $D_{\SL(2)}^{II}$ (Theorem 3.5.17).

Concerning the properness of $D_{\SL(2)}^{II}$ over $\spl(W) \times D_{\SL(2)}(\gr^W)$, we prove a more precise result. 
We define a log modification (3.1.12) $D_{\SL(2)}(\gr^W)^{\sim}$ of $D_{\SL(2)}(\gr^W)$, which is an object of $\cB_\bR(\log)$ and is proper over $D_{\SL(2)}(\gr^W)$, such that the canonical projection $D_{\SL(2)}\to D_{\SL(2)}(\gr^W)$ factors as $D_{\SL(2)}\to D_{\SL(2)}(\gr^W)^{\sim}\to D_{\SL(2)}(\gr^W)$. 
We prove that as an object of $\cB_{\bR}(\log)$, $D_{\SL(2)}^{II}$ is an $\bar L$-bundle over $\spl(W) \times D_{\SL(2)}(\gr^W)^{\sim}$ (Theorem 3.5.15). 
Here $L=\cL(F)$ for any fixed $F\in D(\gr^W)$, and $\bar L$ is the compactified vector space associated to $L$ (3.2.6). 
This is an \lq\lq $\SL$(2)-analogue'' of the fact (\cite{KNU2}, Theorem 8.5) that $D_{\BS}$ is an $\bar L$-bundle over $\spl(W) \times D_{\BS}(\gr^W)$. 
The properness of $D_{\SL(2)}^{II}$ over $\spl(W)\times D_{\SL(2)}(\gr^W)$ follows from this.

\medskip

{\bf 3.5.1.} 
We define the set  $D_{\SL(2)}(\gr^W)^{\sim}$. 

By an {\it $\SL(2)$-orbit on $\gr^W$}, we mean a family $(\rho_w, \varphi_w)_{w\in \bZ}$, where, for some $n\geq 0$, $(\rho_w, \varphi_w)$ is an SL(2)-orbit for $\gr^W_w$ in $n$ variables for any $w\in \bZ$ satisfying the following condition (1). 

\medskip

(1) For $1\leq j\leq n$, there is $w\in \bZ$ such that the $j$-th component of $\rho_w$ is non-trivial.  

\medskip

This $n$ is called the {\it rank of} $(\rho_w, \vf_w)_w$. 

We say two $\SL(2)$-orbits $(\rho_w, \varphi_w)_w$ and $(\rho'_w,\varphi'_w)_w$ on $\gr^W$ are {\it equivalent} if  their ranks coincide, say $n$,  and furthermore there is $t = (t_1,\dots, t_n)\in \bR_{>0}^n$ such that
$$
\rho'_w=\Int(\tilde\rho_w(t))\rho_w, \quad 
\varphi'_w=\tilde\rho_w(t)\varphi_w
$$
for any $w\in \bZ$, where $\tilde\rho_w(t)$ is as in 2.5.1. 

Let $D_{\SL(2)}(\gr^W)^{\sim}$ be the set of all equivalence classes of $\SL(2)$-orbits on $\gr^W$.

  Note that $\SL(2)$-orbits on $\gr^W$ just defined are in fact 
what should be called non-degenerate $\SL(2)$-orbits on $\gr^W$. 
  We omitted this adjective in the above definition, since we use only 
\lq\lq non-degenerate" ones for the study of $D_{\SL(2)}(\gr^W)^{\sim}$.

\medskip

{\bf 3.5.2.}
The canonical map $D_{\SL(2)}\to D_{\SL(2)}(\gr^W)=\tp_{w\in \bZ} D_{\SL(2)}(\gr^W_w)$ factors as 
$$
D_{\SL(2)}\to D_{\SL(2)}(\gr^W)^{\sim}\to D_{\SL(2)}(\gr^W),
$$
where the second arrow is
$$
D_{\SL(2)}(\gr^W)^{\sim} \to D_{\SL(2)}(\gr^W),\;\;(\text{class of}\;(\rho_w,\vf_w)_w) \mapsto (\text{class of}\;(\rho_w,\vf_w))_w,
$$
and the first arrow is defined as follows. 
Let $p\in D_{\SL(2)}$ be the class of an $\SL(2)$-orbit $((\rho_w, \varphi_w)_w, \br)$ in $n$ variables of rank $n$ and let $\Psi$ be the associated set of weight filtrations. 
Then the image $\tilde p$ of $p$ in $D_{\SL(2)}(\gr^W)^{\sim}$ is the class of the following $\SL(2)$-orbit $(\rho_w', \vf_w')_w$ on $\gr^W$. 
If $W\notin \Psi$, $(\rho_w', \vf_w')_w=(\rho_w, \vf_w)_w$ and hence $\tilde p$ is of rank $n$. 
If $W\in \Psi$, then $(\rho_w', \vf'_w)_w$ is an $\SL(2)$-orbit on $\gr^W$ of rank  $n-1$ defined by
$$
\rho_w'(g_1,\dots, g_{n-1})=\rho_w(1, g_1, \dots, g_{n-1}), \quad \vf'_w(z_1,\dots, z_{n-1})=\vf_w(i, z_1, \dots, z_{n-1}),
$$ 
for $w\in \bZ$. 

The map $D_{\SL(2)}\to D_{\SL(2)}(\gr^W)^{\sim}$ is surjective. 

The map $D_{\SL(2)}\to \overline{\cW}, \;p\mapsto \overline{\cW}(p)$, (3.2.2) factors through $D_{\SL(2)}\to D_{\SL(2)}(\gr^W)^{\sim}$. 
For $q\in D_{\SL(2)}(\gr^W)^{\sim}$, we denote by $\overline{\cW}(q)\in \overline{\cW}$ the element $\overline{\cW}(p)$ for $p$ an element of $D_{\SL(2)}$ with image $q$ in $D_{\SL(2)}(\gr^W)^{\sim}$, which is independent of the choice of $p$.

The map $D_{\SL(2)}(\gr^W)^{\sim} \to D_{\SL(2)}(\gr^W)$ is also surjective.
This is shown as follows. 
For each $w\in \bZ$, let $(\rho_w, \vf_w)$ be an $\SL(2)$-orbit on $\gr^W_w$ in $n(w)$ variables of rank $n(w)$. 
Let $n=\max\{n(w)\;|\;w\in \bZ\}$, and let $(\rho'_w, \vf'_w)$ be the $\SL(2)$-orbit on $\gr^W_w$ in $n$ variables defined by $\rho_w'(g_1, \dots, g_n)=\rho_w(g_1, \dots, g_{n(w)})$ and $\vf'_w(z_1, \dots, z_n)=\vf_w(z_1, \dots, z_{n(w)})$. 
Then $(\text{class of}\; (\rho_w, \vf_w))_w\in \prod_w D_{\SL(2)}(\gr^W_w)$ is the image of the element $(\text{class of}\; (\rho'_w, \vf_w')_w)$ in $D_{\SL(2)}(\gr^W)^{\sim}$ (cf.\ 3.5.1). 

The map $D_{\SL(2)}(\gr^W)^{\sim} \to D_{\SL(2)}(\gr^W)$ need not be injective (see Corollary 3.5.12, Example V in 3.5.13). This is explained by two reasons. 
The first reason is as follows. For $\SL(2)$-orbits $(\rho_w, \vf_w)_w$ and $(\rho'_w, \vf'_w)_w$ on $\gr^W$, their images in $D_{\SL(2)}(\gr^W)$ coincide if and only if $(\rho_w, \vf_w)$ and $(\rho'_w, \vf'_w)$ are equivalent for all $w$, and the last equivalences are given by elements of $\bR_{>0}^{n(w)}$ which can depend on $w\in \bZ$ (here $n(w) =\rank (\rho_w,\vf_w)=\rank (\rho'_w, \vf'_w)$) not like the equivalence between $(\rho_w, \vf_w)_w$ and $(\rho'_w, \vf'_w)_w$ defined as in 3.5.1. 
The second reason is as follows. 
For $p\in D_{\SL(2)}$, the image of $p$ in $D_{\SL(2)}(\gr^W)^{\sim}$ still remembers $\overline{\cW}(p)\in \overline{\cW}$, but the image of $p$ in $D_{\SL(2)}(\gr^W)$ remembers only the image of this element of $\overline{\cW}$ in $\prod_w \cW(\gr^W_w)$ (3.3.1).
As in 3.3.2, the map $\overline{\cW}\to \prod_w \cW(\gr^W_w)$ is described as $(\Phi, (\Phi(w))_w)\mapsto (\Phi(w))_w$, and not necessarily injective.

\medskip

{\bf 3.5.3.} 
For $Q=(Q(w))_{w\in \bZ} \in  \tsize\prod_{w\in \bZ} \cW(\gr^W_w)$ (3.3.1), let $D_{\SL(2)}(\gr^W)(Q)$ be the open set of $D_{\SL(2)}(\gr^W)$ defined by 
$$
D_{\SL(2)}(\gr^W)(Q)=\tsize\prod_{w\in \bZ} \; D_{\SL(2)}(\gr^W_w)(Q(w))\sub D_{\SL(2)}(\gr^W),
$$
as in 3.4.21.

Define
$$
D_{\SL(2)}(\gr^W)^{\sim}(Q)\subset D_{\SL(2)}(\gr^W)^{\sim}
$$ 
as the inverse image of $D_{\SL(2)}(\gr^W)(Q)$ in $D_{\SL(2)}(\gr^W)^{\sim}$. For $p\in D_{\SL(2)}(\gr^W)^{\sim}$, $p$ belongs to $D_{\SL(2)}(\gr^W)^{\sim}(Q)$ if and only if $\Phi:=\overline{\cW}(p)$ satisfies $\Phi(w)\subset Q(w)$ for all $w\in \bZ$. 
\medskip

{\bf 3.5.4.} 
Let $Q=(Q(w))_w\in \tsize\prod_{w\in \bZ} \cW(\gr^W_w)$, let $S= D_{\SL(2)}(\gr^W)(Q)$, and let $\cS = \bigoplus_{w\in \bZ} \bN^{Q(w)}$. 
Then we have a canonical surjective homomorphism $\cS\to   M_S/\cO_S^\times$ characterized as follows. For any distance to $Q(w)$-boundary $\b_w=(\b_{w,j})_{j\in Q(w)}:D(\gr^W_w)\to \bR_{>0}^{Q(w)}$ given for each $w\in \bZ$, this homomorphism sends $m = ((m(w, j))_{j\in Q(w)})_w \in \cS$ ($m(w,j)\in \bN$) to  $(\prod_{w\in \bZ,\, j\in Q(w)} \;\b_{w,j}^{m(w,j)}) \modu \cO_S^\times$. 
This homomorphism lifts locally on $S$ to a chart $\cS\to M_{S,>0}$. 

In 3.5.5--3.5.7, we will define and study a finite rational subdivision $\Sig_Q$ of the cone $\Hom(\cS, \bR_{\geq 0}^{\add})= \prod_{w\in \bZ} \bR_{\geq 0}^{Q(w)}$, and in 3.5.9 we will identify $D_{\SL(2)}(\gr^W)^{\sim}(Q)$ with the associated log modification $S(\Sig_Q)$ (3.1.12) of $S$. 
We will see in 3.5.10 that there is a unique structure on $D_{\SL(2)}(\gr^W)^{\sim}$ as an object of $\cB_\bR(\log)$ for which each $D_{\SL(2)}(\gr^W)^{\sim}(Q)$ ($Q\in \tsize\prod_{w\in \bZ} \cW(\gr^W_w)$) is open in $D_{\SL(2)}(\gr^W)^{\sim}$ and the induced structure on it coincides with the structure as the log modification. 
\medskip

{\bf 3.5.5.} 
For  $Q=(Q(w))_{w\in \bZ}\in \tsize\prod_{w\in \bZ} \cW(\gr^W_w)$, we define a finite rational subdivision $\Sig_Q$ of the cone $\prod_{w\in \bZ} \bR_{\geq 0}^{Q(w)}$ as follows.  

First we recall that, for a finite set $\Lam$, the barycentric subdivision $\operatorname{Sd}(\Lambda)$ of the cone $\bR_{\geq 0}^{\Lam}$ is defined as follows (cf.\ \cite{I} 2.8). 
Let $J(\Lam)$ be the set of all pairs $(n, g)$, where $n$ is a non-negative integer and $g$ is a function $\Lam\to \{j\in \bZ\;|\;0\leq j\leq n\}$ such that the image of $g$ contains $\{j\in \bZ\;|\;1\leq j\leq n\}$. 
For $(n , g)\in J(\Lam)$, define the subcone 
$C(n, g)$ of $\bR_{\geq 0}^{\Lam}$ by
$$
C(n, g) =\{(a_\lam)_{\lam\in \Lam}\;|\; a_\lam\leq a_\mu 
\;\text{if $g(\lam)\leq g(\mu)$}, \;a_\lam=0\;\text{if $g(\lam)=0$}\}.
$$
Then the set of cones $\operatorname{Sd}(\Lam):=\{C(n, g)\;|\;(n, g) \in J(\Lam)\}$ is a finite rational subdivision of $\bR_{\geq 0}^{\Lam}$, and is called the {\it barycentric subdivision of $\bR_{\geq 0}^{\Lam}$.} 
The map 
$$
J(\Lam) \to \operatorname{Sd}(\Lam),\quad
(n,g)\mapsto C(n, g),
$$
is bijective. 
For $(n, g)\in J(\Lam)$, the dimension of $C(n, g)$ is equal to $n$. 
\medskip

Let $Q=(Q(w))_w\in \tsize\prod_{w\in \bZ} \cW(\gr^W_w)$. 
For each $w\in \bZ$, we regard $Q(w)$ as a totally ordered set by 2.1.13. 

Let $\Lam=\tsCu_{w\in  \bZ} Q(w)$.
Define a subcone $C$ of $\bR_{\geq 0}^{\Lam} = \tsize\prod_{w\in \bZ} \bR_{\geq 0}^{Q(w)}$ by
$$
C=\{((a_{w,j})_{j \in Q(w)})_w \in \tsize\prod_{w\in \bZ} \bR_{\geq 0}^{Q(w)}\;|\; a_{w,j}\leq a_{w,j'}\;\text{if $w\in \bZ$, 
$j, j'\in Q(w)$ and $j\geq j'$}\}.
$$

Let 
$$
\operatorname{Sd}'(\Lam)=\{\sig\in \operatorname{Sd}(\Lambda)\;|\;\sig\sub C\}\sub \operatorname{Sd}(\Lam),
$$
$$
J'(\Lam) =\{(n, g)\in J(\Lam)\;|\;g(w, j)\leq g(w,j')\;\text{if}\; w\in \bZ, j, j' \in Q(w)\; \text{and}\; j\geq j' \} \sub J(\Lam).
$$
  Here and hereafter, $g(w,-)$ denotes the restriction of the map $g$ 
on $Q(w) \sub \Lambda$ for any $w$. 
Then  
$$\operatorname{Sd}'(\Lam)=\{C(n,g)\;|\; (n, g)\in J'(\Lam)\},$$
and
$\operatorname{Sd}'(\Lam)$ is a subdivision of $C$.

We have an isomorphism of cones 
$$
\bR_{\ge0}^\Lambda = 
\tsize\prod_{w\in \bZ} \bR_{\geq 0}^{Q(w)} @>\sim>> C,\quad b\mapsto c, \tag1
$$
where $c_{w,j}:= \tsize\sum_{k\in Q(w),\, k\geq j} b_{w, k}$ for $w \in \bZ$ and $j\in Q(w)$.

Let $\Sig_Q$ be the subdivision of the cone $\bR_{\ge0}^\Lambda = \tsize\prod_{w\in \bZ} \bR_{\geq 0}^{Q(w)}$
corresponding to the subdivision $\operatorname{Sd}'(\Lam)$ of the cone $C$ via the above isomorphism (1). 

\medskip

{\bf 3.5.6.} 
Let 
$\overline{\cW}\to \prod_{w\in \bZ} \cW(\gr^W_w)$ be the map defined in 3.3.1.

For $Q = (Q(w))_w \in \tsize\prod_{w\in \bZ} \cW(\gr^W_w)$, let $\overline{\cW}(Q)\sub \overline{\cW}$ be the set of all $\Phi\in \overline{\cW}$ such that
 $\Phi(w)\sub Q(w)$ for any $w\in \bZ$.

\proclaim{Proposition 3.5.7} 
Let $Q = (Q(w))_w \in \tsize\prod_{w\in \bZ} \cW(\gr^W_w)$. 
Then we have a bijection $$\overline{\cW}(Q) \to \Sig_Q, \quad \Phi\mapsto \sig_\Phi,$$ where $\sig_{\Phi}$ is the set of all elements $((b_{w,j})_{j\in Q(w)})_{w\in \bZ}$ of $\prod_{w\in \bZ} \bR_{\geq 0}^{Q(w)}$ satisfying the following condition $(1)$.

\medskip

$(1)$ Let $w, w'\in \bZ$, $j\in Q(w)$, $j'\in Q(w')$. Assume that, for any $M\in \Phi$ such that $j\leq M(\gr^W_w)$, we have $j'\leq M(\gr^W_{w'})$ $(2.1.13)$. 
Then 
$$
\tsize\sum_{k \in Q(w),\, k \geq j} b_{w,k} \leq \tsize\sum_{k\in Q(w'),\, k\geq j'} b_{w', k}.
$$ 
\endproclaim 

{\it Remark.}
  The condition (1) is equivalent to the following conditions (1a) and (1b): 
(1a) $b_{w,j}=0$ unless there is an $M \in \Phi$ such that $j = M(\gr^W_w)$.
(1b) $b_{w,j}=b_{w',j'}$ if there is an
$M \in \Phi$ such that $j = M(\gr^W_w)$ and $j' = M(\gr^W_{w'})$.

\medskip

{\it Proof.}
By construction in 3.5.5, we have bijections $J'(\Lambda) \simeq \operatorname{Sd}'(\Lambda) \simeq \Sigma_Q$.
Under these bijections, the above $\sig_{\Phi}$ is equal to the element of $\Sig_Q$ corresponding to the element $C(n, g)\in \operatorname{Sd}'(\Lam)$, where $(n, g)$ is the element of $J'(\Lam)$ ($\Lam=\tsCu_{w\in  \bZ} Q(w)$) defined as follows.   
Let $n$ be the cardinality of $\Phi$, that is, $n =\dim \sig_{\Phi}$. 
Let $M^{(1)}=(M^{(1)}(w))_w, \dots, M^{(n)}=(M^{(n)}(w))_w$ be the all members of $\Phi$ such that $M^{(1)}(w)\leq \dots \leq M^{(n)}(w)$ for any $w\in \bZ$ with respect to the ordering in 2.1.13.
Then, for $w\in \bZ$ and $j\in Q(w)$, define
$$
g(w, j)=\sharp\{k\;|\; 1\leq k\leq n,\;
 M^{(k)}(w)\geq j\}.
$$ 
By 3.3.2, this map $\overline{\cW}(Q) \to J'(\Lam), \;\Phi\mapsto (n, g)$, is  bijective.   
\qed

\medskip

\proclaim{Lemma 3.5.8} 
Let $Q\in \tsize\prod_{w\in \bZ} \cW(\gr^W_w)$, let $p\in S=
D_{\SL(2)}(\gr^W)(Q)$, let $q$ be a point of $D_{\SL(2)}(\gr^W)^{\sim}(Q)$ lying over $p$, let $\Phi=\overline{\cW}(q)$ {\rm(3.2.2)}, and let $\sig_q=\sig_\Phi\in \Sig_Q$ {\rm(3.5.7)}.   
Let $P'(\sig_q)\subset M^{\gp}_{S,>0,p}$ be as in $3.1.13$. 
That is, for $S$ and $\cS$  in {\rm 3.5.4}, let $\cS(\sig_q)$ be the subset of $\cS^{\gp}$ consisting of all elements $m$ of $\cS^{\gp}$ such that 
the homomorphism $\cS^{\gp}\to \bR$ defined by any element of $\sig_q$ sends $m$ into $\bR_{\geq 0}$, let $P(\sig_q)$ be the image of $\cS(\sig_q)$ in 
$(M^{\gp}_S/\cO_S^\times)_p$, and let $P'(\sig_q)$ be the inverse image of $P(\sig_q)$ in $M^{\gp}_{S,>0,p}$.
Then we have  
$$
P'(\sig_q)=\{f\in M^{\gp}_{S,>0,p}\;|\;f(\tau_q(t)\br_q)\;\text{converges in $\bR_{\geq 0}$}\}, \tag 1
$$
$$
P'(\sig_q)^\times=\{f\in M^{\gp}_{S,>0,p}\;|\;f(\tau_q(t)\br_q)\;\text{converges to an element of $\bR_{>0}$}\}. \tag 2
$$ 
Here $\br_q$ is a point on the torus orbit associated to $q$, 
$\tau_q : \bR_{\ge0}^\Phi \to \Aut(\gr^W)$ is $\bar\tau_{q'}$ in $3.2.3$ for a point $q' \in D_{\SL(2)}$ lying over $q$, and $t$ tends to $0^{\Phi}$. 
\endproclaim

{\it Proof.} 
In the notation of 3.5.4, $P'(\sig_q) \sub M^{\gp}_{S,>0,p}$ is written as
$$
P'(\sig_q) = \tsize\bigcup_{m\in\cS(\sigma_q)} \cO_{S,>0,p}^\times\tsize\prod_{w\in\bZ,\,j\in Q(w)}\beta_{w,j}^{m(w,j)},
$$
where $m = ((m(w,j))_{j\in Q(w)})_{w\in\bZ}$.
This coincides with the right hand side of (1) by 3.2.6 (ii).
Since $P'(\sig_q)^\times = P'(\sig_q) \cap P'(\sig_q)^{-1}$, (2) follows.
\qed

\proclaim{Theorem 3.5.9} 
Let $Q\in \tsize\prod_{w\in \bZ} \cW(\gr^W_w)$.

\medskip
{\rm (i)} Let $D_{\SL(2)}(\gr^W)(\Sig_Q)$ be the log modification $(3.1.12)$ 
of $D_{\SL(2)}(\gr^W)(Q)$ corresponding to the subdivision $\Sig_Q$ of the cone $\tsize\prod_{w\in \bZ} \bR_{\ge0}^{Q(w)}$ in $3.5.5$. 
Then we have a  bijection 
$$
D_{\SL(2)}(\gr^W)^{\sim}(Q)\to D_{\SL(2)}(\gr^W)(\Sig_Q)
$$
which sends a point $q$ of $D_{\SL(2)}(\gr^W)^{\sim}(Q)$ lying over $p\in  D_{\SL(2)}(\gr^W)(Q)$ to the point $(p, \sig_q, h_q)$ $(3.1.13)$ of $D_{\SL(2)}(\gr^W)(\Sig_Q)$ lying over $p$, where $\sig_q$ is as in $3.5.8$ and $h_q$ is the homomorphism defined by 
$$
h_q: P'(\sig_q)^\times\to \bR_{>0}, \quad 
f\mapsto \lim_{t\to 0^{\Phi}}\;  f(\tau_q(t)\br_q),
$$
where $\br_q$, $\tau_q$ and $\Phi=\overline{\cW}(q)$ are as in $3.5.8$.

\medskip

{\rm (ii)} Let $\Phi\in \overline{\cW}(Q)$, and let $D_{\SL(2)}(\gr^W)^{\sim}(\Phi)\sub D_{\SL(2)}(\gr^W)^{\sim}$ be the image of $D_{\SL(2)}^{II}(\Phi)$. 
Then $D_{\SL(2)}^{II}(\Phi)$ coincides with the inverse image of $D_{\SL(2)}(\gr^W)^{\sim}(\Phi)$ in $D_{\SL(2)}$. 
Furthermore, let $\sigma_\Phi \in \Sigma_Q$ be as in $3.5.7$, then 
$D_{\SL(2)}(\gr^W)^{\sim}(\Phi)$ coincides with the part of $D_{\SL(2)}(\gr^W)^{\sim}(Q)$ which corresponds to the part $D_{\SL(2)}(\gr^W)(\sig_{\Phi})$ of $D_{\SL(2)}(\gr^W)(\Sig_Q)$ under the bijection in {\rm (i)}.

\endproclaim

{\it Proof.} 
Let $p\in D_{\SL(2)}(\gr^W)$, 
let $A$ be the fiber of $D_{\SL(2)}(\gr^W)^{\sim}\to D_{\SL(2)}(\gr^W)$ on $p$, and let $B$ be the set of all pairs $(\Phi, Z)$, 
where $\Phi$ is an element of $\overline{\cW}$ whose image in $\prod_w \cW(\gr^W_w)$ is $(\cW(p(\gr^W_w)))_w$ 
and $Z$ is an $\bR_{>0}^{\Phi}$-orbit in $D(\gr^W)$ contained in $\prod_w Z_w$, 
 where $Z_w$ is the torus orbit associated to $p(\gr^W_w)$. Then we have a bijection from $A$ to $B$ given by $q\mapsto (\Phi,  Z)$, 
where $\Phi=\overline{\cW}(q)$ and $Z$ is the torus orbit associated to $q$. 

  Assume that $Q(w)=\cW(p(\gr^W_w))$ for all $w$. 
  Then, 
once $\Phi \in \overline \cW (Q)$ 
is fixed, the set $B_{\Phi}$ of all $Z$ such that $(\Phi,Z)\in B$ is a $(\tp_{w\in \bZ} \bR_{>0}^{Q(w)})/\bR_{>0}^{\Phi}$-torsor. On the other hand, let $\sig$ be the cone corresponding to $\Phi$, and let $C_{\Phi}$ be the set of all homomorphisms $P'(\sig)^\times\to \bR_{>0}$ which extend the evaluation $\cO_{>0,p}^\times\to \bR_{>0}$ at $p$. Then $C_{\Phi}$ is also a $(\tp_{w\in \bZ} \bR_{>0}^{Q(w)})/\bR_{>0}^{\Phi}$-torsor with respect to the following action. By the canonical isomorphism $M^{\gp}_{>0,p}/\cO^\times_{>0,p}\simeq \tp_{w\in \bZ} \bZ^{Q(w)}$, we have an isomorphism
$$\Hom(M^{\gp}_{>0,p}/\cO^\times_{>0,p},\bR_{>0})\simeq \tp_{w\in \bZ} \bR_{>0}^{Q(w)}$$
which induces an isomorphism between quotient groups
$$\Hom(P'(\sig)^\times/\cO^\times_{>0,p}, 
\bR_{>0})\simeq (\tp_{w\in \bZ} \bR_{>0}^{Q(w)})/\bR_{>0}^{\Phi}.$$
Since $C_{\Phi}$ is a 
$\Hom(P'(\sig)^\times/\cO^\times_{>0,p}, \bR_{>0})$-torsor in the evident way, it is a 
$(\tp_{w\in \bZ} \bR_{>0}^{Q(w)})/\bR_{>0}^{\Phi}$-torsor. Let $A_{\Phi}$ be the subset of $A$ consisting of all 
$q\in A$ such that $\cW(q)=\Phi$. Then the bijection $A\to B$ induces a bijection $A_{\Phi}\to B_{\Phi}$. The map 
$A_{\Phi}\to C_{\Phi}$ which sends 
$q\in A_{\Phi}$ to the homomorphism
 $P'(\sig)^\times\to 
 \bR_{>0},\;f\mapsto \lim_{t\to 0^{\Phi}} f(\tau_q(t)\br_q)$ (3.5.8) induces a map 
 $B_{\Phi}\to C_{\Phi}$ which is compatible with the action of $(\tp_{w\in \bZ} \bR_{>0}^{Q(w)})/\bR_{>0}^{\Phi}$. Since $B_{\Phi}$ and $C_{\Phi}$ are 
$(\tp_{w\in \bZ} \bR_{>0}^{Q(w)})/\bR_{>0}^{\Phi}$-torsors, this map $B_{\Phi}\to C_{\Phi}$ is bijective.  Hence the map $A_{\Phi}\to C_{\Phi}$ is bijective.

3.5.9 follows from this. 
\qed
\medskip

{\bf 3.5.10.} 
We regard $D_{\SL(2)}(\gr^W)^{\sim}$ as an object of $\cB_\bR(\log)$ as follows. 
For $Q\in \tsize\prod_{w\in \bZ} \cW(\gr^W_w)$, $D_{\SL(2)}(\gr^W)^{\sim}(Q)$ is regarded as an object of $\cB_\bR(\log)$ via the bijection in Theorem 3.5.9. 
If $Q'\in \tsize\prod_{w\in \bZ} \cW(\gr^W_w)$ and $Q'(w)\sub Q(w)$ for all $w\in \bZ$, $D_{\SL(2)}(\gr^W)^{\sim}(Q')$ is open in $D_{\SL(2)}(\gr^W)^{\sim}(Q)$ and the structure of $D_{\SL(2)}(\gr^W)^{\sim}(Q')$ as an object of $\cB_\bR(\log)$ coincides with the one induced from that of $D_{\SL(2)}(\gr^W)^{\sim}(Q)$, as is easily seen. 
Hence there is a unique structure on $D_{\SL(2)}(\gr^W)^{\sim}$ as an object of $\cB_\bR(\log)$ for which $D_{\SL(2)}(\gr^W)^{\sim}(Q)$ are open and which induces on each $D_{\SL(2)}(\gr^W)^{\sim}(Q)$ the above structure as an object of $\cB_\bR(\log)$.

\medskip

\proclaim{Proposition 3.5.11} 
Let $p\in D_{\SL(2)}(\gr^W)$. 
Then the following two conditions are equivalent.
\medskip

{\rm(1)} The fiber of the surjection $D_{\SL(2)}(\gr^W)^{\sim}\to D_{\SL(2)}(\gr^W)$ over $p$ consists of one element. 
\medskip

{\rm(2)} There are at most one $w \in \bZ$ such that the element $p(w)$ of $D_{\SL(2)}(\gr^W_w)$ does not belong to $D(\gr^W_w)$. 
\endproclaim

{\it Proof.} 
This is seen easily by the proof of 3.5.9.
\qed
\medskip

From this it follows:

\proclaim{Corollary 3.5.12} 
The following three conditions are equivalent.
\medskip

{\rm(1)} $D_{\SL(2)}(\gr^W)^{\sim}\to D_{\SL(2)}(\gr^W)$ is bijective.
\medskip

{\rm(2)} $D_{\SL(2)}(\gr^W)^{\sim}\to D_{\SL(2)}(\gr^W)$ is an isomorphism of local ringed spaces over $\bR$. 
\medskip

{\rm(3)} There are at most one $w\in \bZ$ such that $D_{\SL(2)}(\gr^W_w)\neq D(\gr^W_w)$. 
\endproclaim
\medskip

{\bf 3.5.13.}
Consider $D_{\SL(2)}(\gr^W)^{\sim}$ for five examples I--V in 1.1.1.

For Example I--Example IV, we have $D_{\SL(2)}(\gr^W)^{\sim}=D_{\SL(2)}(\gr^W)$ by 3.5.12.
\medskip

{\bf Example V.}
Let $M$ be the increasing filtration on $\gr^W_0$ defined by
$$
M_{-3}=0\sub M_{-2} = M_{-1} = \bR e_1'\sub M_0 = M_1 = M_{-1} + \bR e_2' \sub M_2=\gr^W_0.
$$
Let $M'$ be the increasing filtration on $\gr^W_1$ defined by
$$
M'_{-1}=0 \sub M'_0 = M'_1 = \bR e_4' \sub M'_2=\gr^W_1.
$$
Let $Q = \{Q(w)\}_{w\in\bZ}$ be the following:
$Q(0):=\{M\}$, $Q(1):=\{M'\}$, and $Q(w)$ is the empty set for $w\in \bZ\smallsetminus \{0, 1\}$.
Let $\Lambda:= \{M, M'\}$.

Then the subdivision $\Sigma_Q$ of $\bR_{\ge0}^\Lambda = \prod_{w\in\bZ}\bR_{\ge0}^{Q(w)}$ in 3.5.5 is just the barycentric subdivision of $\bR_{\geq 0}^2$.
In the notation in 3.5.5, $0 \le n \le 2$ and $g$ is a function $\Lambda \to\{0, \dots, n\}$, and hence the fan $\Sigma_Q$ consists of the vertex $\{(0, 0)\}$ and the following cones according to the cases $m = 1, 2, 3, 4, 5$ in 2.3.9.
\medskip

0. $n = 0$, $g(M) = g(M') = 0$, and $C(0, g) = \{(0, 0)\}$.
\medskip

1. $n = 1$, $g(M) = 1$, $g(M') = 0$, and
$C(1, g) = \bR_{\ge0} \times \{0\}$.
\medskip

2. $n = 1$, $g(M) = 0$, $g(M') = 1$, and
$C(1, g) = \{0\} \times \bR_{\ge0}$.
\medskip

3. $n = 1$, $g(M) = g(M') = 1$, and
$C(1, g) = \{(a_\lambda)_\lambda \in \bR_{\geq 0}^2 \;|\; a_M = a_{M'}\}$.
\medskip

4. $n = 2$, $g(M) = 2$, $g(M') = 1$, and
$C(2, g) = \{(a_\lambda)_\lambda \in \bR_{\geq 0}^2 \;|\; a_M \ge a_{M'}\}$.
\medskip

5. $n = 2$, $g(M) = 1$, $g(M') = 2$, and
$C(2, g) = \{(a_\lambda)_\lambda \in \bR_{\geq 0}^2 \;|\; a_M \le a_{M'}\}$.
\medskip

Let $B$ be the corresponding \lq\lq blowing-up of $\bR^2_{\geq 0}$ at $(0, 0)$", i.e., the closure $B$ of $\bR^2_{>0}$ in the corresponding blowing-up of $\bC^2$ at $(0, 0)$.

Let $S = D_{\SL(2)}(\gr^W)(Q)$.
Then the inverse image $D_{\SL(2)}(\gr^W)^{\sim}(Q)$ of $S$ via the projection $D_{\SL(2)}(\gr^W)^{\sim} \to D_{\SL(2)}(\gr^W)$ (3.5.3) is the log modification $S(\Sigma_Q)$ in 3.1.12 (3.5.9 (i)), and we have the following commutative diagram.
$$
\matrix
S(\Sigma_Q) = D_{\SL(2)}(\gr^W)^{\sim}(Q) & \simeq & B \times \bR^2 \times \{\pm 1\}\\
&&\\
@VVV  @VVV \\
&&\\
S = D_{\SL(2)}(\gr^W)(Q) & \simeq & \; \bR_{\ge0}^2 \times \bR^2 \times \{\pm 1\}.
\endmatrix
$$

In the above isomorphism for $D_{\SL(2)}(\gr^W)^{\sim}(Q)$, the class $p_m$ in $D_{\SL(2)}(\gr^W)^{\sim}(Q)$ of the $\SL(2)$-orbit in Case $m$ in 2.3.9 corresponds to the point  $(b_m, (0, 0), 1)$ of $B \times \bR^2\times \{\pm 1\}$, where $b_m$ is the following point of $B$; $b_1$ is the limit of $(t, 1) \in \bR_{>0}^2$ for $t\to 0$,
$b_2$ is the limit of $(1, t)$ for $t\to 0$, $b_3$ is the limit of $(t, t)$ for $t\to 0$,
$b_4$ is the limit of $(t_0t_1, t_1)$ for $t_0, t_1\to 0$, and $b_5$ is the limit of 
$(t_0, t_0t_1)$ for $t_0, t_1\to 0$.

\medskip

\proclaim{Proposition 3.5.14}
The map $D_{\SL(2)}^{II}\to D_{\SL(2)}(\gr^W)^{\sim}$ is a morphism of $\cB_{\bR}(\log)$. 
\endproclaim

The proof is given together with that of Theorem 3.5.15 below.

\proclaim{Theorem 3.5.15} 
Fix any  $F\in D(\gr^W)$, let $L=\cL(F)$, and let $\bar L$ be the compactified vector space associated to the weightened vector space $L$ of weights $\leq -2$.
Then $D_{\SL(2)}^{II}$ is an $\bar L$-bundle over $\spl(W)\times D_{\SL(2)}(\gr^W)^{\sim}$ in $\cB_{\bR}(\log)$. 
\endproclaim

For the definition of the compactified vector space $\bar L$, see the explanation after 3.2.6 (see \cite{KNU2}, \S7 for details). 
\medskip

{\it Proofs of 3.5.14 and 3.5.15.} 
We deduce 3.5.14 and 3.5.15 from 3.4.4.

Let $p\in D_{\SL(2)}^{II}$ and let $p'$ be the image of $p$ in $D_{\SL(2)}(\gr^W)^{\sim}$. 
Let $\br\in D$ be a point on the torus orbit associated to $p$ and let $\bar\br$ be the image of $\br$ in $D(\gr^W)$.
It is sufficient to show that for some open neighborhood $U$ of $p'$ in $D_{\SL(2)}(\gr^W)^{\sim}$, if we denote the inverse image of $U$ in $D_{\SL(2)}^{II}$ by $\tilde U$, then $\tilde U$ is open in $D_{\SL(2)}^{II}$, the projection $\tilde U\to U$ is a morphism of $\cB_\bR(\log)$, and $\tilde U$ is isomorphic to $U\times \spl(W)\times \bar L$ as an object of $\cB_\bR(\log)$ over $U\times \spl(W)$. 

For $w\in \bZ$, let $p_w=p(\gr^W_w)$ and $\br_w=\br(\gr^W_w)$. 
Take $(R_w,S_w)$ for $(p_w,\br_w)$ as a pair in 3.4.1. 
Let $\Phi = \overline\cW(p)$ and $Q(w) = \cW(p_w)$.
Let $R'$ be an $\bR$-subspace of $\prod_w \Lie(\tilde \rho_w(\bR_{>0}^{Q(w)}))$ 
such that $\prod_w \Lie(\tilde \rho_w(\bR_{>0}^{Q(w)}))=
 \Lie(\tilde\rho(\bR_{>0}^{\Phi}))\oplus R'$. Let 
$R=(\prod_w R_w) \oplus R'$ and $S=\prod_w S_w$. Then $(R, S)$ is a pair
 for $(p, \br)$ as in 3.4.1. 

Let ${\bar Y}(p,\br,S)$ (resp.\ ${\bar Y}(p,\br,R,S)$) be the subset of $Z(p)\times S$ (resp.\ $Z(p,R) \times S$) consisting of all elements $(t,f,g,h,k)$ ($(t,f,g,h)\in Z(p)$ (resp.\ $\in Z(p,R)$), $k\in S$) which satisfy the condition (4) in 3.4.2. We define the structure of ${\bar Y}(p,\br,S)$ (resp.\ 
${\bar Y}(p,\br, R,S)$) as an object of $\cB_\bR(\log)$ just in the same way as in the definition for $Y^{II}(p,\br, S)$ (resp.\ $Y^{II}(p,\br, R, S)$) 
in 3.4.2. 
Note that we have evident isomorphisms in $\cB_\bR(\log)$
$$
Y^{II}(p,\br, S)\simeq {\bar Y}(p, \br,S)\times \bar L \times \fg_{\bR,u}, \quad Y^{II}(p, \br, R, S) \simeq {\bar Y}(p, \br, R,S)\times \bar L \times \fg_{\bR,u}.
$$
 Let ${\bar Y}_0(p, \br, S)$ (resp.\ ${\bar Y}_0(p, \br, R, S)$) be the open set of ${\bar Y}(p,\br, S)$ (resp.\ ${\bar Y}(p,\br, R, S)$) consisting of all elements $(t,f,g,h,k)$ such that $t\in \bR_{>0}^{\Phi}$.

For an open neighborhood $U$ of $0$ in $\fg_\bR(\gr^W)\times \fg_\bR(\gr^W)\times \fg_\bR(\gr^W)\times S$ (resp.\ $\fg_\bR(\gr^W)\times R\times \fg_\bR(\gr^W)\times S$), we define ${\bar Y}(p,\br, S,U)$ (resp.\ ${\bar Y}(p,\br,R,S,U)$) as the open set of ${\bar Y}(p,\br,S)$  (resp.\ ${\bar Y}(p,\br, R, S)$) consisting of all elements $(t,f,g,h,k)$ such that $(f, g, h, k)\in U$. Let  ${\bar Y}_0(p,\br,S,U)={\bar Y}_0(p,\br,S)\cap {\bar Y}(p,\br, S,U)$ (resp.\ ${\bar Y}_0(p,\br,R,S,U)={\bar Y}_0(p,\br,R,S)\cap {\bar Y}(p,\br, R, S,U)$).

\medskip

{\bf Claim 1.} 
{\it For a sufficiently small open neighborhood $U$ of $0$ in 
$\fg_\bR(\gr^W)\times R \times \fg_\bR(\gr^W)\times S$, there is an open immersion ${\bar Y}(p,\br,R,S,U) \to D_{\SL(2)}(\gr^W)^\sim$ in 
$\cB_\bR(\log)$ whose restriction to ${\bar Y}_0(p,\br,R,S,U)$
is given as $(t,f,g,h,k) \mapsto \bar \tau_p(t)\exp(g)\exp(k)\bar\br\in D(\gr^W)$ and which sends $(0^{\Phi}, 0,0,0,0)\in {\bar Y}(p,\br, R, S,U)$ to $p'$.} 

\medskip
We give the proof of Claim 1 later. We need one more Claim. 
\medskip

{\bf Claim 2.} 
{\it Let $q\in D_{\SL(2)}$ and let $(q',s)\in D_{\SL(2)}(\gr^W)^\sim\times \spl(W)$ be the image of $q$. Then the fiber on $(q',s)$ in $D_{\SL(2)}$ regarded as a topological subspace of $D_{\SL(2)}^I$ $($resp.\ $D_{\SL(2)}^{II}$$)$ is  homeomorphic to $\bar L$.}

\medskip

Claim 2 is shown easily.

We show that 3.5.14 and 3.5.15 follow from Claim 1 and Claim 2. Let $U$ be a sufficiently small open neighborhood of $0$ in $\fg_\bR(\gr^W)\times 
R\times \fg_\bR(\gr^W)\times S$, let $U'$ be the image of the open immersion $Y^{II}(p,\br, R, S,U)\to D_{\SL(2)}^{II}$ (3.4.4), and let $U''$ be the image of the open immersion ${\bar Y}(p,\br,R,S,U)\to D_{\SL(2)}(\gr^W)^{\sim}$ (Claim 1). Then $U'\to U''$ is a morphism of $\cB_\bR(\log)$ since 
$Y^{II}(p, \br, R, S, U)\to {\bar Y}(p,\br, R,S,U)$, which is identified with the projection ${\bar Y}(p,\br, R,S,U)\times \bar L\times \fg_{\bR,u}\to 
{\bar Y}(p,\br,R,S,U)$, is a morphism of $\cB_\bR(\log)$. 
 The map $U' \to U''
 \times \spl(W)$ is a trivial $\bar L$-bundle  since $Y^{II}(p,\br,R,S,U)\to {\bar Y}(p,\br,R,S,U)\times 
 \spl(W)$ is identified with the projection ${\bar Y}(p,\br,R,S,U)\times \bar L\times \spl(W) \to {\bar Y}(p,\br, R, S,U)\times \spl(W)$. Hence this morphism is proper. Let $V$ be the inverse image of $U'' \times \spl(W)$ under the canonical map $D_{\SL(2)}^{II}\to D_{\SL(2)}(\gr^W)^{\sim}\times \spl(W)$. We prove $V=U'$. Indeed, since $U'$ is proper over $U''\times \spl(W)$, $U'$ is open and closed in $V$. Since all fibers of $V\to U''\times \spl(W)$ is connected by Claim 2, and since $U'\to U''\times \spl(W)$ is surjective, we have $V=U'$. Hence $V$ is open in $D_{\SL(2)}^{II}$,  $V\to U''$ is a morphism of $\cB_\bR(\log)$, and $V\to U''\times \spl(W)$ is a trivial $\bar L$-bundle.

We prove Claim 1. 

For each $w\in \bZ$, let $Q(w)\in \cW(\gr^W_w)$ be the image of $\Phi$.
 For each $w\in \bZ$, by 3.4.4 for the pure case, there is an open neighborhoods $U_w$ of $0$ in $\fg_\bR(\gr^W_w)\times \fg_\bR(\gr^W_w)\times \fg_\bR(\gr^W_w)\times S_w$ 
such that we have a morphism  $Y^{II}(p_w,\br_w,S_w,U_w)\to D_{\SL(2)}(\gr^W_w)$ which sends $(t,f,g,h,k)\in Y_0^{II}(p_w,\br_w,S_w,U_w)$ to
 $\tau_{p_w}(t)\exp(g)\exp(k)\br_w$, which induces an open immersion 
 $Y^{II}(p_w,\br_w,R_w,S_w,U_w')\to D_{\SL(2)}(\gr^W_w)$ 
($U'_w:=U_w\cap(\fg_\bR(\gr^W_w)\times R_w\times \fg_\bR(\gr^W_w)\times S_w))$, and which sends 
 $(0^{Q(w)},0,0,0,0)\in Y^{II}(p_w,\br_w,R_w,S_w,U'_w)$ to $p_w$. By 3.4.13 for the pure case, for some open neighborhood $U''_w\sub U_w$ of $0$ in 
$\fg_\bR(\gr^W_w)\times \fg_\bR(\gr^W_w)\times \fg_\bR(\gr^W_w)\times S_w$, we have a morphism $Y^{II}(p_w,\br_w,S_w,U''_w)\to Y^{II}(p_w,\br_w,R_w,S_w,U'_w)$ which commutes with the morphisms to $D_{\SL(2)}(\gr^W_w)$.
Let ${\bar Y}(p,\br,S)\to Y^{II}(p_w,\br_w,S_w)$ be the morphism
$(t,f,g,h,k)\mapsto (t(\gr^W_w), f(\gr^W_w),g(\gr^W_w),h(\gr^W_w),k(\gr^W_w))$,
where $t(\gr^W_w)$ denotes the image of $t$ under the homomorphism $\bR_{\geq 0}^{\Phi}\to \bR_{\geq 0}^{Q(w)}$ of multiplicative monoids induced by the map  $\Phi\to Q(w)$. Then if $U$ is a 
sufficiently small open neighborhood of $0$ in $\fg_\bR(\gr^W)\times 
R\times \fg_\bR(\gr^W)\times S$, the image of 
${\bar Y}(p,\br, S,U)$ under this morphism is contained in $Y^{II}(p_w,\br_w,S_w,U''_w)$ for any $w$. Hence we have a composite morphism
$$
\xi: {\bar Y}(p,\br, S,U)\to  \tsize\prod_w Y^{II}(p_w,\br_w,S_w,U_w'')\to  \tsize\prod_w Y^{II}(p_w,\br_w, R_w, S_w, U'_w).
$$
 
 Let $P$ be the fiber product of
$$
\tsize\prod_w Y^{II}(p_w, \br_w, R_w, S_w, U'_w) \to 
\tsize\prod_w\bR_{\geq 0}^{Q(w)}\leftarrow \bR_{\geq 0}^{\Phi}\times^{\bR_{>0}^{\Phi}} ( \tsize\prod_w \bR_{>0}^{Q(w)})
$$ 
in $\cB_\bR(\log)$.
Here $\bR_{\geq 0}^{\Phi}\times^{\bR_{>0}^{\Phi}} 
 ( \prod_w \bR_{>0}^{Q(w)})$ is the quotient of 
 $\bR_{\geq 0}^{\Phi}\times ( \prod_w \bR_{>0}^{Q(w)})$ 
 under the action of $\bR_{>0}^{\Phi}$ given by $(x,y)\mapsto (ax, a^{-1}y)$ ($a\in \bR_{>0}^{\Phi}$).  Then $P$ is identified with the fiber product of
$$
\tsize\prod_w Y^{II}(p_w, \br_w, R_w, S_w, U_w') 
\to  \tsize\prod_w D_{\SL(2)}(\gr^W_w)(Q(w))\leftarrow D_{\SL(2)}(\gr^W)^\sim(\Phi).
 $$ 
 Hence we have an open immersion $P\to D_{\SL(2)}(\gr^W)^\sim$.
 
We have a unique morphism $$\xi^{\sim}: {\bar Y}(p,\br, S,U)\to P$$
in $\cB_\bR(\log)$ which is compatible with $\xi$. It is induced from $\xi$ and from the morphism $ {\bar Y}(p,\br, S,U)\to 
\bR_{\geq 0}^{\Phi} \times^{\bR_{>0}^{\Phi}}
 (\prod_w \bR_{>0}^{Q(w)})$ which sends $(t,f,g,h,k)$ to $tt'$, 
where $t'\in \prod_w\bR_{>0}^{Q(w)}$ is the $(\prod_w \bR_{>0}^{Q(w)})$-component of $\xi(1,g,g,g,k)$. 

\medskip

{\bf Claim 3.} 
{\it If $U$ is a sufficiently small open neigborhood of $0$ in 
$T:=\fg_{\bR}(\gr^W)\times R \times \fg_{\bR}(\gr^W)\times S$, the morphism ${\bar Y}(p,\br,R,S,U)\to P$ induced by $\xi^\sim$ is an open immersion.}

\medskip

By Claim 3, the open immersion stated in Claim 1 is obtained as the composite 
${\bar Y}(p,\br,R,S,U)\to P \to D_{\SL(2)}(\gr^W)^\sim$. 
It remains to prove Claim 3. 

For an open neighborhood $U$ of $0$ in $T$, let $P(U)$ be the open set of
$P$ consisting of all elements $(t,f,g,h,k)$ ($t\in \bR_{\geq 0}^{\Phi}\times^{\bR_{>0}^{\Phi}} (\prod_w \bR_{>0}^{Q(w)})$, $f, h\in \fg_\bR(\gr^W)$, $g\in \prod_w R_w$, $k\in \prod_w S_w$) such that 
$t=t'\exp(a)$ for some $t'\in \bR_{\geq 0}^{\Phi}$  and for some $a\in R'$ satisfying  $(f, a+g,h, k)\in U$. 
Then, for a given open neighborhood $U$ of $0$ in $T$, there is an open neighborhood $U'$ of $0$ in $T$ such that the map $\xi^\sim$ induces a morphism ${\bar Y}(p,\br, R, S, U')\to P(U)$. On the other hand, if $U$ is an open neighborhood of $0$ in $T$, then for a sufficiently small open neighborhood $U'$ of $0$ in $T$, we have a morphism $P(U')\to {\bar Y}(p,\br, R, S,U)$. 
This morphism is obtained as the composite 
$P(U') \to {\bar Y}(p,\br, S, U'') \to {\bar Y}(p, \br, R, S, U)$. Here $U''$ is a suitable open neighborhood of $0$ in $T$. The first arrow is
$(t'\exp(a), f,g,h,k) \mapsto (t', f', g', h', k)$, where $f'$, $g'$, $h'$ are near to $f,g,h$, respectively, and defined by $\exp(g')=\exp(a)\exp(g)$, $\exp(f')
=\exp(f)\exp(a)$, $\exp(h')=\exp(2a)\exp(g)\exp(-a)$.  
The second arrow is
a morphism 
constructed in the same way as in the proof of 3.4.13. For an open neighborhood $U$ of $0$ in $T$, the composite 
${\bar Y}(p,\br, R, S,U'')\to P(U')\to {\bar Y}(p,\br, R, S,U)$ and the 
composite 
$P(U'') \to {\bar Y}(p, \br, R, S, U')\to P(U)$ are inclusion maps. Here $U'$ and $U''$ are open neighborhoods of $0$ in $T$, $U'$ is sufficiently small relative to $U$, and $U''$ is sufficiently small relative to $U'$. This proves Claim 3. 
\qed

\proclaim{Theorem 3.5.16}
The canonical map
$$ 
D_{\SL(2)}^{II}\to \spl(W) \times D_{\SL(2)}(\gr^W)
$$ 
is proper. 
\endproclaim

{\it Proof.} 
The map $D_{\SL(2)}^{II}\to \spl(W) \times D_{\SL(2)}(\gr^W)^{\sim}$ is proper by 3.5.15. 
The map $D_{\SL(2)}(\gr^W)^{\sim} \to D_{\SL(2)}(\gr^W)$ is proper (3.5.9--3.5.10).
\qed
\bigskip

\proclaim{Theorem 3.5.17}
Let $\G$ be a subgroup of $G_\bZ$. 
For $* = I, II$, we have the following.
\medskip

{\rm(i)} The action of $\G$ on $D_{\SL(2)}^*$ is proper, and the quotient space $\G\bs D_{\SL(2)}^*$ is Hausdorff. 

\medskip

{\rm (ii)} Assume  that $\G$ is neat. 
Let $\g\in \G$,  $p\in D_{\SL(2)}$, and assume $\g p=p$. Then $\g=1$. 

\medskip

{\rm (iii)} 
Assume that $\G$ is neat. 
Then the quotient $\G \bs D_{\SL(2)}^*$ belongs to $\cB_\bR(\log)$, and the projection $D_{\SL(2)}^* \to \G\bs D_{\SL(2)}^*$ is a local isomorphism of objects of $\cB_\bR(\log)$.

\endproclaim
Here in (iii), we define the sheaf of real analytic functions on 
$\G\bs D^*_{\SL(2)}$ and the log structure with sign on $\G \bs D_{\SL(2)}^*$ in the natural way. 
That is, for an open set $U$ of $\G\bs D_{\SL(2)}^*$, a real valued function $f$ on $U$ is said to be real analytic if the pull-back of $f$ on the inverse image of $U$ in $D_{\SL(2)}^*$ is real analytic. The log structure $M$ of $\G \bs D_{\SL(2)}^*$ is defined to be the sheaf of real analytic functions whose pull-backs on $D_{\SL(2)}^*$ belong to the log structure of $D_{\SL(2)}^*$. The subgroup sheaf $M^{\gp}_{>0}$ of $M^{\gp}$ is defined to be the part of $M^{\gp}$ 
consisting of the local sections whose pull-backs to $D_{\SL(2)}^*$ belong to the $M_{>0}^\gp$ of $D_{\SL(2)}^*$.

Recall that 
a subgroup $\G$ of $G_\bZ$ is said to be {\it neat} if for any $\g \in \G$, the 
subgroup of $\bC^\times$ generated by all eigenvalues of
the action of $\g$ on $H_{0,\bC}$ is torsion free. If $\G$ is neat,
then $\G$ is torsion free. 
There exists a neat subgroup of $G_{\bZ}$ of finite index (cf\. \cite{Bo}).

\medskip

{\it Proof of 3.5.17.} The proof is similar to 
\cite{KNU2}, \S9, where we considered $D_{\BS}$.

(i) 
$D_{\SL(2)}^{II}$ is Hausdorff because $D_{\SL(2)}(\gr^W)$ is Hausdorff (\cite{KU2}) and the map $D_{\SL(2)}^{II}\to \spl(W) \times D_{\SL(2)}(\gr^W)$ 
is proper (3.5.16). It follows that $D_{\SL(2)}^I$ is also Hausdorff. 

Let $\G_u$ be the kernel of 
$\G\to \Aut(\gr^W)$. The properness of the action of $\G$ on $D_{\SL(2)}^{II}$ is reduced to the properness of the action of $\G/\G_u$ on 
$D_{\SL(2)}(\gr^W)$ which is proved in \cite{KU2},
 and to the properness of the action of $\G_u$ on $\spl(W)$. 
  The properness of that on $D_{\SL(2)}^{I}$ follows from this 
because $D_{\SL(2)}^{I}$ is Hausdorff. 
 
Since the action of $\G$ on $D_{\SL(2)}^*$ for $*=I, II$ is proper, 
the quotient space $\G \bs D_{\SL(2)}^*$ is Hausdorff.

\medskip

(ii) The pure case is proved in \cite{KU2}. The general case is reduced to the pure case since the action of $\G_u$ on $\spl(W)$ is fixed point free.

\medskip

(iii) By (i) and (ii), the map $D_{\SL(2)}^*\to \G \bs D_{\SL(2)}^*$ is a local homeomorphism. (iii) follows from this.  
\qed
\medskip

\vskip20pt

\head
\S3.6. Examples
\endhead

\bigskip

We consider $D_{\SL(2)}^I$ and $D_{\SL(2)}^{II}$ for five examples I--V in 1.1.1. 

\medskip

{\bf 3.6.1.} 
We consider  $D_{\SL(2)}^{II}$. 

We use the notation in 1.1.1.
As in 1.2.9, we denote by $L$ the graded vector space $\cL(F)=L^{-1,-1}_\bR(F)\sub \cL$ with $F\in D(\gr^W)$, which is independent of the choice of $F$ for Examples I--V. 
Recall that $D_{\SL(2)}^{II}$ is an $\bar L$-bundle over $\spl(W) \times D_{\SL(2)}(\gr^W)^\sim$ (Theorem 3.5.15), and that for Examples I--IV, $D_{\SL(2)}(\gr^W)^{\sim} = D_{\SL(2)}(\gr^W)$ (3.5.12). 
We will describe the structure of the  open set $D_{\SL(2)}^{II}(\Phi)$ of $D_{\SL(2)}^{II}$ for some $\Phi\in \overline \cW$.

Let $\bar \fh=\{x+iy\;|\;x, y \in \bR, 0<y\leq \infty\}\supset \fh$. 
We regard $\bar \fh$ as an object of $\cB_\bR(\log)$ via $\bar \fh \simeq \bR_{\geq 0}\times \bR, \; x+iy\mapsto (1/\sqrt{y}, x)$ (cf. 3.2.13). 
\medskip

{\bf Example I.}  
We have a commutative diagram in $\cB_{\bR}(\log)$
$$
\matrix D & \simeq & \spl(W) \times L\\ 
&&\\
\cap & &\cap \\ 
&&\\
D_{\SL(2)}^{II} & \simeq & \;\spl(W) \times \bar L, \endmatrix 
$$
where the upper isomorphism is that of 1.2.9. 
Here $\spl(W)\simeq \bR$, $D_{\SL(2)}(\gr^W)=D(\gr^W)$ which is just a one point set, $L\simeq \bR$ with weight $-2$, and $\bar L$ is isomorphic to the interval $[-\infty, \infty]$ endowed with the real analytic structure as in \cite{KNU2}, 
7.5, with $w=-2$ which contains $\bR=L$ in the natural way (1.2.9). 

\medskip

{\bf Example II.} 
Let $Q=\{W'\}\in \cW(\gr^W_{-1})=\prod_w \cW(\gr^W_w)$, where
$$
W'_{-3}=0\sub W'_{-2}=W'_{-1}=\bR e'_1\sub W'_0=\gr^W_{-1}.
$$
The isomorphism $D(\gr^W)=D(\gr^W_{-1})\simeq \fh$ extends to 
$D_{\SL(2)}(\gr^W)(Q) \simeq \bar \fh$.

  Let $\Phi$ be the unique non-empty element of $\overline \cW(Q)$. 
We have a commutative diagram in $\cB_{\bR}(\log)$
$$
\matrix D & \simeq & \spl(W) \times \fh \\ 
&&\\
\cap & &\cap \\
&&\\
D_{\SL(2)}^{II}(\Phi) & \simeq &  \;\spl(W) \times \bar \fh. \endmatrix
$$
Recall that $\spl(W)\simeq \bR^2$ (1.2.9). 
In this diagram, the upper isomorphism is that of 1.2.9. 
The lower isomorphism is induced by the canonical morphisms $D_{\SL(2)}^{II}\to \spl(W)$ and $D_{\SL(2)}^{II}(\Phi) \to D_{\SL(2)}(\gr^W)(Q)\simeq \bar \fh$. 
 
 The specific examples of $\SL(2)$-orbits of rank $1$ in 2.3.9 Example II have classes in $D_{\SL(2)}^{II}(\Phi)$ whose images in $\bar \fh$ are $i\infty$. 
\medskip

{\bf Example III.}  
Let $Q=\{W'\}\in \cW(\gr^W_{-3})=\prod_w \cW(\gr^W_w)$, where
$$
W'_{-5}=0\sub W'_{-4}=W'_{-3}=\bR e'_1\sub W'_{-2}=\gr^W_{-3}.
$$
The isomorphism $D(\gr^W)=D(\gr^W_{-3})\simeq \fh$ extends to 
$D_{\SL(2)}(\gr^W)(Q) \simeq \bar \fh$.

  Let $\Phi$ be the unique non-empty element of $\overline \cW(Q)$. 
We have a commutative diagram in $\cB_{\bR}(\log)$
$$
\matrix D & \simeq & \spl(W) \times \fh \times L&&(s,\, x+iy,\, (d_1,\, d_2)) \\ 
&&\\
\cap & &\downarrow &&\downarrow\\
&&\\
D_{\SL(2)}^{II}(\Phi) & \simeq & \spl(W) \times \bar \fh \times \bar L&&\;
(s,\, x+iy,\, (y^{-2}d_1,\, y^{-1}d_2)). \endmatrix  
$$
Here $\spl(W) \simeq \bR^2$, $L\simeq \bR^2$ with weight $-3$, and $(d_1,d_2)\in \bR^2=L$ (1.2.9). 
In this diagram, the upper isomorphism is that of 1.2.9. 
The lower isomorphism is induced by the canonical morphisms $D_{\SL(2)}^{II} \to \spl(W)$ and  $D_{\SL(2)}^{II}(\Phi)
\to D_{\SL(2)}(\gr^W)(Q) \simeq \bar \fh$, and the following morphism $D_{\SL(2)}^{II}(\Phi)\to \bar L$. 
It is induced by $\nu_{\a,\b}$, where $\a_{-3}:
\bG_{m,\bR}\to \Aut(\gr^W_{-3})$ is 
defined by $\a_{-3}(t)e'_1=t^{-4}e'_1$, $\a_{-3}(t)e'_2=t^{-2}e'_2$, and $\b:D(\gr^W_{-3})=\fh\to \bR_{>0}$ is the distance to $\Phi$-boundary defined by $x+iy\mapsto 1/\sqrt{y}$ (3.2.13). 
Note that the right vertical arrow is {\it not} the evident map, as indicated.
 
The $\SL(2)$-orbits in 2.3.9 Example III, Case 1 (resp. Case 2, resp. Case 3) have classes in $D_{\SL(2)}^{II}(\Phi)$ 
whose images in $\bar \fh \times \bar L$ belong to $\{i\infty\}\times L$ (resp. $\{i\}\times (\bar L\smallsetminus L)$, resp. $\{i\infty\}  \times (\bar L\smallsetminus L)$). 
\medskip

{\bf Example IV.}  
Let $Q=\{W'\}\in \cW(\gr^W_{-1})=\prod_w \cW(\gr^W_w)$, where
$$
W'_{-3}=0\sub W'_{-2}=W'_{-1}=\bR e'_2 \sub W'_0=\gr^W_{-1}.
$$
The isomorphism $D(\gr^W)=D(\gr^W_{-1})\simeq \fh$ extends to 
$D_{\SL(2)}(\gr^W)(Q) \simeq \bar \fh$.

  Let $\Phi$ be the unique non-empty element of $\overline \cW(Q)$. 
We have a commutative diagram in $\cB_{\bR}(\log)$
$$
\matrix D & \simeq & \spl(W) \times \fh \times L &&(s,\, x+iy,\, d) \\ 
&&\\
\cap & &\downarrow&&\downarrow \\ 
&&\\
D_{\SL(2)}^{II}(\Phi) & \simeq &  \spl(W) \times \bar \fh \times \bar L&&
\;(s,\, x+iy,\, y^{-1}d).\endmatrix  
$$
Here  $\spl(W) \simeq \bR^5$, $L\simeq \bR$ with weight $-2$, and $d\in \bR=L$ (1.2.9).
In this diagram, the upper isomorphism is that of 1.2.9.
The lower isomorphism is induced from the canonical morphisms  $D_{\SL(2)}^{II} \to \spl(W)$ and $D_{\SL(2)}^{II}(\Phi)\to D_{\SL(2)}(\gr^W)(Q) \simeq \bar \fh$, and the following morphism $D_{\SL(2)}^{II}(\Phi)\to \bar L$. 
It is induced by $\nu_{\a,\b}$ (3.2.6--3.2.10), where $\a_{-1}:
\bG_{m,\bR}\to \Aut(\gr^W_{-1})$ is 
defined by 
$$
\a_{-1}(t)e_2'=t^{-2}e_2', \;\;
\a_{-1}(t)e_3'=e_3'
$$ 
and $\b:D(\gr^W_{-1})=\fh\to \bR_{>0}$ is the distance to $\Phi$-boundary defined by $x+iy\mapsto 1/\sqrt{y}$ (3.2.13). 
Note that the right vertical arrow is {\it not} the inclusion map, as indicated.

The $\SL(2)$-orbits in 2.3.9 Example IV, Case 1 (resp. Case 2, resp. Case 3) have classes in $D_{\SL(2)}^{II}(\Phi)$ whose images in $\bar \fh \times \bar L$ belong to $\{i\infty\}\times L$ (resp. $\{i\} \times (\bar L\smallsetminus L)$, resp. $\{i\infty\} \times (\bar L\smallsetminus L)$).
\medskip

{\bf Example V.}  
Let $Q\in \prod_w \cW(\gr^W_w)$ and the log modification $B$ of $\bR_{\geq 0}^2$ be as in 3.5.13. 
The isomorphism $D(\gr^W) \simeq \fh^{\pm} \times \fh$ (1.2.9) extends to an isomorphism $D_{\SL(2)}(\gr^W)(Q)\simeq \bar \fh^{\pm} \times \bar \fh$ ($\bar \fh^{\pm}$ is the disjoint union of $\bar \fh^{+}=\bar \fh$ and $\bar \fh^{-}=\{x+iy\;|\;x\in \bR, 0>y\geq -\infty\}$ ($\fh^+\simeq \fh^{-}, \;x+iy\mapsto -x-iy$)), and 
this composite isomorphism is extended to an isomorphism $D_{\SL(2)}(\gr^W)^\sim(Q)\simeq B \times \bR^2\times \{\pm 1\}$ (3.5.13). 
 
  Let $\Phi$ be the maximal element of $\overline \cW(Q)$. 
We have a commutative diagram in $\cB_{\bR}(\log)$
$$
\matrix D & \simeq & \spl(W) \times \fh^{\pm} \times \fh&&(s,\, x+iy,\, x'+iy') \\ 
&&\\
\cap & &\downarrow && \downarrow\\
&&\\
D_{\SL(2)}^{II}(\Phi) & \simeq & \spl(W) \times B \times \bR^2\times \{\pm 1\}&&
\;(s,\, 1/\sqrt{|y|},\, 1/\sqrt{y'},\, x, x',\, \sign(y)).  \endmatrix  
$$
Here $\spl(W) \simeq \bR^6$ (1.2.9). 
In this diagram, the upper isomorphism is that of 1.2.9.
The lower isomorphism is induced from the canonical morphisms $D_{\SL(2)}^{II}\to \spl(W)$ and $D_{\SL(2)}^{II}(\Phi)\to D_{\SL(2)}(\gr^W)^{\sim}(Q)\simeq 
B \times \bR^2\times \{\pm 1\}$. 

The $\SL(2)$-orbits in 2.3.9 Example V have classes in $D_{\SL(2)}^{II}(\Phi)$ whose images in $B$ are described in 3.5.13. 

\medskip

{\bf 3.6.2.} 
We consider $D_{\SL(2)}^I$. 
For Examples I--IV, $D_{\SL(2)}^I=D_{\SL(2)}^{II}$ by Proposition 3.4.29.
\medskip

{\bf Example V.} 
Let $\Psi=\{W'\}\in \cW$, where 
$$
\align
W'_{-3}=0 \sub W'_{-2}=W'_{-1}=\bR e_1&\sub W'_0=W'_{-1}+\bR e_2\\
&\sub W'_1=W'_0+ \bR e_4 +\bR e_5\sub W'_2=H_{0,\bR}.
\endalign
$$
(This $W'$ is $W^{(1)}$ in 2.3.9 V, Case 1.) 
Let $\bar \Psi=\{W'(\gr^W)\}\in \overline{\cW}$. 
Then $D_{\SL(2)}^{II}(\bar \Psi)$ is the open set of $D_{\SL(2)}^{II}(\Phi)$ in 3.6.1 V corresponding to the subcone $\bR_{\geq 0}\times \{0\}$ of $\bR_{\geq 0}^2$. 

We compare $D_{\SL(2)}^I(\Psi)$ and $D_{\SL(2)}^{II}(\bar \Psi)$. For $j=1, 2, 3$, let
$$
A_j=\Hom_\bR(\gr^W_1, \bR e_j).
$$ 
We have an isomorphism of real analytic manifolds 
$$
\spl(W) @>\sim>> \tsize\prod_{j=1}^3 A_j, \quad s \mapsto (a_j)_{1\leq j\leq 3},
$$ 
$$
\text{where}\;\;\; s(v)\equiv \tsize\sum_{j=1}^3 a_j(v) \bmod \bR e_4+\bR e_5 \;\;\;\text{for} \;\; v\in \gr^W_1.
$$ 
Let   
$$
(A_3\times \bar \fh^{\pm})':=\{(v, x+iy)\in A_3\times \bar \fh^{\pm}\;|\;v=0\;\text{if}\;y=\pm \infty\}\sub A_3\times \bar \fh^{\pm}.
$$ 

Then we have a commutative diagram in $\cB_\bR(\log)$
$$
\matrix D & \simeq & (\tsize\prod_{j=1}^3  A_j)\times  \fh^{\pm }\times \fh  \\ 
&&\\
\cap & &\cap \\
&&\\
D_{\SL(2)}^{II}(\bar \Psi) & \simeq & \;(\tsize\prod_{j=1}^3 A_j) \times \bar \fh^{\pm}  \times \fh. \endmatrix  
$$
In this diagram, the upper isomorphism is induced by the isomorphism in 1.2.9 and the above isomorphism $\spl(W)\simeq \prod_{j=1}^3 A_j$. 
On the other hand, we have a commutative diagram in $\cB_\bR(\log)$ 
$$
\matrix D & \simeq & (\prod_{j=1}^3 A_j)\times  \fh^{\pm} \times \fh &&(a_1,\, a_2,\, a_3,\, x+iy, \tau) \\ 
&&\\
\cap & &\downarrow &&\downarrow \\
&&\\
D_{\SL(2)}^I(\Psi) & \simeq & A_1\times A_2\times (A_3 \times \bar \fh^{\pm})' \times \fh   && 
\;(a_1,\, a_2,\, |y|^{1/2}a_3,\, x+iy, \tau). \endmatrix  
$$
In this diagram, the upper isomorphism is the same as in the first diagram. 
The lower isomorphism is induced from the canonical morphisms $D_{\SL(2)}^I \to \spl(W) \to A_1\times A_2$ and $D_{\SL(2)}^I(\Psi) \to D_{\SL(2)}(\gr^W)^\sim(\bar \Psi) \simeq \bar \fh^{\pm} \times \fh$, and the following morphism $D_{\SL(2)}^I(\Psi)
\to A_3$. It is the composite
$$
D_{\SL(2)}^I(\Psi) @>\text{by $\nu_{\a,\b}$}>> D @>\spl_W>> \spl(W) \simeq \tsize\prod_{j=1}^3 A_j \to A_3,
$$
where 
$\nu_{\a,\b}$ is the morphism described in 3.2.6--3.2.10.
Here 
$\a: \bG_{m,\bR}\to \Aut_\bR(H_{0,\bR}, W)$ is the splitting of $\Psi$ defined by $\a(t)e_1=t^{-2}e_1$, $\a(t)e_2=e_2$, $\a(t)e_3=t^2e_3$, $\a(t)e_4=te_4$, $\a(t)e_5=te_5$, and $\b:D \to \bR_{>0}$ is the distance to $\Psi$-boundary defined as the composite $D\to D(\gr^W_0) \simeq \fh^{\pm}\to \bR_{>0}$, where 
the last arrow is $x+iy\mapsto 1/\sqrt{|y|}$. 

Note that the right vertical arrow of the above commutative diagram is {\it not} the inclusion map, as indicated.

The lower isomorphisms in the above two commutative diagrams form a commutative diagram in $\cB_\bR(\log)$ 
$$
\matrix D_{\SL(2)}^I(\Psi) & \simeq &  A_1\times A_2\times  (A_3\times \bar \fh^{\pm})' \times \fh  &\ni&(a_1,\, a_2,\, a_3,\, x+iy,\, \tau)  \\ 
&&\\
\downarrow & &\downarrow&&\downarrow  \\
&&\\
D_{\SL(2)}^{II}(\bar \Psi) & \simeq & (\prod_{j=1}^3 A_j) \times \bar \fh^{\pm} \times \fh&\ni&
\;(a_1,\, a_2,\, |y|^{-1/2}a_3,\, x+iy,\, \tau). \endmatrix  
$$
Here the left vertical arrow is the inclusion map. 
The right vertical arrow is {\it not} the evident map, as indicated. 

The $\SL(2)$-orbits in 2.3.9 Example V, Case 1 have classes in $D_{\SL(2)}^I(\Psi)$ whose images in $\bar \fh^{\pm}\times \fh$ are $(i\infty, i)$. 

\vskip 20pt

\head
\S3.7. $D_{\BS,\val}$ and $D_{\SL(2),\val}$
\endhead

\bigskip

We  outline the definitions of $D_{\SL(2),\val}$ and $D_{\BS,\val}$ in the fundamental diagram in 0.2, which connect $D_{\SL(2)}$ and $D_{\BS}$.  
The detailed studies of these spaces will be given later in this series of papers.
\medskip

{\bf 3.7.1.} 
Let $S$ be an object of $\cB_\bR(\log)$ (3.1). 
Then we have a local ringed space $S_{\val}$ over $S$ with a log structure 
with sign. 
This is the real analytic analogue of the complex analytic theory considered in \cite{KU3}, \S3.6.
In the case when we have a chart $\cS\to M_{S,>0}$ with $\cS$ an fs monoid, 
$$
S_{\val}= \varprojlim_{\Sig} \;S(\Sig),
$$
where $\Sig$ ranges over all finite rational subdivisions of the cone 
$\Hom(\cS, \bR_{\geq0}^{\add})$ (3.1.12). 
The general case is reduced to this case by gluing 
(cf.\ \cite{KU3}, \S3.6). 

\medskip

{\bf 3.7.2.} 
For $* = I, II$, define $D_{\SL(2),\val}^*=(D_{\SL(2)}^*)_{\val}$.   
In the pure case, as topological spaces, they coincide with the topological space $D_{\SL(2),\val}$ in \cite{KU2}.
\medskip

{\bf 3.7.3.} 
$D_{\BS,\val}$ is defined similarly, 
that is, $D_{\BS,\val}=(D_{\BS})_{\val}$. 
  Here we use the log structure with sign of $D_{\BS}$ 
induced by $\bar A_P\simeq \bR^n_{\geq0}$ and 
$\bar B_P \simeq \bR^{n+1}_{\geq0}$ in the notation in \cite{KNU2}, 5.1. 

\medskip

{\bf 3.7.4.}
A canonical injection $D^*_{\SL(2),\val} \to D_{\BS,\val}$ is defined but not 
necessarily continuous (both for $*=I$ and $II$). 
  This is a difference from the pure case, and we try to explain it 
a little more in the next subsection. 

\vskip 20pt

\head
\S3.8. $D_{\BS}$ and $D_{\SL(2)}$
\endhead

\bigskip

  Here in the end of this section, we review some points of our 
constructions and compare it with the construction of $D_{\BS}$ in 
\cite{KNU2}. 

\medskip

{\bf 3.8.1.} 
  First, see 1.2.5, which shows that there are three kinds of coordinate 
functions on $D$, that is, $s$, $F$, and $\delta$. 
  Among these, what are new in the mixed case are $s$ and $\delta$.
  Thus, when we want to endow
a partial compactification like $D_{\SL(2)}$ and 
$D_{\BS}$ with  
a real analytic structure 
by extending coordinate functions, we have to treat 
$s$ and $\d$. 
  Among these, $s$ is more important in applications, and 
the methods to treat $s$ are common to the cases of 
$D_{\SL(2)}$ and $D_{\BS}$.

\medskip

{\bf 3.8.2.}
  On the other hand, the treatment of the $\delta$ coordinate 
for $D_{\SL(2)}$ and that for $D_{\BS}$ 
are 
considerably different.
  See 3.6.1 III and IV, which illustrate the situation of $D_{\SL(2)}$. 
  In there, 
the third components ($\delta$ coordinates) 
of the vertical arrows 
in the diagrams are not the inclusion maps but the twisted ones. 
  In general, the $\bar L$-component of the function which gives the 
real analytic structures on $D_{\SL(2)}$ is not the evident one but 
the one twisted back by torus actions (cf.\ 3.2.6).  
  This twisting 
is natural in view of the relationship with nilpotent orbits and 
crucial in the applications (cf.\ 2.5.7).  

\medskip

{\bf 3.8.3.}
  In the case of $D_{\BS}$, the $\delta$ coordinate 
was also naturally twisted, but there is a difference between 
these two twistings, which explains the discontinuity of 
$D_{\SL(2),\val} \to D_{\BS,\val}$ in 3.7.4. 
  
  More precisely, for example, consider Example III in 3.6.1. 
  Let $p$ be a point of $D_{\SL(2),\val}$. 
  Then, the $\bar L$-component 
of the image of $p$ in $D_{\SL(2)}^{II}$ 
is in the boundary (i.e., belongs to $\bar L \smallsetminus L$) if and 
only if $W \in \cW(p)$, but the $\bar L$-component of 
its image in $D_{\BS}$ is in the boundary 
if and only if $p$ is not split. 
  Hence some arc joining a split point and a 
non-split point in $D_{\SL(2),\val}$ can have 
a disconnected image on $D_{\BS}$.
  These equivalences hold for any Hodge types, and 
we can even prove that for some Hodge types, 
there are no choices of topologies of $D_{\SL(2)}$ 
satisfying both the crucial property 
2.5.7 (ii) and the continuities of the maps 
$D_{\SL(2),\val} \to D_{\BS,\val}$ etc.\ in the fundamental diagram in 0.2.
These topics will be treated later in this series. 

\vskip 20pt

\head
\S4. Applications
\endhead
\medskip

\head
\S4.1. Nilpotent orbits, SL(2)-orbits, and period maps
\endhead
\medskip

In \cite{KNU1}, we generalized the SL(2)-orbit theorem in several variables of Cattani-Kaplan-Schmid for degenerations of polarized Hodge structures, to an SL(2)-orbit theorem in several variables for degenerations of mixed Hodge structures with polarized graded quotients. 
Here we interpret it in the style of a result on the extension of a period map into $D_{\SL(2)}$ defined by a nilpotent orbit. 
\medskip

\proclaim{Theorem 4.1.1}
Assume $(N_1, \dots, N_n, F)$ generates a nilpotent orbit $(2.4.1)$ and the associated $W^{(j)}(\gr^W)$ is rational $(2.2.2)$ for any $j = 1,\dots, n$. 
Then, there is a sufficiently small open neighborhood $U$ of $\bold0:=(0,\dots,0)$ in $\bR_{\ge0}^n$ satisfying the following {\rm (i)} and {\rm (ii)}. 

\medskip

{\rm (i)} 
The real analytic map
$$
p : U\cap \bR^n_{>0}\to D, \;\;t=(t_1,\dots, t_n)\mapsto \exp(\ts_{j=1}^n iy_jN_j)F,
$$
where $y_j=\tp_{k=j}^n t_k^{-2}$, is defined and 
extends to a real analytic map
$$ 
p:U\to D_{\SL(2)}^I.
$$

\medskip

{\rm (ii)} 
For $c \in U$, $p(c)\in D_{\SL(2)}$ is described as follows. 
Let $K = \{j\;|\; 1 \le j \le n,\, c_j=0\}$, and write $K = \{b(1), \dots, b(m)\}$ with $b(1)<\dots<b(m)$. 
Let $b(0)=0$. 
For $1\leq j\leq m$, let $N'_j = \ts_{b(j-1)< k\leq b(j)} (\prod_{k\leq \ell <b(j)} c_{\ell}^{-2})N_k$, where $\prod_{b(j)\leq \ell <b(j)} c_{\ell}^{-2}$ is considered as $1$. 
Let $F'=\exp(i\ts_{b(m) < k \leq n} (\prod_{k \leq \ell \leq n} c_{\ell}^{-2})N_k)F$.
Then $(N'_1, \dots, N'_m, F')$ generates a nilpotent orbit $(2.4.1)$, and $p(c)$ is the class of the $\SL(2)$-orbit associated to $(N'_1, \dots, N'_m, F')$ $(2.4.2)$. 
Hence, when $t \in U$ and $t \to c$, we have the convergence
$$
\exp(\ts_{j=1}^n iy_jN_j)F \to \text{$($class of the $\SL(2)$-orbit associated to $(N'_1, \dots, N'_m, F'))$}
$$
in $D_{\SL(2)}^I$ and hence in $D_{\SL(2)}^{II}$.
  In particular, $p(\bold0)$ is 
the class of the $\SL(2)$-orbit associated to $(N_1, \dots, N_n, F)$.
\endproclaim

{\it Proof.} 
For $(N_1, \dots, N_n, F) \in \cD_{\nilp,n}$, let $\tau$ and $((\rho_w, \vf_w), \br_1, J) \in \cD_{\SL(2),n}$ be as in 2.4.2.
Write $J = \{a(1), \dots, a(r)\}$ with $a(1) < \dots < a(r)$.
Let $W^{(j)}=M(N_1+\dots+N_j, W)$ $(0 \le j \le n)$, where $W^{(0)}:= W$.
Let $\Psi = \{W^{(a(j))}\}_{1\le j\le r}$.
Let $\tau_J$ be the $J$-component of $\tau$.
Take $\a=\tau_J$ as a splitting of $\Psi$ (3.2.3) and take a distance to $\Psi$-boundary $\b$ (3.2.4).

For $t=(t_j)_{1\le j\le n }\in \bR_{>0}^n$, 
let $t'_J=(\prod_{a(j)\leq \ell <a(j+1)} t_{\ell})_{j \in J}$, where 
$a(r+1)$ means $n+1$.
Let $q(t)=\prod_{a(1) \leq \ell}\tau_{\ell}(t_{\ell})^{-1}p(t)$.
Then $q(t)=\tau_J(t'_J)^{-1}p(t)$.

First, we show that $q(t)$ extends to a real analytic map on some open neighborhood $U$ of $\bold0$ in $\bR_{\ge0}^n$.
  To see this, we may assume that $a(1)=1$. 
  Since $\tau(t)$ here coincides with $t(y)$ in \cite{KNU1} 0.5, in the notation there, we have
$$
q(t) = \tau(t)^{-1}p(t)=
{}^eg(y)\exp(\ve(y))\br.
$$
Hence, by loc.\ cit.\  0.5, the assertion follows.
  The extended map, also denoted by $q$, sends $\bold0$ to $\br_1 \in D$ 
in 2.4.2 (ii), that is, $q(\bold0) = \br_1$. 

  In case where $W \in \Psi$, 
since $\br_1 \in D_{\nspl}$, 
shrinking $U$ if necessary, we may assume 
that $p(t) \in D_{\nspl}$ for any $t \in U \cap \bR_{>0}^n$.

\medskip

\noindent
{\bf Claim 1.}
{\it After further replacing $U$, the map
$$
U\cap \bR_{>0}^n \to B:= \bR_{\ge 0}^\Psi \times D \times \spl(W) \times \tp_{W'\in \Psi} \spl(W'(\gr^W)), 
$$
$$
\qquad\qquad
t\mapsto (\b(p(t)),\, \tau_J\b(p(t))^{-1}p(t),\, \spl_W(p(t)),\, (\spl^{\BS}_{W'(\gr^W
)}(p(t)(\gr^W)))_{W'}),
$$
extends to a real analytic map $p':U\to B$ sending $\bold0$ to $(\bold0,\, \tau_J\b(\br_1)^{-1}\br_1,\, s,\, (s^{(W')})_{W'})$. 
  Here $s$ is the limiting splitting of $W$ in \cite{KNU1} {\rm 0.5 (1)}, which coincides with $\spl_W(\br_1)$ $(2.4.2)$, 
and $s^{(W')}$ is the splitting of $W'(\gr^W)$ given by $(\rho_w, \vf_w)_w$ {\rm (}cf. {\rm 3.2.6 (i))}.} 

\medskip

Since $\b(p(t)) = \b(\tau_J(t'_J)q(t)) = t'_J\b(q(t))$ (3.2.4), this extends to a real analytic map on some open neighborhood of $\bold0$ in $\bR_{\ge0}^n$ 
which sends $\bold0$ to $\bold0$. 

Since $\tau_J\b(p(t))^{-1}p(t) = \tau_J\b(q(t))^{-1}q(t)$, this extends to a real analytic map on some open neighborhood of $\bold0$ in $\bR_{\ge0}^n$ which sends $\bold0$ to $\tau_J\b(\br_1)^{-1}\br_1$. 

By \cite{KNU1} 0.5 (2), $\spl_W(p(t))$ extends to a real analytic map on some open neighborhood of $\bold0$ in $\bR^n_{\geq0}$ which sends $\bold0$ to $s$.

Finally, by \cite{KNU1} 8.5, $\spl^{\BS}_{W'(\gr^W)}(p(t)(\gr^W))$  extends to a real analytic map on some open neighborhood of $\bold0$ in $\bR_{\ge0}^n$ which sends $\bold0$ to $(s^{(W')})_{W'}$.

\medskip

  Next, it is easy to see that $(N'_1, \dots, N'_m, F')$ generates a 
nilpotent orbit $(2.4.1)$ for any $c$ in a sufficiently small $U$.
  Since its associated $\SL(2)$-orbits belong to $D_{\SL(2)}^I(\Psi)$, 
once we prove the following claim, the real analytic 
map $p':U \to B$ in Claim 1 factors through the image in $B$ of 
the map $\nu_{\a,\b}$ in 3.2.7 (i). 

\medskip

\noindent
{\bf Claim 2.}
{\it $\exp(\ts_{j=1}^n iy_jN_j)F$ converges to the class of the $\SL(2)$-orbit 
associated to $(N'_1, \dots, N'_m, F')$ in $D_{\SL(2)}^I$ 
when $t \in U$ and $t \to c$.}

\medskip

Thus we reduce both (i) and (ii) to this claim. 

  To prove Claim 2, 
we first consider the case $c = \bold 0$:
  In this case, the image by $\nu_{\a,\b}$ of the class of 
the $\SL(2)$-orbit $((\rho_w, \vf_w)_w, \br_1, J) \in \cD_{\SL(2),n}$ associated to $(N_1, \dots, N_n, F)$ is $\lim\nolimits_{t_J \to \bold 0_J}
(t_J\b(\br_1),\, \tau_J\b(\br_1)^{-1}\br_1,\, s,\, (s^{(W')})_{W'})$ by 
definition of $\nu_{\a, \b}$. 
  On the other hand, $p'(\bold 0)$ is $(\bold0,\, \tau_J\b(\br_1)^{-1}\br_1,\, s,\, (s^{(W')})_{W'})$ by Claim 1. 
  Since $\nu_{\a,\b}$ is injective (3.2.7 (i)), 
the case where $c = \bold 0$ of Claim 2 follows. 

  Now we are in the general case.  Let $c \in U$, $K$ be as in (ii).
Let $t' \in U$ be the element defined by $t'_j = t_j$ if $j \in K$, and $t_j = c_j$ if $j \not\in K$.
Then, by the case where $c = \bold 0$, 
we have the convergence 
$$\exp(\ts_{j\in J} iy'_jN'_j)F' \to \text{$($class of the $\SL(2)$-orbit associated to $(N'_1, \dots, N'_m, F'))$}.
$$
  Together with 
$$
\align
&\nu_{\a,\b}(\lim_{t\to c}p(t))=p'(c) = \lim_{t' \to c}p'(t') \\
&= \nu_{\a, \b}(\lim_{t' \to c}\exp(\ts_{j\in J} iy'_jN'_j)F'),
\endalign
$$
we have the general case of Claim 2.
\qed

\vskip20pt

\head
\S4.2. Hodge metrics at the boundary of $D_{\SL(2)}^I$
\endhead

\medskip

We expect that $D_{\SL(2)}$ plays a role as a natural space in which real analytic asymptotic behaviors of degenerating objects are well described. 
In this subsection we illustrate this by taking the degeneration of the Hodge metric as an example, and explain our previous result on the norm estimate in \cite{KNU1} via $D_{\SL(2)}^I$. 
\medskip

{\bf 4.2.1.}
Let $F \in D$. 
For $c>0$, we define a Hermitian form
$$
(\;,\;)_{F, c}\;:\;H_{0,\bC} \x H_{0,\bC}\to \bC
$$
as follows.

For each $w\in \bZ$, let
$$
(\;,\;)_{F(\gr^W_w)}: \gr^W_{w,\bC} \x \gr^W_{w,\bC}
\to \bC
$$
be the Hodge metric 
$\lan C_{F(\gr^W_w)}(\bullet), \bar\bullet\ran_w$, where $C_{F(\gr^W_w)}$ 
is the Weil operator.
For $v\in H_{0,\bC}$ and for $w\in \bZ$, let $v_{w,F}$ be the image in $\gr^W_{w,\bC}$ of the $w$-component of $v$ with respect to the canonical splitting of $W$ associated to $F$.  
Define
$$
(v, v')_{F,c}=\ts_{w\in \bZ} c^w(v_{w,F},\, v'_{w,F})_{F(\gr^W_w)}\quad (v, v'\in H_{0,\bC}).
$$

\medskip

\proclaim{Proposition 4.2.2}
Let $\Psi$ be an admissible set of weight filtrations on $H_{0,\bR}$ {\rm(3.2.2)}. 
Let $\b$ be a distance to $\Psi$-boundary {\rm(3.2.4--3.2.5)}.
Assume $W \not\in \Psi$ $($resp.\ $W \in \Psi)$. 
For each $W' \in \Psi$, let $\beta_{W'}: D\to \bR_{>0}$ $($resp.\ $D_{\nspl}\to \bR_{>0})$ be the $W'$-component of $\b$. 
For $p\in D$, let 
$$
(\;,\;)_{p, \b}:=(\;,\;)_{p, c}\quad \text{with}\;\;
c= \tp_{W' \in \Psi}\b_{W'}(p)^{-2}.
$$ 

Let $m : \Psi \to \bZ$ be a map, let $V = V_m = \tsize\bigcap_{W'\in\Psi} W'_{m(W'),\bC}$, and let $\Her(V)$ be the space of all Hermitian forms on $V$. 

Let $(\;,\;)_{p,\b,m}\in \Her(V)$ be the restriction of $\tp_{W' \in \Psi}\b_{W'}(p)^{2m(W')}(\;,\;)_{p,\b}$ to $V$.
\medskip 

{\rm(i)} The real analytic map $f : D$ $($resp.\ $D_{\nspl}) \to \Her(V)$, $p \mapsto (\;,\;)_{p,\b,m}$, extends to a real analytic map $f:D_{\SL(2)}^I(\Psi)$ $($resp. $D_{\SL(2)}^I(\Psi)_{\nspl}) \to\Her(V)$.
\medskip

{\rm(ii)} For a point $p \in D_{\SL(2)}^I(\Psi)$ $($resp.\ $p \in D_{\SL(2)}^I(\Psi)_{\nspl})$ such that $\Psi$ is the set of weight filtrations associated to $p$, the limit of $(\;,\;)_{p,\b,m}$ at $p$ induces a positive definite Hermitian form on the quotient space 
$$
V/(\ts_{m'<m}\tsize\bigcap_{W'\in\Psi} W'_{m'(W'),\bC}),
$$ 
where $m'<m$ means $m'(W') \le m(W')$ for all $W'\in\Psi$ and $m' \ne m$. 
\endproclaim

\demo{Proof} 
We prove (i).
Assume $W \not\in \Psi$.
Fix a splitting $\a:(\bR^\x)^\Psi\to\Aut_\bR(H_{0,\bR},W)$ of $\Psi$.
Let $p\in D$.
Let $v, v' \in V$.
Then, we have the weight decompositions $v=\ts_{m' \le m}v_{m'}$, $v'=\ts_{m' \le m}v'_{m'}$ with respect to $\a$.
Since
$$
\align
(v, v')_{p,\b}
&=(\a\b(p)(\a\b(p))^{-1}v, \a\b(p)(\a\b(p))^{-1}v')_{\a\b(p)(\a\b(p))^{-1}p,\b} \\
&=(\a\b(p)^{-1}v,\, \a\b(p)^{-1}v')_{\a\b(p)^{-1}p,1},
\endalign
$$
 we have
$$
\align
&(v, v')_{p,\b,m}\\
=\;&\tp_{W' \in \Psi}\b_{W'}(p)^{2m(W')}
(\a\b(p)^{-1}v,\, \a\b(p)^{-1}v')_{\a\b(p)^{-1}p,1} \tag1\\
=\;&\ts_{m', m'' \le m}
\tp_{W' \in \Psi}\b_{W'}(p)^{(2m-m'-m'')(W')}
(v_{m'},\, v'_{m''})_{\a\b(p)^{-1}p,1}.
\endalign
$$
This extends to a real analytic function on $D_{\SL(2)}^I(\Psi)$, 
because $(2m - m' -m'')(W')\ge0$ for all $W' \in \Psi$, and 
$D\to D, \;p\mapsto \a\b(p)^{-1}p$, extends to a real analytic map
$D_{\SL(2)}^I(\Psi)\to D$ (3.2.10 (i)).

In the case $W \in \Psi$, the argument is analogous.

We prove (ii).
Let $v, v' \in V$ be as above.
Let $\{p_\lambda\}_\lambda$ be a sequence in $D$ which converges to $p$, and let $q = \lim_\lambda \a\b(p_\lambda)^{-1}(p_\lambda) \in D$.
Then, by the result of (i), we have from (1)
$$
\lim_\lambda (v, v')_{p_\lambda,\b,m} 
= (v_m, v'_m)_{q,1}. \tag2
$$
The right-hand side of (2) is nothing but the restriction of the Hermitian metric at $q \in D$ to the $m$-component with respect to $\a$, which is therefore positive definite.
\qed
\enddemo

{\bf 4.2.3.}
As will be shown in a later part of our series, the norm estimate in 
\cite{KNU1} for a given variation of mixed Hodge structure 
$S \to D$ (cf.\ loc.\ cit. \S12) is 
incorporated in the diagram
$$
U\to D_{\SL(2)}^I(\Psi) \; (\text{resp.}\; D_{\SL(2)}^I(\Psi)_{\nspl}) @>f>> \Her(V). 
$$
Here $U$ is an open neighborhood of a point of 
$S_{\val}^\loga$, the first arrow is induced by 
an extension of the period map 
$S_{\val}^\loga \to \Gamma\bs D_{\SL(2)}^I$, where $\Gamma$ is an appropriate 
group (cf\. 4.4.9 below) and $f$ is 
as in 4.2.2.

\medskip

{\bf 4.2.4. Example V.} 
We consider Example V. 
Here the norm estimate is not continuous on $D_{\SL(2)}^{II}$.

Let $\Psi$ and $\bar\Psi$ be as in 3.6.2.  
Fix $u, v\in \bC e_4+ \bC e_5 \sub W'_1$, and let $u'$, $v'$ be their images in $\gr^W_1$, respectively.
Let $\b:D\to \bR_{>0}$ be the distance to $\Psi$-boundary which appears in 3.6.2. 

As in Proposition 4.2.2, the map
$$
f:D\to \bC, \quad 
p\mapsto  \b(p)^2(u, v)_{p, \b},
$$
extends to a real analytic function $f:D_{\SL(2)}^I(\Psi)\to \bC$. We show that however, for some choices of $u$ and $v$, this map $f:D_{\SL(2)}^I(\Psi)\to \bC$
is not continuous with respect to 
the topology of $D_{\SL(2)}^{II}$.
These can be explained by the following commutative diagram at the end of 3.6.2.
$$
\matrix 
(a_1,a_2,\,(a_3,\, x+iy),\, \tau)&\in\;  A_1\times A_2\times (A_3 \times \bar \fh^{\pm})' \times \fh   &\simeq &D_{\SL(2)}^I(\Psi)&@>f>>& \bC\\ 
&&&&&\\
\downarrow &\downarrow && \downarrow && \\
&&&&&\\
(a_1,a_2, \,|y|^{-1/2}a_3,\, x+iy,\, \tau)& \in\;\; (\prod_{j=1}^3 A_j) \times \bar \fh^{\pm} \times \fh  &\simeq & \;D_{\SL(2)}^{II}(\bar\Psi).&& 
\endmatrix
$$ 
Recall that $A_j = \Hom_{\bR}(\gr^W_1, \bR e_j)$ $(j=1,2,3)$. 
The composite $$A_1\times A_2\times (A_3 \times \bar \fh^{\pm})' \times \fh   \simeq D_{\SL(2)}^I(\Psi)@>f>> \bC$$ sends $(a_1,\, a_2,\, (a_3,\, x+iy),\, \tau)$ to  
$$
\align
(&|y|^{-3/2}a_1(u')+|y|^{-1/2}a_2(u')+a_3(u'), \, \\
& |y|^{-3/2}a_1(v')+|y|^{-1/2}a_2(v')+a_3(v'))_{0, (x+iy)/|y|}
+(u', v')_{1,\tau}. 
\endalign
$$ 
Here  $(\;,\;)_{0, (x+iy)/|y|}$ is the Hodge metric on $\gr^W_{0,\bC}$ associated to $(x+iy)/|y|\in \fh^{\pm}=
D(\gr^W_0)$, and $(\;,\;)_{1,\tau}$ is the Hodge metric on $\gr^W_{1,\bC}$ 
associated to $\tau \in \fh=D(\gr^W_1)$. On the other hand, the composition 
$$
\tsize\prod_{j=1}^3 A_j \times  \fh^{\pm} \times \fh \simeq D@>f>> \bC,
$$ 
  where the first arrow is induced by the lower (not upper) horizontal isomorphism of the above diagram, sends $(a_1,\, a_2,\, a_3,\, x+iy,\, \tau)$ to
  $$
\align
(&|y|^{-3/2}a_1(u')+|y|^{-1/2}a_2(u')+|y|^{1/2}a_3(u'), \, \\
& |y|^{-3/2}a_1(v')+|y|^{-1/2}a_2(v')+|y|^{1/2}a_3(v'))_{0, (x+iy)/|y|}
+(u', v')_{1,\tau}. 
\endalign
$$ 
 For some choices of $u$ and $v$, as is explained below precisely, the last map is not extended continuously to the point $(0,0,0, i \infty, i)$ of $\prod_{j=1}^3 A_j \times \bar \fh^{\pm} \times \fh$, for this map has the term $|y|^{1/2}$ which diverges at $i\infty$. 
Since $(0,0,0,i\infty, i)$ is the image of $(0,0,(0,i \infty), i)\in A_1\times A_2\times (A_3 \times \bar \fh^{\pm})' \times \fh$ under the left vertical arrow, this shows that for some choices of $u$ and $v$, $f:D_{\SL(2)}^I(\Psi)\to \bC$
is not continuous for  
the topology of $D_{\SL(2)}^{II}$.

  More precisely, take $u$ and $v$ such that there exists  $b\in A_3$ for which $(b(u'), b(v'))_{0, i}\neq 0$. 
Let $c$ be a real number such that $0 < c <1/2$. 
Then, as $y\to \infty$, $(0,0,y^{c-1/2}b, iy, i)\in \prod_{j=1}^3 A_j \times \fh^{\pm} \times \fh$ converges to $(0,0,0, i\infty, i)\in \prod_{j=1}^3 A_j \times \bar \fh^{\pm} \times \fh$. 
However, $f$ sends the image of $(0,0,y^{c-1/2}b, iy, i)$ in $D$ under the lower isomorphism of the diagram to  $(y^cb(u'), y^cb(v'))_{0,i}+(u',v')_{1,i}$, which diverges. 

\medskip

\vskip20pt

\head
\S4.3. Hodge filtrations at the boundary
\endhead

\medskip

{\bf 4.3.1.} 
In this \S4.3, let $X = D_{\SL(2)}^I$ or $D_{\SL(2)}^{II}$.
\medskip

Let $\cO_X$ be the sheaf of real analytic functions on $X$, and let $\a: M_X \to \cO_X$ be the log structure with sign on $X$.
We define a sheaf of rings $\cO'_X$ on $X$ by $\cO'_X:= \cO_X[q^{-1}\;|\; q \in \a(M_X)] \supset \cO_X$.
Let $\cO'_{X,\bC}=\bC\otimes_{\bR} \cO'_X$. The following theorem shows that the Hodge filtration over $\cO'_{X,\bC}$ extends to the boundary of $X$. 

\proclaim{Theorem 4.3.2} 
Let $X$ be one of $D_{\SL(2)}^I$, $D_{\SL(2)}^{II}$, and let $\cO'_X$ be as in $4.3.1.$ 

Then, for each $p\in \bZ$, there is a unique $\cO'_{X,\bC}$-submodule $F^p$ of $\cO'_{X, \bC} \otimes_{\bZ} H_0$ which is locally a direct summand and whose restriction to $D$ coincides with the filter $F^p$ of
$\cO_{X,\bC} \otimes_{\bZ} H_0$.
\endproclaim

{\it Proof.} 
It is sufficient to prove the case $X=D_{\SL(2)}^{II}$, because the assertion for $X=D_{\SL(2)}^I$ follows from that for $X=D_{\SL(2)}^{II}$ by pulling back.

Assume $X=D_{\SL(2)}^{II}$. 
Let $F$ be the universal Hodge filtration on $D$, and write 
$F = s(\theta(F', \delta))$ ($s \in \spl(W)$, $F' \in D(\gr^W)$, $\delta \in \cL(F')$) as in 1.2.5. 
Let $\Phi$ be an admissible set of weight filtrations on $\gr^W$ (3.2.2), and let $\a$ be a splitting of $\Phi$ and $\b$ a distance to $\Phi$-boundary as in 3.2.5 (ii).
We observe
$$
s(\theta(F', \delta)) = s(\theta(\alpha\beta (F') (\alpha\beta(F'))^{-1}F',\, \Ad(
\alpha\beta(F')) \Ad(\alpha\beta(F'))^{-1}\delta)). \tag1
$$
By 3.2.6 (ii), $(\alpha\beta(F')^{-1}F',\, \Ad(\alpha\beta(F'))^{-1}\delta)$ 
and $s$ extend real analytically over the $\Phi$-boundary. 
Let $G' = \prod_w\Aut(\gr^W_w)$, and we consider the splitting $\a : \bG_m^\Phi \to G'$.
Then, the section $\b(F')$ of $\bG_m^\Phi(\cO_X')$ on $D_{\SL(2)}^{II}(\Phi)$ is sent to a section $\alpha\beta(F')$ of $G'(\cO_X')$ over $D_{\SL(2)}^{II}(\Phi)$.
Thus, $F = s(\theta(F', \delta))$ extends uniquely to a filtration of $\cO'_{X,\bC}\otimes H_0$ consisting of $\cO'_{X,\bC}$-submodules which are locally direct summands.
\qed

\medskip

{\bf 4.3.3.} {\it Remarks.} 
(i) For $D_{\SL(2),\val}$, $D_{\BS}$, $D_{\BS,\val}$, theorems similar to 4.3.2 are analogously proved.

(ii) The Hodge decomposition and the Hodge metric also extend over the 
boundary after tensoring with $\cO'_{X, \bC}$. 
In the pure case, this together 
with the period map $S^{\log}_{\val}\to \Gamma \bs D_{\SL(2)}$ (4.2.3) 
explains the existence of the log $C^{\infty}$ Hodge decomposition in 
\cite{KMN}.

\vskip20pt

\head
\S4.4.  Example IV and height pairing
\endhead
\bigskip

\medskip
We consider Example IV. 
The space $D_{\SL(2)}=D_{\SL(2)}^I=D_{\SL(2)}^{II}$ in this example is related to the asymptotic behavior of the Archimedean height pairing for elliptic curves in degeneration (cf.\ \cite{P2}, \cite{C}, \cite{Si}). 
We describe which kind of $\SL(2)$-orbits appear in such geometric situation of degeneration.

The following observations were obtained in the discussions with Spencer Bloch.

\medskip
{\bf 4.4.1.}
Recall (\cite{A}) that the {\it Archimedean height pairing} for an elliptic curve $E$ over $\bC$ is $\langle Z, W\rangle\in \bR$ defined for divisors $Z$, $W$ on $E$ of degree $0$ such that $|Z|\cap |W|=\emptyset$ ($|Z|$ here denotes the support of $Z$), characterized by the following properties (1)--(4):
\medskip

(1) If $|Z|\cap |W|=|Z'|\cap |W|=\emptyset$, then $\langle Z+Z', W \rangle = \langle Z, W \rangle+\langle Z', W \rangle$.
\medskip

(2) $\langle Z, W\rangle=\langle W, Z\rangle$.
\medskip

(3) If $f$ is a meromorphic function on $E$ such that $|(f)|\cap |W|=\emptyset$ and if $W=\tsize\sum_{w\in |W|} n_w(w)$, then $\langle (f), W\rangle= 
-(2\pi)^{-1}\tsize\sum_{w\in |W|} n_w\log(|f(w)|)$.
\medskip

(4) The map
$(E(\bC) \smallsetminus |W|)  \times (E(\bC) \smallsetminus |W|) \to 
\bR$, $(a,
b) \mapsto \langle (a)-(b), W\rangle$, 
is continuous.
\medskip

{\bf 4.4.2.} 
Consider Example IV.

Let $\tau\in \fh$ and let $E_\tau$ be the elliptic curve $\bC/(\bZ\tau+\bZ)$.

For divisors $Z$, $W$ on $E_\tau$ of degree $0$ such that $|Z|\cap |W|=\emptyset$, we define an element
$$
p(\tau, Z, W)\in G_{\bZ,u}\bs D
$$
as follows. 

For $\tau\in \fh$ and $z\in \bC$, let
$$
\theta(\tau, z) = \tsize\prod_{n=0}^\infty  (1-q^nt) \cdot \tsize\prod_{n=1}^\infty (1-q^nt^{-1}),\quad \text{where}\;\;q=e^{2\pi i\tau}, t=e^{2\pi iz}.
$$
We have
$$
\theta(\tau, z+1)=\theta(\tau, z), \quad \theta(\tau, z+\tau)=-e^{-2\pi iz}\theta(\tau, z).\tag1
$$

Write 
$$
Z=\tsize\sum_{j=1}^r m_j(p_j),  \quad 
W=\tsize\sum_{j=1}^s n_j(q_j)
$$
($p_j, q_j\in E_\tau$, $m_j, n_j\in \bZ$, $\;\tsize\sum_{j=1}^r m_j=0$, $\;\tsize\sum_{j=1}^s n_j=0$), and write 
$$
p_j= (z_j\bmod (\bZ\tau+\bZ)), \quad q_j= (w_j\bmod (\bZ\tau+\bZ)), 
$$
with $z_j, w_j\in \bC$.  
Define 
$$
p(\tau, Z, W) =\text{class of}\; F(\tau, w, \lambda, z) \in G_{\bZ,u}\bs D, \tag2
$$
$$
\text{with}\quad z=\tsize\sum_{j=1}^r m_jz_j, \quad w=\tsize\sum_{j=1}^s n_jw_j, \quad \lambda= (2\pi i)^{-1}\log(\tsize\prod_{j, k} \theta(\tau, z_j-w_k)^{m_jn_k}),
$$ 
\medskip 
\noindent
and with $F(\tau, w, \lam, z)\in D$ as in 1.1.1, Example IV. 
This element $p(\tau, Z, W)$ of $G_{\bZ,u}\bs D$ is well-defined: 
As is easily seen using (1),
 the right hand side of (2) does not change when we replace $((z_j)_j, (w_j)_j)$ by $((z'_j)_j, (w'_j)_j)$ such that $z_j'\equiv z_j\bmod \bZ\tau+\bZ$ and $w'_j\equiv w_j\bmod \bZ\tau+\bZ$ for any $j$. 
  For example, in the case where $z_\ell'=z_\ell+\tau$ for some $\ell$, $z'_j=z_j$ for 
the other $j\neq \ell$, and $w'_j=w_j$ for any $j$, 
by (1), the right hand side of (2) given by $(z'_j)_j,  (w'_j)_j$ is the class of
$F(\tau, w, \lambda+ m_{\ell}w, z+m_{\ell}\tau)=
\gamma F(\tau, w, \lambda, z)$, where $\gamma$ is the element of $G_{\bZ, u}$
which sends $e_j$ ($j=1, 2, 3$) to $e_j$ and $e_4$ to $e_4-m_{\ell}e_3$.

\medskip

{\bf 4.4.3.} 
Let $L=\cL(F)$ with $F\in D(\gr^W)$, which is independent of $F$, and let $\delta: D\to L=\bR$ be the $\delta$-component (1.2.5). 
Note that
$$
\delta(F(\tau, w, \lambda, z))=\Im(\lambda)-\Im(z)\Im(w)/\Im(\tau)
$$
(1.2.9, IV).

\proclaim{Lemma 4.4.4} 
The map $\delta: D\to \bR$ factors through the projection $D\to G_{\bZ,u}\bs D$, and we have 
$$
\delta(p(\tau, Z, W))=\langle Z, W\rangle,
$$
where $\langle Z, W\rangle\in \bR$ is the Archimedean height pairing $(4.4.1)$.
\endproclaim

{\bf 4.4.5.} 
The equality in Lemma 4.4.4 is well known. 
It has also the following geometric (cohomological)  interpretation.

Let $E$ be an elliptic curve over $\bC$, and let $Z$ and $W$ be divisors of degree $0$ on $E$ such that $|Z|\cap |W|=\emptyset$. 
We assume $Z\neq 0$, $W\neq 0$.

Let $U=E\smallsetminus |Z|$, $V=E\smallsetminus (|Z|\cup |W|)$,  and let $j:V\to U$ be the inclusion map. 
Write $Z=\tsize\sum_{z\in |Z|} m_z(z)$, $W=\tsize\sum_{w\in |W|} n_w(w)$. We have exact sequences of mixed Hodge structures
$$
0 \to H^1(E, \bZ)(1)\to H^1(U,\bZ)(1)\to \bZ^{|Z|} \to H^2(E, \bZ)(1)\to 0,
$$
$$
0\to H^0(U, \bZ)(1)\to \bZ(1)^{|W|}\to H^1(U, j_!\bZ)(1)\to H^1(U, \bZ)(1) \to 0.
$$
Note that the map $\bZ^{|Z|}\to H^2(E, \bZ)(1)=\bZ$ is identified with the degree map. 
Let $A\subset B\subset H^1(U, j_!\bZ)(1)$ be sub mixed Hodge structures defined as follows. 
$A$ is the image of $\{x=(x_w)_w 
\in \bZ(1)^{|W|}\;|\;\tsize\sum_w n_wx_w=0\}$ under $\bZ(1)^{|W|}\to H^1(U, j_!\bZ)(1)$. 
$B$ is the inverse image of $\{(m_zx)_z\;|\;x\in \bZ\}$ under the composition $H^1(U, j_!\bZ)(1)\to H^1(U, \bZ)(1) \to \bZ^{|Z|}$. 
Let $H=B/A$. 
Then we have the induced 
injective homomorphism $a:\bZ(1) \to H$, the induced 
surjective homomorphism $b:H\to \bZ$, and $\Ker(b)/\text{Im}(a) = H^1(E,\bZ)(1)$. 
A well-known  cohomological interpretation of the height pairing $\langle Z, W\rangle$ is 
$$
\langle Z, W\rangle = \delta(H).
$$ 
On the other hand, in the case $E=E_\tau$, as is well known,
$$
p(\tau, Z, W) = \text{class}(H).
$$
This explains Lemma 4.4.4.

\medskip

{\bf 4.4.6.}
We consider degeneration.

Let $\Delta=\{q\in \bC\;|\;|q|<1\}$, and let $\Delta^*=\Delta\smallsetminus \{0\}$. 
Fix an integer $c\geq 1$, and consider the family of elliptic curves over $\Delta^*$ whose fiber over $e^{2\pi i\tau/c}$ ($\Im(\tau)>0$) is $\bC/(\bZ\tau +\bZ)$. 
This family has a N\'eron model $E_c$ over $\Delta$ whose fiber over $0\in \Delta$ is canonically isomorphic to $\bC^\times \times \bZ/c\bZ$ as a Lie group. 
If $a\in \bQ$ and $ca\in \bZ$, and if $u$ is a holomorphic function $\Delta\to \bC^\times$, there is a section of $E_c$ over $\Delta$ whose restriction to $\Delta^*$ is given by $e^{2\pi i\tau/c}\mapsto (a\tau + f(e^{2\pi i\tau/c})\modu \bZ\tau +\bZ)$ with $f=(2\pi i)^{-1}\log(u)$ and whose value at $0\in \Delta$ is $(u(0), ca\modu c\bZ)\in \bC^\times \times \bZ/c\bZ$. 
Any section of $E_c$ over $\Delta$ is obtained in this way.

Let $\G\subset G_\bZ$ be the subgroup consisting of all elements $\gamma$ which satisfy $\gamma(e_j)-e_j\in \bigoplus_{1\leq k< j} \bZ e_k$ for $j=1, 2, 3, 4$. 
Note $\G\supset G_{\bZ,u}$. Note also that $\delta:D \to L=\bR$ factors through the projection $D\to \G\bs D$. 

Fix $m_j, n_k\in \bZ$, $a_j, b_k\in \bQ$ ($1\leq j\leq r$, $1\leq k\leq s$) such that $\tsize\sum_j m_j=0$ and $\tsize\sum_k n_k=0$, $ca_j, cb_k\in \bZ$ for any $j, k$, and take holomorphic functions $u_j, v_k : \Delta\to \bC^\times$  ($1\leq j\leq r, 1\leq k\leq s$). 
Assume that, for any $j, k$, the section $p_j$ 
of $E_c$ defined by $(a_j, u_j)$ and the section $q_k$ 
of $E_c$ defined by $(b_k, v_k)$ do not meet over $\Delta$. 
Consider the morphism 
$$
\aligned
p:\Delta^*\to \G \bs D, \quad 
&e^{2\pi i\tau/c} \mapsto  (p(\tau,\, \tsize\sum_j m_j(p_j), \tsize\sum_k n_k(q_k))\modu \G)\\
\text{with}\quad
& p_j:=(a_j\tau +f_j(e^{2\pi i\tau/c})\modu \bZ\tau+\bZ),\\
&q_k:=(b_k\tau+ g_k(e^{2\pi i\tau/c})\modu \bZ\tau+\bZ), \\
\text{where}\quad
&f_j:=(2\pi i)^{-1}\log(u_j), \quad g_k :=(2\pi i)^{-1}\log(v_k).
\endaligned
$$

{\bf 4.4.7.} 
Let $\Delta^{\log}= |\Delta|\times \bS^1$, where
$|\Delta|:=\{r\in \bR\;|\;0\leq r<1\}$, $\bS^1:=\{u\in \bC^\times\;|\;|u|=1\}$.
We have a projection $\Delta^{\log}\to \Delta,\;(r, u) \mapsto ru$ ($r\in |\Delta|$, $u\in \bS^1$), and an embedding $\Delta^*\to \Delta^{\log},\;ru\mapsto (r, u)$ ($r\in |\Delta|$, $r\neq 0$, $u\in \bS^1$). 

We define the sheaf of $C^\infty$ functions on $\Delta^{\log}$ as follows: 
For an open set $U$ of $\Delta^{\log}$ and a real valued function $h$ on $U$,  $h$ is $C^\infty$ if and only if the following (1) holds. 
Let $U'$ be the inverse image of $U$ in $\bR_{\geq 0} \times \bR$ under the surjective map $\bR_{\geq 0}\times \bR \to \Delta^{\log},\;(t, x)\mapsto 
(e^{-1/t^2},\, e^{2\pi ix})$. 
\medskip

(1)  The pull-back of $h$ on $U'$ extends, locally on $U'$, to a $C^\infty$ function  on some open neighborhood of $U'$ in $\bR^2$. 
\medskip

Roughly speaking, a function $h$ on $\Delta^{\log}$ is $C^\infty$ if $h(e^{2\pi i(x+iy)})$ ($x\in \bR,\, 0<y\leq \infty$) is a $C^\infty$ function in $x$ and $1/\sqrt{y}$.

The restriction of this sheaf of $C^\infty$ functions on $\Delta^{\log}$ to the open set $\Delta^*$ coincides with the usual sheaf of $C^\infty$ functions on $\Delta^*$. 

\proclaim{Proposition 4.4.8} 
Let $\Phi\in \overline{\cW}$ be as in $3.6.1$, {\rm{IV}}. 
\medskip

{\rm(i)} The map $p:\Delta^*\to \G\bs D$ in $4.4.6$ extends to a $C^{\infty}$ map $\Delta^{\log}\to \G\bs D_{\SL(2)}^{II}(\Phi)$. 
That is, we have a commutative diagram of local ringed spaces over $\bR$
$$
\matrix \Delta^* & @>p>> & \G\bs D\\
&&\\
\cap & & \cap\\
&&\\
\Delta^{\log} & @>>> & \;\G\bs D_{\SL(2)}^{II}(\Phi).
\endmatrix
$$
\medskip

{\rm(ii)} Let $B_2(x)$ be the second Bernoulli polynomial $x^2-x+1/6$. 
For $x\in \bR$, $\{x\}$ denotes the unique real number such that
$0\leq \{x\}<1$ and $\{x\}\equiv x\bmod \bZ$. 

Then the composite $\Delta^* @>p>> \G \bs D@>\delta>> L=\bR$
has the form 
$$
e^{2\pi i(x+iy)/c}\mapsto \frac{1}{2}(\tsize\sum_{j,k} m_jn_kB_2(\{a_j-b_k\}))y+h(e^{2\pi i(x+iy)/c})
$$ 
for some $C^\infty$ function $h$ on $\Delta^{\log}$. 

\medskip

{\rm (iii)} Let 
$$
D_{\SL(2)}^{II}(\Phi) \simeq \spl(W) \times \bar \fh \times \bar L
$$ 
be the lower isomorphism in the commutative diagram in $3.6.1$, {\rm{IV}}. 
Then the projection $D_{\SL(2)}^{II}(\Phi)\to \bar L$ factors through
$D_{\SL(2)}^{II}(\Phi) \to \G \bs D_{\SL(2)}^{II}(\Phi)$, and the composite
$\Delta^{\log}@>p>> \G \bs D_{\SL(2)}^{II}(\Phi) \to \bar L$ sends  any point of $\Delta^{\log}\smallsetminus \Delta^*$ to 
$$
\frac{1}{2}(\tsize\sum_{j,k} m_jn_k B_2(\{a_j-b_k\})) \in \bR=L\sub \bar L.
$$

\endproclaim

In (ii) and (iii), $B_2(x)$ can be replaced by the polynomial $x^2-x$. 
The constant term of $B_2(x)$ does not play a role, for $\tsize\sum_{j,k} m_jn_k=0$.

Note that the restriction of the map  $D_{\SL(2)}^{II}(\Phi)\to \bar L$ in (iii) to $D$ is {\it not} $\delta: D\to L$ {\it but} is $p\mapsto \Ad(\a\b(p(\gr^W)))^{-1}\d(p)$, where $\a$ and $\b$ are as in 3.6.1, IV.

\medskip

{\it Proof of Proposition 4.4.8.} 
We may and do assume $0\leq a_j<1$ and $0\leq b_k<1$. 
Let $J=\{(j, k)\;|\; 1\leq j\leq r,\, 1\leq k\leq s$,\, $a_j<b_k\}$. 
Then, for each $j$ and $k$, the function
$$
e^{2\pi i\tau/c}\mapsto \theta(\tau,\, (a_j-b_k)\tau+f_j(e^{2\pi i\tau/c})-g_k(e^{2\pi i\tau/c}))
$$
on $\Delta$ is
meromorphic and its order of zero at $0\in \Delta$ is 
$(a_j-b_k)c$ if $(j, k) \in J$ and is $0$ otherwise.
By using this and by using the description of $\spl_W : D\to \spl(W)$ in 1.2.9, IV, we see that the composite
$$
\Delta^* @>p>> \G \bs D @>1.2.9>\simeq> \G\bs (\spl(W) \times \fh)\times L
$$
has the property that the part $\Delta^*\to \G\bs (\spl(W) \times \fh)$ extends to a $C^\infty$ function $\Delta^{\log} \to \G \bs (\spl(W)\times \bar \fh)$, and that the part $\Delta^*\to L=\bR$ has the form
$e^{2\pi i\tau/c}\mapsto (-(\tsize\sum_j m_ja_j)(\tsize\sum_k n_kb_k)+\tsize\sum_{(j,k)\in J} m_jn_k(a_j-b_k))\Im(\tau)+h(e^{2\pi i\tau/c})$, 
where $h$ is a $C^\infty$ function on $\Delta^{\log}$. 
Note that
$$
-({\tsize\sum}_j m_ja_j)({\tsize\sum}_k n_kb_k)+{\tsize\sum}_{(j,k)\in J} m_jn_k(a_j-b_k)
=\frac{1}{2}({\tsize\sum}_{j,k} m_jn_k B_2(\{a_j-b_k\})).
$$
Hence, for the lower isomorphism $D_{\SL(2)}^{II}(\Phi) \cong \spl(W) \times \bar \fh \times \bar L$ in the diagram in 3.6.1, IV, the composite 
$\Delta^*\to \bar L$ is written as $e^{2\pi i\tau/c}\mapsto (1/2)(\tsize\sum_{j,k} m_jn_k B_2(\{a_j-b_k\}))+(\Im(\tau))^{-1}h(e^{2\pi i\tau/c})$, 
where $(\Im(\tau))^{-1}h$ is a $C^\infty$ function on $\Delta^{\log}$ which has value $0$ on $\Delta^{\log}\smallsetminus \Delta^*$. 
These imply the assertions.
\qed

 \medskip

 {\bf 4.4.9.} 
The above proposition 4.4.8 implies 
a special case of the height estimate by Pearlstein (\cite{P2}).

The lower map in the diagram in 4.4.8 (i) is an example of the extended period map (cf.\ 4.2.3).
In a forthcoming part of this series of papers, the existence of the extended period map $X^{\log}_{\val}\to \G\bs D_{\SL(2)}$ ($X$ is a log smooth fs log analytic space) will be proved generally for a variation of mixed Hodge structure on $U=X_{\triv}$ with polarized graded quotients with global monodromy in 
an appropriate group $\G$ which has unipotent local monodromy along 
$D=X\smallsetminus U$ and is admissible at the boundary. 
  This will be accomplished by the CKS map $D^{\sharp}_{\Sig,\val} \to D_{\SL(2)}$ in the fundamental diagram in 0.2 (see \cite{KU3}, 8.4.1 for the pure case), and imply the height estimate of Pearlstein for more general cases.

\vskip40pt

\noindent 
{\bf Correction to Part I.}
  There are some mistakes in calculating examples in Part I (\cite{KNU2}), 
\S10. 
First, the \lq\lq$r^{-2}$''s in 10.2.1 should be \lq\lq$r^{-1}$''. 
  (Note that we gave the real analytic structure on $\bar A_P$ in 
the notation in loc.\ cit.\ 2.6 by using the fundamental roots.) 
  There are similar mistakes also in 10.3, that is, $r$ should be 
replaced by $r^{1/2}$ in the third last line in p.219 of loc.\ cit., which 
should be $(x+ir^{-1}, \ldots)$, in the second line in p.220: 
$(s_1, s_2, x, r, d) \mapsto x+ir^{-1}$, 
and in the second last line in p.220: $t(r)(e_1)=r^{-1/2}e_1$, 
$t(r)(e_2)=r^{1/2}e_2$.

\medskip

\enddocument